\documentclass{amsart} 

\usepackage{amsmath,amssymb}
\usepackage{amscd}
\usepackage{amsthm}
\usepackage{latexsym}

\usepackage{color}

\newtheorem{theorem}{Theorem}[section]
\newtheorem{lemma}[theorem]{Lemma} 

\newtheorem{corollary}[theorem]{Corollary}
\newtheorem{proposition}[theorem]{Proposition}

\theoremstyle{remark}
\newtheorem{remark}[theorem]{Remark} 

\newtheorem{definition}[theorem]{Definition} 

\newtheorem{example}[theorem]{Example}

\numberwithin{equation}{section}

\begin{document}

\title[Endomorphisms and automorphisms of the shift]
{The degrees of onesided resolvingness and the limits 
of onesided resolving directions for  
endomorphisms and automorphisms of the shift} 

\author{Masakazu Nasu}
\address{19-8, 9-ch\={o}me, Takaya-Takamigaoka, Higashi-Hiroshima 
739-2115, Japan}
\curraddr{}
\email{masakazu.nasu@gmail.com}
\thanks{}


\subjclass[2000]{Primary: 37B10; 
Secondary: 37B15, 54H20.}


\dedicatory{}

\begin{abstract} 
\newcommand{\Z}{\mathbf{Z}}
We introduce the notions 
in the title 
for endomorphisms of subshifts, 
and using them we 
characterize various classes of 
``resolving  endomorphisms of subshifts'' in the broad sense 
including onesided and weak ones. Resolving 
endomorphisms of transitive topological Markov shifts 
and resolving automorphisms of topological Markov shifts 
are treated 
particularly in detail. Moreover, 
we understand what the limits of 
onesided resolving directions are  
for ``expansiveness''  
(in the broad sense 
including onesided ones) of endomorphisms of subshifts, and 
we understand the relation  between ``resolvingness'' 
and ``expansiveness''  for endomorphisms and automorphisms 
of subshifts and for $\Z^d$-actions on  zero-dimensional 
compact metric spaces.
\end{abstract}

\maketitle


\newcommand{\Z}{\mathbf{Z}}
\newcommand{\N}{\mathbf{N}}
\newcommand{\R}{\mathbf{R}}
\newcommand{\Q}{\mathbf{Q}}
\newcommand{\bu}{\bar{u}}
\newcommand{\bv}{\bar{v}} 
\newcommand{\G}{\Gamma} 
\newcommand{\D}{\Delta}
\newcommand{\q}{\quad} 
\newcommand{\itl}{\textit} 
\newcommand{\aseq}{(a_j)_{j\in\Z}} 
\newcommand{\bseq}{(b_j)_{j\in\Z}} 
\newcommand{\cseq}{(c_j)_{j\in\Z}} 
\newcommand{\dseq}{(d_j)_{j\in\Z}} 
\newcommand{\vseq}{(a_i)_{i\in\Z}}
\newcommand{\xseq}{(x_i)_{i\in\Z}}
\newcommand{\yseq}{(y_i)_{i\in\Z}}
\newcommand{\bsigma}{\mathbf{\sigma}}
\newcommand{\C}{\mathcal{C}}
\newcommand{\Ll}{\mathcal{L}}
\newcommand{\bell}{\mathbf{\ell}}

\newcommand{\A}{\mathcal{A}}
\newcommand{\cD}{\mathcal{D}}
\newcommand{\E}{\mathcal{E}}

\newcommand{\PL}{\mathcal{PL}}
\newcommand{\QR}{\mathcal{QR}}
\newcommand{\CR}{\mathcal{CR}}
\newcommand{\PR}{\mathcal{PR}}
\newcommand{\QL}{\mathcal{QL}}
\newcommand{\CL}{\mathcal{CL}}

\newcommand{\Oh}{\mathcal{O}}
\newcommand{\plane}{\mathbf{R}^2}
\newcommand{\pplane}{(\mathbf{R}^2)'}
\newcommand{\hf}{\overline} 

\newcommand{\rt}{$\mathcal{P}$ }
\newcommand{\st}{$\mathcal{P}$}

\def\uo{\accent'27u}  
\newcommand{\goesfromto}{\leftrightsquigarrow} 
\newcommand{\pl}{\rightharpoonup}
\newcommand{\qr}{\rightharpoondown}
\newcommand{\lr}{\to}
\newcommand{\mm}{\multimap}
\newcommand{\cG}{\mathcal{G}} 
\newcommand{\M}{\mathcal{M}} 
\newcommand{\cH}{\mathcal{H}} 

\newcommand{\bel}{\bar{\ell}}
\newcommand{\bfk}{\mathbf{k}}
\newcommand{\bfv}{\mathbf{v}}
\newcommand{\bfzero}{\mathbf{0}}
\newcommand{\bfl}{\mathbf{l}}
\newcommand{\bfm}{\mathbf{m}}
\newcommand{\bfn}{\mathbf{n}}

\def\aspace{\vskip 0.5cm}
\def\bspace{\vskip 0.3cm}
\def\cspace{\vskip 0.1cm}
\renewcommand{\labelenumi}{(\theenumi)}
\renewcommand{\labelitemi}{$\circ$}


\tableofcontents

\section{Introduction} 

This paper is a sequel of \cite{Nasu-te}. We continue to study 
the dynamics of endomorphisms and automorphisms of subshifts, 
in particular, the ``overall dynamics'' of an  
endomorphism $\varphi$ of a subshift $(X,\sigma)$, i.e. 
the dynamics of $\varphi^i\sigma^j$ for all $i\in\N,j\in\Z$ 
and the dynamics of $\varphi^i\sigma^j$ 
for all $i,j\in\Z$ when $\varphi$ is an automorphism. 
We are especially interested in
the overall dynamics of  
onto endomorphisms of topological Markov shifts, 
which include onto cellular automata. 
We introduce the notions of 
the degrees of onesided resolvingness, 
which are not topological invariants, 
and the limits of onesided resolving directions, 
which are topological invariants, 
for endomorphisms of subshifts. 
 Using them we characterize various classes of 
``resolving endomorphisms'' (in the broad sense 
including onesided and weak ones) of subshifts 
which were introduced and studied in 
\cite{Nasu-t, Nasu-m, Nasu-te}. 
We also further consider 
left $\tau$-expansiveness and right $\tau$-expansiveness,  
which were introduced  for onto endomorphisms 
of invertible compact dynamical systems $(X,\tau)$ \cite{Nasu-te}. 
Applying geometric methods of Boyle and Lind    
which define expansive lines and expansive components 
for $\Z^d$-actions in \cite{BoyLin}, we define    
left $\tau$-expansive and right $\tau$-expansive analogues 
of these for onto endomorphisms of $(X,\tau)$. 
Using them we understand what the limits of onesided resolving
directions of an onto endomorphism of a subshift are 
and also understand the relation between  
``resolvingness'' and  ``expansiveness'', 
each in the broad sense including onesided ones, for 
endomorphisms and automorphisms of subshifts 
in view of their overall dynamics. We add discussions 
to understand  
the relation between  
``resolvingness'' and  ``expansiveness''
for $\Z^d$-actions on 
zero-dimensional compact metric spaces.

In the remainder of this section, 
we outline the main points 
of this paper. 

Let $(X,\sigma)$ be a subshift over an alphabet $A$. 
For $s\geq 1$, let $L_s(X)$ denote the set of all words 
$a_j\dots a_{j+s+1}$ that appear 
on some point $(a_j)_{j\in\Z}\in X$ 
with $a_j\in A$. For $N\geq 0$, a \itl{local rule 
of neighborhood-size $N$ on $(X,\sigma)$} 
means a mapping $f:L_{N+1}(X)\to L_1(X)$ such that 
$(f(a_j\dots a_{j+N}))_{j\in\Z}\in X$ 
for all points $(a_j)_{j\in\Z} \in X$. 

Let $f:L_{N+1}(X)\to L_1(X)$ be a local rule on $(X,\sigma)$. 
Let $I\geq 0$. We say that $f$ is 
\itl{$I$ left-redundant} if for any points  
$\aseq$ and 
$\bseq$ in $X$ with 
$a_j,b_j\in L_1(X)$, it holds that if
$(a_j)_{j\geq 0}=(b_j)_{j\geq 0}$ then 
$f(a_{-I}\dots a_{-I+N})=f(b_{-I}\dots b_{-I+N})$. 
Symmetrically, 
$f$ is said to be \itl{$I$ right-redundant} if for any points 
$\aseq$ and $\bseq$ in $X$, 
it holds that 
if $(a_j)_{j\leq 0}=(b_j)_{j\leq 0}$ then 
$f(a_{I-N}\dots a_I)=f(b_{I-N}\dots b_I)$.  
We say that $f$ is 
\itl{strictly I left-redundant} 
if it is $I$ left-redundant but 
not $I+1$ left-redundant; if $f$ is $I$ left-redundant for 
all $I\geq 0$, then $f$ is said to be 
\itl{strictly $\infty$ left-redundant}. 
Similarly,  
a \itl{strictly $I$ right-redundant} local rule with $I\geq 0$ and 
a \itl{strictly $\infty$ right-redundant} local rule are defined. 

Let $k\geq0$.
We say that $f$ is \itl{$k$ right-mergible} 
if for any points  
$\aseq$ and 
$\bseq$ in $X$, 
it holds that if
$(a_j)_{j\leq 0}=(b_j)_{j\leq 0}$
and 
$f(a_{j-N}\dots a_{j})
=f(b_{j-N}\dots b_{j})$ for $j=1,\dots, k+1$,  
then $a_1=b_1$. 
We say that $f$ is \itl{$k$ left-mergible} 
if for any points  
$\aseq$ and $\bseq$ in $X$, 
it holds that if 
$(a_j)_{j\geq 0}=(b_j)_{j\geq 0}$
and 
$f(a_{-j}\dots a_{-j+N})
=f(b_{-j}\dots b_{-j+N})$ for $j=1,\dots,k+1$,  
then $a_{-1}=b_{-1}$. 
We say that $f$ is \itl{strictly $0$ right-mergible} if $f$ is 
$0$ right-mergible;
for $k\geq 1$ 
we say that $f$ is \itl{strictly $k$ right-mergible} 
if $f$ is $k$ right-mergible but not $k-1$ right-mergible. 
Similarly, a \itl{strictly $k$ left-mergible} 
local rule is defined for $k\geq 0$. 
We say that $f$ is \itl{strictly $\infty$ right-mergible} 
(respectively, \itl{strictly $\infty$ left-mergible}) if 
$f$ is not $k$ right-mergible (respectively, not $k$ left-mergible) 
for all $k\geq 0$, which is equivalent to the condition that 
an endomorphism of $(X,\sigma)$ given by $f$ (see the next paragraph) is 
not \itl{right-closing} (respectively, not \itl{left-closing}) 
in the well-known sense of Kitchens \cite{Kit-c}.
 
Suppose that $\varphi$ is an endomorphism of 
\itl{$(m,n)$-type} 
of a subshift $(X,\sigma)$ \itl{given} by a local rule 
$f: L_{m+n+1}(X)\to L_1(X)$ on $(X,\sigma)$, i.e., 
$\varphi$ is defined by 
\[\varphi((a_j)_{j\in\Z})= (f(a_{j-m}\dots a_{j+n}))_{j\in\Z}, 
\q (a_j)_{j\in\Z}\in X,\;\, a_j\in L_1(X).\] 
Let $I,J,k,l\in\{0,1,\dots,\infty\}$. 
If $f$ is strictly $I$ left-redundant and strictly $J$ right-redundant, 
then we can define  
the \itl{$p$-L degree $P_L(\varphi)$ of $\varphi$} and 
the \itl{$p$-R degree $P_R(\varphi)$ of $\varphi$} by  
\[P_L(\varphi)=I-m,\q\q\q\q P_R(\varphi)=J-n\] 
(Proposition 6.3). We can define 
\itl{$q$-R degree $Q_R(\varphi)$ of $\varphi$} 
and 
\itl{$q$-L degree $Q_L(\varphi)$ of $\varphi$} as follows
(Proposition 3.3):
if $X$ is infinite and  $f$ is strictly $k$ right-mergible and  
strictly $l$ left-mergible, then
\[Q_R(\varphi)=n-k,\q\q\q\q Q_L(\varphi)=m-l\] 
and if $X$ is finite, then 
$Q_R(\varphi)=\infty$ and $Q_L(\varphi)=\infty$. 
 $P_L(\varphi)$, $P_R(\varphi)$, $Q_L(\varphi)$ and 
$Q_R(\varphi)$ are generically called the 
\itl{degrees of onesided resolvingness of $\varphi$}. 

Moreover, for any endomorphism $\varphi$ of any subshift 
$(X,\sigma)$, 
we can define $p_L(\varphi)$, $p_R(\varphi)$ and 
$q_R(\varphi),q_L(\varphi)$ by 
\begin{align*}
p_L(\varphi)&=\lim_{s\to\infty} P_L(\varphi^s)/s,\q\q\q\q
p_R(\varphi)=\lim_{s\to\infty} P_R(\varphi^s)/s,\\
q_R(\varphi)&=\lim_{s\to\infty} Q_R(\varphi^s)/s,\q\q\q\q
q_L(\varphi)=\lim_{s\to\infty} Q_L(\varphi^s)/s. 
\end{align*} 
which equal $\sup_s P_L(\varphi^s)/s$, $\sup_s P_R(\varphi^s)/s$, 
$\sup_s Q_R(\varphi^s)/s$, $\sup_s Q_L(\varphi^s)/s$, 
respectively (Theorem 9.2). We call $-p_L(\varphi)$ 
the \itl{limit of $p$-L directions of $\varphi$} 
or the \itl{$p$-L limit of $\varphi$} for short, 
and call
$-q_R(\varphi)$, $p_R(\varphi)$ and $q_L(\varphi)$ 
in the same way. 
They are generically called the 
\itl{limits of onesided resolving directions of $\varphi$}. 

For an endomorphism $\varphi$ of an infinite subshift 
$(X,\sigma)$, the following hold: 
if $\varphi^i(X)$ is infinite  for all $i\geq 1$ then 
$p_L(\varphi), p_R(\varphi)\in\R$, and 
otherwise, $p_L(\varphi)=\infty$ and $p_R(\varphi)=\infty$;
if $\varphi$ is right-closing, then $q_R(\varphi)\in\R$, 
and otherwise, $q_R(\varphi)=-\infty$, and  
if $\varphi$ is left-closing, then $q_L(\varphi)\in\R$, 
and otherwise, $q_L(\varphi)=-\infty$ 
(see the facts remarked after Definition 9.3). 
(For an endomorphism $\varphi$ of a finite subshift, 
$p_L(\varphi)=p_R(\varphi)=q_R(\varphi)=q_L(\varphi)=\infty$.)

As stated above the degrees of onesided resolvingness of an 
endomorphism of a subshift are not a topological invariant (i.e., 
an invariant of topological conjugacy 
between endomorphisms of dynamical systems). 
However if $\varphi$ and 
$\psi$ are topologically-conjugate endomorphisms of subshifts 
such that $P_L(\varphi)$ (respectively, $P_R(\varphi)$,
$Q_R(\varphi)$, $Q_L(\varphi)$) is an integer, then 
$\{P_L(\varphi^i)-P_L(\psi^i)\,|\,i\in\N\}$ 
(respectively, $\{P_R(\varphi^i)-P_R(\psi^i)\,|\,i\in\N\}$, 
$\{Q_R(\varphi^i)-Q_R(\psi^i)\,|\,i\in\N\}$, 
$\{Q_L(\varphi^i)-Q_L(\psi^i)\,|\,i\in\N\}$) is 
a finite set (Propositions 6.15, 5.4).

The limits $-p_L(\varphi)$, $-q_R(\varphi)$, $p_R(\varphi)$,
$q_L(\varphi)$ of onesided resolving directions are an 
topological invariant (Theorem 9.4). Moreover, 
for any endomorphisms 
$\varphi$ and $\psi$ of a subshift $(X,\sigma)$, the following hold: 
if $\varphi^i(X)$ and $\psi^i(X)$ are infinite 
for all $i\geq 0$, then $-p_L(\varphi\psi)=-p_L(\psi\varphi)$ and 
$p_R(\varphi\psi)=p_R(\psi\varphi)$; if
$\varphi$ and $\psi$ are right-closing 
then $-q_R(\varphi\psi)=-q_R(\psi\varphi)$, 
and if $\varphi$ and $\psi$ are left-closing 
then $q_L(\varphi\psi)=q_L(\psi\varphi)$ 
(Theorem 9.5). 

The limits 
$-p_L(\varphi)$, $-q_R(\varphi)$, $p_R(\varphi)$ 
and $q_L(\varphi)$ of onesided resolving directions 
as well as the degrees $P_L(\varphi)$, $Q_R(\varphi)$,  
$P_R(\varphi)$ and $Q_L(\varphi)$ of onesided resolvingness 
of an endomorphism $\varphi$ of a subshift 
have been introduced  in relation to 
the ``resolvingness'' (in the broad sense   
including onesided and weak ones) 
which means certain regular structures of 
the textiles (woven by a textile system \cite{Nasu-t}) 
in which the orbits of the endomorphism can be embedded. 
Though our main concern in this paper is onto endomorphisms of
subshifts, they are defined for any, not necessarily onto, 
endomorphism $\varphi$ of any subshift $(X,\sigma)$. 
It will turn out that $p_L(\varphi)$ 
and $p_R(\varphi)$ are closely related in definition to 
the right and left Lyapunov 
exponents for a cellular automaton $\varphi$ on $X$   
which were defined and treated from 
a very different viewpoint 
by M. A. Shereshevsky \cite{Sher}
and further studied by P. Tisseur \cite{Tiss} and, in particular, 
M. Hochman's treatment of them in \cite{Hoch} (Proposition 9.13). 

Here some remarks are in order. We follow 
the terminology of \cite{Nasu-te} and \cite{Nasu-t} 
on ``resolvingness'' for onto endomorphisms of subshifts 
and generalize it 
to not necessarily onto endomorphisms of subshifts further. 
(The definitions of the notions on ``resolvingness'' 
will be described in Subsection 2.3.) 
Each one of the ten terms ``$p$-L'', ``$p$-R'', 
``$q$-R'', ``$q$-L'', 
``LR'', ``RL'', ``LL'', ``RR'', ``$p$-biresolving'' 
and ``$q$-biresolving'', which will 
be called a \itl{resolving term},  
has in itself  a meaning only for onto endomorphisms of 
topological Markov shifts, 
and a corresponding notion for 
onto endomorphisms of SFTs (subshifts of finite type)  
can be given  as 
``\rt up to higher-block conjugacy 
between endomorphisms of subshifts''  
or more generally  ``essentially \st'', 
for each resolving term \st.  
Here and throughout the remainder of this paper, 
the term ``essentially''  means 
``up to topological conjugacy 
between endomorphisms of dynamical systems''. 
The term ``weakly \st'' has a meaning generally 
for all (not necessarily onto) endomorphisms 
of subshifts for every resolving term \st. 
In the remainder 
of this section, 
we shall outline our main results 
on weak resolving properties for 
onto endomorphisms $\varphi$ of
general infinite subshifts $(X,\sigma)$. 
However, all ``weakly \st'' appearing in them, 
where \rt is any resolving term, can be replaced by ``\st''
for the important special case 
(including onto cellular automata)   
that $\varphi$ is an onto endomorphism  
of a topological Markov shift $(X,\sigma)$ with 
$\varphi$ invertible 
or $\sigma$ topologically transitive.
We treat this case  
in particular in detail providing direct proofs for many results  
for the reader who is interested in the theory within the case. 
Some results particular to onto endomorphisms of 
topological Markov shifts and to those of full shifts 
are also given in this paper.  

Another main subject of this paper 
is ``expansiveness'' in the broad sense 
including certain types of onesided ones. Definitions 
for ``expansiveness''  
are found in Subsection 2.1.  The notions 
of ``left $\tau$-expansiveness'' 
and ``right $\tau$-expansiveness'' 
are defined for onto endomorphisms of an
invertible (compact) dynamical system $(X,\tau)$ 
\cite{Nasu-te}. (Equivalent notions 
for automorphisms of subshifts were defined with a directional 
dynamics treatment of them in \cite{Sab}). 
Also for onto endomorphisms, the notions of
``left $\tau$-expansiveness 
on the upper side'' and 
``right $\tau$-expansiveness 
on the upper side'' will be defined; 
they are special cases of 
``left $\tau$-expansiveness'' and 
``right $\tau$-expansiveness''. 
For not necessarily onto endomorphisms 
of an invertible dynamical system $(X,\tau)$, 
we shall define  
the notions of ``positively left $\tau$-expansiveness''
and ``positively right $\tau$-expansiveness'', 
which for onto endomorphisms are 
also special cases of ``left $\tau$-expansiveness'' 
and ``right $\tau$-expansiveness''. 
Equivalent notions to ``positively left $\tau$-expansiveness''
and ``positively right $\tau$-expansiveness''
for endomorphisms of subshifts with a directional dynamics 
treatment of them were defined and 
studied by Sablik \cite{Sab} and also treated 
by K\uo rka \cite{Kur2} (see the explanation given
after the proof of Theorem 11.8 for detail).  

Throughout the remainder of this section, 
we assume that $\varphi$ is an onto endomorphism of 
an infinite subshift $(X,\sigma)$, unless otherwise stated. 
However, many results hold without the assumption ``onto'' or 
with a weaker one, as seen in the referred results whose indexes 
are given in parentheses. 

We begin by stating  
the following two results, which are key results 
for understanding the relation between ``resolvingness'' and 
``expansiveness'': $\varphi$ is essentially weakly $p$-L and 
right $\sigma$-expansive (respectively, 
essentially weakly $p$-R and 
left $\sigma$-expansive) if and only if $\varphi$ is 
right $\sigma$-expansive on the upper side 
(respectively, left $\sigma$-expansive 
on the upper side) (Theorem 8.10); 
$\varphi$ is essentially weakly $q$-R and 
left $\sigma$-expansive (respectively, essentially weakly $q$-L and 
right $\sigma$-expansive)
 if and only if $\varphi$ is positively
left $\sigma$-expansive (respectively, 
positively right $\sigma$-expansive) (Theorem 8.9).

Let $s\in\Z$. 
We have the following basic equalities (Propositions 6.4 and 3.4):
\begin{align*}
P_L(\varphi\sigma^s)&=P_L(\varphi)+s,\q\q\q\q 
P_R(\varphi\sigma^s)=P_R(\varphi)-s,\\ 
Q_R(\varphi\sigma^s)&=Q_R(\varphi)+s,\q\q\q\q 
Q_L(\varphi\sigma^s)=Q_L(\varphi)-s.    
\end{align*} 
Hence $P_L(\varphi)+P_R(\varphi)$, $Q_R(\varphi)+Q_L(\varphi)$, 
$P_L(\varphi)+Q_L(\varphi)$ and $P_R(\varphi)+Q_R(\varphi)$ 
are shift-invariant. 
Except for $Q_R(\varphi)+Q_L(\varphi)$, these are 
nonpositive, and hence we have the following 
(Propositions 6.13 and 7.9):
if $\varphi$ is of $(m,n)$-type then
\[P_R(\varphi)\leq -P_L(\varphi), \q\q 
-n\leq P_R(\varphi)\leq -Q_R(\varphi),
\q\q Q_L(\varphi)\leq -P_L(\varphi)\leq m.\] 

We have the following results (Proposition 6.12, Theorem 5.2): 
$\varphi\sigma^s$ is weakly $p$-L if and only if $s\geq -P_L(\varphi)$, 
and if $s> -P_L(\varphi)$ then $\varphi\sigma^s$ is 
right $\sigma$-expansive; 
$\varphi\sigma^s$ is weakly $p$-R if and only if $s\leq P_R(\varphi)$,  
and if $s< P_R(\varphi)$ then $\varphi\sigma^s$ is 
left $\sigma$-expansive;
$\varphi\sigma^s$ is weakly $q$-R if and only if $s\geq -Q_R(\varphi)$, 
and if $s>-Q_R(\varphi)$ then $\varphi\sigma^s$ is left $\sigma$-expansive;
$\varphi\sigma^s$ is weakly $q$-L if and only if $s\leq Q_L(\varphi)$, 
and if $s<Q_L(\varphi)$ then $\varphi\sigma^s$ is right $\sigma$-expansive.

Let $C_R(\varphi)= \max\{-P_L(\varphi), -Q_R(\varphi)\}$ 
and $C_L(\varphi)=\min\{P_R(\varphi), Q_L(\varphi)\}$.
Then we have the following results (Theorem 7.8): 
$\varphi\sigma^s$ is weakly LR 
up to higher block conjugacy 
if and only if 
$s\geq C_R(\varphi)$, and if 
$s>C_R(\varphi)$ then $\varphi\sigma^s$ is expansive; 
$\varphi\sigma^s$ is weakly RL up to  
higher block conjugacy 
if and only if $s\leq C_L(\varphi)$, and if $s<C_L(\varphi)$ 
then $\varphi\sigma^s$ is expansive. 

$Q_R(\varphi)+Q_L(\varphi)$ 
can be negative, zero and positive.  
If $\varphi$ is invertible, then 
$Q_R(\varphi)+Q_L(\varphi)\leq 0$ (Proposition 7.11(2)), 
but the converse does not hold (Example 9.19). 
We have the following result 
(Theorem 5.3):  
$\varphi\sigma^s$ is weakly $q$-biresolving 
up to higher-block conjugacy if and only if
$-Q_R(\varphi)\leq s\leq Q_L(\varphi)$; if 
$-Q_R(\varphi)< s< Q_L(\varphi)$ 
then $\varphi\sigma^s$ is positively expansive. 

Let $i\in\N$ and $j\in\Z$. Then we have the following equations 
(Proposition 9.6):
\begin{align*}
p_L(\varphi^i\sigma^j)&=ip_L(\varphi)+j,\q\q\q\q 
p_R(\varphi^i\sigma^j)=ip_R(\varphi)-j,\\
q_R(\varphi^i\sigma^j)&=iq_R(\varphi)+j,\q\q\q\q 
q_L(\varphi^i\sigma^j)=iq_L(\varphi)-j. 
\end{align*} 
Hence
$p_L(\varphi)+p_R(\varphi)$, $q_R(\varphi)+q_L(\varphi)$, 
$p_L(\varphi)+q_L(\varphi)$ and $p_R(\varphi)+q_R(\varphi)$ 
are shift-invariant. Except for $q_R(\varphi)+q_L(\varphi)$, 
these are nonpositive, and hence 
we have the following (Proposition 9.9):
if $\varphi$ is 
of $(m,n)$-type then 
\[p_R(\varphi)\leq -p_L(\varphi),\q\q 
-n\leq p_R(\varphi)\leq -q_R(\varphi),\q\q 
q_L(\varphi)\leq -p_L(\varphi)\leq m.\] 
We have the following results (Theorem 9.7): 
$\varphi^i\sigma^j$ is essentially weakly $p$-L 
and right $\sigma$-expansive if and only if $j/i> -p_L(\varphi)$;
$\varphi^i\sigma^j$ is essentially weakly $q$-R 
and left $\sigma$-expansive if and only if $j/i> -q_R(\varphi)$;
$\varphi^i\sigma^j$ is essentially weakly $p$-R 
and left $\sigma$-expansive if and only if $j/i< p_R(\varphi)$;
$\varphi^i\sigma^j$ is essentially weakly $q$-L 
and right $\sigma$-expansive if and only if $j/i< q_L(\varphi)$. 

Let $c_R(\varphi)=\max\{-p_L(\varphi), -q_R(\varphi)\}$ 
and let 
$c_L(\varphi)=\min\{p_R(\varphi), q_L(\varphi)\}$.  
We also have the results (Theorem 9.10) that 
$\varphi^i\sigma^j$ is essentially weakly LR  and expansive 
if and only if $j/i> c_R(\varphi)$ and that 
$\varphi^i\sigma^j$ is essentially weakly RL and expansive 
if and only if $j/i<c_L(\varphi)$.  

We find that $q_R(\varphi)+q_L(\varphi)$ 
can be negative, zero and positive. 
If $\varphi$ is an  automorphism of $(X,\sigma)$, 
then it is nonpositive, 
but the converse does not hold (Example 9.19). 
We obtain the following result (Theorem 9.11): 
$\varphi^i\sigma^j$ is positively expansive if and only if 
$-q_R(\varphi)<j/i<q_L(\varphi)$. 

In Subsection 9.2, we discuss types of the limits 
$-p_L(\varphi)$, $-q_R(\varphi)$, $p_R(\varphi)$ and $q_L(\varphi)$ 
in $\R$. 
For example, $-p_L(\varphi)$ is said to be of \itl{type I} 
if there exist $i\in\N,j\in\Z$  
such that $j/i=-p_L(\varphi)$ and 
$\varphi^i\sigma^j$ is an essentially weakly $p$-L endomorphism 
of $(X,\sigma)$; the definition of being of \itl{type I} 
are given similarly for each of the limits $-q_R(\varphi)$, 
$p_R(\varphi)$ 
and $q_L(\varphi)$. 
An example due to Lind and Boyle \cite{BoyLin} shows 
that there exists an automorphism having \itl{type III}
(i.e. irrational) limits and 
an example due to Hochman \cite{Hoch} shows 
that there exists an automorphism having \itl{type II} 
(i.e. rational but not type I) limits (Example 9.14). However 
we cannot answer the question whether or not 
an onto endomorphism of an SFT 
can have only type I limits of onesided resolving directions. 

We know no general method for calculating 
the limits of onesided resolving directions 
even for onto endomorphisms of SFTs. 
For onto endomorphisms of SFTs in some specific cases, 
we can calculate type I limits (using Propositions 9.15, 9.16). 
Examples of the calculation  
are found in Examples 9.17--9.20. 

Section 10 treats onesided resolving endomorphisms from 
a much more general viewpoint than that of the overall dynamics 
of a single endomorphism of a subshift. Let $X$ be a zero-dimensional 
compact metric space. Let $S(X)$ 
denote the monoid of all surjective continuous maps of $X$ 
onto itself and 
$H(X)$ the group of all homeomorphisms of $X$ onto itself. 
Let $J$ be a submonoid of $S(X)$ and let $\tau\in H(X)\cap J$ 
be expansive with $\tau^{-1}\in J$. We call the set $C_J(\tau)$ of 
all essentially weakly LR endomorphisms of $(X,\tau)$ in $J$
the \itl{essentially-weakly-LR cone}, or \itl{cone, of $\tau$ 
in $J$}. 

If $(X,\tau)$ is conjugate to an SFT, then the set of all 
essentially LR endomorphisms of $(X,\tau)$ in $J$ will be 
called the \itl{essentially-LR cone}, 
or \itl{ELR cone}, \itl{of $\tau$ in $J$} and 
denoted by $(C_0)_J(\tau)$. If in addition,  
$\tau$ is topologically transitive 
then $C_J(\tau)=(C_0)_J(\tau)$, 
and if $J=K$ for a subgroup $K$ of $H(X)$
then $C_K(\tau)=(C_0)_K(\tau)$. 
We emphasize that when $(X,\tau)$ is conjugate to 
an SFT, the dynamics of the elements of $(C_0)_J(\tau)$ 
of $\tau$ in $J$ is 
lucid  and regular, even when $J=S(X)$ (Theorem 10.1). 

Let $\varphi,\tau\in H(X)$. We write 
$\varphi\pl\tau$ if 
$\tau$ is expansive and $\varphi$ is 
an essentially weakly $p$-L automorphism of $(X,\tau)$, and
if in addition, $\varphi$ is right $\tau$-expansive, then 
we write $\varphi\pl^\circ\tau$, which is 
equivalent to the condition that $\tau$ is expansive and 
$\varphi$ is left $\tau$-expansive on 
the upper side (by Theorem 8.10(1)); 
we write  
$\varphi\qr\tau$ if 
$\tau$ is expansive and $\varphi$ is 
an essentially weakly $q$-R automorphism of $(X,\tau)$, and
if in addition, $\varphi$ is left $\tau$-expansive, then 
we write $\varphi\qr^\circ\tau$, which is equivalent to 
the condition that $\tau$ is expansive and 
 $\varphi$ is positively left $\tau$-expansive 
(by Theorem 8.9(2)). 
We write $\varphi\lr\tau$ to mean that 
$\tau$ is expansive and $\varphi$ is 
an essentially weakly LR 
automorphism of $(X,\tau)$, and 
$\varphi\lr^\circ\tau$ to mean that $\varphi\lr\tau$ 
with $\varphi$ expansive, 
which is equivalent to the 
condition that 
$\varphi\pl^\circ\tau$ and $\varphi\qr^\circ\tau$ 
(by \cite[Proposition 8.1]{Nasu-te})(see Theorem 8.12(3) 
for other equivalent 
conditions to this).
Let $K$ be any subgroup of $H(X)$ and let $E(K)$ denote the 
set of all expansive elements in $K$. If $\tau\in E(K)$, 
then $C_K(\tau)=\{\varphi\in K\,|\,\varphi\lr\tau\}$. We define 
$C_K^\circ(\tau)=
\{\varphi\in K\,|\,\varphi\lr^\circ\tau\}=
\{\varphi\in E(K)\,|\,\varphi\lr\tau\}$, 
which is called the \itl{interior} of $C_K(\tau)$. 
   
If $K$ is commutative, then 
the relation $\lr^{\circ}$ 
is an equivalence relation on $E(K)$, 
$C_K^\circ(\tau)$ is the equivalence class containing $\tau$ 
with respect to $\lr^\circ$ and $C_K(\tau)=C_K(\tau')$ 
for every $\tau'\in C^\circ_K(\tau)$ (Proposition 10.3); 
hence 
we call $C_K(\tau)$ \itl{an essentially-weakly-LR cone in $K$} 
or \itl{a cone in $K$}. 
We know that
if a cone $C$ in a commutative subgroup $K$ of $H(X)$ 
contains an element with POTP (the pseudo orbit tracing property) 
in $C^\circ$ (i.e., $C$ is an ELR cone in $K$), then 
all elements in $C$ have POTP (see Proposition 10.2). 

For $\varphi,\tau\in H(X)$, 
write $\varphi\mm\tau$ to mean 
that $\varphi\pl\tau$ or $\varphi\qr\tau$, write 
$\varphi\mm^\circ\tau$ to mean 
that $\varphi\pl^\circ\tau$ or $\varphi\qr^\circ\tau$.
For any (not necessarily commutative) subgroup $K$ of $H(X)$ 
and $\varphi,\tau\in K$
we write $\varphi\mm^\ast_K\tau$ to mean 
that there exist $r\geq 1$ and elements 
$\tau_i, i=1,\dots, r$, in $E(K)$ such that 
$\varphi\mm\tau_1 \mm\cdots\mm\tau_r=\tau$ 
(or $\varphi\mm\tau_1 \mm^\circ\cdots\mm^\circ\tau_r=\tau$), and write 
$\varphi\mm^{\circ,\ast}_K\tau$ to mean that 
there exist $r\geq 1$ and elements 
$\tau_i, i=1,\dots, r$, in $E(K)$ such that 
$\varphi\mm^\circ\tau_1 \mm^\circ\cdots\mm^\circ\tau_r=\tau$.
For any (not necessarily commutative) subgroup 
$K$ of $H(X)$, 
the relation $\mm^\ast_K$ ( $\mm^{\circ,\ast}_K$) restricted on $E(K)$ 
is an equivalence relation on $E(K)$ (Proposition 10.4). 
Define 
\[D_K^\ast(\tau)
=\{\varphi\in K\,|\,\varphi\mm^\ast_K\tau\},\q\q
(D_K^\ast)^\circ(\tau)
=\{\varphi\in K\,|\,\varphi\mm^{\circ,\ast}_K\tau\}.\]
Then $D_K^\ast(\tau)$ is called 
the \itl{extended district of $\tau$ in $K$}
or an \itl{extended district in $K$} 
(for the reason that for every $\tau'$ in $D_K^\ast(\tau)\cap E(K)$, 
$D^\ast_K(\tau)=D^\ast_K(\tau')$ and  
$(D^\ast_K)^\circ(\tau)=(D^\ast_K)^\circ(\tau')$), and
$(D_K^\ast)^\circ(\tau)$ is called the \itl{interior} of 
an extended district $D_K^\ast(\tau)$. 

Let $D^\ast$ be  
any extended district in $K=H(X)$. 
If there exists $\tau\in D^\ast$ with $(X,\tau)$
conjugate to a mixing SFT, then 
for every $\varphi$ in $D^\ast\cap E(K)$, 
$(X,\varphi)$ is conjugate to 
a mixing SFT (Theorem 10.5). 

In Section 11,  we apply geometric methods of 
Boyle and Lind \cite{BoyLin} defining expansive lines 
and expansive components for $\Z^d$-actions to 
``expansiveness'' of  
endomorphisms and automorphisms of subshifts. 
First, in Subsection 11.1 we make a general treatment. 
For an onto endomorphism $\varphi$ of 
an invertible (compact) dynamical system $(X,\tau)$, 
we define the notions of left $\tau$-expansive 
and right $\tau$-expansive (non-horizontal) lines 
(in the plane $\R^2$) and directions (in $\R$) for $\varphi$, 
and those of left $\tau$-expansive
and right $\tau$-expansive lines 
and directions on the upper side for $\varphi$, and 
for a not necessarily onto endomorphism $\varphi$ of $(X,\tau)$, 
the notions of  positively left $\tau$-expansive and 
positively right $\tau$-expansive lines and directions 
for $\varphi$. 
For an onto endomorphism $\varphi$ of 
an invertible dynamical system $(X,\tau)$, 
let $E_L(\varphi)$ and $E_R(\varphi)$ denote the sets of 
all left $\tau$-expansive directions and 
the set of all right $\tau$-expansive directions, 
respectively, for $\varphi$ in $\R$ 
(hence, if $r$ is a rational number and $r=j/i$ 
with $i\in\N,j\in\Z$, then $r\in E_L(\varphi)$ 
(respectively, $r\in E_R(\varphi)$) 
if and only if $\varphi^i\tau^j$ is 
left $\tau$-expansive 
(respectively, right $\tau$-expansive)). 
Each of 
$E_L(\varphi)$ and $E_R(\varphi)$
is an open subset of $\R$ (Proposition 11.2). 
Define $E(\varphi)=E_L(\varphi)\cap E_R(\varphi)$. 
Then $E(\varphi)$ is the set of 
all expansive directions for $\varphi$ 
(Proposition 11.1). 
The main result of 
Subsection 11.1 is given as follows (Proposition 11.6): 
for an onto endomorphism $\varphi$ 
of invertible, expansive dynamical system 
$(X,\tau)$,  
the set of right $\tau$-expansive 
(respectively, left $\tau$-expansive)
directions on the upper side for $\varphi$ is a right-unbounded 
(respectively, left-unbounded) open interval
and is a (connected) component of $E_R(\varphi)$ 
(respectively, $E_L(\varphi)$), and so is
the set of positively left $\tau$-expansive 
(respectively, positively right $\tau$-expansive)
directions for $\varphi$. 

In Subsection 11.2, we combine the results in Subsections 9.1 
and 11.1 to obtain the following results (Theorem 11.8) 
which are described with the convention that  
for $\alpha,\beta\in\R\cup\{\infty,-\infty\}$, the interval
$(\alpha,\beta)$ with $\alpha\geq\beta$ means the empty set. 

Let $\varphi$ be an endomorphism of a subshift $(X,\sigma)$. 
If $\varphi$ is onto, then the following 
(1), (2), (3) and (4) hold.  
\begin{enumerate}
\item 
The interval $(-p_L(\varphi),\infty)$  
(respectively, $(-\infty,p_R(\varphi))$) is the set of all
right $\sigma$-expansive (respectively, left $\sigma$-expansive) 
directions on the upper side for $\varphi$ and the right-unbounded  
(connected) component of $E_R(\varphi)$ (respectively, 
the left-unbounded component of $E_L(\varphi)$).  
\item 
The interval $(-q_R(\varphi),\infty)$  
(respectively, $(-\infty,q_L(\varphi))$) is 
the set of all positively
left $\sigma$-expansive 
(respectively, positively right $\sigma$-expansive) 
directions for $\varphi$, and it is the right-unbounded 
component of $E_L(\varphi)$ (respectively, 
the left-unbounded component of $E_R(\varphi)$) if it is nonempty. 
\item Particularly, 
neither $-p_L(\varphi)$ nor $q_L(\varphi)$ is a 
right $\sigma$-expansive direction for $\varphi$, and neither
$p_R(\varphi)$ nor $-q_R(\varphi)$ is a
left $\sigma$-expansive direction for $\varphi$. 
\item 
The interval  
$(c_R(\varphi),\infty)$ (respectively, 
$(-\infty, c_L(\varphi))$) is the right-unbounded 
(respectively, left-unbounded) component of $E(\varphi)$ 
if it is nonempty. Further, the interval 
$(-q_R(\varphi),q_L(\varphi))$ is a 
component of $E(\varphi)$ if it is nonempty. 
\end{enumerate} 
If $\varphi$ is not necessarily onto, then the following (5) holds.
\begin{enumerate}
\item[(5)]
The interval $(-q_R(\varphi),\infty)$  
(respectively, $(-\infty,q_L(\varphi))$) is 
the set of all positively
left $\sigma$-expansive 
(respectively, positively right $\sigma$-expansive) 
directions for $\varphi$, and the interval 
$(-q_R(\varphi),q_L(\varphi))$ is the set of all positively 
expansive directions for $\varphi$. 
\end{enumerate} 

This result (5) refines a result 
of Sablik \cite[Theorem 5.2]{Sab} 
and that of K\uo rka \cite[Theorem 3.3]{Kur2} 
(see the explanation given
after the proof of Theorem 11.8 for detail). 

Suppose that $\varphi$ is an automorphism of $(X,\sigma)$. 
Let $K$ be the subgroup of $H(X)$ generated by $\{\sigma,\varphi\}$. 
We define $\E_L(\varphi)$ 
(respectively, $\E_R(\varphi)$, $\E(\varphi)$)
to be the set-union of all left $\sigma$-expansive 
(respectively, right $\sigma$-expansive, expansive) lines 
passing through $(0,0)$ (including the horizontal one)
for $\varphi$, with $\{(0,0)\}$ subtracted. 
Then $\E_L(\varphi)$ (respectively, $\E_R(\varphi)$, $\E(\varphi)$) 
is the open subset of $\R^2$ 
which is the disjoint union of 
open cones with apex $(0,0)$ in $\R^2$ 
and the set of lattice points in which 
equals the set of all $(i,j)\in \Z^2$ such that 
$\varphi^i\sigma^j$ is a left $\sigma$-expansive 
(respectively, right $\sigma$-expansive, 
expansive) homeomorphism. 
Furthermore, for any $(k,l)\in \E(\varphi)\cap\Z^2$, 
$(\E_L(\varphi)\cup\E_R(\varphi))\cap\Z^2$ is 
the set of all lattice points $(i,j)$ such that 
$\varphi^i\sigma^j$ is a \itl{onesided $\varphi^k\sigma^l$-expansive} 
(i.e, left $\varphi^k\sigma^l$-expansive or 
right $\varphi^k\sigma^l$-expansive) 
homeomorphism (Proposition 11.10), and hence 
if $(i,j)\in \E_L(\varphi)\cup\E_R(\varphi)$, then
we can call $\varphi^i\sigma^j$  
a \itl{onesided-expansive automorphism of $(X,\sigma)$}. 
We prove the result (Theorem 11.12) that 
for any $(k,l)\in \E(\varphi)\cap\Z^2$ it 
holds that if $\C^L$, $\C^R$ and $\C$ are the components of 
$\E_L(\varphi)$, $\E_R(\varphi)$ and $\E(\varphi)$, respectively, 
each containing $(k,l)$, then 
\begin{align*}
\C\cap\Z^2=
(\C^L\cap\C^R)\cap\Z^2
=&\{(i,j)\in\Z^2\,|\,
\varphi^i\sigma^j\lr^\circ\varphi^k\sigma^l\},\\
(\C^L\cup\C^R)\cap\Z^2
=&\{(i,j)\in\Z^2\,|\,
\varphi^i\sigma^j\mm^\circ\varphi^k\sigma^l\}.
\end{align*}
From these we obtain the following (1) and (2) (Theorem 11.13).  

(1) A subset of $\Z^2$ is given as 
$\C\cap\Z^2$, where $\C$ is 
an \itl{expansive component for $\varphi$}, 
i.e. a component of   
$\E(\varphi)=\E_R(\varphi)\cap\E_L(\varphi)$,  
if and only if it is given as 
$\{(i,j)\in\Z^2 \,|\, \varphi^i\sigma^j\in C^\circ\}$, 
where $C^\circ$ is the interior of  
an essentially-weakly-LR cone $C$ in $K$. 

(2) A subset of $\Z^2$ is given as
$\cD\cap\Z^2$, where $\cD$ is 
a \itl{onesided-expansive component for $\varphi$}, i.e. 
 a component  
of $\E_R(\varphi)\cup\E_L(\varphi)$,  with 
$\cD\cap\E(\varphi)\neq\emptyset$,  
if and only if it is given as 
$\{(i,j)\in\Z^2 \,|\, \varphi^i\sigma^j\in (D^\ast)^\circ\}$, 
where $(D^\ast)^\circ$ is the interior of an extended district 
$D^\ast$ in $K$ (Theorem 11.13). 

Since an expansive component for the automorphism $\varphi$ 
is the same as an ``expansive component of 1-frames'' 
for the $\Z^2$-action $(i,j)\mapsto \varphi^i\sigma^j$
in the sense of Boyle and Lind \cite{BoyLin}, 
the statement (1) above 
is the same as \cite[Remark 9.6]{Nasu-te} 
(in view of Theorem 
12.2). A complete, self-contained proof for this  
is provided (by Theorem 11.12),
and moreover it is generalized 
(using Theorem 11.12) to 
$\Z^d$-actions on infinite zero-dimensional 
compact metric spaces (Theorem 11.18). 
Boyle and Lind \cite{BoyLin} proved that 
any ``expansive component  
of 1-frames'' $\C$ for a $\Z^d$-action $\alpha$
on an infinite compact metric space $X$ is 
an open cone in $\R^d$ such 
that $\C\cap(-\C)=\emptyset$ and that 
if one vector in $\C$ is ``Markov'' 
then so is every vector in $\C$. 
They said that a component 
containing  a Markov vector is ``Markov''. 
By the above 
we can understand that for a 
$\Z^d$-action $\alpha$ on 
an infinite zero-dimensional compact 
metric space $X$,  
an expansive component of 1-frames for $\alpha$ 
is the ``same'' 
as the interior of an essentially-weakly-LR cone 
in the subgroup $K=\{\alpha^\bfv\,|\,\bfv\in \Z^d\}$ of $H(X)$ 
(up to the natural correspondence between them as in (1)). 
We can also understand that a Markov component 
for a $\Z^d$-action in the zero-dimensional case  
is the ``same'' as the interior of 
an ELR (essentially-LR) cone in $K$.  
This is presented for $d=2$ 
as a result of automorphisms of 
subshifts (Corollary 11.15) and then 
given generally as that of  
$\Z^d$-actions (Theorem 11.18). 
The dynamics of elements of ELR cones is lucid and 
regular (see Theorem 10.1 together with Proposition 10.2). 

For an automorphism $\varphi$ of a subshift $(X,\sigma)$ 
(or a $\Z^2$-action $(i,j)\mapsto \varphi^i\sigma^j$), 
by the above (2) we understand that 
a onesided-expansive component 
containing an expansive element 
for the automorphism $\varphi$ 
is the ``same'' as the interior of an extended district 
in the subgroup $K$ of $H(X)$ generated 
by $\{\sigma,\varphi\}$. Furthermore we see the following.
If $\cD$ is a component of 
$\E_L(\varphi)\cup\E_R(\varphi)$, 
then $\cD$ is either an open cone 
in $\R^2$ with $\cD\cap(-\cD)=\emptyset$ or 
equals $\R^2/\{(0,0)\}$. 
There exists an automorphism $\varphi$ of a full-shift such that 
$\R^2/\{(0,0)\}$ is a unique onesided-expansive 
component for $\varphi$ (Example 11.17). 
If a component $\cD$ of $\E_L(\varphi)\cup\E_R(\varphi)$ 
contains a lattice point $(m,n)$ with 
$(X,\varphi^m\sigma^n)$ conjugate to a mixing SFT, 
then for every $(k,l)\in\cD\cap\Z^2$ with $\varphi^k\sigma^l$ 
expansive, $(X,\varphi^k\sigma^l)$ is conjugate to a 
mixing SFT, and hence every expansive component for 
$\varphi$ included 
in $\cD$ is the interior of an ELR cone in $K$ (Corollary 11.15). 

In Subsection 11.4, we discuss about the relation 
between  ``resolvingness'' and  ``expansiveness''
for $\Z^d$-actions on a zero-dimensional 
compact metric space along the framework provided by 
Boyle and Lind \cite{BoyLin}. We prove not only the result 
stated above but also the following (Theorem 11.24) : for a 
$\Z^d$-action $\alpha$ on a infinite zero-dimensional 
compact metric space and an expansive integral vector $\bfk\in\R^d$ 
(i.e., $\alpha^\bfk$ is expansive), 
an integral vector
$\bfm$ belongs to the component of $\E_R(\alpha,\bfk)$ 
(respectively, $\E_L(\alpha,\bfk)$)
containing $\bfk$ if and only if 
$\alpha^\bfm\pl^\circ\alpha^\bfk$ 
(respectively, $\alpha^\bfm\qr^\circ\alpha^\bfk$ , 
Here $\E_R(\alpha,\bfk)$ (respectively, $\E_L(\alpha,\bfk)$)
is the set-union of all ``right $\bfk$-expansive'' 
(respectively, ``left
$\bfk$-expansive'') one-dimensional subspaces of $\R^d$, 
with $\{\bfzero\}$ subtracted, each component $\C$ of which is proved to be 
an open cone in $\R^d$ with $\C\cap(-\C)=\emptyset$. 

Mike Boyle \cite{Boy-p} proved that  
if $\varphi^s$ is weakly $p$-L (respectively, weakly $p$-R) 
with $s\in\N$, 
then the ``dual higher block 
presentation of order $s$''  $\varphi^{[^*s]}$ is weakly $p$-L 
(respectively, weakly $p$-R) 
and hence $\varphi$ is essentially weakly $p$-L (respectively, 
essentially weakly $p$-R), 
together with more for the case that $\varphi$ is an automorphism 
(see Theorem 8.1), and furthermore he suggested the possibility
of proving the main results of \cite [Section 7]{Nasu-t}
(on resolving endomorphisms of topological Markov shifts) 
without using the long theory of ``resolvable textile systems'' 
presented there. In this paper, using $q$-R and $q$-L degrees 
we prove that if $\varphi^s$ is weakly $q$-R 
(respectively, weakly $q$-L) with $s\in\N$, then 
$\varphi^{[^*s]}$ is weakly $q$-R (respectively, weakly $q$-L) 
(Theorem 8.6). Boyle's results and these are essential fundamental 
results in this paper; in fact, using them 
we further obtain the result that   
for every resolving term \st, 
it holds that if $\varphi^s$ is weakly \rt for some $s\in\N$, 
then $\varphi$ is essentially weakly \rt 
(Theorems 8.7 together with Theorems 8.1 and 8.6) 
and the result that for every resolving term \st, 
the condition that 
$\varphi^i\sigma^j$ is essentially weakly \rt depends only on 
the direction $j/i$ (Corollary 8.8). 

In Section 12, using another method we prove that if $\varphi$ 
is an essentially weakly \rt endomorphism of $(X,\sigma^s)$ 
for some $s\in\N$, then $\varphi$ is 
essentially weakly \st, where \rt is 
any resolving term. 
Therefore we obtain the result that for every  
resolving term \st, it holds that if $\varphi$ is 
\itl{directionally} essentially weakly \rt 
(i.e., $\varphi^r$ is an essentially weakly \rt 
endomorphism of $(X,\sigma^s)$ 
for some $r, s\in\N$), then $\varphi$ is essentially weakly \rt 
(Theorem 12.2). 
This is a generalization of 
the above-mentioned important special case of 
the main results of \cite [Section 7]{Nasu-t}, 
to onto endomorphisms of general subshifts, and more.  
(However, the theory of resolvable textile systems 
of \cite [Section 7]{Nasu-t} has some significant results 
which we cannot cover in this paper.) 

We prove that if $(X,\sigma)$ is a topological Markov shift with 
defining matrix $M$ and $\varphi$ is 
its endomorphism which is LR up to higher block conjugacy, 
then there exists a nonnegative integral matrix $N$ 
with $MN=NM$ such that for all $i\geq 0, j\geq 1$, 
the inverse limit system of $\varphi^i\sigma^j$
is topologically conjugate to the 
topological Markov shift 
with defining matrix $N^iM^j$ and such that 
$(\tilde{X},\tilde{\varphi}^i\tilde{\sigma}^j)$ is 
topologically conjugate to  the 
onesided topological Markov shift 
with defining matrix $N^iM^j$,  
where $(\tilde{X},\tilde{\sigma})$ is 
the induced onesided subshift of 
$(X,\sigma)$ and $\tilde{\varphi}$ is 
its endomorphism induced by $\varphi$ 
(Theorem 7.6). (Note that this is a new 
result only because it contains 
``up to higher block conjugacy''; see 
\cite[Corollaries 6.5, 6.7]{Nasu-t}.) 
We also give a construction method for obtaining 
positively expansive endomorphisms of any full-shift from  
right-closing and left-closing endomorphisms of the full-shift 
(Theorem 4.8).   

\textbf{Standing convention I for the description of proofs.}  
For each of the theorems, propositions, lemmas, and 
remarks in this paper, if it has the second version  
together with the first, 
we describe a proof only 
for the first version, because the proof of the second one 
is analogous, or 
the second one follows from the first one \itl{by symmetry} 
(i.e. by reversing the direction of the shift).


\section{Preliminaries} 

In this section, we give preliminaries to the subsequent sections. 
For more information, 
the reader is referred to \cite{Kit-s} or \cite{LinMar} on 
symbolic dynamics, and \cite{AokHir} on topological dynamics. 

\subsection{Some basic definitions} 

A \itl{dynamical system}  means an ordered pair $(X,\tau)$ of a 
compact metric space $X$ and a continuous map $\tau:X\to X$. 
When $\tau$ is a homeomorphism, we say that 
the dynamical system $(X,\tau)$ is \itl{invertible}.

A \itl{commuting system} $(X,\tau,\varphi)$ means an ordered 
pair of commuting continuous maps $\tau:X\to X$ and
$\varphi:X\to X$  of a compact metric space.
If $(X,\tau,\varphi)$ is a commuting system, then 
$\varphi$ is an \itl{endomorphism} of the dynamical system 
$(X,\tau)$. If $\varphi$ is a homeomorphism in this definition, 
then $\varphi$ is an \itl{automorphism} of $(X,\tau)$. 

Two commuting systems $(X,\tau,\varphi)$ and 
$(X',\tau',\varphi')$ are said to 
be \itl{topologically conjugate} 
if there exists 
a \itl{a topological conjugacy} 
$\theta:(X,\tau,\varphi)\to (X',\tau',\varphi')$, 
i.e., 
a homeomorphism $\theta:X\to X'$ which gives topological 
conjugacies $\theta:(X,\tau)\to (X,\tau')$ and
$\theta:(X,\varphi)\to (X,\varphi')$ between 
the dynamical systems 
at the same time. We call 
the topological conjugacy $\theta:(X,\tau,\varphi)\to (X',\tau',\varphi')$
a \itl{topological conjugacy of the endomorphism $\varphi$ of $(X,\tau)$
to the endomorphism $\varphi'$ of $(X',\tau')$.} 

For an onto continuous map $\varphi: X\to X$ 
of a compact metric space, 
a sequence $\xseq$ of points in $X$ is called a \itl{$\varphi$-orbit} 
if $\varphi(x_i)=x_{i+1}$ for all $i\in\Z$. We say that
$\varphi$ is \itl{expansive}, if there exists 
an \itl{expansive constant for $\varphi$}, i.e. 
$\delta>0$ such that for any $\varphi$-orbits $\xseq,\yseq$, 
it holds that if 
$d_X(x_i,y_i)\leq\delta$ for all $i\in\Z$ then $\xseq=\yseq$. 
(This notion appears 
in \cite[p.57]{AokHir} with the term ``c-expansive''.) 

Let $X$ be a compact metric space, 
$\varphi:X\to X$ an onto continuous map and 
$\tau:X\to X$ a homeomorphism.
When $\varphi\tau=\tau\varphi$ (or $\varphi$ is 
an endomorphism of $(X,\tau)$),
we say that $\varphi$ is \itl{left $\tau$-expansive} 
(respectively, \itl{right $\tau$-expansive})
if there is $\delta>0$ such that 
for any $\varphi$-orbits $(x_i)_{i\in\Z}, (y_i)_{i\in\Z}$
it holds that if $d_X(\tau^j(x_i),\tau^j(y_i))\leq\delta$ 
for all $i\in\Z$ and $j\leq 0$ (respectively, $j\geq 0$) 
then $\xseq=\yseq$ (\cite[Section 6]{Nasu-te}). 

As will be shown by Proposition 11.1, 
$\varphi$ is left $\tau$-expansive 
and right $\tau$-expansive if and only if $\varphi$ is 
expansive. (It is natural to define that $\varphi$ is 
\itl{$\tau$-expansive} if there exist $\delta>0$ and $t\geq 0$ 
such that if for any $\varphi$-orbits $\xseq, \yseq$ it holds that 
if $d_X(\tau^j(x_i),\tau^j(y_i))\leq\delta$ for all 
$i\in\Z, -t\leq j\leq t$ then $\xseq=\yseq$. However $\varphi$ 
is $\tau$-expansive if and only if it is expansive.) 
We shall say that $\varphi$ is \itl{onesided $\tau$-expansive} if 
$\varphi$ is right $\tau$-expansive or left $\tau$-expansive. 

We note that a homeomorphism $\tau: X\to X$ is 
expansive if and only if $\tau$ 
is left $\tau$-expansive (respectively, right $\tau$-expansive). 
When $X$ is infinite, the identity map $i_X$ 
is not onesided $\tau$-expansive (i.e. neither left $\tau$-expansive 
nor right $\tau$-expansive). This follows from 
Schwartsman's theorem cited as \cite[Theorem 3.9]{BoyLin}. 

An automorphism $\varphi$ of an invertible dynamical system 
$(X,\tau)$ is said to be \itl{onesided-expansive},  
if $\varphi$ is onesided $\varphi^k\tau^l$-expansive 
for some $(k,l)\in\Z^2$ such 
that $\varphi^k\tau^l$ is expansive. 
If an automorphism $\varphi$ of an invertible dynamical system 
$(X,\tau)$ is onesided-expansive, then $\varphi$ is onesided 
$\varphi^k\tau^l$-expansive for all $(k,l)\in \Z^2$ such that 
$\varphi^k\tau^l$ is expansive. (This will be proved after 
the proof of Proposition 11.10.) 

An onto endomorphism $\varphi$ of an invertible dynamical 
system $(X,\tau)$ is said to be 
\itl{left $\tau$-expansive on the upper side} (respectively, 
\itl{right $\tau$-expansive on the upper side} if there 
exists $\delta>0$ such that 
for any $\varphi$-orbits $(x_i)_{i\in\Z},(y_i)_{i\in\Z}$ 
it holds that if $d_X(\tau^j(x_i),\tau^j(y_i))\leq \delta$
for all $i\leq 0$ and $j\leq 0$ (respectively, $j\geq 0$) then 
$(x_i)_{i\in\Z}=(y_i)_{i\in\Z}$. 

It is clear that if an onto endomorphism $\varphi$ 
of an invertible dynamical system $(X,\tau)$ is 
left $\tau$-expansive on the upper side (respectively, 
right $\tau$-expansive on the upper side), then it is 
left $\tau$-expansive (respectively, right $\tau$-expansive). 

We note that all right or left $\tau$-expansiveness notions above 
have been defined only for onto endomorphisms of 
invertible dynamical systems $(X,\tau)$. However the following 
notions are defined for not necessarily onto endomorphisms 
of $(X,\tau)$. 

A (not necessarily onto) endomorphism $\varphi$ of 
an invertible dynamical system $(X,\tau)$ 
is said to be \itl{positively left $\tau$-expansive} 
(respectively, \itl{positively right $\tau$-expansive}) 
if there exists $\delta>0$ such that 
for any points $x,y$ of $X$ it holds that if 
$d_X(\tau^j\varphi^i(x),\tau^j\varphi^i(y)\leq\delta$ 
for all $i\geq 0$ and $j\leq 0$ (respectively, $j\geq 0$) 
then $x=y$. 

Let $\varphi$ be an endomorphism of 
an invertible dynamical system $(X,\tau)$. 
It is clear that if $\varphi$ is onto and 
positively left $\tau$-expansive 
(respectively, positively right $\tau$-expansive)  
then $\varphi$ is left $\tau$-expansive 
(respectively, right $\tau$-expansive). 
As will be seen in Proposition 11.2, $\varphi$ is 
positively expansive (in the usual sense) 
if and only if 
$\varphi$ is positively left $\tau$-expansive 
and positively right $\tau$-expansive. 

We note that M. Sablik \cite{Sab} 
already defined notions 
equivalent to 
``positively left $\sigma$-expansiveness''
and ``positively right $\sigma$-expansiveness'' 
for (not necessarily onto) endomorphisms 
of subshifts $(X,\sigma)$. 
(He also defined notions
equivalent to ``left $\sigma$-expansiveness''
and ``right $\sigma$-expansiveness'' 
for automorphisms of subshifts $(X,\sigma)$ 
independently of the definition of those 
for onto endomorphisms in \cite{Nasu-te}.) 
See the explanation 
after the proof of Theorem 11.8 for detail. 

Let $X$ be a compact metric space. 
Let $\varphi:X\to X$ be an onto continuous map. 
Let $\mathcal{O}_{\varphi}$ 
be the metric space consisting of  
all $\varphi$-orbits 
endowed with the 
metric $d_{\mathcal{O}_{\varphi}}$ such that 
\[d_{\mathcal{O}_{\varphi}}((x_i)_{i \in \Z},(y_i)_{i \in \Z})= 
\sup\{2^{-|i|}d_X(x_i,y_i) \bigm| i \in \Z \},\q 
(x_i)_{i \in \Z}, (y_i)_{i \in \Z} 
\in \mathcal{O}_{\varphi}. 
\] 
Let $\bsigma_\varphi:\Oh_\varphi\to\Oh_\varphi$
be the homeomorphisms\ defined by
\begin{align*}
\bsigma_\varphi((x_i)_{i\in\Z})=(x_{i+1})_{i\in\Z}.
\end{align*}
The dynamical system $(\mathcal{O}_{\varphi},\bsigma_\varphi)$
is the \itl{inverse limit system of $\varphi$}. 
(See, e.g., \cite{AokHir} for dynamical properties 
of inverse limit systems.) If $\varphi$ is a 
homeomorphism then the inverse limit system of $\varphi$ is 
conjugate to the dynamical system $(X,\varphi)$. 

Let $A$ be an alphabet 
(i.e. a finite nonempty set of symbols). Throughout this paper, 
we assume that an alphabet is not a singleton. 
A finite sequence $a_1\dots a_n$ with $a_j\in A$ and $n\geq 1$ 
is called a \itl{word} or \itl{block} of \itl{length} $n$ over $A$. 

Let $A^{\Z}$ be endowed with 
the metric $d$
such that 
for $x=(a_j)_{j \in \Z}$ and 
$y=(b_j)_{j \in \Z}$ with $a_j, b_j \in A,$ 
$d(x,y)=0$ if $x=y$, and otherwise 
$d(x,y)= 1/(1+k)$, where 
$k= \min \{\lvert j \rvert \bigm| a_j \neq b_j \}$. 
The metric $d$ is 
compatible with the product topology of 
the discrete topology on $A$. 
Let $\sigma_A: A^{\Z} \to A^{\Z}$ be defined by 
$\sigma_A((a_j)_{j \in \Z})=(a_{j+1})_{j \in \Z}$.
The dynamical system $(A^{\Z},\sigma_A)$ is called 
the \itl{full shift} over $A$, This is also called 
the \itl{full $n$-shift} if the cardinality of $A$ is $n$. 
For a closed subset 
$X$ of $A^{\Z}$ with $\sigma_{A}(X) = X$, 
the dynamical system
$(X,\sigma)$ or $(X,\sigma_{X})$ is called a \itl{subshift} 
over $A$, where $\sigma$ or $\sigma_{X}$   is 
the restriction of $\sigma_A$ on $X$. 

For a subshift $(X,\sigma)$ over an alphabet $A$ 
and $k\geq 1$, let $L_k(X)$ denote the set of 
all words $a_j\dots a_{j+k-1}$ 
that appear on some point $\aseq\in X$ with $a_j\in A$. 

Let $(X,\sigma_X)$ and $(Y,\sigma_Y)$ be subshifts. 
Let $N$ be a nonnegative integer. 
A mapping $f:L_{N+1}(X)\to L_1(Y)$ is called 
a \itl{local rule} of \itl{neighborhood-size} $N$ 
\itl{on $(X,\sigma_X)$ to $(Y,\sigma_Y)$}   
if $(f(a_j\dots a_{j+N}))_{j\in\Z}\in Y$ 
for all $(a_j)_{j\in\Z} \in X$ with $a_j\in L_1(X)$. 
If $X=Y$ in the above, then 
 $f$ is called a \itl{local rule on $(X,\sigma_X)$}.

Let $f:L_{N+1}(X)\to L_1(Y)$ be a local rule 
on $(X,\sigma_X)$ to $(Y,\sigma_Y)$. 
Let $m$ and $n$ be nonnegative integers with $m+n=N$. 
A mapping $\phi:X\to Y$ is called 
a \itl{block map of $(m,n)$-type given by $f$}, if 
\begin{align*}
\phi((a_j)_{j\in\Z})=(b_j)_{j\in\Z}, 
\q\q\text{where}\;\> 
b_j= f(a_{j-m}\dots a_{j+n}) \;\>\text{for all}\;\> j\in \Z. 
\end{align*} 

For $m,n\geq 0$, a mapping $\phi:X\to Y$ is called a 
\itl{block map of $(m,n)$-type}, if there exists a
local rule $f:L_{m+n+1}(X)\to L_1(Y)$ such that 
$\phi$ is a block map of $(m,n)$-type given by $f$. 
A block map of $(m,n)$ type is said to have 
\itl{memory} $m$ and \itl{anticipation} $n$. 
A block map of $(m,n)$-type for some $m,n\geq 0$ is 
simply called a \itl{block map}. 

For subshifts $(X,\sigma_X)$ and $(Y,\sigma_Y)$,  
a mapping $\phi:X\to Y$ is a homomorphism of 
$(X,\sigma_X)$ into $(Y,\sigma_Y)$
(i.e. a continuous map
with $\phi\sigma_X=\sigma_Y\phi$) if and only if 
$\phi$ is a block map 
(Curtis-Hedlund-Lyndon Theorem \cite{Hedlund}).  Hence we use the 
same terminology as above for 
homomorphisms between subshifts. 

Let $A^\N=\{(a_j)_{j\in\N}\,|\, a_j\in A\}$ 
be endowed with a metric compatible with  the product topology of 
the discrete topology on $A$. Let
$\tilde \sigma_A: A^\N\to A^\N$ be defined by 
$\tilde \sigma_A((a_j)_{j\in\N}) = (a_{j+1})_{j\in\N}$. 
The dynamical system 
$(A^\N,\tilde \sigma_A)$ is called the 
\textit{onesided full shift} over $A$. 
For a subshift $(X,\sigma)$ 
over $A$, let 
$\tilde{X} = \{(a_j)_{j \in\N}\,|\, \exists (a_j)_{j\in\Z}\in X\}$.  
Then with the onto continuous map 
$\tilde{\sigma}=\tilde{\sigma}_A|\tilde{X}$ 
we have a dynamical system $(\tilde{X},\tilde{\sigma})$, 
which is called
a \itl{onesided subshift} over $A$ and is said to be \itl{induced by} 
$(X,\sigma)$.  From  
$(\tilde{X},\tilde{\sigma})$ we can uniquely recover $(X,\sigma)$, 
and hence $(X,\sigma)$ is said to be \itl{induced by} 
$(\tilde{X},\tilde{\sigma})$. 
Let $s_X: X\to \tilde{X}$ be the continuous map which maps  
$(a_j)_{j\in\Z}$ to $(a_j)_{j\in\N}$. 

If a homomorphism $\phi$ of a 
subshift $(X,\sigma_X)$ into another $(Y,\sigma_Y)$ has memory zero, 
then it uniquely 
\itl{induces} a homomorphism 
$\tilde{\phi}$ between the induced onesided subshifts 
$(\tilde{X},\tilde{\sigma}_X)$ and $(\tilde{Y},\tilde{\sigma}_Y)$
such that 
$s_Y\phi=\tilde{\phi}s_X$, and conversely 
each homomorphism $\tilde{\phi}$ of a onesided subshift 
$(\tilde{X},\tilde{\sigma}_X)$ into another 
$(\tilde{Y},\tilde{\sigma}_Y)$ uniquely \itl{induces}  
a homomorphism $\phi$ between the induced subshifts $(X,\sigma_X)$ 
and $(Y,\sigma_Y)$
such that $\phi$ has memory zero and 
$s_Y\phi=\tilde{\phi}s_X$. 

For a subshift $(X,\sigma)$ and $n\geq 1$, 
we define the \itl{higher block system of order $n$} 
of $(X,\sigma)$
to be the subshift $(X^{[n]},\sigma^{[n]})$ over the 
alphabet $L_n(X)$ with 
$X^{[n]}=\{ (a_j\dots a_{j+n-1})_{j\in\Z}\bigm|
(a_j)_{j\in\Z}\in X,\;a_j\in L_1(X)\}$. 

For a homomorphism $\phi:(X,\sigma_X)\to (Y,\sigma_Y)$ 
between subshifts and $n\geq 1$, we define the homomorphism 
$\phi^{[n]}:(X^{[n]},\sigma_{X^{[n]}})\to (Y^{[n]},\sigma_{Y^{[n]}})$, 
which is called the 
\itl{higher block presentation of order $n$ of $\phi$}, 
as follows:
if $\phi$ maps $\aseq\in X$ to $\bseq\in Y$ 
with $a_j\in L_1(X)$ and $b_j\in L_1(Y)$, then 
$\phi^{[n]}$ maps $(a_j\dots a_{j+n-1})_{j\in\Z}$
to $(b_j\dots b_{j+n-1})_{j\in\Z}$. 

For a subshift $(X,\sigma)$ and $m,n\geq 0$, the block-map 
$\rho_{X,m,n}$ which maps $(a_j)_{j\in\Z}$ 
to $(a_{j-m}\dots a_{j+n})_{j\in\Z}$ gives a topological conjugacy 
of $(X,\sigma)$ onto $(X^{[N+1]},\sigma^{[N+1]})$ 
with $N=m+n$,
which is called the \itl{higher-block conjugacy of 
$(m,n)$-type on $X$} or a \itl{higher-block conjugacy}.  
If $\varphi$ is an endomorphism of a subshift $(X,\sigma)$ 
then $\rho_{X,m,n}$ is also a topological conjugacy of 
$(X,\sigma,\varphi)$ onto 
$(X^{[N+1]},\sigma^{[N+1]},\varphi^{[N+1]})$, which is 
a \itl{higher-block conjugacy} between the commuting systems and 
considered to be a \itl{higher-block conjugacy} of
the endomorphism $\varphi$ onto 
its \itl{higher-block presentation} $\varphi^{[N+1]}$. 

\subsection{Graph-homomorphisms}

Let $G$ be a graph. Here a \itl{graph} means 
a directed graph which may  
have multiple arcs and loops.
Let $A_G$ and $V_G$ denote the arc-set and 
the vertex-set, respectively, of $G$.
Let $i_G: A_G \rightarrow V_G$ 
and $t_G : A_G \rightarrow V_G$ 
be the mappings
such that for each arc $a \in A_G$, 
$i_G(a)$ and $t_G(a)$ are 
its initial and terminal vertices. 
Hence the graph $G$ is represented by 
\[V_G\stackrel{i_G}{\longleftarrow}A_G
\stackrel{t_G}{\longrightarrow}V_G.\]
We say that $G$ is \itl{nondegenerate} 
if both $i_G$ and $t_G$ are onto. 

Let $X_G$ be the set of all points 
$(a_j)_{j \in \Z}$ in $A_G^{\Z}$ 
such that $t_G(a_j) = i_G(a_{j+1})$ 
for all $j\in\Z$. 
Then we have a subshift $(X_G,\sigma_G)$ 
over $A_G$, 
which is called the 
\itl{topological Markov shift} 
defined by $G$. We call $G$  the
\itl{defining graph} of $(X_G,\sigma_G)$. 
The topological Markov shift $(X_G,\sigma_G)$ 
is also denoted by $(X_M,\sigma_M)$, 
where $M$ is the \itl{adjacency matrix of $G$}, i.e.,  
the square matrix
$M_G=(m_{u,v})_{u,v\in V_G}$ 
such that the $(u,v)$-component $m_{u,v}$ is 
the number of arcs going from 
vertex $u$ to vertex $v$. 
We call $M$ the \itl{defining matrix} of the 
topological Markov shift 
$(X_M,\sigma_M)=(X_G,\sigma_G)$. 
In what follows, 
we assume, without loss of generality, that 
the defining graph of
any topological Markov shift 
is nondegenerate. 
(Hence, when we write $(X_G,\sigma_G)$ or $X_G$, 
$G$ is always assumed to be nondegenerate.)

For an alphabet $A$, let $G_A$ denote the one-vertex graph with 
arc-set $A$. Then the full-shift over $A$ is the topological 
Markov shift defined by $G_A$. 

A \itl{subshift finite type}, abbreviated \itl{SFT}, 
is a subshift 
which is topologically conjugate to a 
topological Markov shift. 
A \itl{sofic system} is a subshift 
which is the image of a topological Markov shift 
under a block map. 
 
For graphs $\G$ and $G$, 
a \itl{graph-homomorphism} $h$ of $\G$ into $G$, 
written by 
$h: \G \rightarrow G$, is a pair $(h_A,h_V)$ 
of mappings $h_A: A_\G \rightarrow A_G$ 
(\itl {arc-map}) and 
$h_V: V_\G\rightarrow V_G$ (\itl {vertex-map}) 
such that 
the following diagram is commutative.  
\begin{align*} \begin{CD}
V_\G @<i_{\G}<< A_\G @>{t_\G}>> V_\G  \\
@V{h_V}VV                    @V{h_A}VV           
@VV{h_V}V  \\
V_G @<{i_G}<<  A_G       
@>{t_G}>>   V_G 
\end{CD}
\end{align*} 
A graph-homomorphism is said to be \itl{onto} if 
$h_A$ and $h_V$ are onto. 

We call a block map of $(0,0)$ type between subshift spaces a 
\itl{1-block map}. For a graph-homomorphism $h:\G\to G$, we define 
$\phi_h:X_\G\to X_G$ to be the 1-block map given by the local rule 
$h_A:A_\G\to A_G$. 

A graph-homomorphism $h:\G \to G$ is said to be
\textit{right-resolving} if for each pair $(u,a)$ of $u \in V_{\G}$ 
and $a\in A_G$ with $i_G(a)=h_V(u)$, there exists  
unique $\alpha\in A_\G$ with $i_\G(\alpha)=u$ and $h_A(\alpha)=a$.
It is said to be \itl{left-resolving} if for each pair $(u,a)$ 
of $u \in V_{\G}$ 
and $a\in A_G$ with $t_G(a)=h_V(u)$, there exists 
unique $\alpha\in A_\G$ with $t_\G(\alpha)=u$ and $h_A(\alpha)=a$. 
We say that $h$ is \itl{biresolving} if $h$ is right-resolving and 
left-resolving. 
Note that if $h$ is a right-resolving or left-resolving, onto 
graph-homomorphism, then $\phi_h$ is onto. 

A graph-homomorphism $h:\G\to G$ is said to be 
\itl{weakly right-resolving} 
if each $\alpha \in A_{\G}$ is uniquely determined by 
$i_\G(\alpha)$ and $h_A(\alpha)$. We say that $h$ is 
\textit{weakly left-resolving} if 
each $\alpha \in A_{\G}$ is uniquely determined by 
$t_\G(\alpha)$ and $h_A(\alpha)$. We say that $h$ is 
\textit{weakly biresolving} if $h$ is weakly 
right-resolving and weakly left-resolving. 

By using the following well-known lemma,  
the definitions and results concerning 
weakly right-resolving 
and weakly left resolving graph-homomorphisms are 
interpreted as those concerning 
right-resolving and left-resolving graph-homomorphisms  
for important special cases. 
\begin{lemma} Let $\G$ and $G$ be nondegenerate graphs. 
Let $h:\G\to G$ be a weakly right-resolving 
(respectively, weakly left-resolving) graph-homomorphism. 
\begin{enumerate} 
\item  If $\phi_h$ is bijective, then $h$ is 
right-resolving (respectively, left-resolving)
\item  If $G$ and $\G$ are irreducible and have 
the same spectral radius, then $h$ is 
right-resolving (respectively, left-resolving) . 
\end{enumerate} 
\end{lemma} 

Let $G$ be a graph. Let $L_0(G)=V_G$, let $L_1(G)=A_G$ and 
for $k\geq 2$ let $L_k(G)$ be the set of all words $a_1\dots a_k$
with $a_j\in A_G$
such that  $t_G(a_j)=i_G(a_{j+1})$ for $j=1,\dots,k-1$. For $k\geq 0$, 
we call an element of $L_k(G)$ a \itl{path} of \itl{length $k$} in $G$. 
Let $L(G)$ denote the set of all paths of any length in $G$. 
We extend $i_G:A_G\to V_G$ 
and $t_G:A_G\to V_G$ to
$i_G:L(G)\to V_G$ and $t_G:L(G)\to V_G$, respectively, as follows: 
$i_G(v)=t_G(v)=v$ for $v\in L_0(G)$; for $w=a_1\dots a_k\in L_k(G)$,  
$i_G(w)=i_G(a_1)$ and $t_G(w)=t_G(a_k)$. A path $w$ is said to 
\itl{go from $i_G(w)$ to $t_G(w)$}, and $i_G(w)$ and $t_G(w)$ are 
called the \itl{initial vertex} and 
\itl{terminal vertex}, respectively, of the path $w$. 
For a graph-homomorphism $h:\G\to G$, we also extend 
$h_V$ and $h_A$ to $h:L(\G)\to L(G)$ as follows:
$h(v)=v$ for $v\in L_0(G)$; for $\mu=\alpha_1\dots \alpha_k\in L_k(\G)$ 
with $k\geq 1$ and $\alpha_j\in A_\G$, 
$h(\mu)=h(\alpha_1)\dots h(\alpha_k)$. We say that a path 
$\mu$ \itl{generates} $h(\mu)$ under $h$. 

Let $h:\G\to G$ be a graph-homomorphism between nondegenerate graphs.  
For $U\subset V_{\G}$ and $w\in L(G)$, define 
\begin{gather*} 
S^+_h(U,w)=\{t_{\G}(\mu)\bigm| \mu\in L(\G), i_{\G}(\mu)\in U, h(\mu)=w\}, \\
S^-_h(w,U)=\{i_{\G}(\mu)\bigm| \mu\in L(\G), t_{\G}(\mu)\in U, h(\mu)=w\}.
\end{gather*}
and define 
\begin{gather*} 
B^+_h(U,w)=\{\mu\in L(\G)\bigm| i_{\G}(\mu)\in U, h(\mu)=w\}, \\
B^-_h(w,U)=\{\mu\in L(\G)\bigm|  h(\mu)=w, t_{\G}(\mu)\in U\}.
\end{gather*}
We call $S^+_h(U,w)$  the
\itl{$w$-successor} or a \itl{successor} \itl{of $U$ under $h$}  
and call $S^-_h(w,U)$ the \itl{w-predecessor} or 
a \itl{predecessor} \itl{of $U$ 
under $h$}. 
A \itl{right-compatible set for $h$} and
a \itl{left-compatible set for $h$} are defined to be a 
nonempty successor of a singleton of $V_\G$ under $h$ and a nonempty 
predecessor of a singleton of $V_\G$ under $h$, respectively. 
Let $C^+_h$ (respectively, $C^-_h$) be the family consisting of all 
maximal right-compatible (respectively, left-compatible) sets 
and their nonempty successors (respectively, predecessors). 
(Throughout this paper, we shall often abuse the element 
of a singleton to denote the singleton.) 
We define $h^+:\G^+_h\to G$ to be the graph-homomorphism such that 
the vertex-set of    
$\G^+_h$ is $\C_h^+$, the arc-set of $\G^+_h$ is 
$\{(U,a)\,|\, U\in\C^+_h,\, a\in A_G,\, h(U)=i_{G}(a)\}$
and each arc $(U,a)$ goes from $U$ to $S^+_h(U,a)$ with 
$h^+((U,a))=a$. We also define
$h^-:\G^-_h\to G$ to be the graph-homomorphism such that 
the vertex-set of    
$\G^-_h$ is $\C_h^-$, the arc-set of $\G^-_h$ is 
$\{(a,U)\,|\,  a\in A_G,\, U\in\C^-_h,\, h(U)=t_{G}(a)\}$ 
and each arc $(a,U)$ goes from  $S^-_h(a,U)$ to $U$ with 
$h^-((a,U))=a$. 
We call $h^+$ (respectively, $h^-$)
the \itl{induced right-resolving}
(respectively, \itl{induced left-resolving}) graph-homomorphism 
of $h$, though generally $h^+$ (respectively, $h^-$) is a 
weakly right-resolving (respectively weakly left-resolving) 
graph-homomorphism.  
(For the case that $G=G_A$ for an alphabet $A$, 
the induced right-resolving 
and left-resolving graph-homomorphisms
are the same as the ``induced right-resolving and induced left-resolving 
$\lambda$-graphs'' introduced in \cite{Nasu-t-ex}. 
These were first introduced in \cite{Nasu-l} for a special type of 
$\lambda$-graphs.) 
However, we note the following:
\begin{lemma} 
\begin{enumerate} 
\item If $h$ is a graph-homomorphism between nondegenerate 
graphs with $\phi_h$ bijective, then 
$h^+$ is right-resolving and $h^-$ is left-resolving; 
\item \cite{Nasu-c} if 
$h$ is a graph-homomorphism between irreducible 
graphs having the same spectral radius with $\phi_h$ onto, 
then $h^+$ is right-resolving and $h^-$ is left-resolving. 
\end{enumerate}
\end{lemma} 
\begin{proof}
(1) Since $\phi_{h^+}$ and $\phi_{h^-}$ are bijections (see 
\cite[the proof of Lemma 5.5]{Nasu-c}), 
the conclusion follows by Lemma 2.1(1). 
 
(2) By \cite[Proposition 4.1]{Nasu-c}. 
\end{proof} 

For a graph-homomorphism $h:\G\to G$ and an integer $k\geq 0$, 
we say that $h$ is $k$  right-mergible 
(respectively, $k$ left-mergible) 
if for any pair of paths 
$\mu_1=\alpha_1\dots\alpha_{k+1}$ and 
$\mu_2=\beta_1\dots\beta_{k+1}$ in $L_{k+1}(X_\G)$ 
with $\alpha_j,\beta_j\in A_\G$ it holds 
that if $i_\G(\mu_1)=i_\G(\mu_2)$  
(respectively, $t_\G(\mu_1)=t_\G(\mu_2)$) and
$h(\mu_1)=h(\mu_2)$, then $\alpha_1=\beta_1$. 
(These notions correspond to ``nonexistence of  
right (respectively, left) $f$-branch of length exceeding $k$'' 
in \cite[Section 16]{Hedlund} and appear 
for a special class of graph-homomorphisms in \cite{Nasu-l} 
and for the general class of graph-homomorphisms  
with somewhat different terms in \cite[p.400]{Nasu-c}.) 

For a graph $G$ and $n\geq 1$, we define the 
\itl{higher-block presentation}
$G^{[n]}$ of order $n$ of $G$ as follows: $G^{[1]}=G$;
if $n\geq 2$, $G^{[n]}$ is the graph 
such that $A_{G^{[n]}}=L_n(G),V_{G^{[n]}}=L_{n-1}(G)$ and 
for $\alpha=a_1\dots a_n\in A_{G^{[n]}}$
with $a_j\in A_G$, $i_{G^{[n]}}(\alpha)=a_1\dots a_{n-1}$ and 
$t_{G^{[n]}}(\alpha)=a_2\dots a_n$. 
For a graph-homomorphism $h:\G\to G$ and $n\geq 1$, we define the 
\itl{higher-block presentation} $h^{[n]}:\G^{[n]}\to G^{[n]}$
\itl{of order} $n$ of $h$ by 
$h_A^{[n]}(\alpha)=h(\alpha), \alpha\in L_{n}(\G)$. 
\begin{remark} Let $h:\G\to G$ be a graph-homomorphism. Let 
$k\geq 0$ and $n\geq 1$.
\begin{enumerate}
\item \cite[Lemma 5.6]{Nasu-c}  
If $h$ is $k$ right-mergible, then so is $h^{[n]}$, and 
if $h$ is $k$ left-mergible, then so is $h^{[n]}$. 
\item  
If $h:\G\to G$ is $k$ right-mergible (respectively, 
$k$ left-mergible), then  
$\{B^+_h(i_\G(\mu),h(\mu))\,|\,\mu\in L_k(\G)\}=\C^+_{h^{[k+1]}}$ 
(respectively, 
$\{B^-_h(h(\mu),t_\G(\mu))\,|\,\mu\in L_k(\G)\}=\C^-_{h^{[k+1]}}$), 
and hence
any two distinct sets in $\C^+_{h^{[k+1]}}$ 
(respectively, in $\C^-_{h^{[k+1]}}$)
are disjoint, and hence all sets 
in $\C^+_{h^{[k+1]}}$ (respectively, in $\C^-_{h^{[k+1]}}$)
are maximal right-compatible sets 
(respectively, maximal left-compatible sets) for $h^{[k+1]}$.
\end{enumerate}
\end{remark}  

Let $(X_G,\sigma_G)$ and $(X_H,\sigma_H)$ be topological 
Markov shifts. 
Let $f:L_{N+1}(G)\to A_H$ be 
a local rule (on $(X_G,\sigma_G)$ to $(X_H,\sigma_H)$).  
Then we naturally define the graph-homomorphism $q_f:G^{[N+1]}\to H$
such that $q_f(w)=f(w)$ for $w\in A_{G^{[N+1]}}=L_{N+1}(G)$. 
For $k\geq 0$, we say that $f$ is \itl{$k$ right-mergible}
(respectively, \itl{$k$ left-mergible}) if $q_f$ is 
\itl{$k$ right-mergible}
(respectively, \itl{$k$ left-mergible}), 
(These definitions are compatible 
with the general definition of mergibility for a local rule 
on a subshift to another presented in the next section.) 

For convenience's sake, we extend the local rule
$f:L_{N+1}(G) \to A_H$ to the mapping
$f:\cup_{s\geq 1}L_{N+s}(G)\to \cup_{s\geq 1}L_s(H)$ such that 
for $s\geq 1$ and for $a_1\dots a_{N+s}\in L_{N+s}(G)$ with $a_j\in A_G$,
\[f(a_1\dots a_{N+s})=f(a_1\dots a_{N+1})\dots f(a_s\dots a_{N+s}).\] 

Suppose that the local rule $f$ is $k$ right-mergible.   
For $w\in L_{N+k}(G)\cup L_{N+k+1}(G)$,
we define $D^+_{f;k}(w)$ as follows: 
if $w=w_0w_1$ with 
$w_0\in L_N(G)\cup L_{N+1}(G)$ and $w_1\in L_k(G)$, then 
\[D^+_{f;k}(w)=\{w_0w'_1\in L(G)\,\bigm|\, 
w'_1\in L_k(G), f(w_0w'_1)=f(w)\}.\] 
Since $f$ is $k$ right-mergible, we can define a graph $G^+_{f;k}$  
and a graph-homomorphism 
\[q^+_{f;k}:G^+_{f;k}\to H^{[k+1]}\]
as follows: 
the vertex-set of $G^+_{f;k}$ is 
$\{D^+_{f;k}(w)\,|\,w\in L_{N+k}(G)\}$; 
the arc-set of $G^+_{f;k}$ is 
$\{D^+_{f;k}(w)\,|\,w\in L_{N+k+1}(G)\}$; 
each arc $D^+_{f;k}(w)$ with $w\in L_{N+k+1}(G)$ goes from vertex
$D^+_{f;k}(w')$ to vertex $D^+_{f;k}(w'')$, 
where $w'$ and $w''$ are the initial and terminal subpaths, 
respectively, of length $N+k$ of $w$, and 
$q^+_{f;k}(D^+_{f;k}(w))=f(w)$.

Suppose that $f$ is 
$l$ left-mergible.    
For $w\in L_{N+l}(G)\cup L_{N+l+1}(G)$,
we define $D^-_{f;,l}(w)$ as follows: 
if $w=w_{-1}w_0$ with 
$w_{-1}\in L_l(G)$ and $w_0\in L_N(G)\cup L_{N+1}(G)$, then 
\[D^-_{f;l}(w)=\{w'_{-1}w_0\in L(G)\,\bigm|\, 
w'_{-1}\in L_l(G), f(w'_{-1}w_0)=f(w)\}.\]
Since $f$ is $l$ left-mergible, we can define a graph $G^-_{f;l}$ 
and a graph-homomorphism
\[q^-_{f;l}:G^-_{f;l}\to H^{[l+1]}\] 
in the way symmetric to the above. 

Suppose that $f$ is $k$ right-mergible and $l$ left-mergible. 
For $w\in L_{N+k+l}(G)\cup L_{N+k+l+1}(G)$,
we define $D^{-+}_{f;l,k}(w)$ as follows:
if $w=w_{-1}w_0w_1$ with 
$w_{-1}\in L_l(G)$, $w_0\in L_N(G)\cup L_{N+1}(G)$ 
and $w_1\in L_k(G)$, then  
\[D^{-+}_{f;l,k}(w)=\{w'_{-1}w_0w'_1\in L(G)\,\bigm|\, 
w'_{-1}\in L_l(G), w'_{1}\in L_k(G), f(w'_{-1}w_0w'_1)=f(w)\}.\]  
($D^{-+}_{f;l,k}(w)$ corresponds to a``maximal 
connected $f$-covering'' of $w$ in \cite[p.360]{Hedlund}.)
Since $f$ is $k$ right-mergible and $l$ left-mergible, 
we can define 
a graph $G^{-+}_{f;l,k}$ and a graph-homomorphism 
\[q^{-+}_{f;l,k}:G^{-+}_{f;l,k}\to H^{[l+k+1]}\]
as follows: 
the vertex-set of $G^{-+}_{f;l,k}$ is 
$\{D^{-+}_{f;l,k}(w)\,|\,w\in L_{N+k+l}(G)\}$; 
the arc-set of $G^{-+}_{f;l,k}$ is 
$\{D^{-+}_{f;l,k}(w)\,|\,w\in L_{N+k+l+1}(G)\}$; 
each arc $D^{-+}_{f;l,k}(w)$ with $w\in L_{N+k+l+1}(G)$ goes from vertex
$D^{-+}_{f;l,k}(w')$ to vertex $D^{-+}_{f;l,k}(w'')$, 
where $w'$ and $w''$ are the initial and terminal subpaths, 
respectively, of length $N+k+l$ of $w$, and 
$q^{-+}_{f;l,k}(D^{-+}_{f;l,k}(w))=f(w)$.

These definitions are due to Bruce Kitchens \cite{Kit-c}, 
\cite[Section 4.3]{Kit-s}.   
We easily see that  
$q^+_{f;k}=(q_f^{[k+1]})^{^+}$ 
(up to isomorphism between graph-homomorphisms) 
when $f$ is $k$ right-mergible, and that 
$q^-_{f;k}=(q_f^{[k+1]})^{-}$ 
when $f$ is $k$ left-mergible. One can see that 
$q^{-+}_{f;l,k}=((q_f^{[k+l+1]})^+)^-=((q_f^{[k+l+1]})^-)^+$ 
when $f$ is $k$ right-mergible and $l$ left-mergible 
(see \cite[Section 5]{Nasu-c}). 

\begin{proposition}[Kitchens]
Let $(X_G,\sigma_G)$ and $(X_H,\sigma_H)$ be topological 
Markov shifts. 
Let $f:L_{N+1}(G)\to A_H$ be 
a local rule (on $(X_G,\sigma_G)$ to $(X_H,\sigma_H)$). 
\begin{enumerate}
\item\cite{Kit-c} 
When $f$ is $k$ right-mergible,  
$q^+_{f;k}:G^+_{f;k}\to H^{[k+1]}$ 
is weakly right-resolving.  When $f$ is $l$ left-mergible,  
$q^-_{f;l}:G^-_{f;l}\to H^{[l+1]}$ 
is weakly left-resolving. 
\item\cite{Kit-s}
When $f$ is $k$ right-mergible and $l$ left-mergible, 
$q^{-+}_{f;l,k}:G^{-+}_{f;l,k}\to H^{[l+k+1]}$ is weakly biresolving. 
\end{enumerate}
\end{proposition}

\subsection{Textile systems and textile-subsystems}

A \textit{textile system} $T$ 
\textit{over} 
a graph $G$ is an ordered pair 
of graph-homomorphisms $p: \G \to G$ and 
$q: \G \to G$. 
We write
\begin{align*}
T=(p,q: \G \to G).
\end{align*}
We have the following commutative diagram.
\begin{align*}
\begin{CD}
V_{G}@< i_{G} << A_{G}@> t_{G} >> V_{G} \\
@AA{p_V}A    @A{p_A}AA         @A{p_V}AA \\
V_{\G}@< i_{\G}<< A_{\G}@> t_{\G} >> V_{\G}\\ 
@VV{q_V}V    @V{q_A}VV         @V{q_V}VV \\ 
V_{G}@<i_{G} << A_{G}@> t_{G} >> V_{G} 
\end{CD}
\end{align*}
If we observe this diagram vertically, 
then we have the ordered pair of
graph-homomorphisms
\begin{align*}
\begin{CD}
V_{G}@< i_{G} << A_{G} \\
@AA{p_V}A    @A{p_A}AA        \\
V_{\G}@< i_{\G}<< A_{\G}\\ 
@VV{q_V}V    @V{q_A}VV         \\ 
V_{G}@<i_{G} << A_{G} 
\end{CD}
\qquad
\mathrm{and}
\qquad
\begin{CD}
 A_{G}@> t_{G} >> V_{G} \\
    @A{p_A}AA         @A{p_V}AA \\
 A_{\G}@> t_{\G} >> V_{\G}\\ 
   @V{q_A}VV         @V{q_V}VV \\ 
A_{G}@> t_{G} >> V_{G} 
\end{CD}.
\end{align*}
This defines another textile system 
\begin{align*}T^*=(p^*,q^*:\G^*\to G^*)\end{align*}
called the \textit{dual} of $T$, 
where 
$i_{\G^*}=p_A, t_{\G^*}=q_A$, $i_{G^*}=p_V$ 
and $t_{G^*}=q_V$.

Let $T=(p,q:\G\to G)$ be a textile system. 
Let $\xi=\phi_p$ and let $\eta=\phi_q$.
A two-dimensional configuration 
$(\alpha_{ij})_{i,j\in \Z}, \alpha_{ij} \in A_{\G}$, 
is called a \textit{textile} woven by $T$ 
if $(\alpha_{ij})_{j\in \Z}\in X_{\G}$ 
and $\eta((\alpha_{i-1,j})_{j\in \Z})=\xi((\alpha_{ij})_{j\in \Z})$
for all $i\in \Z $. 
Let $U_T$ denote the set of all textiles woven by $T$. 
Then $U_T$ is 
a closed subset of the space  
$A_\G^{\Z^2}$ 
equipped with  the product topology of 
the discrete topology on $A_\G$ and is
\itl{shift invariant} 
(i.e., if $(\alpha_{i,j})_{i,j\in\Z}\in U_T$ 
with $\alpha_{i,j}\in A_{\G}$, then 
$(\alpha_{i+k,j+l})_{i,j\in\Z}\in U_T$ for all $k,l\in\Z$). 
Define
\[
X_T = \{\xi((\alpha_{0j})_{j \in \Z}) \bigm| 
(\alpha_{ij})_{i,j \in \Z} \in U_T\},\q 
Z_T = \{(\alpha_{0j})_{j \in \Z} \bigm| 
(\alpha_{ij})_{i,j \in \Z} \in U_T\}. 
\] 
Then we have subshifts $(X_T,\sigma_T)$ and $(Z_T,\varsigma_T)$ .
We call $(X_T,\sigma_T)$  the \itl{woof shift} 
of $T$ and $(X_{T^*},\sigma_{T^*})$ the \itl{warp shift} 
of $T$. 
We say that $T$ is \itl{nondegenerate} 
if $(X_T,\sigma_T) = (X_G,\sigma_G)$ (with $G$ nondegenerate).  
We assume that a nondegenerate textile system is always defined 
over a nondegenerate graph. 
We define onto maps 
$\xi_T: Z_T \rightarrow X_T$ and $\eta_T: Z_T \rightarrow X_T$ 
to be the restrictions of $\xi$ and $\eta$, respectively. 
If $T$ is \textit{onesided 1-1}, i.e., $\xi_T$ is 1-1, 
then an onto endomorphism $\varphi_T$ of $(X_T,\sigma_T)$ 
is defined by 
\[
\varphi_T = \eta_T\xi_T^{-1}. 
\]
If $T$ is \textit{1-1}, i.e., both $\xi_T$ and $\eta_T$ are 1-1, 
then $\varphi_T$ is an automorphism of $(X_T,\sigma_T)$. 
We also have the onesided subshifts $(\tilde Z_T,\tilde \varsigma_T)$ 
and $(\tilde X_T,\tilde \sigma_T)$ induced by 
$(Z_T,\varsigma_T)$ and $(X_T,\sigma_T)$, respectively. 

For a textile system $T=(p,q:\G\to G)$ and $n\geq 1$, 
we define the
\textit{higher block system of order $n$ of $T$} to 
be the textile system 
$T^{[n]}=(p^{[n]},q^{[n]}:\Gamma^{[n]}\to G^{[n]})$. 
($p^{[n]}$ and $q^{[n]}$ are the higher block 
presentations of order $n$ of $p$ and $q$, respectively.) 

Here we make the following remarks. The definition of
a textile system was simplified 
with no essential change in \cite{Nasu-te} so that 
it is a pair of graph-homomorphisms 
$T=(p,q:\G\to G)$. 
The additional condition 
that the correspondence between each arc
$\alpha$ in $\G$ and its quadruple
$(i_\G(\alpha),t_\G(\alpha), p_A(\alpha),q_A(\alpha))$ 
should be one-one was removed from 
the original definition in \cite{Nasu-t} with no 
impact on the statements of the results in the 
preceding papers 
\cite{Nasu-note,Nasu-t,Nasu-m,Nasu-e,Nasu-d,Nasu-n}. 
(The removal makes easier to define the composition 
and product of textile systems; see the remarks on them 
at the end of this subsection.)
However it is useful and natural by the duality of $T$ and $T^*$ 
to imagine that 
each arc $\alpha$ in the graph $\G$ is a ``square'' 
whose left, right, upper and lower sides 
are $i_\G(\alpha)$, $t_\G(\alpha)$, $p_A(\alpha)$ 
and $q_A(\alpha)$, respectively, 
(see  \cite[pp. 16,17]{Nasu-t}), 
which may be regarded as a Wang tile of a Wang tiling system
by unifying different ``squares'' 
with the same quadruple if any. For
onesided 1-1 textile systems and ``resolving'' textile 
systems, which were treated in the preceding papers 
and will be treated in this paper, we need not to consider 
multiple ``squares'' with the same quadruple. 

For a textile system $T$, 
a closed, shift-invariant subset $U$ of $U_T$ is called
a \itl{textile-subsystem of $T$} or simply  called 
a \itl{textile-subsystem}. Hence $U_T$ is 
a \itl{textile full-system of $T$} or simply  
a \itl{textile full-system}. 

Let $U$ be a textile-subsystem of a textile system $T$. 
The \itl{dual} $U^*$ of $U$ is defined by 
$U^*=\{(\alpha_{j,i})_{i,j\in\Z}\bigm| 
(\alpha_{i,j})_{i,j\in\Z}\in U, \alpha_{i,j}\in A_\G\}$. 
Clearly $U^*$ is a textile-subsystem of $T^*$. 
The \itl{woof shift} $(X_U,\sigma_U)$ of $U$ 
is defined by 
$X_U=\{\xi_T((\alpha_{0,j})_{j\in\Z})\bigm| 
(\alpha_{i,j})_{i,j\in\Z}\in U\}$ and 
$\sigma_U=\sigma_T| X_U$. We also define 
another subshift $(Z_U,\varsigma_U)$ by
$Z_U=\{(\alpha_{0,j})_{j\in\Z}\bigm| 
(\alpha_{i,j})_{i,j\in\Z}\in U\}$ and 
$\varsigma_U=\varsigma_T| Z_U$. 
We define onto maps $\xi_U: Z_U\to X_U$ and 
$\eta_U: Z_U\to X_U$ restricting 
$\xi_T$ and $\eta_T$, respectively.
If $U$ is \itl{onesided 1-1}, 
i.e., $\xi_U$ is one-to-one, then we have 
an onto endomorphism $\varphi_U$ of $(X_U,\sigma_U)$
by $\varphi_U=\eta_U\xi_U^{-1}$. 
We say that $U$ is \itl{1-1} if both $\xi_U$ 
and $\eta_U$ are one-to-one. 
For $n\geq 1$ 
we define a textile-subsystem $U^{[n]}$, which is called the 
\itl{higher block system of order $n$} of $U$, by 
$U^{[n]}=\{(\alpha_{i,j}\dots\alpha_{i,j+n-1})_{i,j\in\Z} 
\bigm| (\alpha_{i,j})_{i,j\in\Z}\in U\}$. Clearly 
$U^{[n]}$ is a textile-subsystem of $T^{[n]}$. 

If $U$ is a onesided 1-1 textile-subsystem, then $Z_U$ 
determines $U$. If $T$ is a textile system, then 
for any subshift space 
$Z\subset Z_T$ with $\xi_T(Z)=\eta_T(Z)$ and 
$\xi_T|Z$ one-to-one, there exists a unique 
onesided 1-1 textile-subsystem $U$ of $T$ with $Z_U=Z$.  

Let $T=(p,q:\G\to G)$ be a textile system. 
Let $\xi=\phi_p$ and $\eta=\phi_q$ as above. 
A two-dimensional configuration 
$(\alpha_{ij})_{i\in\N,j\in \Z}, \alpha_{ij} \in A_{\G}$, 
is called a \textit{half-textile} woven by $T$ 
if $(\alpha_{ij})_{j\in \Z}\in X_{\G}$ 
and $\eta((\alpha_{i-1,j})_{j\in \Z})=\xi((\alpha_{ij})_{j\in \Z})$
for all $i\geq 2 $. Let $\hf{U}_T$ 
denote the set of all half-textiles woven by $T$. 
Then $\hf{U}_T$ is 
a closed subset of the space  
$A_\G^{\N\times\Z}$ 
equipped with the product topology of 
the discrete topology on $A_\G$ and is
\itl{shift invariant} 
(i.e., if $(\alpha_{i,j})_{i\in\N,j\in\Z}\in \hf{U}_T$ 
with $\alpha_{i,j}\in A_{\G}$, then 
$(\alpha_{i+k,j+l})_{i\in\N,j\in\Z}\in \hf{U}_T$ 
for all $k\geq 0,l\in\Z$). 
Define
\[
\hf{X}_T = \{\xi((\alpha_{1j})_{j\in\Z})\bigm|
(\alpha_{ij})_{i\in\N,j\in\Z}\in\hf{U}_T\},\q 
\hf{Z}_T = \{(\alpha_{1j})_{j\in\Z}\bigm|
(\alpha_{ij})_{i\in\N,j\in \Z}\in\hf{U}_T\}. 
\] 
Then we have subshifts $(\hf{X}_T,\hf{\sigma}_T)$ and 
$(\hf{Z}_T,\hf{\varsigma}_T)$. 
We define an onto map 
$\hf{\xi}_T: \hf{Z}_T \rightarrow \hf{X}_T$ and 
a (not necessarily onto) map
$\hf{\eta}_T:\hf{Z}_T \rightarrow \hf{X}_T$ 
to be the restrictions of $\xi$ and $\eta$, respectively. 
If $\hf{\xi}_T$ is 1-1, 
then a (not necessarily onto) endomorphism $\hf{\varphi}_T$ 
of $(\hf{X}_T,\hf{\sigma}_T)$ is defined by 
\[
\hf{\varphi}_T = \hf{\eta}_T\hf{\xi}_T^{-1}. 
\] 
Clearly $\hf{Z}_T\supset Z_T$ and $\hf{X}_T\supset X_T$. 
If $\hf{\eta}_T$ is onto, then $\hf{Z}_T= Z_T$ and 
$\hf{X}_T= X_T$, and hence $\hf{\xi}_T=\xi_T$ and 
$\hf{\eta}_T=\eta_T$.

For a textile system $T$, 
a closed, shift-invariant subset $\hf{U}$ of $\hf{U}_T$ is called
a \itl{half-textile-subsystem of $T$} or simply  called 
a \itl{half-textile-subsystem}. 

Let $\hf{U}$ be a half-textile-subsystem of a textile system $T$. 
We define a subshift $(X_{\hf{U}},\sigma_{\hf{U}})$ by 
$X_{\hf{U}}=\{\hf{\xi}_T((\alpha_{1,j})_{j\in\Z})\bigm| 
(\alpha_{i,j})_{i\in\N,j\in\Z}\in \hf{U}\}$ and 
$\sigma_{\hf{U}}=\hf{\sigma}_T| X_{\hf{U}}$. We also define 
another subshift $(Z_{\hf{U}},\varsigma_{\hf{U}})$ by
$Z_{\hf{U}}=\{(\alpha_{1,j})_{j\in\Z}\bigm| 
(\alpha_{i,j})_{i\in\N,j\in\Z}\in \hf{U}\}$ and 
$\varsigma_{\hf{U}}=\hf{\varsigma}_T| Z_{\hf{U}}$. 
We define an onto map $\xi_{\hf{U}}: Z_{\hf{U}}\to X_{\hf{U}}$ and 
a (not necessarily onto) map 
$\eta_{\hf{U}}: Z_{\hf{U}}\to X_{\hf{U}}$ restricting 
$\hf{\xi}_T$ and $\hf{\eta}_T$, respectively.
If $\hf{U}$ is \itl{onesided 1-1}, i.e., $\xi_{\hf{U}}$ is one-to-one, 
then we have 
a (not necessarily onto) endomorphism 
$\varphi_{\hf{U}}$ of $(X_{\hf{U}},\sigma_{\hf{U}})$
by $\varphi_{\hf{U}}=\eta_{\hf{U}}\xi_{\hf{U}}^{-1}$. 

If $\hf{U}$ is a onesided 1-1 half-textile-subsystem, 
then $Z_{\hf{U}}$ determines $\hf{U}$. 
If $T$ is a textile system, then 
for any subshift space 
$Z\subset \hf{Z}_T$ with $\hf{\xi}_T(Z)\supset\hf{\eta}_T(Z)$ and
$\hf{\xi}_T|Z$ one-to-one, there exists a unique 
onesided 1-1 half-textile-subsystem 
$\hf{U}$ of $T$ with $Z_{\hf{U}}=Z$. 

\begin{remark}  
Suppose that $\hf{U}$ is a onesided 1-1 half-textile-subsystem of 
a textile system $T$
with $\varphi_{\hf{U}}$ onto. Then there exists a unique  
textile-subsystem $U$ of $T$ such that $Z_U=Z_{\hf{U}}$, which 
we call the \itl{extension of $\hf{U}$}. If $U$ is the extension 
of $\hf{U}$ then $\varphi_U=\varphi_{\hf{U}}$.
\end{remark} 
\begin{proof} Suppose that 
$\hf{U}$ is a half-textile-subsystem of $T=(p,q:\G\to G)$
with $\varphi_{\hf{U}}$ onto. Then $\eta_{\hf{U}}$ is onto and so is 
$\Phi_{\hf{U}}=\xi_{\hf{U}}^{-1}\eta_{\hf{U}}$. 
Let $U$ be the set of all $\Phi_{\hf{U}}$-orbits. Then 
$U$ is a closed subset of the product space $A_\G^{\Z\times\Z}$ 
and shift-invariant. If $(z_i)_{i\in\Z}$ is a $\Phi_{\hf{U}}$-orbit 
with $z_i\in Z_{\hf{U}}$, then
$\eta_{\hf{U}}(z_{i-1})=\xi_{\hf{U}}(z_i)$ for all $i\in\Z$, 
hence $\eta(z_{i-1})=\xi(z_i)$ for all $i\in\Z$ with 
$\xi=\phi_p$ and $\eta=\phi_q$,  
and hence $(z_i)_{i\in\Z}\in Z_T$. 
Therefore $U\subset U_T$. 
Hence $U$ is a textile-subsystem of $T$ with $Z_U=Z_{\hf{U}}$. 
Such $U$ is unique for $\hf{U}$ because  $Z_U=Z_{\hf{U}}$. 
Since $Z_U=Z_{\hf{U}}$, 
$\xi_U=\xi|Z_U=\xi|Z_{\hf{U}}=\xi_{\hf{U}}$; similarly 
$\eta_U=\eta_{\hf{U}}$.  Hence $\varphi_U=\varphi_{\hf{U}}$. 
\end{proof}

The following remark is a generalization of \cite[Remark 7.1]{Nasu-te}. 

\begin{remark} 
Let $\varphi$ be an endomorphism of a subshift 
$(X,\sigma)$. Then there exists a onesided 1-1 
half-textile-subsystem $\hf{U}$ with 
$(X_{\hf{U}},\sigma_{\hf{U}},\varphi_{\hf{U}})=(X,\sigma,\varphi)$. 
\end{remark} 
\begin{proof} Suppose that 
$\varphi$ is of $(m,n)$ type and
given by a local rule $f:L_{m+n+1}(X)\to L_1(X)$. 
Let
$T_{f,m,n}=(p,q:\Gamma\to G)$ be the textile system 
defined as follows: 
$G=G_A$ with $A=L_1(X)$; 
if $m=n=0$, then $\Gamma=G$ with $p_A(a)=a$ and 
$q_A(a)=f(a)$ for $a\in A_\Gamma=L_1(X)$; if $m+n\geq 1$, 
then 
$A_\Gamma=L_{m+n+1}(X)$, 
$V_\Gamma=L_{m+n}(X)$, and 
for $w=a_1\dots a_{m+n+1} \in A_\Gamma$ with 
$a_j\in L_1(X)$, $i_\Gamma(w)=a_1\dots a_{m+n}$, 
$t_\Gamma(w)=a_2\dots a_{m+n+1}$, $p_A(w)=a_{m+1}$ and 
$q_A(w)=f(w)$. Then $T_{f,m,n}$ is onesided 1-1. 
Since $X^{[m+n+1]}\subset\hf{Z}_{T_{f,m,n}}$ with 
$\hf{\xi}_{T_{f,m,n}}(X^{[m+n+1]})=X\supset\varphi(X)= 
\hf{\eta}_{T_{f,m,n}}(X^{[m+n+1]})$
and $T_{f,m,n}$ is onesided 1-1,
there exists a unique 
onesided 1-1 half-textile-subsystem 
$\hf{U}_{f,m,n}$ of $T_{f,m,n}$ with 
$Z_{\hf{U}_{f,m,n}}=X^{[m+n+1]}$. 
Clearly $\hf{U}=\hf{U}_{f,m,n}$ has the properties required. 
\end{proof}

A textile system $T=(p,q:\G\to G)$ is said to be 
\itl{$p$-L} if $p$ is left resolving, 
\itl{$p$-R} if $p$ is right resolving, 
\itl{$q$-L} if $q$ is left resolving, and 
\itl{$q$-R} if $q$ is right resolving. 
An onto endomorphism $\varphi$ of a topological 
Markov shift $(X,\sigma)$ is said to be  
\itl{$p$-L} if there exists a 
onesided 1-1, nondegenerate, $p$-L textile system $T$ 
with 
$(X,\sigma,\varphi)=(X_T,\sigma_T,\varphi_T)$. 
We similarly define 
\itl{$p$-R, 
$q$-L} and \itl{$q$-R endomorphisms} 
of a topological Markov shift. 

A textile system is said to be  
\itl{LR} if it is both $p$-L and $q$-R, 
\itl{RL} if it is both $p$-R and $q$-L, 
\itl{LL} if it is both $p$-L and $q$-L, 
\itl{RR} if it is both $p$-R and $q$-R, 
\itl{$p$-biresolving} if it is both $p$-L and $p$-R 
and
\itl{$q$-biresolving} if it is both $q$-L and $q$-R. 
An onto endomorphism $\varphi$ of a topological 
Markov shift $(X,\sigma)$ is said to be \itl{LR} if there exists a 
onesided 1-1, nondegenerate, LR textile system $T$ with  
$(X_T,\sigma_T,\varphi_T)=(X,\sigma,\varphi)$. 
The definitions of \itl{RL, LL, RR}, 
\itl{$p$-biresolving} and \itl{$q$-biresolving 
endomorphisms} of a topological Markov shift 
are similarly given. 
Note that if
a textile system $T=(p,q: \G\to G)$ with $p_V$ and $q_V$ onto 
(or equivalently, with $T^*$ defined over a nondegenerate graph) is 
LR, RL, LL or RR, then $T$ is necessarily nondegenerate. 

We note that by definition,  
only an onto endomorphism of a topological Markov shift 
can be \rt for every resolving term \rt. Here 
we recall that each of 
the ten terms ``$p$-L'', ``$p$-R'', 
``$q$-R'', ``$q$-L'', 
``LR'', ``RL'', ``LL'', ``RR'', ``$p$-biresolving'' 
and ``$q$-biresolving'' is called a \itl{resolving term}. 

A textile system $T=(p,q:\G\to G)$ is said to be  
\itl{weakly $p$-L} 
if $p$ is weakly left resolving. We similarly 
define \itl{weakly $p$-R, weakly $q$-L} and 
\itl{weakly $q$-R} \itl{textile systems}. 
A textile system is said to be \itl{weakly LR}
if it is both weakly $p$-L and 
weakly $q$-R. Similarly \itl{weakly RL}
\itl{weakly LL} and \itl{weakly RR} \itl{textile systems} are 
defined. A textile system is said to be
\itl{weakly $p$-biresolving}, if it is both 
weakly $p$-R and weakly $p$-L, and similarly a 
\itl{weakly $q$-biresolving textile system} is defined. 

For each resolving term \st, a (not necessarily onto) endomorphism 
$\varphi$ of a subshift $(X,\sigma)$ 
is said to be 
\itl{weakly \st} if there 
exists a onesided 1-1 half-textile-subsystem $\hf{U}$ of 
a weakly \rt textile system such that 
$(X_{\hf{U}},\sigma_{\hf{U}},\varphi_{\hf{U}})=(X,\sigma,\varphi)$. 

If an endomorphism $\varphi$ of a subshift $(X,\sigma)$ is 
onto, then we can replace  every 
``half-textile-system $\hf{U}$'' by 
``textile-system $U$'' in the above definitions. In \cite{Nasu-te} 
we defined  the ``weakly resolving'' notions 
only for onto endomorphisms of subshifts.  
However in the above we have defined these notions for not necessarily 
onto endomorphisms of subshifts without losing consistency. 
In fact, using Remark 2.5 we easily see the following. 
\begin{remark} Let \rt be any resolving term. Then 
an onto endomorphism $\varphi$ of a subshift $(X,\sigma)$ 
is weakly \rt if and only if there 
exists a onesided 1-1 textile-subsystem $U$ of 
a weakly \rt textile system such that 
$(X_U,\sigma_U,\varphi_U)=(X,\sigma,\varphi)$. 
\end{remark} 

Here we add a refinement of \cite[Proposition 7.6]{Nasu-te}. 
A version of it for an onto endomorphism of a topological Markov shift 
will appear as Proposition 2.11. 
\begin{proposition} 
Let $\varphi$ be a (not necessarily onto) endomorphism of 
a subshift $(X,\sigma)$. 
\begin{enumerate}
\item 
If $\varphi$ is weakly $p$-L (respectively, weakly $p$-R), 
then $\varphi\sigma$ (respectively, $\varphi\sigma^{-1}$) is 
weakly $p$-L (respectively, weakly $p$-R); 
if in addition $\varphi$ is 
onto, then $\varphi\sigma$ (respectively, $\varphi\sigma^{-1}$)
is right $\sigma$-expansive (respectively, 
left $\sigma$-expansive) on the upper side. 
\item
If $\varphi$ is weakly $q$-R (respectively, weakly $q$-L), 
then $\varphi\sigma$ (respectively, $\varphi\sigma^{-1}$) is 
weakly $q$-R (respectively, weakly $q$-L) and 
positively left $\sigma$-expansive (respectively, 
positively right $\sigma$-expansive). 
\item 
If $\varphi$ is weakly LR (respectively, weakly RL), 
then $\varphi\sigma$ (respectively, $\varphi\sigma^{-1}$) is 
weakly LR (respectively, weakly RL); if in addition 
$\varphi$ is onto then $\varphi\sigma$ (respectively,
$\varphi\sigma^{-1}$) is 
right $\sigma$-expansive  
(respectively, left $\sigma$-expansive) on the upper side 
and positively left $\sigma$-expansive  
(respectively, positively right $\sigma$-expansive)
and hence expansive. 
\end{enumerate}
\end{proposition}
\begin{proof} 
Let $T=(p,q:\G\to G)$ be a textile system. Let 
$T'=(p',q':\G'\to G)$ be the textile system 
such that  
$\G'=\G^{[2]}$ and for 
$\alpha\beta\in A_{\G'}=L_2(X_{\G})$ 
with $\alpha,\beta\in A_\G$, 
$p'_A(\alpha\beta)=p_A(\alpha)$ and
$q'_A(\alpha\beta)=q_A(\beta)$. Then we see that if  $p$ is 
weakly left-resolving then so is $p'$, and that
if $q$ is weakly right-resolving then so is $q'$. 

Suppose that $\hf{U}$ is a onesided 1-1 
half-textile-subsystem of $T$
with $(X_{\hf{U}},\sigma_{\hf{U}},\varphi_{\hf{U}})
=(X,\sigma,\varphi)$. Put $Z_{\hf{U}}^{[2]}=Z$.  Then, 
since $Z$ is 
a subshift-space such that 
$Z\subset \hf{Z}_{T'}$ with 
$\hf{\xi}_{T'}(Z)
\supset\hf{\eta}_{T'}(Z)$ and
$\hf{\xi}_{T'}| Z$ one-to-one, we have a 
onesided 1-1 half-textile-subsystem $\hf{U}'$ of $T'$ 
such that $Z_{\hf{U}'}=Z$, and see that 
$(X_{\hf{U}'},\sigma_{\hf{U}'},\varphi_{\hf{U}'})
=(X,\sigma,\varphi\sigma)$. 

First we prove (2). 

(2) Suppose that $\varphi$ is weakly $q$-R with $q$ 
(in the above) 
weakly right-resolving. 
Then by the above $\varphi\sigma$ is weakly $q$-R. 

Let $(\beta_{i,j})_{i\in\N,j\in\Z}\in\hf{U}'$ with 
$\beta_{i,j}\in A_{\G'}$. Then 
there exists $(\alpha_{i,j})_{i\in\N,j\in\Z}\in\hf{U}$ 
with $\alpha_{i,j}\in A_\G$  such that 
for any $j\in\Z$ 
$(\beta_{i,j})_{i\in\N}=(\alpha_{i,j+i-1}\alpha_{i,j+i})_{i\in\N}$. 
Since $(\alpha_{i,j})_{i\in\N,j\in\Z}$ is a half-textile woven by a 
weakly $q$-R textile system $T$, it follows that $\alpha_{i,j}$ and 
$\alpha_{i+1,j+1}$ uniquely determine $\alpha_{i,j+1}$, which implies 
that $(\beta_{i,j})_{i\in\N}$ uniquely determines 
 $(\beta_{i,j+1})_{i\in\N}$ for $j\in\Z$. Therefore it follows 
that the endomorphism $\Phi_{\hf{U}'}=\xi_{\hf{U}'}^{-1}\eta_{\hf{U}'}$ 
of $(Z_{\hf{U}'},\varsigma_{\hf{U}'})$ is positively 
left $\varsigma_{\hf{U}'}$-expansive, and hence 
the endomorphism $\varphi\sigma$ of $(X,\sigma)$ is positively 
left $\sigma$-expansive (because the two endomorphisms 
are topologically conjugate). Hence (2) is proved. 

(1) Suppose that $\varphi$ is weakly $p$-L with $p$ 
weakly left-resolving. Then by the above $\varphi\sigma$ is 
weakly $p$-L.

To prove the remainder of (1), suppose further that 
$\varphi$ is onto. 
Since $\varphi_{\hf{U}}$ and $\varphi_{\hf{U}'}$ are onto, 
we can consider the textile-subsystem $U$ of $T$ with 
$Z_U=Z_{\hf{U}}$ and the textile-subsystem $U'$ of $T'$ with 
$Z_{U'}=Z_{\hf{U}'}$. 
Let $(\beta_{i,j})_{i,j\in\Z}\in U'$ with 
$\beta_{i,j}\in A_{\G'}$. Then 
there exists $(\alpha_{i,j})_{i,j\in\Z}\in U$ 
with $\alpha_{i,j}\in A_\G$  such that 
for any $j\in\Z$ 
$(\beta_{i,j})_{i\in\Z}=(\alpha_{i,j+i-1}\alpha_{i,j+i})_{i\in\Z}$. 
Since $(\alpha_{i,j})_{i,j\in\Z}$ is a textile woven by a 
weakly $p$-L textile system $T$, it follows that $\alpha_{i-1,j-1}$ and 
$\alpha_{i,j}$ uniquely determine $\alpha_{i,j-1}$, which implies 
that $(\beta_{i,j})_{i\leq 0}$ uniquely determines 
 $(\beta_{i,j-1})_{i\leq 0}$ for $j\in\Z$. Therefore it follows 
that the endomorphism $\Phi_{U'}=\xi_{U'}^{-1}\eta_{U'}$ 
of $(Z_{U'},\varsigma_{U'})$ is  
right $\varsigma_{U'}$-expansive, and hence 
the endomorphism $\varphi\sigma$ of $(X,\sigma)$ is positively 
right $\sigma$-expansive (because the two endomorphisms 
are topologically conjugate). Hence (1) is proved.

(3) Suppose that $\varphi$ is weakly LR with $p$ 
weakly left-resolving and $q$ weakly right-resolving. 
Then by the above $\varphi\sigma$ is weakly LR.

To prove the remainder of (3), suppose further that 
$\varphi$ is onto. Let $U$ and $U'$ be the same as in the 
proof of (1). By the proofs of 
(1) and (2) the endomorphism $\Phi_{U'}=\xi_{U'}^{-1}\eta_{U'}$ 
of $(Z_{U'},\varsigma_{U'})$ is 
left $\varsigma_{U'}$-expansive on the upper side 
and positively right $\varsigma_{U'}$-expansive
with $(Z_{U'},\varsigma_{U'},\Phi_{U'})$ conjugate to 
$(X,\sigma,\varphi\sigma)$. Hence (3) is proved.
\end{proof}

We define a \itl{textile relation-system} to be 
an ordered pair of graph-homomorphisms 
(\cite[p.191]{Nasu-te}). 
Let $T=(p:\G\to G,q:\G\to H)$ be a textile relation-system. 
Let $X_T=\phi_p(X_\G)$, let $Y_T=\phi_q(X_\G)$ and let 
$\xi_T:X_\G\to X_T$ and $\eta_T:X_\G\to Y_T$ be the onto maps 
naturally induced by $\phi_p$ and $\phi_q$, respectively.
$T$ is said to be \itl{onesided 1-1} if $\xi_T$ is one-to-one. 
For a onesided 1-1 textile-relation system $T$, we 
define $\phi_T=\eta_T\xi_T^{-1}$.  

Note that a textile system is 
considered to be a special case of a textile relation-system. 
The definition of  \itl{weakly \rt textile relation system} 
(and its ``non-weakly'' version) should be clear 
for every resolving term \st. 

For every resolving term \st, 
a factor map (onto homomorphism) $\phi$ of an SFT
$(X,\sigma_X)$ onto a sofic system $(Y,\sigma_Y)$ is 
said to be \itl{weakly \st} 
if there exists 
a onesided 1-1, weakly \rt textile relation system 
$T$ such that $\phi_T=\phi$ with $X_T=X$ and $Y_T=Y$. 

We admit the following inconsistency 
in our terminology without confusion: 
for the case where  
\rt is ``LR'', ``RL'', ``$p$-biresolving'' 
or ``$q$-biresolving'', we cannot prove that 
a weakly \st, ``onto endomorphism'' of a sofic system 
 is a weakly \rt ``factor map'' of the sofic system. 
(This can be proved for the case where 
\rt is ``$p$-L'', ``$q$-L'', ``$q$-R'', ``$q$-L'', ``LL'', or 
$RR$.) 
However they are essentially the same, 
more exactly, the same up to higher block conjugacy, 
as seen by the following: 
 
\begin{remark} Let \rt be any resolving term. 
A weakly \rt factor map of an SFT onto itself 
is a weakly \st, onto endomorphism of the SFT, 
and a weakly \rt onto endomorphism of 
an SFT is a weakly \rt factor map of the SFT 
onto itself up to higher-block conjugacy. 
\end{remark} 
\begin{proof} 
Let $\varphi$ be a weakly \rt factor map 
of an SFT $(X,\sigma)$ onto itself. Then there exists 
a onesided 1-1, weakly \rt textile relation-system 
$T=(p:\G\to G,q:\G\to G)$ such that 
$\phi_p(X_\G)=\phi_q(X_\G)=X$ and $\phi_T=\varphi$. 
If we regard $T$ as a weakly \rt textile system, then 
$U_T$ is a textile-subsystem of $T$ 
with $\varphi_{U_T}=\varphi$. 

Suppose that $U$ is a onesided 1-1 textile-subsystem of 
a weakly \rt textile system $T_0$, 
such that $(X_U,\sigma_U)$ is an SFT. Then $(Z_U,\varsigma_U)$ 
is also an SFT, because $U$ is onesided 1-1. 
There exist $n\geq 1$ and 
graphs $G_1$ and $\G_1$
such that $(X_{G_1},\sigma_{G_1})=(X_U^{[n]},\sigma_U^{[n]})$ and
$(X_{\G_1},\sigma_{\G_1})=(Z_U^{[n]},\varsigma_U^{[n]})$. 
Since $T_0^{[n]}$ is weakly \rt and 
$U^{[n]}$ is a textile-subsystem of $T_0^{[n]}$,   
we can define 
a weakly \rt textile relation-system 
$T_1=(p_1:\G_1\to G_1,q_1:\G_1\to G_1)$
such that $\phi_{T_1}=\varphi_{U^{[n]}}=(\varphi_U)^{[n]}$. 
\end{proof} 

\begin{remark} 
Let $\varphi$ be an onto endomorphism of an SFT $(X,\sigma)$. 
\begin{enumerate}
\item  Let \rt be any one of ``$p$-L'' 
``$p$-R'' and  ``$p$-biresolving''.
If $\varphi$ is weakly \st, then 
$\varphi$ is \rt up to higher-block conjugacy. 

\item  Let \rt be any one of 
``$q$-R'', ``$q$-L'', ``LR'', ``RL'', ``LL'', 
``RR'' and ``$q$-biresolving''. 
Suppose that $\varphi$ is one-to-one or $\sigma$ is topologically 
transitive. If $\varphi$ is weakly \st, then 
$\varphi$ is \rt up to higher-block conjugacy. 
\item 
If $(X,\sigma)$ is a topological Markov shift, then 
$\varphi$ is weakly $p$-L (respectively, weakly $p$-R) if and 
only if $\varphi$ is $p$-L (respectively, p-R). 
\end{enumerate}
\end{remark}
\begin{proof} 
(1) We consider the case that \rt is ``$p$-L''. 
Since $\varphi$ is weakly $p$-L,
there exists a onesided 1-1 textile-subsystem $U$ of 
a weakly $p$-L textile system $T$ such that 
$(X_U,\sigma_U,\varphi_U)=(X,\sigma,\varphi)$. 
Since $(X,\sigma)=(X_U,\sigma_U)$ is an SFT, 
$(Z_U,\varsigma_U)$ is an SFT because $U$ is onesided 1-1. 
There exists $n\geq 1$ such that both
$(Z_U^{[n]},\varsigma_U^{[n]})$ and $(X_U^{[n]},\sigma_U^{[n]})$ 
are topological Markov shifts and hence 
there exists a onesided 1-1, nondegenerate 
textile system $T'=(p',q':\G'\to G')$ with $U_{T'}=U^{[n]}$ 
(with $G'$ and $\G'$ nondegenerate). 
Since $T^{[n]}$ is weakly $p$-L, 
so is $T'$. Since $\phi_{p'}$ is 
a conjugacy, it follows from Lemma 2.1(1) that $p'$ is left resolving. 
Since $(X_{T'},\sigma_{T'},\varphi_{T'})
=(X^{[n]},\sigma^{[n]},\varphi^{[n]})$, 
$\varphi$ is $p$-L up to 
higher block conjugacy.  

The proofs for the other cases are similar. 

(2) If $\varphi$ is one-to-one, 
then in the proof of (1),
not only $\phi_{p'}$ but also $\phi_{q'}$ is a conjugacy and 
hence the proof of (2) is similarly given by using Lemma 2.1(1). 

If $\sigma$ is transitive, then 
the proof of (2) is similarly given by using Lemma 2.1(2). 

(3) By \cite[Propositions 5.1(1) and 7.5(1) ]{Nasu-te}. 
\end{proof}

In this paper, we give not only results on 
``weakly resolving'' (i.e. weakly $p$-L, weakly $p$-R, 
weakly $q$-R, weakly $q$-L, weakly LR, weakly RL, 
weakly LL, weakly RR, weakly $p$-biresolving 
and weakly $q$-biresolving) 
endomorphisms of general subshifts, but also 
those on ``resolving'' 
(i.e.  $p$-L, $p$-R, $q$-R, $q$-L,  LR, RL, LL,  RR, 
$p$-biresolving, $q$-biresolving) 
endomorphisms of transitive topological Markov shifts and 
automorphisms of topological Markov shifts.  By using 
Remark 2.10, many of the latter can be 
derived from the former. However, 
since this method is logically very roundabout 
(except for in Sections 11 and 12), 
we shall provide direct proofs
for the latter (except for in Sections 11 and 12).  
In order to do so, 
we need the ``non-weakly'' versions of 
\cite[Propositions 7.6, 7.7 and 7.11]{Nasu-te}.  
The following propositions provide them with some refinements. 

\begin{proposition} 
Let $\varphi$ be an onto endomorphism of a 
topological Markov shift $(X_G,\sigma_G)$. 
\begin{enumerate}
\item  If $\varphi$ is $p$-L (respectively, $p$-R) , then 
$\varphi\sigma_G$ (respectively, $\varphi\sigma_G^{-1}$) is $p$-L 
(respectively, $p$-R) 
and right (respectively, left) $\sigma_G$-expansive 
on the upper side. 
\item  If $\varphi$ is $q$-R (respectively, $q$-L), then 
$\varphi\sigma_G$ (respectively, $\varphi\sigma_G^{-1}$) is 
$q$-R (respectively, $q$-L) and positively 
left $\sigma_G$-expansive (respectively, right $\sigma_G$-expansive). 
\item  If $\varphi$ is LR (respectively, RL), then 
$\varphi\sigma_G$ (respectively, $\varphi\sigma_G^{-1}$) is 
LR (respectively, RL) and 
right $\sigma_G$-expansive  
(respectively, left $\sigma_G$-expansive) on the upper side 
and positively left $\sigma_G$-expansive  
(respectively, positively right $\sigma_G$-expansive)
and hence expansive.  
\end{enumerate}
\end{proposition} 
\begin{proof} 
Suppose that $T=(p,q:\G\to G)$ be
a onesided 1-1, nondegenerate textile system 
with $(X_T,\sigma_T,\varphi_T)=(X_G,\sigma_G,\varphi)$. 
Let $T'=(p',q':\G'\to G)$ be 
the textile system such that  
$\G'=\G^{[2]}$ and for 
$\alpha\beta\in A_{\G'}=L_2(X_{\G})$ with $\alpha,\beta\in A_\G$, 
$p'_A(\alpha\beta)=p_A(\alpha)$ and
$q'_A(\alpha\beta)=q_A(\beta)$. 
Then $T'$ is onesided 1-1, nondegenerate 
with $\varphi_{T'}=\varphi_T\sigma_T=\varphi\sigma_G$. 
We easily see that if $T$ is $p$-L then so is $T'$, 
if  $T$ is $q$-R then so is $T'$, and if $T$ is LR then 
so is $T'$. 

Let $T'^*=(p'^*,q'^*:\G'^*\to G'^*)$ be the dual of $T'$. 
Any point in $X_{G'^*}$ is of the form 
$(\beta_i)_{i\in\Z}$ with $\beta_i\in A_\G$ 
and $q_A(\beta_{i-1})p_A(\beta_{i})\in L_2(X_G)$ 
for all $i\in\Z$. 

Suppose that $p$ is left-resolving. Then 
for each $\beta_0\beta_1\in L_2(G'^*)$ 
with $\beta_0,\beta_1\in A_\G$, 
an arc $\alpha_1\in A_\G$ with
$p_A(\alpha_1)=q_A(\beta_0)$ and 
$t_\G(\alpha_1)=i_\G(\beta_1)$ is unique and hence 
for any $(\alpha_i)_{i\in\Z}, (\beta_i)_{i\in\Z}\in X_{G'^*}$  
such that $(\alpha_i\beta_i)\in L_2(\G)$ for all $i\in\Z$, 
$(\alpha_i)_{i\leq 0}$ is uniquely determined by $(\beta_i)_{i\leq 0}$. 
Hence it follows that $\varphi_{T'}$ 
is right $\sigma_G$-expansive on the upper side. 
By this and the above, (1) is proved. 

Suppose that $q$ is right-resolving. Then 
for each $\alpha_0\alpha_1\in L_2(G'^*)$ 
with $\alpha_0,\alpha_1\in A_\G$,  
an arc $\beta_0\in A_\G$ 
with $t_\G(\alpha_0)=i_\G(\beta_0)$ and 
$q_A(\beta_0)=p_A(\alpha_1)$ 
is unique and hence 
for any $(\alpha_i)_{i\in\Z}, (\beta_i)_{i\in\Z}\in X_{G'^*}$  
such that $(\alpha_i\beta_i)\in L_2(\G)$ for all $i\in\Z$, 
$(\beta_i)_{i\geq 0}$ is uniquely determined by $(\alpha_i)_{i\geq 0}$. 
Hence it follows that $\varphi_{T'}$ 
is positively left $\sigma_G$-expansive. 
By this and the above, (2) is proved. 

Suppose that $p$ is left-resolving and 
$q$ is right-resolving. 
For any $(\alpha_i)_{i\in\Z}, (\beta_i)_{i\in\Z}\in X_{G'^*}$  
such that $(\alpha_i\beta_i)\in L_2(\G)$ for all $i\in\Z$, 
$(\alpha_i)_{i\leq 0}$ is uniquely  
determined by $(\beta_i)_{i\leq 0}$, 
and $(\beta_i)_{i\geq 0}$ is uniquely determined by
 $(\alpha_i)_{i\geq 0}$. 
Hence it follows that $\varphi_{T'}$ 
is right $\sigma_G$expansive on the upper side 
and positively left $\sigma_G$-expansive. 
By this and the above, (3) is proved.   
\end{proof} 

\begin{proposition} 
Let $\varphi$ be an onto endomorphism of 
a topological Markov shift $(X_G,\sigma_G)$. 
\begin{enumerate} 
\item 
If $\varphi$ is $p$-L (respectively, $p$-R) and 
right (respectively, left) $\sigma_G$-expansive,   
then for some $k\geq 1$, 
$\varphi^k\sigma_G^{-1}$ (respectively, $\varphi^k\sigma_G$) is 
$p$-L (respectively, $p$-R). 
\item 
If $\varphi$ is $q$-R (respectively, $q$-L) and 
left (respectively, right) $\sigma_{G}$-expansive, 
then for some $k\geq 1$, 
$\varphi^k \sigma_G^{-1}$ (respectively, $\varphi^k \sigma_{G}$) is 
$q$-R (respectively, $q$-L). 
\item 
If $\varphi$ is LR (respectively, RL) and expansive, 
then for some $k\geq 1$, 
$(\varphi^k \sigma_G^{-1})^{[2]}$ (respectively, 
$(\varphi^k \sigma_G)^{[2]}$) is LR (respectively, RL). 
\end{enumerate} 
\end{proposition}
\begin{proof} 
(1) Assume that $\varphi$ is $p$-L and 
right $\sigma_G$-expansive. Then there exists  
a onesided 1-1, nondegenerate textile system 
$T=(p,q:\G\to G)$ with 
$(X_T,\sigma_T,\varphi_T)=(X_G,\sigma_G,\varphi)$, $p$ 
left-resolving and $\eta_{T^*}$ one-to-one 
(by \cite[Proposition 6.2]{Nasu-te}). 
Let $T^*=(p^*,q^*:\G^*\to G^*)$. Then $q^*$ is 
right resolving and 
the bijective block map $\eta_{T^*}:Z_{T^*}\to X_{T^*}$ is 
the restriction of $\phi_{q^*}$. 
Therefore it follows that $\eta_{T^*}^{-1}$ 
is a block map of
$(k-1,0)$ type for some $k\geq 1$ (see 
\cite[Proof of Proposition 5.1]{Nasu-te}).
Hence for any 
$(\alpha_i)_{1\leq i\leq k}\in L_k(Z_{T^*})$, 
$(q^*_A(\alpha_i))_{1\leq i\leq k}$ uniquely 
determines $\alpha_k$ and hence 
$q_A(\alpha_k)$, and consequently we can define 
a mapping $F_1:L_k(X_{T^*})\to L_1(X_T)$ by 
\begin{align*} F_1((q^*_A(\alpha_i))_{1\leq i\leq k})=q_A(\alpha_k) 
\q\text{with}\; 
(\alpha_i)_{1\leq i\leq k}\in L_k(Z_{T^*}).
\end{align*}
Let $\bar{T}=(\bar{p},\bar{q}:\bar{\G}\to G)$ be the $k$-th power 
of $T$, i.e., the textile system such that 
$\bar{\G}$, $\bar{p}$ and $\bar{q}$ are defined as follows: 
$A_{\bar{\G}}=L_k(Z_{T^*})$, $V_{\bar{\G}}=L_k(X_{T^*})$, each arc 
$\bar{\alpha}=(\alpha_i)_{1\leq i\leq k}$ goes from 
$(p^*_A(\alpha_i))_{1\leq i\leq k}$ to $(q^*_A(\alpha_i))_{1\leq i\leq k}$ 
with $\bar{p}_A(\bar{\alpha})=p_A(\alpha_1)$ and 
$\bar{q}_A(\bar{\alpha})=q_A(\alpha_k)$. 
Then $\varphi_{\bar{T}}=\varphi_T^k=\varphi^k$. Since $T$ is $p$-L, 
it follows that so is $\bar{T}$. 
Let $q':\bar{\G}\to G$ be the graph-homomorphism 
such that $q'_A(\bar{\alpha})=F_1((p^*_A(\alpha_i))_{1\leq i\leq k})$ 
for $\bar{\alpha}=(\alpha_i)_{1\leq i\leq k}\in A_{\bar{\G}}=L_k(Z_{T^*})$ 
(accordingly, for 
$u=(a_i)_{1\leq i\leq k}\in V_{\bar{p}}=V_{q'}=L_k(X_{T^*})$ 
with $a_i\in L_1(X_{T^*})$, 
$q'_V(u)=i_G(F_1(u))$ though $\bar{p}_V(u)=p_V(a_1)$).
Then we have a onesided 1-1, nondegenerate, $p$-L textile system 
$T'=(\bar{p},q':\bar{\G}\to G)$ with 
$\varphi_{T'}=\varphi_{\bar{T}}\sigma_T^{-1}=\varphi^k\sigma_G^{-1}$.  

(2) The proof is given by 
straightforward modifications in the arguments above. 
Instead of $F_1$, we have, for some $k\geq 1$, 
 $F_2:L_k(X_{T^*})\to L_1(X_T)$ such that  
\begin{align*} F_2((p^*_A(\alpha_i))_{1\leq i\leq k})=p_A(\alpha_1) 
\q\text{with}\; 
(\alpha_i)_{1\leq i\leq k}\in L_k(Z_{T^*}).
\end{align*} 
Let $p':\bar{\G}\to G$ be the graph-homomorphism 
such that $p'(\bar{\alpha})=F_2((q^*_A(\alpha_i))_{1\leq i\leq k})$ 
for $\bar{\alpha}=(\alpha_i)_{1\leq i\leq k}\in A_{\bar{\G}}=L_k(Z_{T^*})$.
Then we have a onesided 1-1, nondegenerate, $q$-R textile system 
$T'=(p',\bar{q}:\bar{\G}\to G)$ with 
$\varphi_{T'}=\varphi_{\bar{T}}\sigma_T^{-1}=\varphi^k\sigma_G^{-1}$. 

(3) The proof is also given by straightforward modifications 
in the arguments above. 
We have both $F_1$ and $F_2$ in the above for some $k\geq 1$. 
Let $T''=(p'',q'':\bar{\G}\to G^{[2]})$ be the textile system
such that 
$p''_A(\bar{\alpha})=p_A(\alpha_1)F_2((q^*_A(\alpha_i))_{1\leq i\leq k})$ and 
$q''_A(\alpha')=F_1((p^*_A(\alpha_i))_{1\leq i\leq k})q_A(\alpha_k)$ 
for $\bar{\alpha}=(\alpha_i)_{1\leq i\leq k}$. 
It follows that $T''$ is LR. (Note that for 
$u=(a_i)_{1\leq i\leq k}\in V_{p''}=V_{q''}=L_k(X_{T^*})$ with $a_i\in 
L_1(X_{T^*})$, $p''_V(u)=F_2(u)$ and $q''_V(u)=F_1(u)$). Therefore, since
$T''$ is onesided 1-1 and nondegenerate and 
$\varphi_{T''}=(\varphi_T^k\sigma_T^{-1})^{[2]}
=(\varphi^k\sigma_G^{-1})^{[2]}$, 
(3) is proved. 
\end{proof} 

\begin{proposition} 
Let $\varphi$ be an onto endomorphism of a 
topological Markov shift $(X,\sigma)$. 
\begin{enumerate} 
\item 
If $\varphi$ is essentially $p$-L (respectively, $p$-R)
and right (respectively, left) $\sigma$-expansive, 
then there exists $m\geq 1$ such 
that for all $n\geq m$, $\varphi^n$ is $p$-L (respectively, $p$-R). 
\item 
If $\varphi$ is essentially  $q$-R (respectively, $q$-L)
and left (respectively, right) $\sigma$-expansive,  
then there exists $m\geq 1$ such 
that for all $n\geq m$, $\varphi^n$ is $q$-R (respectively, $q$-L). 
\item 
If $\varphi$ is essentially LR (respectively, RL) 
and expansive, 
then there exists $m\geq 1$ such 
that for all $n\geq m$, $\varphi^n$ is  LR (respectively, RL). 
\end{enumerate}
\end{proposition} 
\begin{proof} 
To prove (1),
assume that $\varphi$ is essentially $p$-L 
and right $\sigma$-expansive. There exist 
a topological Markov shift $(X_1,\sigma_1)$ and a conjugacy 
$\theta_1:(X,\sigma)\to (X_1,\sigma_1)$ such that 
$\varphi_1=\theta_1\varphi\theta_1^{-1}$ is a  
$p$-L endomorphism of $(X_1,\sigma_1)$. 
We may assume 
that $\theta_1^{-1}$ is an LR conjugacy and so is 
$\theta_1\sigma^{l_1}$ with some $l_1\geq 0$, 
because any ``forward conjugacy'' between topological 
Markov shifts is LR (see \cite[p.102]{Nasu-t}) and 
for any conjugacy $\theta:(X,\sigma)\to (X',\sigma')$, 
$\theta\sigma^s$ is forward for sufficiently large $s$.
Since $\varphi_1$ is  
$p$-L and right $\sigma_1$-expansive, it follows from Proposition 2.11 
that there exists $k\geq 1$ such 
that $\varphi_1^k\sigma_1^{-1}$ is 
an essentially $p$-L endomorphism 
of $(X_1,\sigma_1)$. Hence there exist a topological Markov shift 
$(X_2,\sigma_2)$ and a conjugacy 
$\theta_2:(X_1,\sigma_1)\to (X_2,\sigma_2)$ such that 
$\varphi_2=\theta_2(\varphi_1^k\sigma_1^{-1})\theta_2^{-1}$ is 
a $p$-L endomorphism of $(X_2,\sigma_2)$. Also we may 
assume that $\theta_2^{-1}$ is an LR conjugacy and so is 
$\theta_2\sigma_1^{l_2}$ with some $l_2\geq 0$. 

Let $m=k(l_1+l_2)$ and $i$ be any nonnegative integer. 
Then we have 
\[\theta_1^{-1}\theta_2^{-1}\varphi_2^{l_1+l_2}\theta_2\sigma_1^{l_2}
\varphi_1^i\theta_1\sigma^{l_1}=\varphi^{m+i}.\] 
Therefore, 
since each of $\theta_1^{-1}$, $\theta_2^{-1}$, $\varphi_2$, 
$\theta_2\sigma_1^{l_2}$, 
$\varphi_1$ and $\theta_1\sigma^{l_1}$ is $p$-L, (1) follows from 
a trivial generalization of \cite[Fact 3.16]{Nasu-t} to 
conjugacies.

The proofs of (2) and (3) are similar to the above. 
\end{proof} 

Proposition 2.13(3) was given as \cite[Proposition 8.8]{Nasu-t} 
for the case that $\varphi$ is an automorphism. 

A \itl{symbolic endomorphism} of a subshift means 
a 1-block map of the subshift into itself. 
Hence a \itl{symbolic automorphism} is 
a 1-block map of the subshift onto itself 
given by a local rule 
which is just a permutation of the symbols. 

\begin{remark} 
\begin{enumerate} 
\item 
A symbolic endomorphism of a topological Markov shift is 
equal to $\varphi_T$ for 
a onesided 1-1, $p$-biresolving textile system $T$;
a symbolic automorphism of a topological Markov shift 
is equal to $\varphi_T$ for 
a 1-1, $p$-L, $p$-R, $q$-R, $q$-L textile system $T$.
\item 
A symbolic endomorphism of a subshift is equal to 
$\varphi_{\hf{U}}$ for 
a onesided 1-1 half-textile-subsystem $\hf{U}$ of a 
weakly $p$-biresolving textile system; 
a symbolic automorphism of a subshift is equal to $\varphi_U$ for 
a 1-1 textile-subsystem $U$ of a 
weakly $p$-L, weakly $p$-R, weakly $q$-R, weakly $q$-L textile system.
\item 
If $\varphi$ is an onto endomorphism of a finite subshift, then 
it is expansive and for some 
$t\geq 1$, $\varphi^{[t]}$ is
a symbolic automorphism of 
a finite topological Markov shift and hence $\varphi^{[t]}=\varphi_T$
for a 1-1, $p$-L, $p$-R, $q$-R, $q$-L textile system $T$. 
\end{enumerate} 
\end{remark} 
\begin{proof} 
(1) A symbolic endomorphism $\varphi$ of a topological
Markov shift $(X_G,\sigma_G)$ is a 1-block map given 
by a local rule $f$ such that $q_f$ is a graph-endomorphism of $G$. 
It is equal to $\varphi_T$ for the textile system 
$T=(\imath_G,q_f:G\to G)$, where $\imath_G$ is 
the identity graph-automorphism of $G$. Clearly, 
$T$ is onesided 1-1, $p$-L and $p$-R. 
If $\varphi$ is onto, then $q_f$ is a graph-automorphism of $G$
and hence $T$ is 1-1, $p$-L, $p$-R, $q$-R and $q$-L. 

(2) A symbolic endomorphism $\varphi$ of 
a subshift $(X,\sigma)$ with $A=L_1(X)$ 
is a 1-block map given by a local rule $f:A\to A$. Let 
$T=(p,q:G_A\to G_A)$ be the textile system such that $p(a)=a$ and 
$q(a)=f(a)$ for $a\in A_{G_A}=A$. Then $T$ is onesided 1-1, 
$p$-L and $p$-R. There 
exists a (onesided 1-1) half-textile-subsystem $\hf{U}$ of $T$ such that 
$Z_{\hf{U}}=X_{\hf{U}}=X$ and $\varphi_{\hf{U}}=\varphi$. 
If $\varphi$ is onto, then $f$ is a permutation on $A$ 
and hence $T$ is 1-1, $p$-L, $p$-R, $q$-R and $q$-L, 
and there exists a (1-1) textile-subsystem $U$ of $T$ such that 
$Z_U=X_U=X$ and $\varphi_U=\varphi$. 

(3) A finite subshift is, up to higher block conjugacy, 
a topological Markov shift. Hence for some $t\geq 1$ 
$\varphi^{[t]}$ is an onto endomorphism of a finite topological 
Markov shift, which is a symbolic automorphism of 
the shift and is expansive. 
\end{proof} 

Finally we make some remarks on the definitions of 
the composition $T_1\circ T_2$ and the product $T_1T_2$ of 
textile systems $T_1$ and $T_2$ given in \cite{Nasu-t} 
and give some properties of them, which will be needed later.

According to the simplification 
of the definition of a textile system made in \cite[p.172]{Nasu-te}
(with no influence on the statements of the results of 
\cite{Nasu-note,Nasu-t,Nasu-m,Nasu-d,Nasu-n}), 
we have the following simplified definitions of them. 

For textile systems $T_1=(p_1,q_1:\G_1\to G)$ 
and $T_2=(p_2,q_2:\G_2\to G)$ 
defined over a common graph $G$, we 
define the \itl{composition} $T_1\circ T_2$  
to be the textile system 
$\Hat{T}=(\Hat{p},\Hat{q}:\Hat{\Gamma}\to G)$, 
where the arcs of $\Hat{\Gamma}$ are 
the elements of the form 
$\alpha/\beta$
with $\alpha\in A_{\Gamma_1}, \beta\in A_{\Gamma_2}$ and 
$(q_1)_A(\alpha)=(p_2)_A(\beta)$, and for each 
$\gamma
=\alpha/\beta\in A_{\Hat{\Gamma}}$, 
$i_{\Hat{\Gamma}}(\gamma)
=i_{\Gamma_1}(\alpha)/ i_{\Gamma_2}(\beta), 
t_{\Hat{\Gamma}}(\gamma)
=t_{\Gamma_1}(\alpha)/ t_{\Gamma_2}(\beta), 
 \Hat{p}_{A}(\gamma)=(p_1)_A(\alpha)$ and 
$\Hat{q}_{A}(\gamma)=(q_2)_A(\beta)$. 
Then we easily see the following facts: 
if $T_1$ and $T_2$ are onesided 1-1 then 
so is $T_1\circ T_2$ and 
$\varphi_{T_1\circ T_2}=\varphi_{T_2}\varphi_{T_1}$ 
\cite[p.56]{Nasu-t}; if $T_1$ and $T_2$ are nondegenerate then 
so is $T_1\circ T_2$; for every resolving term \st, 
if $T_1$ and $T_2$ are \rt (respectively, weakly \st) 
then so is $T_1\circ T_2$  
(\cite[Fact 3.16, Corollary 3.17]{Nasu-t} with obvious 
supplements). 

For any textile system $T$ over a graph $G$ and $n\geq 0$,
we define the \itl{$n$-th composition power} $T^n$ of $T$ 
by the following:
$T^0=(\imath_G,\imath_G:G\to G)$, 
where $\imath_G$ is the identity graph-automorphism of $G$, 
and $T^n=T^{n-1}\circ T$ for $n\geq 1$. 

Let $G_1$ and $G_2$ be graphs with $V_{G_1}=V_{G_2}$. 
Then we define the \itl{product $G_1G_2$ of $G_1$ and $G_2$} 
to be the graph such that $A_{G_1G_2}=
\{ab\,|\,a\in A_{G_1}, b\in A_{G_2}, t_{G_1}(a)=i_{G_2}(b)\}$, 
$V_{G_1G_2}=V_{G_1}$ and $i_{G_1G_2}(ab)=i_{G_1}(a)$ and 
$t_{G_1G_2}(ab)=t_{G_2}(b)$ 
for $ab\in A_{G_1G_2}$. Clearly we have $M_{G_1G_2}=M_{G_1}M_{G_2}$. 

For textile systems $T_1=(p_1,q_1:\G_1\to G_1)$ 
and $T_2=(p_2,q_2:\G_2\to G_2)$ such that $T_1^*$ and $T_2^*$ are 
defined over a common graph 
we define the \itl{product} 
$T_1 T_2$ to be the textile system 
$T_1T_2=(p,q:\G_1\G_2\to G_1G_2)$ such that 
$p_A(\alpha_1\alpha_2)=(p_1)_A(\alpha_1)(p_2)_A(\alpha_2)$ 
and $q_A(\alpha_1\alpha_2)=(q_1)_A(\alpha_1)(q_2)_A(\alpha_2)$ 
for $\alpha_1\alpha_2\in A_{\G_1\G_2}$ with $\alpha_1\in A_{\G_1}$
$\alpha_2\in A_{\G_2}$. For any textile system $T$  
and for $s\geq 1$,
we define the \itl{$s$-th product power} $T^{(s)}$ of $T$ by the 
following: $T^{(1)}=T$ and $T^{(s)}=T^{(s-1)} T$ for $s\geq 2$. 

If $T_1$ and $T_2$ are textile systems 
with $T_1^*$ and $T_2^*$  
defined over a common graph, 
then $T_1T_2=(T_1^*\circ T_2^*)^*$. 
For any textile and $s\geq 1$, $T^{(s)}=((T^*)^s)^*$; if 
$T$ is nondegenerate then so is $T^{(s)}$; for every resolving 
term \st, if $T$ is \rt (respectively, weakly \st) 
then so is $T^{(s)}$; if $T$ is onesided 1-1 then so is 
$T^{(s)}$ and $\varphi_{T^{(s)}}=\varphi^{(s)}$, where 
the definition of $\varphi^{(s)}$ is given below.

For a subshift $(X,\sigma_X)$ and an integer $s\geq 1$, 
we define the \itl{$s$-th power system} 
of $(X,\sigma_X)$
to be the subshift $(X^{(s)},\sigma_{X^{(s)}})$ over the 
alphabet $L_s(X)$ with 
$X^{(s)}=\{(a_{(j-1)s+1}\dots a_{js})_{j\in\Z}\,\bigm|
\, (a_j)_{j\in\Z}\in X,\;a_j\in L_1(X)\}$.

For a homomorphism $\phi:(X,\sigma_X)\to (Y,\sigma_Y)$ 
between subshifts and $s\geq 1$, define 
$(\phi)^{(s)}:(X^{(s)},\sigma_{X^{(s)}})\to (Y^{(s)},\sigma_{Y^{(s)}})$ 
to be the homomorphism  such that 
if $\phi$ maps $\aseq$ to $\bseq$ 
with $a_j,b_j\in L_1(X)$, then 
$(\phi)^{(s)}$ maps $(a_{(j-1)s+1}\dots a_{js})_{j\in\Z}$
to $(b_{(j-1)s+1}\dots b_{js})_{j\in\Z}$ (see \cite[p.12]{Nasu-t}). 

In the following remark, (1) 
is \cite[Corollary 3.18]{Nasu-t} with 
obvious supplements, and (2) is a slight extension of 
\cite[Remark 7.9]{Nasu-te}. 

\begin{remark} Let \rt be any resolving term.
\begin{enumerate} 
\item[(1)]  
If $\varphi_1$ and $\varphi_2$ are \rt 
endomorphisms of a topological Markov shift $(X,\sigma)$, then so is
$\varphi_2\varphi_1$. 
If $\varphi$ is a \rt 
endomorphism of a topological Markov shift $(X,\sigma)$, then 
$(\varphi^r)^{(s)}$ is a  \rt 
endomorphism of $(X^{(s)},\sigma^{(s)})$, 
and hence $\varphi^r$ is an essentially  \rt endomorphism 
of $(X,\sigma^s)$.
\item[(2)]
If $\varphi_1$ and $\varphi_2$ are weakly \rt 
endomorphisms of a subshift $(X,\sigma)$, then so is
$\varphi_2\varphi_1$. 
If $\varphi$ is a weakly \rt 
endomorphism of a subshift $(X,\sigma)$, then 
$(\varphi^r)^{(s)}$ is a weakly \rt 
endomorphism of a subshift $(X^{(s)},\sigma^{(s)})$, 
and hence $\varphi^r$ is an essentially weakly \rt endomorphism 
of $(X,\sigma^s)$.
\end{enumerate}
\end{remark} 
\begin{proof}  
(1) By the facts on the composition and product 
of textile systems above. 

(2) Suppose that $\varphi_i=\varphi_{\hf{U}_i}$ for a onesided 1-1
half-textile-subsystem $\hf{U}_i$ of a weakly \rt textile 
system $T_i$ for $i=1,2$. Let $\Hat{T}=T_1\circ T_2$. 
Then $\Hat{T}$ is weakly \rt.

Let $Z$ be the the set of 
all bi-sequences $(\alpha_j/\beta_j)_{j\in\Z}$ 
with $\alpha_j\in A_{\G_1}$ and $\beta_j\in A_{\G_2}$ such that 
$(\alpha_j)_{j\in\Z}\in X_{\hf{U}_1}$, 
$(\beta_j)_{j\in\Z}\in X_{\hf{U}_2}$ and
$\eta_{\hf{U}_1}((\alpha_j)_{j\in\Z})=
\xi_{\hf{U}_2}((\beta_j)_{j\in\Z})$. 
Then, since $Z$ is a subshift space 
such that  $Z\subset\hf{Z}_{\Hat{T}}$ and 
$\hf{\xi}_{\Hat{T}}(Z)\supset\hf{\eta}_{\Hat{T}}(Z)$ with
$\hf{\xi}_{\Hat{T}}|Z$ one-to-one, there exists a unique 
onesided 1-1 half-textile-subsystem 
$\hf{U}$ of $\Hat{T}$ with $Z_{\hf{U}}=Z$, and 
we have $\varphi_{\hf{U}}=\varphi_2\varphi_1$. 
Hence $\varphi_2\varphi_1$ is weakly \rt. 

The proof of the remainder is straightforwardly given by using 
the facts on the composition and product of textile systems above, 
and hence omitted.
\end{proof}

\begin{proposition} 
Let $T_1$ and $T_2$ be onesided 1-1, 
nondegenerate LR textile systems over 
a common (nondegenerate) graph, and suppose that $T_1^*$ and $T_2^*$
are defined over graphs $H_1$ and $H_2$, respectively. 
\begin{enumerate}
\item 
$T_1\circ T_2$ is onesided 1-1, 
nondegenerate and LR, and 
its dual $(T_1\circ T_2)^*$ 
is LR, defined over the graph $H_1H_2$, and nondegenerate.   
\item 
If $\varphi_{T_1}\varphi_{T_2}=\varphi_{T_2}\varphi_{T_1}$
then $M_{H_1}M_{H_2}=M_{H_2}M_{H_1}$.
\end{enumerate}
\end{proposition}
\begin{proof} 
(1) It is straightforward to check that 
$T_1\circ T_2$ is onesided 1-1, 
nondegenerate and LR, and 
its dual $(T_1\circ T_2)^*$ 
is LR and defined over the graph $H_1H_2$. 

Let $T_i=(p_i,q_i:\G_i\to G)$ for $i=1,2$. Since 
$T_i$ is nondegenerate, $G$ is nondegenerate and  
$(p_i)_V, (q_V)_i$ are onto.  This implies that  
$H_i$ is nondegenerate and hence $H_1H_2$ is 
a nondegenerate graph. 
Since $(T_1\circ T_2)^*$ is LR and 
its dual  $T_1\circ T_2$  is defined 
over a nondegenerate graph $G$,  
$(T_1\circ T_2)^*$ are nondegenerate.   

(2) For an LR endomorphism $\varphi$ of a topological Markov shift 
$(X_G,\sigma_G)$, a onesided 1-1, nondegenerate LR textile system $T$ 
with $\varphi_T=\varphi$ is unique \cite[Corollary 7.25]{Nasu-t};
more precisely, the ``representation $\lambda$-graph $\cG_T$'' of $T$ is 
the ``canonical $\lambda$-graph'' of the SFT  whose space is 
$\{(a_j/b_j)_{j\in\Z}\,|\, (a_j)_{j\in\Z} \in X_G, (b_j)_{j\in\Z}=
\varphi((a_j)_{j\in\Z})\}$ (see 
\cite[Proof of Corollary 7.25]{Nasu-t}). This implies that 
if $T=(p,q:\G\to G)$ and $T'=(p',q':\G'\to G)$ are 
onesided 1-1, nondegenerate LR textile systems with 
$\varphi_T=\varphi_{T'}$, 
then there exists a graph-isomorphism 
$k:\G\to \G'$ such that $p=p'k$ and $q=q'k$, which implies that 
for each vertex $u$ in $\G$, 
$p_V(u)=p'_V(k_V(u))$ and $q_V(u)=q'_V(k_V(u))$. Therefore, 
if $T^*$ and $T'^*$ is defined over graphs $H$, $H'$ 
respectively, then $V_H=V_{H'}=V_G$ and we have 
a graph-isomorphism $l:H\to H'$ such that $l_A=k_V$ and 
$l_V:V_H\to V_{H'}$ is the identity mapping of $V_G$. 
Hence $M_H=M_{H'}$ with indexing set $V_G\times V_G$.

By (1), $(T_1\circ T_2)^*$ is defined over $H_1H_2$ 
and $(T_2\circ T_1)^*$ is defined over $H_2H_1$. 
Since $\varphi_{T_1\circ T_2}=
\varphi_{T_2}\varphi_{T_1}=\varphi_{T_1}\varphi_{T_2} 
=\varphi_{T_2\circ T_1}$, it follows from the above that 
$M_{H_1H_2}=M_{H_2H_1}$. Hence $M_{H_1}M_{H_2}=M_{H_2}M_{H_1}$. 
\end{proof}

\section{$q$-R and $q$-L degrees}
Let $(X,\sigma_X)$ and $(Y,\sigma_Y)$ be subshifts. 
Let $f:L_{N+1}(X)\to L_1(Y)$ be a local rule 
on $(X,\sigma_X)$ to $(Y,\sigma_Y)$. 
Let $k$ be a nonnegative integer.
We say that $f$ is \itl{$k$ right-mergible} 
if for any points  
$\aseq$ and 
$\bseq$ in $X$ with 
$a_j,b_j\in L_1(X)$, it holds that if
$(a_j)_{j\leq 0}=(b_j)_{j\leq 0}$
and 
$f(a_{j-N}\dots a_{j})
=f(b_{j-N}\dots b_{j})$ for $j=1,\dots, k+1$,  
then $a_1=b_1$. 
We say that $f$ is \itl{$k$ left-mergible} 
if for any points  
$\aseq$ and 
$\bseq$ in $X$ with 
$a_j,b_j\in L_1(X)$, it holds that if
$(a_j)_{j\geq 0}=(b_j)_{j\geq 0}$
and 
$f(a_{-j}\dots a_{-j+N})
=f(b_{-j}\dots b_{-j+N})$ for $j=1,\dots,k+1$,  
then $a_{-1}=b_{-1}$.

We say that $f$ is \itl{strictly $0$ right-mergible} 
(respectively,  \itl{strictly $0$ left-mergible}) if 
$f$ is $0$ right-mergible (respectively, $0$ left-mergible).
If $k\geq 1$, then we say that 
$f$ is \itl{strictly $k$ right-mergible} 
(respectively,  \itl{strictly $k$ left-mergible}) if 
$f$ is $k$ right-mergible but not $k-1$ right-mergible
(respectively, $k$ left-mergible but not $k-1$ left-mergible). 
If
$f$ is not $k$ right-mergible 
(respectively, not $k$ left-mergible) for any $k\geq 0$, 
then we say that 
$f$ is \itl{strictly $\infty$ right-mergible} 
(respectively, \itl{strictly $\infty$ left-mergible}). 

(Note that we have never defined above ``$\infty$ right-mergible 
local rules'' and  ``$\infty$ left-mergible 
local rules''.)

Note that any local rule $f$ on a finite subshift to a subshift 
is zero right-mergible and  zero left-mergible. 

The notion of being $k$ right-mergible (left-mergible) for 
local rules defined above is a generalization of that 
of the ``non-existence of  right (respectively, left) 
$f$-branch of length exceeding $k$'' for local rules of the
endomorphisms of the full shifts
in \cite[Section 16]{Hedlund},  notions corresponding to 
which were applied to 
cellular automata in \cite{Nasu-l}, to 
graph-homomorphisms in \cite{Nasu-c} and 
to $\lambda$-graphs in \cite{Nasu-t-ex}.   

A homomorphism $\phi:(X,\sigma_X)\to (Y,\sigma_Y)$ 
between subshifts is said to be 
\itl{right-closing} (respectively, \itl{left-closing}) 
if it never collapses 
distinct left (respectively, right) 
$\sigma_X$-asymptotic points 
(Kitchens \cite{Kit-c}).
It is well known and easily seen that 
a homomorphism $\phi$ between subshifts is 
right-closing (respectively, left-closing) 
if and only if $\phi$ is given by a 
$k$ right-mergible (respectively, $k$ left-mergible) 
local rule for some $k\geq 0$. 
A local rule $f$ does not give a right-closing 
(respectively, left-closing) homomorphism 
if and only if $f$ is 
strictly $\infty$ right-mergible
(respectively, strictly $\infty$ left-mergible). 
\begin{definition}
Let $m,n\geq 0$ and $N=m+n$. Let $k,l\in\{0,1,2,\dots,\infty\}$. 
Let $\phi: (X,\sigma_X)\to (Y,\sigma_Y)$ be 
a homomorphism of between subshifts. Suppose that $\phi$ is of  
$(m,n)$-type and given by a local rule $f:L_{N+1}(X)\to L_1(Y)$ 
which is strictly $k$ right-mergible and strictly $l$ left-mergible.  
Define the \itl{$q$-R degree $Q_R(\phi)$ of $\phi$} 
and the \itl{$q$-L degree $Q_L(\phi)$ of $\phi$} as follows:
if $X$ is infinite, then 
\[Q_R(\phi)=n-k\q\text{and}\q Q_L(\phi)=m-l.\] 
if $X$ is finite, then define $Q_R(\phi)=\infty$ and 
$Q_L(\phi)=\infty$. 
\end{definition} 

We follow the usual convention about the  
arithmetic calculations and inequalities containing $\pm\infty$. 
Hence, if $\phi$ is not right-closing (respectively, 
not left-closing) if and only if 
$Q_R(\phi)=-\infty$ (respectively, $Q_L(\phi)=-\infty$). 
For example, if $Q_R(\phi)=-\infty$, 
then $Q_R(\phi)<I$ for all $I\in\Z$. 

We must prove that $Q_R(\phi)$ and $Q_L(\phi)$ are uniquely 
determined by $\phi$
for any homomorphism $\phi$ between subshifts. 
To do this we need the following lemma. 

Let $N,s,t\geq 0$. 
Let $f:L_{N+1}(X)\to L_1(Y)$ be a local rule on a 
subshift $(X,\sigma_X)$ to another $(Y,\sigma_Y)$.  Let 
$g:L_{N+s+t+1}(X)\to L_1(Y)$ 
be the mapping such that for $a_{-t}\dots a_{N+s}\in L_{N+s+t+1}(X)$ 
with $a_j\in L_1(X)$, 
\[g(a_{-t}a_{-t+1}\dots a_{N+s})
=f(a_0\dots a_N).\]  Then, clearly $g$ is a local rule 
on $(X,\sigma_X)$ to $(Y,\sigma_Y)$.  
We call $g$ the \itl{local rule obtained from $f$ by adding 
right redundancy by $s$ and left redundancy by $t$}. 

\begin{lemma} 
Let $N,s,t\geq 0$. Let
$f:L_{N+1}(X)\to L_1(Y)$ be a local rule on  
a subshift $(X,\sigma_X)$ to another $(Y,\sigma_Y)$ with $X$ infinite. 
Let $g$ be the local rule obtained from $f$ by adding 
right redundancy by $s$ and left redundancy by $t$.  
Then, if $f$ is strictly $k$ right-mergible 
and strictly $l$ left-mergible then
$g$ is strictly $k+s$ right-mergible and 
strictly $l+t$ left-mergible. 
\end{lemma}
\begin{proof}  
Suppose that $f$ is strictly $k$ right-mergible. If $k=\infty$ then 
it is clear that $g$ is strictly $k+s$ right-mergible. 
Hence we assume that $k\in\Z$. 

Suppose that 
$\aseq,\bseq$ are points of $X$ with 
$a_j,b_j\in L_1(X)$ such that 
$(a_j)_{j\leq 0}=(b_j)_{j\leq 0}$ and 
\[g(a_{j-s-t-N}\dots a_{j})
=g(b_{j-s-t-N}\dots b_{j})\q\text{for}\; j=1,\dots, k+s+1.\]
Then, by the definition of $g$ we have 
\[f(a_{j-s-N}\dots a_{j-s})=f(b_{j-s-N}\dots b_{j-s})\]
for $j=1,\dots,k+s+1$ and hence for $j=s+1,\dots, k+s+1$. 
Hence 
\[f(a_{j-N}\dots a_{j})=f(b_{j-N}\dots b_{j}) 
\q\text{for}\; j=1,\dots,k+1.\]
Therefore, since $f$ is
$k$ right-mergible, it follows that $a_1=b_1$. This implies that
$g$ is $k+s$ right-mergible. 

To see that $g$ is strictly $k+s$ right-mergible, 
suppose first that $k=0$. 
If $s=0$ then $g$ is $0$ right-mergible by the above and hence 
strictly $k+s$ right-mergible with $k+s=0$ by definition. Therefore
we assume that $s\geq 1$.
Since $(X,\sigma_X)$ is an infinite subshift, 
there 
exist
$\aseq,\bseq\in X$ such that 
$(a_j)_{j\leq 0}=(b_j)_{j\leq 0}$ and $a_1\neq b_1$. 
Since $a_j=b_j$ for $j\leq 0$, for $j=1,\dots,s$ 
\begin{align*}
g(a_{j-s-t-N}\dots a_{j})&=f(a_{j-s-N}\dots a_{j-s})\\
&=f(b_{j-s-N}\dots b_{j-s})=g(b_{j-s-t-N}\dots b_{j}).
\end{align*}
Therefore $g$ is not $s-1$ right-mergible. 

Suppose next that $k\geq 1$. 
Since $f$ is not $k-1$ right-mergible, there 
exist points 
$\aseq,\bseq\in X$ such that 
$(a_j)_{j\leq 0}=(b_j)_{j\leq 0}$, $a_1\neq b_1$ and 
\[f(a_{j-N}\dots a_{j})
=f(b_{j-N}\dots b_{j})\q\text{for}\; j=1,\dots, k.\]
Therefore, since
$a_j=b_j$ for $j\leq 0$, 
it follows that 
\[f(a_{j-N}\dots a_{j})
=f(b_{j-N}\dots b_{j})\q\text{for}\; j=-s+1,\dots, k.\]
Hence we have 
\[g(a_{j-t-N}\dots a_{j+s})=g(b_{j-t-N}\dots b_{j+s})
\q\text{for}\; j=-s+1,\dots, k,\] 
and hence we have 
\[g(a_{j-s-t-N}\dots a_{j})=g(b_{j-s-t-N}\dots b_{j})
\q\text{for}\; j=1,\dots, k+s.\]
Therefore, since $(a_j)_{j\leq 0}=(b_j)_{j\leq 0}$ and $a_1\neq b_1$,  
$g$ is not $k+s-1$ right-mergible. 

Hence we have proved that if $f$ is strictly $k$ right-mergible 
then $g$ is strictly $k+s$ right-mergible. 

By symmetry (i.e. by reversing the direction of the shift), 
it follows from the result proved above that 
if $f$ is strictly $l$ mergible then 
$g$ is strictly $l+t$ left-mergible. 
\end{proof} 

\begin{proposition}
Let $\phi$ be a homomorphism between subshifts. 
Then $Q_R(\phi)$ and $Q_L(\phi)$ are uniquely determined by $\phi$. 
\end{proposition} 
\begin{proof}   
Suppose that for $i=1,2$,  
$\phi$ is a homomorphism 
of a subshift $(X,\sigma_X)$ into another $(Y,\sigma_Y)$. 
We may assume that $X$ is infinite. Suppose that for $i=1,2$,
$\phi_i$ is of $(m_i,n_i)$-type and 
given by a strictly $k_i$ right-mergible 
local-rule $f_i:L_{m_i+n_i+1}(X)\to L_1(Y)$. 
Let $m=\max\{m_1,m_2\}$ and $n=\max\{n_1,n_2\}$. For $i=1,2$, let 
$g_i:L_{m+n}(X)\to L_1(X)$ be the local rule such that
\[g_i(a_{-m}\dots a_n)=f_i(a_{-m_i}\dots a_{n_i}),\q 
a_{-m}\dots a_n\in L_{m+n+1}(X),\,\> a_j\in L_1(X).\]
Then, since $\phi$ is of $(m_i,n_i)$-type and given by $f_i$, 
$\phi$ is of $(m,n)$-type and given by $g_i$, for $i=1,2$. 
Since $f_i$ is strictly $k_i$ right-mergible,  
it follows from Lemma 3.2 that $g_i$ is 
strictly $k_i+(n-n_i)$ right-mergible for $i=1,2$. 
Since $g_1:L_{m+n+1}(X)\to L_1(Y)$ and $g_2:L_{m+n+1}(X)\to L_1(Y)$
give the same block-map $\phi$ of $(m,n)$-type, 
it follows that $g_1=g_2$. Hence 
$k_1+(n-n_1)=k_2+(n-n_2)$, and hence $n_1-k_1=n_2-k_2$. 

Similarly, it is proved that if $f_i$ is strictly $l_i$ left-mergible 
for i=1,2, then $m_1-l_1=m_2-l_2$. 
\end{proof} 

Here we give a little more compact description 
of the definition of
the $q$-R, $q$-L degrees of a homomorphism 
$\phi:(X,\sigma_X)\to (Y,\sigma_Y)$ between subshifts.

Let $f:L_{N+1}(X)\to L_1(Y)$ be a local rule 
on $(X,\sigma_X)$ to $(Y,\sigma_Y)$. 
A pair of points $x,y$ in $X$ such that 
$x=\aseq, y=\bseq$ with $a_j,b_j\in L_1(X)$, 
$(a_j)_{j\leq 0}=(b_j)_{j\leq 0}$ and $a_1\neq b_1$ 
(respectively, $(a_j)_{j\geq 0}=(b_j)_{j\geq 0}$ and $a_{-1}\neq b_{-1}$) 
will be called a  
\itl{standard pair of left} (respectively, 
\itl{right}) \itl{$\sigma_X$-asymptotic points}. 
For any such pair of 
$x, y$, 
define $\mu_{f,R}(x,y)$ (respectively, $\mu_{f,L}(x,y)$)
to be the supremum of the integers $s\geq 0$ 
such that for $j=0,\dots,s$, 
$f(a_{j-N}\dots a_j)=f(b_{j-N}\dots b_j)$  
(respectively, $f(a_{-j}\dots a_{-j+N})=f(b_{-j}\dots b_{-j+N})$).  
Define $\mu_{f,R}=\sup\mu_{f,R}(x,y) $ 
(respectively, $\mu_{f,L}=\sup\mu_{f,L}(x,y)$), where 
$(x,y)$ runs over all standard pairs of  
left (respectively, right) $\sigma_X$-asymptotic points. 
We call $\mu_{f,R}$ (respectively, 
$\mu_{f,L}$) the \itl{maximum length of right mergibility} 
(respectively, the \itl{maximum length of left mergibility})
of $f$. It is clear that 
$f$ is strictly $k$ right-mergible
(respectively, strictly $k$ left-mergible if and only if 
$\mu_{f,R}=k$ (respectively, $\mu_{f,L}=k$). 

We can define 
$Q_R(\phi)$ (respectively, $Q_L(\phi)$) 
to be the supremum of $n-\mu_{f,R}$
(respectively, $m-\mu_{f,L}$) 
such that $\phi$ is of $(m,n)$-type and given by a
local rule $f:L_{m+n+1}(X)\to L_1(Y)$. If $X$ is finite, then
$\mu_{f,R}=\mu_{f,L}=0$ (because $f$ is strictly
$0$ right-mergible and strictly $0$ left-mergible), 
and hence $Q_R(\phi)=\infty$ and $Q_L(\phi)=\infty$. 

\textbf{Standing convention II for the description of proofs.} 
Hereafter, we shall mention nothing on proofs for 
the case when the domains of homomorphisms 
(including endomorphisms) between subshifts 
appearing in the results are finite, 
whenever the proofs are trivial or given by Remark 2.14(3). 

\begin{proposition}
Let $\phi:(X,\sigma_X)\to (Y,\sigma_Y)$ be a homomorphism 
between subshifts. Let $s\in\Z$. Then 
\[Q_R(\phi\sigma_X^s)=Q_R(\phi)+s,\q\q 
Q_L(\phi\sigma_X^s)=Q_L(\phi)-s\] 
and hence $Q_R(\phi)+Q_L(\phi)$ is shift-invariant. 
\end{proposition}
\begin{proof}
Suppose that $\phi$ is of $(m,n)$-type given by a local rule 
$f:L_{m+n+1}(X)\to L_1(Y)$. Increasing the redundancies of $f$ 
as in Lemma 3.2 if necessary, we may assume that 
$m, n\geq |s|$, by Proposition 3.3. 
Since $\phi\sigma_X^s$ is of $(m-s,n+s)$-type and 
given by $f$, the proposition follows. 
\end{proof} 

Let $(X,\sigma)$ be a subshift and let
$i_X$ denote the identity automorphism of the shift.
If $(X,\sigma)$ is an infinite subshift, then by definition
\begin{equation} 
Q_R(i_X)= 0\q\q\text{and}\q\q Q_L(i_X)=0,\end{equation}
hence by Proposition 3.4, for all $s\in\Z$
\begin{equation}
Q_R(\sigma^s)= s\q\q\text{and}\q\q Q_L(\sigma^s)=-s, 
\end{equation} 
and for any homomorphism $\phi$ of $(m,n)$-type of $(X,\sigma)$ 
into another subshift, 
\begin{equation}
Q_R(\phi)\leq n\q\q\text{and}\q\q Q_L(\phi)\leq m. 
\end{equation} 
However if $(X,\sigma)$ is finite, then for any 
homomorphism $\phi$ of $(X,\sigma)$ into another subshift, 
 $Q_R(\phi)=\infty$ 
and $Q_L(\phi)=\infty$ by definition. 

Let $N\geq 0$ and $s\geq 1$. 
Let $f:L_{N+1}(X)\to L_1(Y)$ be a local rule  
on a subshift $(X,\sigma_X)$ to another $(Y,\sigma_Y)$. We define 
the \itl{higher block presentation 
$f^{[s]}:L_{N+1}(X^{[s]})\to L_1(Y^{[s]})$  
of order $s$ of $f$} by 
\begin{align*}
f^{[s]}&((a_1\dots a_s)(a_2\dots a_{1+s})\dots (a_{N+1}\dots a_{N+s}))\\
&=f(a_1\dots a_{N+1})f(a_2\dots a_{N+2})\dots f(a_s\dots a_{N+s}),
\q\q a_j\in L_1(X).
\end{align*}
Clearly, $f^{[s]}$ is a local rule on $(X^{[s]},\sigma_{X^{[s]}})$ 
to $(Y^{[s]},\sigma_{Y^{[s]}})$. 

As is easily seen, the following proposition holds:
\begin{proposition} 
Let $f:L_{N+1}(X)\to L_1(Y)$ be a local rule 
on a subshift $(X,\sigma_X)$ to 
another $(Y,\sigma_Y)$ and let $s\geq 1$. 
\begin{enumerate}
\item If  $f$ is strictly 
$k$ right-mergible (respectively, strictly $l$ left-mergible), 
then $f^{[s]}$ is strictly
$k$ right-mergible (respectively, strictly $l$ left-mergible). 
\item If a homomorphism $\phi:(X,\sigma_X)\to (Y,\sigma_Y)$ is  
of $(m,n)$-type with $m+n=N$ and given by $f$, then 
$\phi^{[s]}$ is of $(m,n)$-type and given by  
$f^{[s]}$. 
\item For any homomorphism between subshifts,  
$Q_R(\phi)=Q_R(\phi^{[s]})$ and
$Q_L(\phi)=Q_L(\phi^{[s]})$.
\end{enumerate}
\end{proposition}

Let $(X,\sigma_X), (Y,\sigma_Y)$, and $(Z,\sigma_Z)$ be subshifts. 
Let $N,N'\geq 0$.
Let $f:L_{N+1}(X)\to L_{1}(Y)$ and $g:L_{N'+1}(Y)\to L_{1}(Z)$
be local rules on $(X,\sigma_X)$ to $(Y,\sigma_Y)$ and 
on $(Y,\sigma_Y)$ to $(Z,\sigma_Z)$, respectively. 
We define the \itl{composition} 
$gf:L_{N+N'+1}(X)\to L_1(Z)$ of $f$ and $g$ by 
\[gf(a_1\dots a_{N+N'+1})=g(f(v_1)f(v_2)\dots f(v_{N'+1})),\]
where $v_j=a_j\dots a_{j+N}$ for $j=1\dots N'+1$. 
Clearly, $gf$ is a local rule on $(X,\sigma_X)$ to $(Z,\sigma_Z)$. 

\begin{lemma} Let $(X,\sigma_X), (Y,\sigma_Y)$, 
and $(Z,\sigma_Z)$ be subshifts. Let
$f:L_{N+1}(X)\to L_{1}(Y)$ and $g:L_{N'+1}(Y)\to L_{1}(Z)$
be local rules on $(X,\sigma_X)$ to $(Y,\sigma_Y)$ and 
on $(Y,\sigma_Y)$ to $(Z,\sigma_Z)$, respectively. 
Then if $f$ is $k$ right-mergible (respectively, $k$ left-mergible) 
and $g$ is 
$k'$ right-mergible 
(respectively, $k'$ left-mergible), 
then $gf$ is $k+k'$ 
right-mergible (respectively, $k+k'$ left-mergible). 
\end{lemma} 
\begin{proof} 
Let $\aseq$ and $\bseq$ are points of $X$ such that 
$(a_j)_{j\leq 0}=(b_j)_{j\leq 0}$ and 
$gf(a_{j-N-N'}\dots a_j)=gf(b_{j-N-N'}\dots b_j)$ 
for $j=1,\dots,k+k'+1$. 
Then $f(v_j)=f(w_j)$ for $j\leq 0$ and 
$g(f(v_{j-N'})\dots f(v_j))=g(f(w_{j-N'})\dots f(w_j))$ 
for $j=1,\dots,k+k'+1$, 
where $v_j=a_{j-N}\dots a_j$ and $w_j=b_{j-N}\dots b_j$ 
for $j\in\Z$. Therefore, 
since $g$ is $k'$ right-mergible, 
it follows that the prefixes of length 
$(k+k'+1)-k'$ of 
$f(v_1)\dots f(v_{k+k'+1})$ and 
$f(w_1)\dots f(w_{k+k'+1})$ are the same.
Hence $f(a_{j-N}\dots a_j)=f(b_{j-N}\dots b_j)$ 
for $j=1,\dots, k+1$. Therefore, 
since $(a_j)_{j\leq 0}=(b_j)_{j\leq 0}$ and $f$ is 
$k$ right-mergible, we have $a_1=b_1$. 
This implies that $gf$ is 
$k+k'$ right-mergible. 
\end{proof}
\begin{remark} Let 
$f$ and $g$ be  the same as in Lemma 3.6. Then the following hold:
\begin{enumerate} 
\item 
If $gf$ is $k$ right-mergible 
(respectively, $k$ left-mergible),  
then so is $f$. 
\item 
If $g$ is 0 right-mergible (respectively, 0 left-mergible), 
then $f$ is strictly $k$ right-mergible 
(respectively, strictly $k$ left-mergible) 
if and only if so is $gf$. 
\end{enumerate}
\end{remark} 
\begin{proof} 
(1) Let $\aseq$ and $\bseq$ be points of $X$ such that 
$(a_j)_{j\leq 0}=(b_j)_{j\leq 0}$ and 
$f(a_{j-N}\dots a_j)=f(b_{j-N}\dots b_j)$ for 
$j=1,\dots, k+1$. Then, putting
$v_j=a_{j-N}\dots a_j$ and $w_j=b_{j-N}\dots b_j$ 
for $j\in\Z$, we have
$f(v_j)=f(w_j)$ for $j\leq k+1$ and hence
$g(f(v_{j-N'})\dots f(v_j))=g(f(w_{j-N'})\dots f(w_j))$ 
for $j=1,\dots,k+1$, which implies that 
$gf(a_{j-N-N'}\dots a_j)=gf(b_{j-N-N'}\dots b_j)$ 
for $j=1,\dots,k+1$. Therefore, since $gf$ is $k$ right-mergible 
and $(a_j)_{j\leq 0}=(b_j)_{j\leq 0}$, 
we have $a_1=b_1$. 

(2) It follows from Lemma 3.6 and (1) that 
if $g$ is $0$ right-mergible, then 
$f$ is $k$ right-mergible if and only if so is $gf$. 
\end{proof}

\begin{proposition} Let $(X,\sigma_X)$, $(Y,\sigma_Y)$ 
and $(Z,\sigma_Z)$ be subshifts 
such that if $X$ is infinite then so is $Y$. 
Let $\phi:(X,\sigma_X)\to (Y,\sigma_Y)$ and 
$\psi:(Y,\sigma_Y)\to (Z,\sigma_Z)$ be homomorphisms. 
Then 
\[Q_R(\psi\phi)\geq Q_R(\phi)+Q_R(\psi)\q\;\text{and}\q\;  
Q_L(\psi\phi)\geq Q_L(\phi)+Q_L(\psi).\] 
\end{proposition} 
\begin{proof} 
If $\phi$ is of $(m,n)$-type and given by 
a local rule $f:L_{m+n+1}(X)\to L_1(Y)$ and 
$\psi$ is of $(m',n')$-type and given by 
a local rule $g:L_{m'+n'+1}(Y)\to L_1(Z)$,  
then $\psi\phi$ 
is of $(m+m',n+n')$-type and given by $gf$. Suppose that 
$f$ is strictly $k$ right-mergible, $g$ is 
strictly $k'$ right-mergible and 
$fg$ is strictly $k''$ right mergible. If $k\neq\infty$
and $k'\neq\infty$ 
then by Lemma 3.6 $k''\leq k+k'$, 
and if $k=\infty$ or $k'=\infty$ then $k''\leq k+k'$.
Hence we have 
$Q_R(\psi\phi)= n+n'-k''\geq n-k + n'-k'=Q_R(\phi)+Q_R(\psi)$. 

The proof of the second inequality is similar. 
\end{proof}

\begin{remark}
Let $\phi:(X,\sigma_X)\to (Y,\sigma_Y)$ be a homomorphism 
between subshifts. Let $s\geq 1$. 
Then the following hold:
\[Q_R(\phi)\geq sQ_R((\phi)^{(s)}),\q\q 
Q_L(\phi)\geq sQ_L((\phi)^{(s)}).\] 
\end{remark}
\begin{proof} We prove only the first inequality, because the 
second one follows from the first by symmetry. 
Suppose that 
$\phi$ is of $(m,n)$-type and given by a 
local rule $f:L_{N+1}(X)\to L_1(Y)$ with $N=m+n$. It follows from 
Proposition 3.3 that 
the addition of any prefix or suffix redundancy to 
the local rule does not change $Q_R(\phi)$. Hence
we may assume that both $m$ and $n$ are divisible 
by $s$. Let $m'=m/s$, $n'=n/s$ and $N'=m'+n'$. 
Define $g:L_{N'+1}(X^{(s)})\to L_1(Y^{(s)})$
as follows: for $a_1\dots a_{N+s}\in L_{N+s}(X)$ 
with $a_j\in L_1(X)$, 
\[g(c_1\dots c_{N'+1})
=f(a_1\dots a_{N+1})f(a_2\dots a_{N+2})\dots f(a_s\dots a_{N+s}),\]
where $c_j=a_{(j-1)s+1}\dots a_{js}$ for $j=1,\dots,N'+1$. 
Then $(\phi)^{(s)}$ is of $(m',n')$-type and given by the local rule $g$. 

Suppose that $f$ is $k$ right-mergible.  
Let $k'=\lceil k/s\rceil$. 
Let $\aseq, \bseq\in X$  with $a_j,b_j\in L_1(X)$. For $j\in\Z$, 
let $c_j=a_{(j-1)s+1}\dots a_{js}$ and let 
$d_j=b_{(j-1)s+1}\dots b_{js}$.
Then $(c_j)_{j\in\Z}, (d_j)_{j\in\Z}\in X^{(s)}$. 
Assume that $(c_j)_{j\leq 0}=(d_j)_{j\leq 0}$ and  
$g(c_{j-N'}\dots c_j)=g(d_{j-N'}\dots d_j)$ 
for $j=1,\dots, k'+1$, then it follows that 
$(a_j)_{j\leq 0}=(b_j)_{j\leq 0}$ and 
$f(a_{j-N}\dots a_{j})=f(b_{j-N}\dots b_{j})$ for 
$j=1,\dots, (k'+1)s$. 
Since $f$ is $k$ right-mergible, 
the prefixes of length $(k'+1)s-k$ of $a_1\dots a_{(k'+1)s}$ 
and $b_1\dots b_{(k'+1)s}$ must be the same. 
Since 
$(k'+1)s-k=\lceil k/s\rceil s +s-k\geq s$, we have 
$c_1=a_1\dots a_s=b_1\dots b_s=d_1$. 
Therefore, $g$ is $k'$ right-mergible. 

By the above,  
for $k=0,\infty$, 
if $f$ is strictly $k$ right-mergible, then so is $g$. 
the first inequality is proved when $k=0,\infty$. 

Assume that  $f$ is strictly $k$ right-mergible 
with $k\geq 1$ and $k\neq \infty$. 
Then it occurs that 
$(a_j)_{j\leq 0}=(b_j)_{j\leq 0}$, $a_1\neq b_1$ and 
$f(a_{j-N}\dots a_j)=f(b_{j-N}\dots b_j)$ for $j=1,\dots,k$. 
Let $\bar{c}_j=a_{(j-2)s+2}\dots a_{(j-1)s+1}$ and let   
$\bar{d}_j=b_{(j-2)s+2}\dots b_{(j-1)s+1}$, for $j\in\Z$.
Then $(\bar{c}_j)_{j\in\Z}$ and $(\bar{d}_j)_{j\in\Z}$ 
are points of $X^{(s)}$, 
$(\bar{c}_j)_{j\leq 0}=(\bar{d}_j)_{j\leq 0}$, 
$\bar{c}_1\neq \bar{d}_1$ and 
$g(\bar{c}_{j-N'}\dots \bar{c}_j)=
g(\bar{d}_{j-N'}\dots \bar{d}_j)$ for $j=1,2,\dots,k'$, 
because $1+\lfloor (k-1)/s\rfloor=\lceil k/s\rceil=k'$. 
Therefore, it is proved that $g$ is 
strictly $k'$ right-mergible. Hence 
$Q_R(\phi)=n-k\geq n-sk'=(n'-k')s=sQ_R((\phi)^{(s)})$,
and the first inequality is proved. 
\end{proof}

\section{$q$-R and $q$-L endomorphisms of topological Markov shifts}

We shall see what the $q$-R and $q$-L degrees 
are for factor maps of SFTs onto sofic systems 
for onto endomorphisms of transitive 
topological Markov shifts and 
for automorphisms of topological Markov shifts.

\begin{proposition} 
Let $\phi: (X,\sigma_X)\to (Y,\sigma_Y)$ be a factor map of 
an SFT onto a 
sofic system. Let $s\in\Z$.
\begin{enumerate}
\item $\phi$ is weakly $q$-R if and only if $Q_R(\phi)\geq 0$; moreover
$\phi\sigma_X^s$ is weakly $q$-R if and only if $s\geq -Q_R(\phi)$. 
\item $\phi$ is weakly $q$-L if and only if $Q_L(\phi)\geq 0$; moreover
$\phi\sigma_X^s$ is weakly $q$-L if and only if $s\leq Q_L(\phi)$.
\end{enumerate}
\end{proposition}
\begin{proof} 
(1) 
Suppose that $\phi$ 
is weakly $q$-R. Then there exists 
a textile relation system $T=(p:\G\to G,q:\G\to H)$ such that 
$q$ is weakly right-resolving, 
$\xi_T: (X_\G,\sigma_\G)\to (X,\sigma_X)$ is a conjugacy and 
$\phi=\eta_T\xi_T^{-1}$. 
There exist 
$m,n\geq 0$ and a local rule $g:L_{m+n+1}(X)\to L_1(\G)$ such that 
$\xi_T^{-1}:X\to X_\G$ is a block-map of $(m,n)$-type given
by $g$. Let $f=q_Ag$. Then $\phi$ is a block-map 
of $(m,n)$-type given by the local rule $f$. 
Let $\aseq$ and $\bseq$ be points in $X$ such that 
$(a_j)_{j\leq 0}=(b_j)_{j\leq 0}$ and 
$f(a_{j-m-n}\dots a_j)=f(b_{j-m-n}\dots b_j)$ for $j=1,\dots,n+1$. 
Let $\alpha_j=g(a_{j-m}\dots a_{j+n})$ 
and let $\beta_j=g(b_{j-m}\dots b_{j+n})$ for $j\in \Z$. Then 
$\alpha_{-n}=g(a_{-m-n}\dots a_0)=g(b_{-m-n}\dots b_0)=\beta_{-n}$ 
and 
$q(\alpha_{j-n})=f(a_{j-m-n}\dots a_j)
=f(b_{j-m-n}\dots b_j)=q(\beta_{j-n})$ for $j=1,\dots, n+1$. 
Therefore, since $q$ is weakly right-resolving, we have 
$\alpha_{j-n}=\beta_{j-n}$ for all $j=1,\dots,n+1$, and hence, 
in particular, $\alpha_1=\beta_1$. 
Hence $a_1=q(\alpha_1)=q(\beta_1)=b_1$, which implies that 
$f$ is $n$ right-mergible. 
Therefore, 
$Q_R(\phi)\geq n-n=0$ if $X$ is infinite. Otherwise, 
$Q_R(\phi)=\infty>0$. 

Conversely, assume that $Q_R(\phi)\geq 0$. 
Let $A=L_1(X)$ and $B=L_1(Y)$. 
Since $(X,\sigma_X)$ is an SFT, for all sufficiently large 
integer $I$ 
$(X^{[I]},\sigma_X^{[I]})$ is the topological Markov shift whose 
defining graph, say $G_I$, is defined as follows: 
$A_{G_I}=L_I(X)$; $V_{G_I}=L_{I-1}(X)$; for $w=a_1\dots a_I \in L_I(X)$ 
with $a_j\in A$, $i_{G_I}(w)=a_1\dots a_{I-1}$ 
and $t_{G_I}(w)=a_2\dots a_I$. By this 
and Proposition 3.3,  
we may assume,
by adding necessary right and left redundancies to 
the local rule giving $\phi$, 
that $\phi$ is of $(m,n)$ type and given by a local rule 
$f:L_{m+n+1}(X)\to B$  with $m,n\geq 0$ such that 
$(X^{[m+n+1]},\sigma_X^{[m+n+1]})$ is the topological Markov shift 
whose defining graph is $G_{m+n+1}$. Suppose that 
$f$ is strictly $k$ right-mergible. Then $n-k\geq 0$ 
(because $n-k=Q_R(\phi)\geq 0$ if $Q_R(\phi)\in\Z$, and otherwise, 
$X$ is finite and hence $k=0$). 
We have the graph-homomorphism 
$q_f:\G\to G_B$ with $\G=G_{m+n+1}$ naturally defined by $f$. 
Since $f$ is $k$ right-mergible, so is $q_f$.  
We consider the induced right-resolving graph-homomorphism
$(q_f)^{+}:\G^+_{q_f}\to G_B$. 
Each arc in $\G^+_{q_f}$ is of the form 
$(U,b)$, 
where $U\in \C_{q_f}$ and $b$ is an arc in $G_B$
(see Subsection 2.2 for the notation). 
By definition, $B^+_{q_f}(U,b)$ is the set of  
all arcs $w$ in $G_{m+n+1}$ 
such that $i_{G_{m+n+1}}(w)\in U$ and $q_f(w)=b$. We claim that 
all arcs in $B^+_{q_f}(U,b)$ are words 
in $L_{m+n+1}(X)$ having the same 
prefix of length $m+n+1-k$. 
We prove this claim in the following. 

If $n=0$, then $q_f$ is $0$ right-mergible (because $k\leq n$), 
and hence every set $U\in\C^+_h$ is a singleton 
and the claim is valid. 
Therefore, we assume that $n\geq 1$. 
Suppose that $w=a_1\dots a_{m+n+1}$ 
and $\bar{w}=\bar{a}_1\dots \bar{a}_{m+n+1}$ are in 
$B_{q_f}(U,b)$ with $a_j,\bar{a}_j\in A$. Since 
the vertices $a_1\dots a_{m+n}$ 
and $\bar{a}_1\dots \bar{a}_{m+n}$ in $G_{m+n+1}$ are in the 
right-compatible set $U$, there exist
$s\geq m+n$ and paths 
$w_1\dots w_s$ and $\bar{w}_1\dots \bar{w}_s$ 
of length $s$ in $G_{m+n+1}$
with $w_j,\bar{w}_j\in A_{G_{m+n+1}}=L_{m+n+1}(X)$ 
starting from the same vertex and ending in 
$a_1\dots a_{m+n}$ and $\bar{a}_1\dots \bar{a}_{m+n}$, respectively, 
in $G_{m+n+1}$ 
and generating the same path (of length $s$) in $G_B$ 
under $q_f$. 
Clearly $P=w_1\dots w_sw$ 
and $\bar{P}=\bar{w}_1\dots \bar{w}_s\bar{w}$ 
are paths of length $s+1$ going from the same vertex 
and generating the same path. 
Since $q_f$ is  $k$ right-mergible, 
the initial subpaths of length 
$s+1-k$ of $P$ and $\bar{P}$ are the same. Hence 
$w_{s+1-k}=\bar{w}_{s+1-k}$. Since 
$w_s=c_0a_1\dots a_{m+n}$ and 
$\bar{w}_s=\bar{c}_0\bar{a}_1\dots \bar{a}_{m+n}$ with some 
$c_0,\bar{c}_0\in A$ and $n\geq k$, it follows that 
$w_{s+1-k}=c_{-k+1}\dots c_0a_1\dots a_{m+n+1-k}$ 
and $\bar{w}_{s+1-k}
=\bar{c}_{-k+1}\dots \bar{c}_0\bar{a}_1\dots \bar{a}_{m+n+1-k}$ 
with some $c_j,\bar{c}_j\in A$. Hence  
$a_1\dots a_{m+n+1-k}=\bar{a}_1\dots \bar{a}_{m+n+1-k}$, 
and the claim is proved. 

Since $n-k\geq 0$, it follows from
the claim above that 
we can define a textile relation system
$T_0=(p_0:\G^+_{q_f}\to G_A,(q_f)^+:\G^+_{q_f}\to G_B)$ such that 
$(p_0)_A((U,b))$ is the $(m+1)$-st arc (in $G_A$) of 
any path (in $L_{m+n+1}(G_A)$) in $B^+_{q_f}(U,b)$.  
Then, $T_0$ is onesided 1-1 and weakly $q$-R and 
$\phi_{T_0}=\phi$. Hence $\phi$ is weakly $q$-R. 

We have proved that $\phi$ is weakly $q$-R 
if and only if $Q_R(\phi)\geq 0$. 
It follows from Proposition 3.4 that 
$Q_R(\phi\sigma_X^s)\geq 0$ if and only if 
$s\geq -Q_R(\phi)$. Therefore, 
$\phi\sigma_X^s$ is weakly $q$-R if and only if 
$s\geq -Q_R(\phi)$. 

(2) By symmetry, (2) follows from (1). 
\end{proof} 

\begin{theorem} Let $\varphi$ be 
an onto endomorphism of a topological Markov shift 
$(X_G,\sigma_G)$ with $\varphi$ one-to-one or $G$  
irreducible. Let $s\in\Z$.  Then  
\begin{enumerate}
\item $\varphi$ is $q$-R if and only if $Q_R(\varphi)\geq 0$; moreover
$\varphi\sigma^s$ is $q$-R if and only if 
$s\geq -Q_R(\varphi)$; if $s>-Q_R(\varphi)$, then $\varphi\sigma^s$ 
is positively left $\sigma$-expansive. 
\item $\varphi$ is $q$-L if and only if $Q_L(\varphi)\geq 0$; moreover
$\varphi\sigma^s$ is $q$-L if and only if 
$s\leq Q_L(\varphi)$; if $s<Q_L(\varphi)$, then $\varphi\sigma^s$ 
is positively right $\sigma$-expansive. 
\end{enumerate}
\end{theorem} 
\begin{proof} (1) The proof that if $\varphi$ is $q$-R 
then $Q_R(\varphi)\geq 0$, is similar to the proof of the corresponding 
part of Proposition 4.1.  

Suppose that $Q_R(\varphi)\geq 0$.  Suppose that 
$\varphi$ is of $(m,n)$-type and given 
by a local rule 
$f:L_{m+n+1}(G)\to A_G$. 
If we modify the proof of Proposition 4.1 (the part 
of the proof that if $Q_R(\phi)\geq 0$ then $\phi$ is weakly $q$-R) 
letting $\G=G^{[m+n+1]}$, replacing $G_A$ and $G_B$ by $G$ 
and using Lemma 2.2 (note that this can be used because  
$\varphi$ is one-to-one or $G$ is  
irreducible), then 
we directly have a onesided 1-1, nondegenerate, $q$-R textile system 
$T_0=(p_0,(q_f)^+:\G^+_{q_f}\to G)$ such that  
$\varphi_{T_0}=\varphi$.  Hence it is proved that 
$\varphi$ is $q$-R. 

We have proved that $\varphi$ is 
$q$-R if and only if $Q_R(\varphi)\geq 0$. From this and 
Proposition 3.4 it follows that 
$\varphi\sigma^s$ is $q$-R if and only 
if $s\geq -Q_R(\varphi)$. If $s\geq -Q_R(\varphi)+1$, then 
$Q_R(\varphi\sigma^{s-1})\geq 0$. Hence $\varphi\sigma^{s-1}$ is 
$q$-R, and hence by Proposition 2.11(2) 
$\varphi\sigma^s$ is positively left $\sigma$-expansive. 

(2) By symmetry, (2) follows from (1). 
\end{proof}

\begin{proposition}  
Let $\phi:(X,\sigma_X)\to (Y,\sigma_Y)$ be a factor map 
of an SFT onto a sofic system. Let $s\in\Z$.
\begin{enumerate}
\item If $\phi^{[t]}\sigma_{X^{[t]}}^s$ is 
weakly $q$-biresolving with some $t\geq 1$, then  
\[-Q_R(\phi)\leq s\leq Q_L(\phi).\] 
\item If $Q_R(\phi)+Q_L(\phi)\geq 0$, then there exists 
$t\geq 1$ such that $\phi^{[t]}\sigma_{X^{[t]}}^s$ 
is weakly $q$-biresolving 
for all $-Q_R(\phi)\leq s\leq Q_L(\phi)$; 
more precisely, if $\phi$ is given by a strictly $k$ right-mergible, 
strictly $l$ left-mergible local rule, then for $t=k+l+1$ 
we can construct a 
onesided 1-1, weakly $q$-biresolving textile-relation system 
\[\q\q\q\q\q\q T_s=(p_s:\G_0\to G_A^{{[t]}}, q:\G_0\to G_B^{[t]})\q
\text{with}\;\,A=L_1(X), B=L_1(Y)\]
such that $\phi_{T_s}=\phi^{[t]}\sigma_{X^{[t]}}^s$ 
for all $-Q_R(\phi)\leq s\leq Q_L(\phi)$. 
\end{enumerate}
\end{proposition} 
\begin{proof} (1)
If $\phi^{[t]}\sigma_{X^{[t]}}^s$ is 
weakly $q$-biresolving, then 
by Propositions 3.4, 3.5 and 4.1 
$Q_R(\phi)+s
=Q_R(\phi^{[t]}\sigma_{X^{[t]}}^s)\geq 0$ and 
$Q_L(\phi)-s=Q_L(\phi^{[t]}\sigma_{X^{[t]}}^s)\geq 0$,  
and hence $-Q_R(\phi)\leq s\leq Q_L(\phi)$. 

(2) Suppose that 
$Q_R(\phi)+Q_L(\phi)\geq 0$. 
For the same reason as in the second paragraph of the proof of 
Proposition 4.1, 
we may assume that 
$\phi$ is of $(m,n)$ type given by a local rule $f: L_{m+n+1}(X)\to B$
with $m,n\geq 0$ such that   
$(X^{[m+n+1]},\sigma_X^{[m+n+1]})$ is the topological Markov shift 
with the defining graph $G_{m+n+1}$ as defined there. 
Since $Q_R(\phi)+Q_L(\phi)\geq 0$, it cannot be the case that 
$Q_R(\phi)=-\infty$ or $Q_L(\phi)=-\infty$. If $Q_R(\phi)=\infty$ then 
let $k=l=0$. If $Q_R(\phi)\in\Z$ then let
$k=n-Q_R(\phi)$ and $l=m-Q_L(\phi)$. 
Then $f$ is strictly $k$ right-mergible and 
strictly $l$ left-mergible with 
\[m+n-k-l\geq 0.\]  To prove (2),  we construct a 
onesided 1-1, weakly $q$-biresolving textile-relation system 
\[T_s=(p_s:\G_0\to G_A^{{[k+l+1]}}, q:\G_0\to G_B^{[k+l+1]})\]
such that $\phi_{T_s}=\phi^{[k+l+1]}\sigma_{X^{[k+l+1]}}^s$ 
for $-Q_R(\phi)\leq s\leq Q_L(\phi)$. 

Let $\G=G_{m+n+1}$. Let $F:L_1(\G)\to G_B$ be the local rule on 
$(X_\G,\sigma_\G)$ to $(X_{G_B},\sigma_{G_B})$ such that 
$F(w)=f(w)$ for $w\in A_\G=L_{m+n+1}(X)$. 
Since $f$ is $k$ right-mergible and $l$ left-mergible, so is $F$.
We use the graph-homomorphism
\[q^{-+}_{F;l,k}:\G^{-+}_{F;l,k}\to G_B^{[l+k+1]}\] 
due to Kitchens \cite{Kit-s}
(see Subsection 2.2). Each arc $\bar{\alpha}$ in $\G^{-+}_{F;l,k}$ 
is written in the form 
\[\bar{\alpha}=D^{-+}_{F;l,k}(\mu_{-1}\alpha\mu_1),\] 
where $\mu_{-1}\alpha\mu_1\in L_{l+k+1}(\G), \mu_{-1}\in L_l(\G), 
\alpha\in A_\G$ and $\mu_1\in L_k(\G)$, and we have  
$q^{-+}_{F;l,k}(\bar{\alpha})=F(\mu_{-1}\alpha\mu_1)$. 
By definition 
$\alpha$ is uniquely determined by $\bar{\alpha}$. 
Since $\G=G_{m+n+1}$, the arc $\alpha$ is a word in 
$L_{m+n+1}(X)$ and hence is written  
\[\alpha=a_1\dots a_{m+n+1}\] and 
the path $\mu_{-1}\alpha\mu_1$ is written
$\mu_{-1}\alpha\mu_1=\alpha_{-l}\alpha_{-l+1}\dots\alpha_k$, 
where $\alpha_0=\alpha$ and 
$\alpha_j=a_{j+1}\dots a_{j+m+n+1}$ for 
$j=-l,\dots,k$ with 
\[a_{-l+1}\dots a_{k+m+n+1}\in L_{l+m+n+1+k}(X),\q a_j\in A,\] 
and hence we have
\begin{equation}
q^{-+}_{F;l,k}(\bar{\alpha})=f(a_{-l+1}\dots a_{k+m+n+1}).
\end{equation} 
Noting that $-(n-k)\leq m-l$, 
we define, for $-(n-k)\leq s\leq m-l$, 
a graph-homomorphism $p_s:\G^{-+}_{F;l,k}\to G_A^{[l+k+1]}$ by 
\begin{equation}
p_s(\bar{\alpha})=a_{m+1-l-s}\dots a_{m+1+k-s}.
\end{equation}
Note that the word of the right side of (4.2) 
is a subword of $\alpha$ for 
for all $-(n-k)\leq s\leq m-l$, 
because its length is $k+l+1$,  
$p_{-(n-k)}(\bar{\alpha})=a_{m+n+1-k-l}\dots a_{m+n+1}$,
$p_{m-l}(\bar{\alpha})=a_1\dots a_{k+l+1}$ and
$m+n+1-(k+l+1)\geq 0$. 
Therefore, the definition of $p_s$ 
is possible for $-(n-k)\leq s\leq m-l$. 
For $-(n-k)\leq s\leq m-l$, 
define the textile-relation system 
\[T_s=(p_s:\G^{-+}_{F;l,k}\to G_A^{{[l+k+1]}}, 
q^{-+}_{F;l,k}:\G^{-+}_{F;l,k}\to G_B^{[l+k+1]}).\]
By Lemma 2.4(2) 
$q^{-+}_{F;l,k}:\G^{-+}_{F;l,k}\to G_B^{[l+k+1]}$
is weakly biresolving, and hence 
$T_s$ is weakly $q$-biresolving. Clearly $T_s$ is onesided 1-1. 
It follows from (4.1) and (4.2) that 
\[\phi_{T_s}
=\phi^{[k+l+1]}\sigma_{X^{[k+l+1]}}^{s}, \q\; -(n-k)\leq s\leq m-l.\]
Therefore (2) is proved (note that if $X$ is finite, then 
$k=l=0$ and we can take $m,n$ arbitrarily large). 
\end{proof} 

\begin{theorem} 
Suppose that $\varphi$ is an onto endomorphism of a 
topological Markov shift $(X_G,\sigma_G)$ with $G$ irreducible and 
$Q_R(\varphi)+Q_L(\varphi)\geq 0$. Let $s\in\Z$. 
If $\varphi$ is of $(m,n)$-type 
and given by 
a strictly $k$ right-mergible, strictly $l$ left-mergible
local rule $f:L_{m+n+1}(G)\to A_G$, 
then 
for all $-Q_R(\varphi)\leq s\leq Q_L(\varphi)$, we can construct
the textile system 
\[T_s=(p_s,q_{f;l,k}:G^{-+}_{f;l,k}\to G^{[k+l+1]})\] 
such that  
\[p_s(D^{-+}_{f;l,k}(w_1a_1\dots a_{m+n+1}w_2))
= a_{m+1-l-s}\dots a_{m+1+k-s},\] 
where $w_1a_1\dots a_{m+n+1}w_2\in L_{m+n+k+l+1}(G)$ 
with $a_j\in A_G$, $w_1\in L_l(G)$ and $w_2\in L_k(G)$. 
$T_s$ is onesided 1-1, nondegenerate and 
$q$-biresolving with
\[\varphi _{T_s}=\varphi^{[k+l+1]}(\sigma_G^{[k+l+1]})^s,\]
and hence $\varphi^{[k+l+1]}(\sigma_G^{[k+l+1]})^s$ is $q$-biresolving 
for all $-Q_R(\varphi)\leq s\leq Q_L(\varphi)$.  
\end{theorem} 
\begin{proof} 
By a proof which is similar to (but more direct than) the proof of 
Proposition 4.3 (take $G$ instead of $G_A$ and $G_B$, and 
take $G^{[m+n+1]}$ instead of $G_{m+n+1}$),   
we know that the textile system $T_s$ in the theorem
can be constructed and is onesided 1-1 and weakly $q$-biresolving with
\[\varphi _{T_s}=\varphi^{[k+l+1]}(\sigma_G^{[k+l+1]})^s\]
for $-Q_R(\varphi)\leq s\leq Q_L(\varphi)$. 
Since $\phi_{p_0}$ is onto and so is
$\phi_{q^{-+}_{f;l,k}}$ 
(because $f$ gives $\varphi$ which is onto), $T_s$ is nondegenerate. 
Since $G$ is irreducible and $\phi_{p_0}$ 
is a topological conjugacy, $G^{-+}_{f;l,k}$ is irreducible and has the 
same spectral radius as $G$. Therefore, since $q^{-+}_{f;l,k}$ is 
weakly $q$-biresolving, it follows from Lemma 2.1(2) that 
$q^{-+}_{f;l,k}$ is $q$-biresolving, and hence 
$T_s$ is $q$-biresolving. 
\end{proof}

Let $\varphi$ be an onto endomorphism of a subshift $(X,\sigma)$. 
\itl{Deciding the $\sigma$-expansiveness situation of $\varphi$} 
means deciding which one of the following is the case for $\varphi$:
\begin{enumerate}
\item $\varphi$ is expansive; 
\item $\varphi$ is left $\sigma$-expansive but not right $\sigma$-expansive;
\item $\varphi$ is right $\sigma$-expansive but not left $\sigma$-expansive;
\item $\varphi$ is not onesided $\sigma$-expansive (i.e. 
neither left $\sigma$-expansive nor right $\sigma$-expansive).
\end{enumerate} 

Here we note that \cite[the proof of Proposition 8.3(2)]{Nasu-te} 
proves: 
\begin{theorem}[\cite{Nasu-te}] 
An onto endomorphism of a subshift is positively expansive 
if and only if it is expansive and essentially weakly $q$-biresolving. 
\end{theorem}

We also recall the following. 

If $\varphi$ is 
a positively expansive endomorphism of
an irreducible topological Markov shift $(X,\sigma)$, then 
$\varphi$ is onto and
$(X,\varphi)$ is conjugate to a onesided topological 
Markov shift \cite{Kur1}. Hence by 
\cite[Theorem 3.9(1),(3), Theorem 2.5]{Nasu-t}, 
an onto endomorphism of 
an irreducible topological Markov shift is positively expansive 
if and only if it is expansive and essentially $q$-biresolving. 
(This also follows from Theorem 4.5 and Lemma 2.1(2).)

If $\varphi$ is 
a positively expansive endomorphism of
an irreducible topological Markov shift $(X,\sigma)$, then 
$\varphi$ is biclosing \cite{Kur1} and hence exactly $c$-to-one 
with some positive integer $c$ \cite{Nasu-c}; 
if $\sigma$ is topologically mixing in addition, $(X,\varphi)$ 
is topologically-conjugate to the onesided full 
$c$-shift (see \cite[p.4083]{Nasu-d}). 

\begin{theorem}
Let $\varphi$ be an onto endomorphism of 
an irreducible topological Markov shift $(X,\sigma)$.  
\begin{enumerate}
\item 
$Q_R(\varphi)+Q_L(\varphi)\geq 0$ if and only if 
there exists $s\in\Z$ such that $(\varphi\sigma^s)^{[t]}$ is 
$q$-biresolving with some $t\geq 1$; for $s\in\Z$
it holds that 
$(\varphi\sigma^s)^{[t]}$ is $q$-biresolving with some $t\geq 1$ 
if and only if $-Q_R(\varphi)\leq s\leq Q_L(\varphi)$; if 
$-Q_R(\varphi)<s<Q_L(\varphi)$, then $\varphi\sigma^s$ is 
positively expansive. 
\item 
If $Q_R(\varphi)+Q_L(\varphi)\geq 0$ with $X$ infinite, 
then for 
$s=-Q_R(\varphi), Q_L(\varphi)$, we can 
decide the $\sigma$-expansiveness situation of $\varphi\sigma^s$ 
and hence we can decide whether $\varphi\sigma^s$ is 
positively expansive or not. 
\item 
If $Q_R(\varphi)+Q_L(\varphi)\geq 0$ and $\varphi$ is 
$c$-to-one with $c\in\N$, then the topological entropy 
$h(\varphi\sigma^s)$ of $\varphi\sigma^s$ is equal to $\log c$ 
for each integer $-Q_R(\varphi)\leq s\leq Q_L(\varphi)$, 
in particular, when $X$ is infinite we have 
\[h(\varphi\sigma^{-Q_R(\varphi)})
=h(\varphi\sigma^{Q_L(\varphi)})=\log c.\]
\end{enumerate} 
\end{theorem}
\begin{proof} 
(1) By Proposition 3.5 and Theorems 4.2 and 4.4, 
it follows that for $s\in\Z$, 
$(\varphi\sigma^s)^{[t]}$ is $q$-biresolving with some $t\geq 1$ 
if and only if $-Q_R(\varphi)\leq s\leq Q_L(\varphi)$. 
If $-Q_R(\varphi)<s<Q_L(\varphi)$, then $\varphi\sigma^{s-1}$ 
and $\varphi\sigma^{s+1}$
is $q$-biresolving, and hence $\varphi\sigma^s$ right $\sigma$-expansive 
and left $\sigma$-expansive (by Proposition 2.11(2)) and hence 
expansive (by \cite[Corollary 7.3]{Nasu-te}). 
Therefore, by the fact recalled above, 
if $-Q_R(\varphi)<s<Q_L(\varphi)$ then $\varphi\sigma^s$ is 
positively expansive. 

(2) Let $-Q_R(\varphi)\leq s\leq Q_L(\varphi)$. 
By Theorem 4.4, 
we can construct a onesided 1-1, nondegenerate, 
$q$-biresolving textile system 
$T_s=(p_s,q :\G_0\to G^{[t]})$ 
such that  
$\varphi _{T_s}=(\varphi\sigma^s)^{[t]}$ for 
$t\geq 1$ which we can effectively find. 
Since the dual textile system $T^*_s$  is LL, we can decide 
whether $\xi_{T^*}$ is one-to-one or not and 
whether $\eta_{T^*}$ is one-to-one or not. 
(This follows from \cite[Lemma 6.25]{Nasu-t}.) 
Hence the result follows from  
\cite[Proposition 6.2]{Nasu-te} 
and \cite[Theorem 2.5]{Nasu-t}. 

(3) Notation being the same as above, let $G^*_s$ be 
the graph such that $T^*_s$ is defined over it. Since 
$T^*_s$ is LL, we have $\sigma_{T^*_s}=\sigma_{G^*_s}$. Therefore,
using \cite[Theorem 2.13, Corollary 2.14]{Nasu-t} we have 
\[h(\varphi\sigma^s)=h((\varphi\sigma^s)^{[t]})=h(\varphi_{T_s})
=h(\sigma_{T^*_s})=h(\sigma_{G^*_s}).\] 
It follows from \cite[Proposition 3.1(1)]{Nasu-n} that the spectral 
radius of $G^*_s$ is equal to $c$. Therefore (3) is proved. 
\end{proof} 

Next we shall see a method for constructing positively expansive 
endomorphisms of full shifts. ``Positively expansive endomorphisms 
of full-shifts'' is not a very confined subject. 
Boyle and Maass \cite{BoyMaa} proved that 
any mixing topological Markov shift having a positively 
expansive $c$-to-one endomorphism is shift equivalent to some 
full $d$-shift with $c$ and $d$ divisible by the same primes; 
they also conjectured that in this, ``shift equivalent'' 
can be replaced by 
``topologically conjugate'' (see \cite[p.4083]{Nasu-d}).

Let $A$ be an alphabet. 
Let $f:A^{N+1}\to A$ and $f':A^{N'+1}\to A$ be local rules 
with $N,N'\geq 0$. Let $t$ be an integer with 
$t\leq \min\{N,N'\}$.  Let $\pi: A\times A\to A$ 
be a \itl{bipermutation} 
(i.e., for any $a,a',b,b'\in A$, 
if $a\neq a'$ then $\pi(a,b)\neq \pi(a',b))$, 
and if $b\neq b'$ then $\pi(a,b)\neq \pi(a,b')$). 
Let $g:A^{N+N'-t+1}\to A$ be the local rule defined by 
\[g(a_1\dots a_{N+N'-t+1})=\pi(f'(w'),f(w)),\]
where $w'$ and $w$ are the $(N'+1)$-prefix and the $(N+1)$-suffix, 
respectively, of $a_1\dots a_{N+N'-t+1}$ with $a_j\in A$.  
(Note that $N+N'-t\geq \max\{N,N'\}$
by our hypothesis for $t$.)
We call $g$ the 
\itl{local rule 
defined by $(f',f)$, $t$ and $\pi$}. 

Let $f:A^{N+1}\to A$ be a local rule which gives an onto block-map. 
Let $R(f)$ (respectively, $L(f)$) be the cardinality of 
a maximal right (respectively, left) compatible-set for 
the graph-homomorphism $q_f:G_A^{[N+1]}\to G_A$ 
(see Subsection 2.2 for the terminology and notation). 
$R(f)$ and $L(f)$ are the same as those of 
L. R. Welch which and whose properties 
are found in \cite[Section 14]{Hedlund}.

\begin{lemma} 
Let $A$ be an alphabet.  
Let $f:A^{N+1}\to A$ and $f':A^{N'+1}\to A$ be local rules. 
Let $g:A^{N+N'-t+1}\to A$ be the local rule  
defined by $(f',f)$, an integer $t\leq \min\{N,N'\}$ and 
a bipermutation $\pi:A\times A\to A$. 
\begin{enumerate}
\item 
If $f$ is strictly $k$ right-mergible and $t<N-k$, 
then $g$ is strictly $k$ right-mergible and $R(g)=R(f)$. 
\item 
If $f'$ is strictly $l'$ left-mergible and $t<N'-l'$, 
then $g$ is strictly $l'$ left-mergible and L(g)=L(f').
\end{enumerate}
\end{lemma} 
\begin{proof} (1) Let $\widehat{N}=N+N'-t$. 
Let $\aseq$ and $\bseq$ be any two points in $A^\Z$ with 
$a_j, b_j\in A$ such that $(a_j)_{j\leq 0}=(b_j)_{j\leq 0}$. Then 
for $j\in\Z$, 
\begin{align*}
g(a_{j-\widehat{N}}\dots a_j)
&=\pi(f'(a_{j-\widehat{N}}\dots a_{j-\widehat{N}+N'}),f(a_{j-N}\dots a_j)),\\
g(b_{j-\widehat{N}}\dots b_j)
&=\pi(f'(b_{j-\widehat{N}}\dots b_{j-\widehat{N}+N'}),f(b_{j-N}\dots b_j)).
\end{align*}
Since $a_j=b_j$ for all $j\leq 0$ and 
\[k+1-\widehat{N}+N'=t-(N-k)+1\leq 0\q\text{by assumption,}\]
we have $a_{j-\widehat{N}}\dots a_{j-\widehat{N}+N'}=
b_{j-\widehat{N}}\dots b_{j-\widehat{N}+N'}$ 
for all $j=1,\dots, k+1$. Therefore 
it holds for all $j=1,\dots, k+1$ that 
$f(a_{j-N}\dots a_j)=f(b_{j-N}\dots b_j)$ if and only if 
$g(a_{j-\widehat{N}}\dots a_j)=g(b_{j-\widehat{N}}\dots b_j)$. 

Therefore, 
if $g(a_{j-\widehat{N}}\dots a_j)=g(b_{j-\widehat{N}}\dots b_j)$ 
for $j=1,\dots, k+1$, then we have 
$f(a_{j-N}\dots a_j)=f(b_{j-N}\dots b_j)$ 
for $j=1,\dots, k+1$, and hence we have $a_1=b_1$, 
because $f$ is $k$ right-mergible. 
If $a_1\neq b_1$ and
$f(a_{j-N}\dots a_j)=f(b_{j-N}\dots b_j)$ 
for $j=1,\dots, k$, then  $a_1\neq b_1$ and
$g(a_{j-\widehat{N}}\dots a_j)=g(b_{j-\widehat{N}}\dots b_j)$ 
for $j=1,\dots, k$. Therefore, since
$f$ is strictly $k$ right-mergible, so is $g$.

Fix $a_{-\widehat{N}+1}\dots a_k \in A^{\widehat{N}+k}$ arbitrarily. 
Let $C=S^+_{q_f}(a_{-N+1}\dots a_0, f(a_{-N+1}\dots a_k))$ 
and let 
$D=S^+_{q_g}(a_{-\widehat{N}+1}\dots a_0, g(a_{-\widehat{N}+1}\dots a_k))$.
Since $f$ and $g$ are $k$ right-mergible, 
it follows from \cite[Lemma 5.3]{Nasu-c} that $C$ is 
a maximal right compatible set for $q_f$ and $D$ is 
a maximal right compatible set for $q_g$. 
Let 
\begin{align*}
E&=\{a_{-\widehat{N}+1}\dots a_{0}w\,|\, w\in A^k, 
g(a_{-\widehat{N}+1}\dots a_0w)=g(a_{-\widehat{N}+1}\dots a_k)\},\\
F&=\{a_{-\widehat{N}+1}\dots a_{0}w\,|\, w\in A^k, 
f(a_{-\widehat{N}+1}\dots a_0w)=f(a_{-\widehat{N}+1}\dots a_k)\}.
\end{align*}
For the same reason as in the above, we see that $E=F$. 
Since the cardinality of $E$ is equal to that of $D$ and 
the cardinality of $F$ is equal to that of $C$, we have 
$R(g)=R(f)$.

(2) By symmetry, (2) follows from (1). 
\end{proof}

\begin{theorem} 
Let $A$ be an alphabet with cardinality $|A|$. 
Suppose that $f:A^{N+1}\to A$ is a strictly $k$ right-mergible 
local rule and $f':A^{N'+1}\to A$ is a strictly 
$l'$ left-mergible local rule. 
Let $t$ be an integer such that  
\[t<\min\{N-k,N'-l'\}.\] 
Let $g:A^{N+N'-t+1}\to A$ be the local rule 
defined by $(f',f)$, $t$ and a 
bipermutation $\pi:A\times A\to A$. 
Let $\varphi$ be the endomorphism of $(0, N+N'-t)$-type of 
the full-shift $(A^\Z,\sigma_A)$ given 
by $g$. Then $g$ is strictly $k$ right-mergible and 
strictly $l'$ left-mergible, 
\[Q_R(\varphi)+Q_L(\varphi)>0\] 
and $\varphi$ is $c$-to-one with $c=|A|^{N+N'-t}/(L(f')R(f))$. 
\end{theorem}
\begin{proof} By Lemma 4.7, $g$ is strictly $k$ right-mergible 
and strictly $l'$ left-mergible. Therefore $Q_L(\varphi)=0-l'$
and $Q_R(\varphi)=N+N'-t-k$, and hence 
\[Q_R(\varphi)+Q_L(\varphi)=N+N'-k-l'-t> 0. \] 

The remainder follows from Lemma 4.7 and 
\cite[Theorem 14,9, Section 16]{Hedlund} 
(see \cite[Section 5]{Nasu-l}). 
\end{proof} 

Recalling Theorem 4.6(1), we know that the theorem above presents 
a method for obtaining  
positively expansive endomorphisms of any full shift 
from any pair of a left-closing endomorphism and a right-closing 
endomorphism of the shift.  In particular, we can obtain positively 
expansive endomorphism of any full shift 
from any biclosing endomorphism 
(in particular, any automorphism) of the shift 
by expanding its neighborhood-size sufficiently 
large and using a bipermutation. 
An application of the method to permutation cellular automata 
is found in \cite[Subsection 4.3]{Jad-Nasu-Yaz}. 

\section{Weakly $q$-R and weakly $q$-L endomorphisms of subshifts} 

Let $(X,\sigma)$ be a subshift.  For $s\geq 1$, we define 
a graph $G[X,s]$ as follows. $G[X,1]=G_A$,  
where $A=L_1(X)$; if $s\geq 2$, then 
$A_{G[X,s]}=L_s(X)$, $V_{G[X,s]}=L_{s-1}(X)$, and 
for each $w=a_1\dots a_s\in A_{G[X,s]}$ with $a_j\in A$, 
$i_{G[X,s]}$ and $t_{G[X,s]}$ map $w$ to $a_1\dots a_{s-1}$ and 
$a_2\dots a_s$, respectively. 
Let $(X^{\langle s\rangle},\sigma_{X^{\langle s\rangle}})$ be the 
SFT over $A$ such that 
$((X^{\langle s\rangle})^{[s]},(\sigma_{X^{\langle s\rangle}})^{[s]})
=(X_{G[X,s]},\sigma_{G[X,s]})$. That is, 
$(X^{\langle s\rangle},\sigma_{X^{\langle s\rangle}})$ is the SFT 
with $A^s\setminus L_s(X)$ as ``the set of forbidden words'' that defines it.  
We call $(X^{\langle s\rangle},\sigma_{X^{\langle s\rangle}})$ 
the \itl{approximation SFT of order $s$ of} 
$(X,\sigma)$. If $f:L_{N+1}(X)\to A$ is a local rule on 
$(X,\sigma)$ 
and $s\geq N+1$, then $f$ is a local rule on
$(X^{\langle s\rangle},\sigma_{X^{\langle s\rangle}})$ 
to $(X_A,\sigma_A)$ (because $L_{N+1}(X^{\langle s\rangle})=L_{N+1}(X)$). 
If $\varphi$ is the endomorphism of $(m,n)$-type of $(X,\sigma)$ given 
by a local rule $f:L_{N+1}(X)\to A$ with 
$m+n=N$, 
then for each $s\geq N+1$, $f$ gives a homomorphism of $(m,n)$-type
$\phi_0^{\langle s,\varphi\rangle}:
(X^{\langle s\rangle},\sigma^{\langle s\rangle})\to (X_A,\sigma_A)$, 
which is an extension of $\varphi$. 
Let $\phi^{\langle s,\varphi\rangle}$ 
be the factor map induced by $\phi_0^{\langle s,\varphi\rangle}$. 
Then $\phi^{\langle s,\varphi\rangle}$ is an extension of $\varphi$, 
of $(m,n)$-type and given by $f$.
We call $\phi^{\langle s,\varphi\rangle}$ 
the \itl{approximation factor map 
of order $s$ of $\varphi$}. 
\begin{lemma}
Let $f:L_{N+1}(X)\to L_1(X)$ 
be a local rule on a subshift $(X,\sigma)$ 
and $\varphi$ an endomorphism of 
$(X,\sigma)$ given by $f$. 
\begin{enumerate}
\item If $f$ is a strictly 
$k$ right-mergible (respectively, strictly $k$ left-mergible) 
local rule on $(X,\sigma)$, 
then there exists an integer $J\geq N+1$
such that for all $J'\geq J$, $f$ is a strictly 
$k$ right-mergible (respectively, strictly $k$ left-mergible) 
local rule on  
$(X^{\langle J'\rangle},\sigma_{X^{\langle J'\rangle}})$ 
to its image
under $\phi^{\langle J',\varphi\rangle}$. 
\item There exists an integer $J\geq N+1$
such that for all $J'\geq J$, 
$Q_R(\phi^{\langle J',\varphi\rangle})=Q_R(\varphi)$
(respectively, 
$Q_L(\phi^{\langle J',\varphi\rangle})=Q_L(\varphi)$). 
\end{enumerate}
\end{lemma}
\begin{proof}
(1) First we prove that if $f$ is $k$ right-mergible 
on $(X,\sigma)$, then there exists $J\geq N+1$ such that 
$f$ is a $k$ right-mergible local rule on  
$(X^{\langle J\rangle},\sigma_{X^{\langle J\rangle}})$ 
to its image
under $\phi^{\langle J,\varphi\rangle}$. 

To do this, assume the contrary. Then there 
would exist no $J\geq N+1$ 
such that $f$ is a 
$k$ right-mergible local rule on 
$(X^{\langle J\rangle},\sigma_{X^{\langle J\rangle}})$ 
to its image 
under $\phi^{\langle J,\varphi\rangle}$, and hence
for each $s\geq N+1$, there would exist a pair of points 
$(a^{(s)}_j)_{j\in\Z}$ and $(b^{(s)}_j)_{j\in\Z}$ 
in $X^{\langle s\rangle}$ such that 
$(a^{(s)}_j)_{j\leq 0}=(b^{(s)}_j)_{j\leq 0}$, 
$a_1^{(s)}\neq b_1^{(s)}$ and 
$f(a^{(s)}_{j-N}\dots a^{(s)}_j)=f(b^{(s)}_{j-N}\dots b^{(s)}_j)$ 
for $j=1,\dots, k+1$. 
Since $X^{\langle s+1\rangle}\subset X^{\langle s\rangle}$ 
and $X=\cap_{s\geq N+1}X^{\langle s\rangle}$
with each $X^{\langle s\rangle}$ closed, 
a standard compactness argument shows 
that there would exist points 
$(a_j)_{j\in\Z}$ and $(b_j)_{j\in\Z}$ 
in $X$ such that 
$(a_j)_{j\leq 0}=(b_j)_{j\leq 0}$, 
$a_1\neq b_1$ and 
$f(a_{j-N}\dots a_j)=f(b_{j-N}\dots b_j)$ 
for $j=1,\dots, k+1$. This contradicts the hypothesis that $f$ is 
$k$ right mergible on $(X,\sigma)$. 

Secondly we note that if $f$ is not $ k$ right-mergible on   
$(X,\sigma)$, then for each $s\geq N+1$,  
$f$ is not a $k$ right-mergible local rule on  
$(X^{\langle s\rangle},\sigma_{X^{\langle s\rangle}})$ to 
its image under $\phi^{\langle s,\varphi\rangle}$. 
(This follows because the assumption of this implies that 
there exist points $\cseq,\dseq$ in $X$ such that  
$(c_j)_{j\leq 0}=(d_j)_{j\leq 0}$, $c_1\neq d_1$, and 
$f(c_{j-N}\dots c_j)=f(c_{j-N}\dots d_j)$ for $j=1,\dots,k+1$, 
and because $X^{\langle s\rangle}\supset X$.) 

By the fact proved first, (1) follows when $k=0$. 
By the facts proved and noted above, we see that 
if $f$ is $k$ right-mergible but not $k-1$ right-mergible on
$(X,\sigma)$ with $k\geq 1$, then there  exists $J\geq N+1$ such that 
$f$ is a $k$ right-mergible but not $k-1$ right-mergible 
local rule on  $(X^{\langle J\rangle},\sigma_{X^{\langle J\rangle}})$ 
to its image under $\phi^{\langle J,\varphi\rangle}$. 
In addition, for all $J'\geq J$, $f$ is 
a $k$ right-mergible but not $k-1$ right-mergible local rule on  
$(X^{\langle J'\rangle},\sigma_{X^{\langle J'\rangle}})$ 
to its image under $\phi^{\langle J',\varphi\rangle}$, 
because $X\subset X^{\langle J'\rangle}\subset X^{\langle J\rangle}$. 
Therefore (1) follows when $k\in\Z$. 
By the fact secondly noted, (1) follows when $k=\infty$.

(2) Suppose that $\varphi$ is of $(m,n)$-type with 
$m+n=N$. 
Then for all $J\geq N+1$, 
$\phi^{\langle J,\varphi\rangle}$ is of $(m,n)$-type 
and given by $f$. Hence (2) follows from (1). 
\end{proof}

\begin{theorem}
Let $\varphi$ be an endomorphism of 
a subshift $(X,\sigma)$. Let $s\in\Z$. 
\begin{enumerate}
\item
$\varphi$ is 
weakly $q$-R if and only if $Q_R(\varphi)\geq 0$; furthermore,  
$\varphi\sigma^s$ is weakly $q$-R if and only 
if $s\geq -Q_R(\varphi)$; if $s> -Q_R(\varphi)$, then 
$\varphi\sigma^s$ is positively left $\sigma$-expansive. 
\item
$\varphi$ is 
weakly $q$-L if and only if $Q_L(\varphi)\geq 0$; furthermore,  
$\varphi\sigma^s$ is weakly $q$-L if and only 
if $s\leq Q_L(\varphi)$; if $s<Q_L(\varphi)$, then 
$\varphi\sigma^s$ is positively right $\sigma$-expansive.
\end{enumerate}
\end{theorem}
\begin{proof}
(1) Suppose that $\varphi$ is weakly $q$-R. Then there 
exists a onesided 1-1 half-textile-subsystem $\hf{U}$ of 
a weakly $q$-R textile system $T=(p,q:\G\to G)$ such that 
$(X_{\hf{U}},\sigma_{\hf{U}},\varphi_{\hf{U}})=(X,\sigma,\varphi)$. 
There exist 
$m,n \geq 0$ and a local rule 
$g:L_{m+n+1}(X_{\hf{U}})\to L_1(Z_{\hf{U}})$ such that 
$\xi_{\hf{U}}^{-1}:X_{\hf{U}}\to Z_{\hf{U}}$ 
is a block-map of $(m,n)$-type given by $g$. 
Let $f=q_Ag$. Then $\varphi$ is a block-map 
of $(m,n)$-type given by the local rule $f$. 
By the same proof as in the proof 
of Proposition 4.1, we see that $f$ is $n$ right-mergible 
and hence $Q_R(\varphi)\geq 0$. 

Conversely suppose that $Q_R(\varphi)\geq 0$. 
By Lemma 5.1, 
there exists $J\in\N$
such that     
$Q_R(\phi^{\langle J,\varphi\rangle})=Q_R(\varphi)$ . 
Since $Q_R(\phi^{\langle J,\varphi\rangle})\geq 0$, 
it follows from Proposition 4.1 that the factor map 
$\phi^{\langle J,\varphi\rangle}$ is weakly $q$-R. 
Therefore, there exists a onesided 1-1 textile relation system 
$T=(p:\G\to G_A,q:\G\to G_A)$ such that 
$q$ is weakly right-resolving and 
$\phi_T=\phi^{\langle J,\varphi\rangle}$, 
where $A=L_1(X^{\langle J\rangle})=L_1(X)$. 
We can regard $T$ as a textile system $T=(p,q:\G\to G_A)$, 
which is onesided 1-1 and weakly $q$-R. 
Let $Z=\hf{\xi}_T^{-1}(X)$. Then, since 
$Z\subset\hf{Z}_T$ with 
$\hf{\xi}_T(Z)\supset\hf{\eta}_T(Z)$
with $\hf{\xi}_T(Z)\supset\hf{\eta}_T(Z)$ 
and $\hf{\xi}_T|Z$ is one-to-one,  
there exists a unique onesided 1-1 
half-textile-subsystem $\hf{U}$ of $T$ with $Z_{\hf{U}}=Z$. 
Since $\varphi=\varphi_{\hf{U}}$, 
we conclude that $\varphi$ is weakly $q$-R. 

We have proved that $\varphi$ is 
weakly $q$-R if and only if $Q_R(\varphi)\geq 0$. From this and 
Proposition 3.4 it follows that 
$\varphi\sigma^s$ is weakly $q$-R if and only 
if $s\geq -Q_R(\varphi)$. If $s\geq -Q_R(\varphi)+1$, then 
$Q_R(\varphi\sigma^{s-1})\geq 0$. Hence $\varphi\sigma^{s-1}$ is 
weakly $q$-R, and hence 
$\varphi\sigma^s$ is positively left $\sigma$-expansive, 
by Proposition 2.8(2). 

(2) By symmetry, (2) follows from (1). 
\end{proof}

\begin{theorem} 
Let $\varphi$ be an endomorphism of a 
subshift $(X,\sigma)$. 
$Q_R(\varphi)+Q_L(\varphi)\geq 0$ if and only if 
there exists $s\in\Z$ such that $(\varphi\sigma^s)^{[t]}$ is 
weakly $q$-biresolving with some $t\geq 1$; 
it holds for $s\in\Z$ that 
$(\varphi\sigma^s)^{[t]}$ is weakly $q$-biresolving with some $t\geq 1$
if and only if $-Q_R(\varphi)\leq s\leq Q_L(\varphi)$; if 
$-Q_R(\varphi)<s<Q_L(\varphi)$, then $\varphi\sigma^s$ is 
positively expansive. 
\end{theorem} 
\begin{proof} 
It follows from Theorem 5.2 and Proposition 3.5 that 
if there exists $t\geq 1$ such that $\varphi^{[t]}$ is 
weakly $q$-biresolving, then 
$Q_R(\varphi)
=Q_R(\varphi^{[t]}) \geq 0$ and 
$Q_L(\varphi)=Q_L(\varphi^{[t]}) \geq 0$. 
Therefore it follows from Proposition 3.4 that  
if $\varphi^{[t]}(\sigma^{[t]})^s$ 
is weakly $q$-biresolving with some $t\geq 1$, 
then $-Q_R(\varphi)\leq s\leq Q_L(\varphi)$. 

Suppose that $-Q_R(\varphi)\leq s\leq Q_L(\varphi)$. 
By Lemma 5.1, there exists $J\in\N$ such that 
$Q_R(\phi^{\langle J,\varphi\rangle})=Q_R(\varphi)$ and 
$Q_L(\phi^{\langle J,\varphi\rangle})=Q_L(\varphi)$. 
Since $-Q_R(\varphi^{\langle J,\varphi\rangle})\leq s 
\leq Q_L(\varphi^{\langle J,\varphi\rangle})$,  
it follows from Proposition 4.3 that for some $t\geq 1$, 
we have 
a onesided 1-1, weakly $q$-biresolving textile-relation system 
\[\q\q\q T_s=(p_s:\G_0\to G_A^{[t]}, q:\G_0\to G_A^{[t]}) 
\q \text{with}\;\, A=L_1(X)\]
such that $\phi_{T_s}
=(\phi^{\langle J,\varphi\rangle})^{[t]}
\sigma_{(X^{\langle J \rangle})^{[t]}}^s$. 
We can regard $T_s$ as a textile system 
$T_s=(p_s,q:\G_0\to G_A^{[t]})$, 
which is onesided 1-1 and weakly $q$-biresolving. 
Let $Z_s=\hf{\xi}_{T_s}^{-1}(X^{[t]})$. 
Then, since $Z_s\subset\hf{Z}_{T_s}$ with  
$\hf{\xi}_{T_s}(Z_s)\supset\hf{\eta}_{T_s}(Z_s)$ and
$\hf{\xi}_{T_s}|Z_s$ one-to-one, 
there exists a unique onesided 1-1 
half-textile-subsystem $\hf{U}_s$ of $T_s$ with $Z_{\hf{U}_s}=Z_s$. 
Since $\varphi^{[t]}(\sigma^{[t]})^s=\varphi_{\hf{U}_s}$, 
we conclude that $\varphi^{[t]}(\sigma^{[t]})^s$ 
is weakly $q$-biresolving. 

Suppose that $-Q_R(\varphi)<s< Q_L(\varphi)$. Then by Theorem 5.2
$\varphi\sigma^s$ is positively left $\sigma$-expansive 
and positively right $\sigma$-expansive.  
Therefore, by Proposition 11.3(3) (appearing in Subsection 11.1), 
$\varphi\sigma^s$ is positively expansive. 
\end{proof} 

\begin{proposition} 
If right-closing (respectively, left-closing) 
endomorphisms 
$\varphi$ and $\psi$ of infinite subshifts 
are topologically-conjugate, then 
there exists an integer 
$s\geq 0$ such that for all $i\geq 0$, 
$|Q_R(\varphi^i)-Q_R(\psi^i)|\leq s$ (respectively,
$|Q_L(\varphi^i)-Q_L(\psi^i)|\leq s$). 
\end{proposition} 
\begin{proof} Suppose that $\varphi$ and $\psi$ are 
 right-closing
endomorphisms of subshifts
$(X,\sigma)$ and $(X',\sigma')$, respectively, 
and they are topologically conjugate. Then for 
all $i\geq 0$ 
$Q_R(\varphi^i),Q_R(\psi^i)\in\Z$ and 
there exists a conjugacy 
$\theta:(X,\sigma)\to (X',\sigma')$ such 
that $\psi=\theta\varphi\theta^{-1}$. 
Since $X$ is infinite, $Q_R(i_X)=0$ by (3.1). Hence 
by Proposition 3.8 
$Q_R(\theta)+Q_R(\theta^{-1})\leq 0$. 
Since $\theta$ and $\theta^{-1}$ are right-closing, 
$Q_R(\theta), Q_R(\theta^{-1}) \in\Z$. 
Let 
\[s=-Q_R(\theta)-Q_R(\theta^{-1}).\] 
Then $s$ is a nonnegative integer. 
By Proposition 3.8, for all $i\geq0$ 
\begin{align*}
Q_R(\psi^i)&=Q_R(\theta\varphi^i\theta^{-1})
\geq Q_R(\theta)+Q_R(\varphi^i)+ Q_R(\theta^{-1}) 
=Q_R(\varphi^i)-s, \q and \\
Q_R(\varphi^i)&=Q_R(\theta^{-1}\psi^i\theta) 
\geq Q_R(\theta^{-1})+Q_R(\psi^i)+ Q_R(\theta) 
=Q_R(\psi^i)-s. 
\end{align*}
Hence it follows that $|Q_R(\varphi^i)-Q_R(\psi^i)|\leq s$ 
for all $i\geq 0$. 
\end{proof} 

\begin{proposition} Let $\varphi$ be an endomorphism 
of a subshift $(X,\sigma)$. 
\begin{enumerate} 
\item 
If $\varphi$ essentially has nonnegative $q$-R degree 
(respectively, nonnegative $q$-L degree), 
then there exists an integer $s\geq 0$ 
such that for all $i\geq 0$, $Q_R(\varphi^i)\geq -s$ 
(respectively, $Q_L(\varphi^i)\geq -s$).
\item
If $\varphi$ is essentially weakly $q$-R (respectively, 
essentially weakly $q$-L), 
then there exists an integer $s\geq 0$ such that for all $i\geq 0$,
$\varphi^i\sigma^s$ is weakly $q$-R 
(respectively, $\varphi^i\sigma^{-s}$ is weakly $q$-L). 
\item Suppose in addition that $(X,\sigma)$ is 
a topological Markov shift and that 
$\varphi$ is invertible or $\sigma$ is topologically transitive, 
(2) with all ``weakly'' in it deleted hold. 
\end{enumerate}
\end{proposition}
\begin{proof} 
(1) If $X$ is finite, then (1) is true by definition. 
Hence we assume that $X$ is infinite.  
Since $\varphi$ essentially has nonnegative $q$-R degree, 
there exists an endomorphism $\psi$ 
of a subshift $(X',\sigma')$ such that 
$(X,\sigma,\varphi)$ is conjugate to $(X',\sigma',\psi)$ 
and $Q_R(\psi)\geq 0$, which implies that 
$\psi$ is right-closing and hence so is $\varphi$. 
By Proposition 5.4, there exists 
an integer $s\geq 0$ such that 
$Q_R(\varphi^i)\geq Q_R(\psi^i)-s$ for all $i\geq 0$. 
Since by Proposition 3.8 $Q_R(\psi^i)\geq iQ_R(\psi)\geq 0$, we have 
$Q_R(\varphi^i)\geq -s$ for all $i\geq 0$. 

(2) By (1), Theorem 5.2 and Proposition 3.4, (2) follows. 

(3) By (1), Theorem 4.2 and Proposition 3.4, (3) follows. 
\end{proof}

\section{$p$-L and $p$-R degrees}

Let $N\geq 0$. Let 
$f:L_{N+1}(X)\to L_1(Y)$  be a local rule on a subshift 
$(X,\sigma_X)$ to another $(Y,\sigma_Y)$. 
Let $I\geq 0$.
We say that $f$ is 
\itl{$I$ left-redundant} if for any points  
$\aseq$ and 
$\bseq$ in $X$ with 
$a_j,b_j\in L_1(X)$, it holds that if
$(a_j)_{j\geq 0}=(b_j)_{j\geq 0}$, then 
$f(a_{-I}\dots a_{-I+N})=f(b_{-I}\dots b_{-I+N})$. 
Symmetrically,  
$f$ is said to be \itl{$I$ right-redundant} 
if for any points 
$\aseq$ and 
$\bseq$ in $X$, 
it holds that 
if $(a_j)_{j\leq 0}=(b_j)_{j\leq 0}$, 
then $f(a_{I-N}\dots a_I)=f(b_{I-N}\dots b_I)$. 

We remark that if $f$ is $I$ left-redundant (respectively, 
$I$ right-redundant), then for all $0\leq I'\leq I$, 
$f$ is $I'$ left-redundant (respectively, 
$I'$ right-redundant). 

For $I\geq 0$ we say that $f$ is 
\itl{strictly I left-redundant} 
if it is $I$ left-redundant but 
not $I+1$ left-redundant; we say that $f$ is strictly 
\itl{$\infty$ left-redundant} if $f$ is $I$ left-redundant 
for all $I\geq 0$. 
Similarly,  
a \itl{strictly $I$ right-redundant} local rule with $I\geq 0$ and 
a \itl{strictly $\infty$ right-redundant} local rule are defined. 

(Note that we have never defined above ``$\infty$ left-redundant 
local rules'' and ``$\infty$ right-redundant 
local rules''.)

A homomorphism $\phi:(X,\sigma_X)\to (Y,\sigma_Y)$ 
between subshifts is said to be 
\itl{left-recognizing} 
(respectively, \itl{right-recognizing}) 
if there exists at least one pair of 
distinct right (respectively, left) $\sigma_X$-asymptotic points 
which $\phi$ does not collapse.

Let $\phi:(X,\sigma_X)\to (Y,\sigma_Y)$ be a homomorphism
between subshifts. Then if $\phi$ is left-recognizing or 
right-recognizing, then the set $X$ is infinite. 
The homomorphism 
$\phi$ is left-recognizing (respectively, right-recognizing) 
if and only if $\phi$ is given by a 
strictly $I$ left-redundant (respectively, 
strictly $I$ right-redundant) local rule 
for some integer $I\geq 0$. 
If the set $\phi(X)$ is finite, 
then $\phi$ is neither left-recognizing 
nor right recognizing. (For assume that there were
distinct right (respectively, left) 
$\sigma_X$-asymptotic points $x,x'\in X$ with 
$\phi(x)\neq\phi(x')$. Then $\phi(x)$ and $\phi(x')$ 
would be distinct right (respectively, left) 
$\sigma_Y$-asymptotic points in 
the finite subshift-space $\phi(X)$, which cannot be the case.) 
\begin{definition}
Let $m,n\geq 0$ and $N=m+n$. Let $I,J\in\{0,1,2,\dots,\infty\}$. 
Let $\phi: (X,\sigma_X)\to (Y,\sigma_Y)$ be 
a homomorphism 
between subshifts. Let $\phi$ be of  
$(m,n)$-type given by a local rule $f:L_{N+1}(X)\to L_1(Y)$ 
which is strictly $I$ left-redundant and strictly $J$ right-redundant.  
Define the \itl{$p$-L degree $P_L(\phi)$ of $\phi$} 
and the \itl{$p$-R degree $P_R(\phi)$ of $\phi$} by
\[P_L(\phi)=I-m\q\text{and}\q P_R(\phi)=J-n.\] 
\end{definition}

We follow the usual convention about the  
arithmetic calculations and inequalities containing $\infty$. 
Hence $P_L(\phi)<\infty$ (respectively, $P_R(\phi)<\infty$)
if and only if $\phi$ is left-recognizing (respectively, 
right-recognizing). 
If $\phi(X)$ is finite, 
then $P_L(\phi)=\infty$ and $P_R(\phi)=\infty$.  
Note that for example, 
$P_L(\phi)\geq 0$ means that $P_L(\phi)$ is some nonnegative integer 
or $P_L(\phi)=\infty$. 

The following lemma is obvious. 
\begin{lemma}  
Let $N,s,t\geq 0$. Let
$f:L_{N+1}(X)\to L_1(Y)$ be a local rule on  
a subshift $(X,\sigma_X)$ to another $(Y,\sigma_Y)$. 
Let $g$ be the local rule obtained from $f$ by adding 
\itl{right redundancy} by $s$ and \itl{left redundancy} by $t$ 
(as defined before Lemma 3.2). 
If $f$ is strictly $I$ left-redundant and strictly $J$ right-redundant,
then $g$ is strictly $I+s$ left-redundant and 
strictly $J+t$ right-redundant. 
\end{lemma}

\begin{proposition}
If $\phi$ is a homomorphism  between 
subshifts, then 
$P_L(\phi)$ and $P_R(\phi)$ are 
uniquely determined by $\phi$. 
\end{proposition} 
\begin{proof}
Suppose that for $i=1,2$, $\phi$ is a homomorphism of $(m_i,n_i)$-type
of a subshift $(X,\sigma_X)$ into another $(Y,\sigma_Y)$ 
given by a local-rule 
$f_i:L_{m_i+n_i+1}(X)\to L_1(Y)$ which is 
strictly $I_i$ left-redundant and strictly $J_i$ right-redundant 
with $I_i,J_i\in\{0,1,2,\dots,\infty\}$. 
Let $m=\max\{m_1,m_2\}$ and let $n=\max\{n_1,n_2\}$. For $i=1,2$, let 
$g_i:L_{m+n+1}(X)\to L_1(Y)$ be the local rule such that
$g_i(a_{-m}\dots a_n)=f_i(a_{-m_i}\dots a_{n_i})$, where 
$a_{-m}\dots a_n\in L_{m+n+1}(X)$ with $a_j\in L_1(X)$.  
Then, since $\phi$ is of $(m_i,n_i)$-type and given by $f_i$, 
$\phi$ is of $(m,n)$-type and given by $g_i$, for $i=1,2$. 
Since $f_i$ is strictly $I_i$ left-redundant and strictly 
$J_i$ right-redundant, 
it follows from Lemma 6.2 that $g_i$ is 
strictly $I_i+(m-m_i)$ left-redundant and strictly 
$J_i+(n-n_i)$ right-redundant, for $i=1,2$. 
Since $g_1:L_{m+n+1}(X)\to L_1(X)$ and $g_2:L_{m+n+1}(X)\to L_1(X)$
give the same block-map $\phi$ of $(m,n)$-type, 
it follows that $g_1=g_2$. Hence 
$I_1+(m-m_1)=I_2+(m-m_2)$ and $J_1+(n-n_1)=J_2+(n-n_2)$.
Thus we have 
$I_1-m_1=I_2-m_2$ and $J_1-n_1=J_2-n_2$. 
\end{proof} 

Here we give a little more compact description of the definition of
the $p$-R, $p$-L degrees of a homomorphism 
$\phi:(X,\sigma_X)\to (Y,\sigma_Y)$ between subshifts. 

Let $f:L_{N+1}(X)\to L_1(Y)$ be a local rule 
on $(X,\sigma_X)$ to $(Y,\sigma_Y)$. For $s\in\Z$ and 
$x=\aseq\in X$ with $a_j\in L_1(X)$, 
let $W_s^+(x)=\{\bseq\in X\,|\, (b_j)_{j\geq s}=(a_j)_{j\geq s}\}$
and $W_s^-(x)=\{\bseq\in X\,|\, (b_j)_{j\leq s}=(a_j)_{j\leq s}\}$. 
For $x=\aseq\in X$, define $\rho_{f,L}(x)$ 
(respectively, $\rho_{f,R}(x)$) to be the supremum of 
the integers $s\geq 0$ such that for all $\bseq\in W_0^+(x)$  
(respectively, $\bseq\in W_0^-(x)$), 
it holds that
for $j=0,\dots,s$, $f(a_{-j}\dots a_{-j+N})=f(b_{-j}\dots b_{-j+N})$
(respectively, $f(a_{j-N}\dots a_j)=f(b_{j-N}\dots b_j)$). 
Let $\rho_{f,L}=\min_{x\in X}\rho_{f,L}(x)$ 
(respectively, $\rho_{f,R}=\min_{x\in X}\rho_{f,R}(x)$).
We call $\rho_{f,L}$ (respectively, 
$\rho_{f,R}$) the \itl{minimum length of left redundancy} 
(respectively, \itl{minimum length of right redundancy})
of $f$. It is clear that 
$f$ is strictly $I$ left-redundant
(respectively, strictly $I$ right-redundant) if and only if 
$\rho_{f,L}=I$ (respectively, $\rho_{f,R}=I$).  

We can define 
$P_L(\phi)=\rho_{f,L}-m$ 
(respectively, $P_R(\phi)=\rho_{f,R}-n$) 
if $\phi$ is of $(m,n)$ type and given by a
local rule $f:L_{m+n+1}(X)\to L_1(Y)$. 

\begin{proposition}
Let $\phi$ be a homomorphism of a subshift 
$(X,\sigma_X)$ into another $(Y,\sigma_Y)$. Let $s\in\Z$. Then 
\[
P_L(\phi\sigma_X^s)=P_L(\phi)+s,\q\q  
P_R(\phi\sigma_X^s)=P_R(\phi)-s, 
\] 
and hence $P_R(\phi)+P_L(\phi)$ is shift-invariant.
\end{proposition} 
\begin{proof}
Suppose that $\phi$ is of $(m,n)$-type given by a local rule 
$f:L_{m+n+1}(X)\to L_1(Y)$. Increasing the redundancies of $f$ 
if necessary, we may assume that 
$m, n\geq |s|$ by Proposition 6.3. 
Since $\phi\sigma_X^s$ is of $(m-s,n+s)$-type and 
given by $f$, we have the first two equations, 
and hence the remainder follows.
\end{proof} 

We remark the following: 
for any infinite subshift $(X,\sigma_X)$  
\begin{equation}
P_L(i_X)= 0\q\q\text{and}\q\q P_R(i_X)=0
\end{equation}
and hence by Proposition 6.4, for all $s\in\Z$
\begin{equation}
P_L(\sigma_X^s)= s\q\q\text{and}\q\q P_R(\sigma_X^s)=-s.
\end{equation} 
For any homomorphism $\phi$ of $(m,n)$-type between subshifts, 
\begin{equation}
P_L(\phi)\geq -m\q\q\text{and}\q\q P_R(\phi)\geq -n.
\end{equation} 

As is easily seen, we have: 
\begin{proposition} Let $(X,\sigma_X)$, $(Y,\sigma_Y)$ be subshifts. 
Let $t\geq 1$. 
\begin{enumerate}
\item  
If a local rule $f:L_{N+1}(X)\to L_1(Y)$ on $(X,\sigma_X)$ 
to $(Y,\sigma_Y)$ is strictly
$I$ left-redundant (respectively, strictly $J$ right-redundant), 
then the higher block presentation $f^{[t]}$ is strictly
$I$ left-redundant (respectively, strictly $J$ right-redundant). 
\item 
If $\phi:(X,\sigma_X)\to (Y,\sigma_Y)$ is a 
homomorphism, then  $P_L(\phi)=P_L(\phi^{[t]})$ and 
$P_R(\phi)=P_R(\phi^{[t]})$. 
\end{enumerate}
\end{proposition} 

\begin{lemma} Let $(X,\sigma_X), (Y,\sigma_Y)$, 
and $(Z,\sigma_Z)$ be subshifts. Let
$f:L_{N+1}(X)\to L_{1}(Y)$ and $g:L_{N'+1}(Y)\to L_{1}(Z)$
be local rules on $(X,\sigma_X)$ to $(Y,\sigma_Y)$ and 
on $(Y,\sigma_Y)$ to $(Z,\sigma_Z)$, respectively. 
Then if $f$ is $I$ left-redundant 
(respectively, $I$ right-redundant) 
and $g$ is 
$I'$ left-redundant 
(respectively, $I'$ right-redundant), 
then $gf$ is $I+I'$ 
left-redundant (respectively, $I+I'$ right-redundant). 
\end{lemma} 
\begin{proof} 
Let $\aseq$ and $\bseq$ be points in $X$ 
with $a_j,b_j\in L_1(X)$
such that  $(a_j)_{j\geq 0}=(b_j)_{j\geq 0}$. 
Put $a_{-I+j}\dots a_{-I+N+j}=v_j$ and 
$b_{-I+j}\dots b_{-I+N+j}=w_j$ for $j\geq 0$. Then, 
since $f$ is $I$ left-redundant, 
$(f(v_j))_{j\geq 0}=(f(w_j))_{j\geq 0}$. Since  
$g$ is $I'$ left-redundant, 
$g(f(v_{-I'})\dots f(v_{-I'+N'}))=g(f(w_{-I'})\dots f(w_{-I'+N'}))$ 
and hence 
$gf(a_{-I-I'}\dots a_{-I-I'+N+N'})=
gf(b_{-I-I'}\dots b_{-I-I'+N+N'})$. 
Therefore $gf$ is $I+I'$ left-redundant. 
\end{proof}

\begin{proposition}
Let $\phi:(X,\sigma_X)\to (Y,\sigma_Y)$ and 
$\psi:(Y,\sigma_Y)\to (Z,\sigma_Z)$ be homomorphisms between subshifts. 
Then 
\[P_L(\psi\phi)\geq P_L(\phi)+P_L(\psi)\q\text{and}\q  
P_R(\psi\phi)\geq P_R(\phi)+P_R(\psi).\] 
\end{proposition} 
\begin{proof} 
If $\phi$ is of $(m,n)$-type and given by 
a local rule $f:L_{m+n+1}(X)\to L_1(Y)$ and 
$\psi$ is of $(m',n')$-type and given by 
a local rule $g:L_{m'+n'+1}(Y)\to L_1(Z)$,   
then $\psi\phi$
is of $(m+m',n+n')$-type and given by $gf$. Suppose that 
$f$ is strictly $I$ left-redundant, $g$ is 
strictly $I'$ left-redundant and 
$fg$ is strictly $I''$ left-redundant. If $I\neq\infty$
and $I'\neq\infty$ 
then by Lemma 6.6  $I''\geq I+I'$. 
If $I=\infty$, then 
for every integer $J\geq 0$, $f$ is $J$ left-redundant, and 
for some integer $0\leq J'\leq I'$, $g$ is $J'$ left-redundant
and hence by Lemma 6.6 $gf$ is $J+J'$ left-redundant, 
which implies that
$gf$ is strictly $\infty$ left-redundant, i.e., $I''=\infty$. 
Similarly, if $I'=\infty$ then $I''=\infty$.   
Therefore we see that in any case, $I''\geq I+I'$. 
Hence we have 
$P_L(\psi\phi)= I''-(m+m')\geq I-m + I'-m'=P_L(\phi)+P_L(\psi)$. 

The proof of the second inequality is similar. 
\end{proof}

The following lemma is clear. 
\begin{lemma} Let $f:L_{N+1}(X_G)\to L_1(Y)$ be a local rule  
on a topological Markov shift $(X_G,\sigma_G)$ 
to a subshift $(Y,\sigma_Y)$. 
Let $0\leq I\leq N$.  Then $f$ is
$I$ left-redundant (respectively, $I$ right-redundant)
if and only if $f(w_1)=f(w_2)$ for any $w_1,w_2\in L_{N+1}(G)$ 
with the same terminal subpath 
(respectively, the same initial subpath) of length $N+1-I$ in $G$. 
\end{lemma} 

\begin{proposition} 
Let $\varphi$ be an endomorphism of a 
topological Markov shift $(X,\sigma)$. 
\begin{enumerate} 
\item  $P_L(\varphi)\geq 0$ (respectively, $P_R(\varphi)\geq 0$) 
if and only if  
$\varphi$ has memory zero (respectively, anticipation zero). 
\item  If $\varphi$ is onto, then 
$P_L(\varphi)\geq 0$ (respectively, $P_R(\varphi)\geq 0$) 
if and only if $\varphi$ is $p$-L (respectively, $p$-R). 
\end{enumerate}
\end{proposition}
\begin{proof} 
(1) Suppose that 
$\varphi$ is of $(m,n)$-type 
and given by a strictly $I$ left-redundant 
local rule $f: L_{m+n+1}(X)\to L_1(X)$ 
(recall that $I$ may be $\infty$).
If $m=0$, then $P_L(\varphi)=I-0\geq 0$. Conversely, if 
$P_L(\varphi)\geq 0$, then $I-m\geq 0$ and hence $f$ is 
$m$ left-redundant. By Lemma 6.8, 
we can decrease the left-redundancy of $f$ by $m$ 
to obtain a local rule $g: L_{n+1}(X)\to L_1(X)$ so that 
$\varphi$ may be of $(0,n)$-type and given by $g$. 
Therefore (1) is proved.

(2) By \cite[Proposition 5.1(1)]{Nasu-te} and (1).
\end{proof}

\begin{lemma} Let $(X,\sigma)$ be a subshift and 
$f: L_{N+1}(X)\to L_1(X)$ a local rule on it.  
Let $\varphi$ be 
an endomorphism of $(X,\sigma)$ given by $f$.  
\begin{enumerate}
\item 
If $f$ is an 
$I$ left-redundant (respectively, $I$ right-redundant) 
local rule on $(X,\sigma)$, 
then there exists $J\geq N+1$
such that for all $J'\geq J$, 
$f$ is an $I$ left-redundant 
(respectively, $I$ right-redundant) 
local rule on  
$(X^{\langle J'\rangle},\sigma_{X^{\langle J'\rangle}})$ to 
its image under $\phi^{\langle J',\varphi\rangle}$. 
\item 
If $P_L(\varphi)\geq 0$ (respectively, $P_R(\varphi)\geq 0$), 
then there exists $J\geq N+1$
such that $P_L(\phi^{\langle J',\varphi\rangle})\geq 0$ 
(respectively, 
$P_R(\phi^{\langle J',\varphi\rangle})\geq 0$) for all $J'\geq J$. 
\end{enumerate}
\end{lemma}
\begin{proof} 
(1) If $f$ is an $I$ left-redundant local rule on  
$(X^{\langle J\rangle},\sigma_{X^{\langle J\rangle}})$ to 
its image under $\phi^{\langle J,\varphi\rangle}$, then 
for all $J'\geq J$  
it is an $I$ left-redundant local rule on  
$(X^{\langle J'\rangle},\sigma_{X^{\langle J'\rangle}})$ to 
its image under $\phi^{\langle J',\varphi\rangle}$, because 
$X^{\langle J'\rangle}\subset X^{\langle J\rangle}$ if $J'\geq J$.
Therefore it suffices to show that there exists $J\geq N+1$
such that $f$ is an $I$ left-redundant local rule on  
$(X^{\langle J\rangle},\sigma_{X^{\langle J\rangle}})$ 
to its image under $\phi^{\langle J,\varphi\rangle}$. To do this,
assume the contrary. Then
there would exist no $J\geq N+1$ 
such that $f$ is 
$I$ left-redundant on 
$(X^{\langle J\rangle},\sigma_{X^{\langle J\rangle}})$ to its image 
under $\phi^{\langle J,\varphi\rangle}$.  
Then for each $s\geq N+1$, there would exist a pair of points 
$(c^{(s)}_j)_{j\in\Z}$ and $(d^{(s)}_j)_{j\in\Z}$ 
in $X^{\langle s\rangle}$ such that 
$(c^{(s)}_j)_{j\geq 0}=(d^{(s)}_j)_{j\geq 0}$
and 
$f(c^{(s)}_{-I}\dots c^{(s)}_{-I+N})
\neq f(d^{(s)}_{-I}\dots d^{(s)}_{-I+N})$. 
Since $X^{\langle s+1\rangle}\subset X^{\langle s\rangle}$ 
and $X=\cap_{s\geq N+1}X^{\langle s\rangle}$
with each $X^{\langle s\rangle}$ closed, 
a standard compactness argument shows that 
there would exist points 
$(c_j)_{j\in\Z}$ and $(d_j)_{j\in\Z}$ 
in $X$ such that 
$(c_j)_{j\geq 0}=(d_j)_{j\geq 0}$  and 
$f(c_{-I}\dots c_{-I+n})\neq f(d_{-I}\dots d_{-I+N})$. 
This contradicts the hypothesis that $f$ is $I$ left-redundant 
on $(X,\sigma)$. 

(2) Suppose that $\varphi$ is of $(m,n)$-type 
with $m+n=N$. Then for all $J\geq N+1$, 
$\phi^{\langle J,\varphi\rangle}$ is of $(m,n)$-type 
and given by $f$. Since $P_L(\varphi)\geq 0$, 
there exists $I\geq m$ such that $f$ is $I$ left-redundant 
on $(X,\sigma)$. 
Therefore, by (1) there exists $J\geq N+1$ 
such that for all $J'\geq J$, 
$f$ is an $I$ left-redundant local rule on  
$(X^{\langle J'\rangle},\sigma_{X^{\langle J'\rangle}})$ to 
its image under $\phi^{\langle J',\varphi\rangle}$ with $I\geq m$. 
From this, (2) follows. 
\end{proof} 
\begin{theorem} 
Let $\varphi$ be an endomorphism of a subshift $(X,\sigma)$. 
The following three conditions are equivalent:
\begin{enumerate} 
\item  $P_L(\varphi)\geq 0$ (respectively, $P_R(\varphi)\geq 0$);
\item  $\varphi$ is weakly $p$-L (respectively, weakly $p$-R);
\item  $\varphi$ has memory zero (respectively, anticipation zero). 
\end{enumerate}
\end{theorem}
\begin{proof}
By the definition of $P_L(\varphi)$ (3) implies (1) . 

To prove that (1) implies (3), 
suppose that $P_L(\varphi)\geq 0$. Let $\varphi$ be 
of $(m,n)$-type and be given by a local rule 
$f: L_{N+1}(X)\to L_1(X)$ on $(X,\sigma)$ with $N=m+n$. 
Then by Lemma 6.10(2), there exists $t\geq N+1$ such that 
$P_L(\phi^{\langle t,\varphi\rangle})\geq 0$. 
Put $\psi=\phi^{\langle t,\varphi\rangle}$ 
and put $(X_1,\sigma_1)=
(X^{\langle t\rangle},\sigma^{\langle t\rangle})$. 
Then by definition 
(see the first paragraph of Section 5)
$\psi$ is a factor map of  
$(X_1,\sigma_1)$ onto 
some sofic system $(Y,\sigma_Y)$
such that $\psi$ is an extension of 
$\varphi$, of $(m,n)$-type and given by 
the local rule $f_1$ on $(X_1,\sigma_1)$
to $(Y,\sigma_Y)$ that is identical with $f$ as a mapping
(note that $L_{N+1}(X_1)=L_{N+1}(X))$. Further 
$\psi^{[t]}$ is a factor map of $(X_1^{[t]},\sigma_1^{[t]})$ 
onto $(Y^{[t]}, \sigma_{Y^{[t]}})$ 
such that $\psi^{[t]}$ is an extension of $\varphi^{[t]}$, 
of $(m,n)$-type and given by the local rule 
$f_1^{[t]}$ on $(X_1^{[t]},\sigma_1^{[t]})$ to 
$(Y^{[t]},\sigma_{Y^{[t]}})$. 
By definition
$(X_1^{[t]},\sigma_1^{[t]})$ is a topological Markov shift (in fact,
the topological Markov shift with defining graph 
$G[X,t]$, which was defined in the first paragraph of Section 5). 
Therefore, since by Proposition 6.5 
$P_L(\psi^{[t]})=P_L(\psi)\geq 0$, 
it follows that $\psi^{[t]}$ is 
$m$ left-redundant. 
Hence by Lemma 6.8 
we can decrease the left-redundancy of $f_1^{[t]}$ by $m$ to 
have a local rule $g_1:L_{n+1}(X_1^{[t]})\to L_1(Y^{[t]})$ 
so that $\psi^{[t]}$ may be of $(0,n)$-type and given by $g_1$. 
Let $g:L_{n+1}(X^{[t]})\to g_1(L_{n+1}(X^{[t]}))$ be the restriction 
of $g_1$. Then, 
since $\varphi^{[t]}$ is a restriction of $\psi^{[t]}$, 
it follows that $g_1(L_{n+1}(X^{[t]}))=L_1(X^{[t]})$ and 
$\varphi^{[t]}$ is of $(0,n)$-type and $given$ by $g$. 
Let $g':L_{n+t}(X)\to L_1(X)$ be defined as follows: 
if $g((a_0\dots a_{t-1})(a_1\dots a_t)\dots(a_n\dots a_{n+t}))
= (b_0\dots b_{t-1})$, 
where $(a_0\dots a_{t-1})(a_1\dots a_t)\dots(a_n\dots a_{n+t})
\in L_{n+1}(X^{[t]})$
and $(b_0\dots b_{t-1})\in L_1(X^{[t]})$, with $a_i,b_i \in L_1(X)$, then 
$g'(a_0\dots a_{n+t})=b_0$. Then 
$\varphi$ is of $(0,n+t)$-type and given by $g'$. Therefore 
$\varphi$ has memory $0$ and (3) is proved. 

The proof of the equivalence of (2) and (3) is 
given by a straightforward 
modification in \cite[Proof of Proposition 7.5(1)]{Nasu-te}. 
We describe it for completeness. 
Suppose (3) and $\varphi$  
is of $(0,l)$ type 
and given by a 
local rule $f: L_{l+1}(X)\to L_1(X)$, then 
$\varphi=\varphi_{\hf{U}_{f,0,l}}$, where $\hf{U}_{f,0,l}$ 
is the half-textile-subsystem of the textile system
$T_{f,0,l}$ defined in the proof of Remark 2.6. 
Since $T_{f,0,l}$ is weakly $p$-L, $\varphi$ is 
weakly $p$-L, and hence (2) follows. 

To show that (2) implies (3), 
suppose $\varphi=\varphi_{\hf{U}}$, where $\hf{U}$ 
is a onesided 1-1 half-textile-subsystem of a weakly $p$-L
textile system $T=(p,q:\Gamma\to G)$.   
Since $\hf{U}$ is onesided 1-1, 
there exist 
$m,n\geq 1$ and a local rule 
$g:L_{m+n+1}(X_{\hf{U}})\to L_1(Z_{\hf{U}})$ such that  
$\xi_{\hf{U}}^{-1}((a_j)_{\in\Z})= 
(g(a_{j-m}\dots a_{j+n}))_{j\in\Z}$. 
Since $\xi_{\hf{U}}$ is a restriction of 
$\phi_p$ with $p$ weakly left resolving, for any 
$\alpha_{-m}\dots \alpha_n\in L_{m+n+1}(Z_{\hf{U}})$ 
with $\alpha_j\in A_\G$,
$\alpha_{-m}\dots\alpha_{-1}$ is  
determined by $p_A(\alpha_{-m})\dots p_A(\alpha_{-1})$ 
and $t_\Gamma(\alpha_{-1})=i_\Gamma(\alpha_0)=
i_\Gamma(g(p_A(\alpha_{-m})\dots p_A(\alpha_{n})))$. 
Hence we can define a local rule 
$g':L_{m+n+1}(X_{\hf{U}})\to L_1(Z_{\hf{U}})$ such that  
$g'(p_A(\alpha_{-m})\dots p_A(\alpha_n))=\alpha_{-m}$, 
with which $\xi_{\hf{U}}^{-1}$ is of $(0,m+n+1)$ type.  
Since $\varphi_{\hf{U}}$ is of $(0,m+n+1)$ type given by 
the local rule $q_Ag'$, (3) follows. 
\end{proof} 

\begin{proposition} 
Let $\varphi$ be an endomorphism of a subshift $(X,\sigma)$. 
Let $s\in\Z$.  
\begin{enumerate} 
\item 
$\varphi\sigma^s$ has memory (respectively anticipation) zero 
if and only if $s\geq -P_L(\varphi)$ 
(respectively, $s\leq P_R(\varphi)$)
\item Suppose that $\varphi$ is onto. 
$\varphi\sigma^s$ is weakly $p$-L (respectively, weakly $p$-R)
if and only if $s\geq -P_L(\varphi)$ 
(respectively, $s\leq P_R(\varphi)$)
; if $s> -P_L(\varphi)$ (respectively, $s< P_R(\varphi)$), then 
$\varphi\sigma^s$ is right $\sigma$-expansive 
(respectively, left $\sigma$-expansive) on the upper side. 
\item 
If $(X,\sigma)$ is a topological Markov shift, then 
the statements in (2) with all ``weakly'' deleted in them hold. 
\end{enumerate}
\end{proposition} 
\begin{proof} 
(1) By Theorem 6.11 and Proposition 6.4, (1) holds. 

(2) By Theorem 6.11 and Proposition 6.4 
$\varphi\sigma^s$ is weakly $p$-L 
if and only if $s\geq -P_L(\varphi)$. 
If $s> -P_L(\varphi)$, then 
$P_L(\varphi\sigma^{s-1})\geq 0$, and hence 
$\varphi\sigma^{s-1}$ is weakly $p$-L. Therefore 
$\varphi\sigma^s$ is right $\sigma$-expansive 
on the upper side, 
by Proposition 2.8. 

(3) By a similar proof to that of (1), but use 
Propositions 6.9 and 2.11 instead of Theorem 6.11 and 
Proposition 2.8.
\end{proof} 

\begin{proposition} Let $\varphi$ be an endomorphism of 
a subshift $(X,\sigma)$. 
\begin{enumerate}
\item[(1)] If $\varphi^i(X)$ is infinite for all $i\geq 0$, then 
the following hold:
\begin{enumerate} 
\item[(a)]
$P_L(\varphi)+P_R(\varphi)$ is nonpositive, 
and hence if in addition, $\varphi$ is 
of $(m,n)$-type then 
\[-n\leq P_R(\varphi)\leq -P_L(\varphi)\leq m;\] 
\item[(b)]
$\varphi$ is bi-recognizing (i.e. left-recognizing and right-recognizing); 
\item[(c)]
if $\varphi$ is given by a strictly $I$ left-redundant, 
strictly $J$ right-redundant local rule $f:L_{N+1}(X)\to L_1(X)$, 
then $I+J\leq N$. 
\end{enumerate}
\item[(2)]
If $\varphi^k(X)$ is finite with $k\geq 0$, then for all $i\geq k$ 
$P_L(\varphi^i)=\infty$ and $P_R(\varphi^i)=\infty$. 
\end{enumerate}
\end{proposition} 
\begin{proof}   
(1)(a) Assume that $P_L(\varphi)+P_R(\varphi)>0$. 

If $P_L(\varphi)\in\Z$, then let 
$\psi=\varphi\sigma^{-P_L(\varphi)}$. Then by Proposition 6.4, 
$P_L(\psi)=0$  and $P_R(\psi)=P_L(\varphi)+P_R(\varphi)>0$. 
Hence, by Theorem 6.11 $\psi$ is an endomorphism 
of $(0,n_1)$-type and of $(m_1,0)$-type of $(X,\sigma)$ with some 
$m_1,n_1\geq 0$. 
Therefore $\psi^{[m_1+n_1]}$ is an endomorphism of $(0,0)$-type 
of $(X^{[m_1+n_1]},\sigma^{[m_1+n_1]})$. Let 
$g:L_1(X^{[m_1+n_1]})\to L_1(X^{[m_1+n_1]})$ be the local rule giving $\psi$. 
Then there exists $k\geq 1$ such that 
$g^k(L_1(X^{[m_1+n_1]}))=g^{k+1}(L_1(X^{[m_1+n_1]}))$. 
We see that
$(\psi^{[m_1+n_1]})^k(X^{[m_1+n_1]})=
(\psi^{[m_1+n_1]})^{k+1}(X^{[m_1+n_1]})$. 
Let $B=g^k(L_1(X^{[m_1+n_1]}))$ and $g'=g|B$. Let 
$Y=(\psi^{[m_1+n_1]})^k(X^{[m_1+n_1]})$ and $\psi'=\psi^{[m_1+n_1]}|Y$. 
Then $\psi'$ is an onto symbolic endomorphism  
and hence a symbolic automorphism of $(Y',\sigma_{Y'})$ 
given by the local rule $g':L_1(Y')\to L_1(Y')$. 
It follows by hypothesis that $Y'$ is infinite. Therefore, since 
$g'$ is a permutation on $L_1(Y')$ 
it follows that $g'$ is strictly 
zero right-redundant and hence so is $g$. 
Therefore $P_R(\psi^{[m_1+n_1]})=0$, 
and hence by Proposition 6.5 
$P_L(\varphi)+P_R(\varphi)=P_R(\psi)=0$, which is  
contrary to the assumption. 

If $P_R(\varphi)\in\Z$, then let $\psi=\varphi\sigma^{P_R(\varphi)}$. 
Then $P_R(\psi)=0$ and $P_L(\psi)=P_L(\varphi)+P_R(\varphi)>0$. 
Hence, by a similar argument to the above, 
we are led to a contradiction.

If $P_L(\varphi)=\infty$ and $P_R(\varphi)=\infty$, then since
$P_L(\varphi)>0$ and $P_R(\varphi)>0$, we also see using Theorem 6.11 
that $\varphi$ would be of $(0,0)$-type up to 
higher block conjugacy, which leads us to a contradiction. 

Hence $P_L(\varphi)+P_R(\varphi)\leq 0$, which with (6.3) proves (1). 

(1)(b)  By (1)(a). 

(1)(c) Suppose that $\varphi$ is 
of $(m,n)$-type and given by $f$ with $m+n=N$. 
It follows from (1)(a) that 
$(I-m)+(J-n)=P_L(\varphi)+P_R(\varphi)\leq 0$, and hence 
$I+J\leq N$. 

(2) By the fact stated in the paragraph before Definition 6.1. 
\end{proof} 

A direct proof of Proposition 6.13(1)(a) for an onto endomorphism of 
an infinite topological Markov shift is given 
by using Proposition 6.9 instead of Theorem 6.11.

\begin{remark} Let $\varphi$ be an endomorphism of 
a subshift. 
\begin{enumerate}
\item If $X$ is finite, then $P_L(\varphi)=\infty$, 
$P_R(\varphi)=\infty$, $Q_R(\varphi)=\infty$ and $Q_L(\varphi)=\infty$. 
\item 
If $X$ is infinite, then
$P_L(\varphi)=\infty$ (respectively, $P_R(\varphi)=\infty$) 
if and only if $\varphi$ is not left-recognizing 
(respectively, not right-recognizing) 
or equivalently, $\varphi$ collapses 
every pair of distinct right 
(respectively, left) $\sigma$-asymptotic points. 
\item 
$Q_L(\varphi)=-\infty$ (respectively, $Q_R(\varphi)=-\infty$) 
if and only if $X$ is infinite and $\varphi$ is not left-closing 
(respectively, not right-closing) 
or equivalently, $\varphi$ collapses 
some pair of distinct right 
(respectively, left) $\sigma$-asymptotic points. 
\item 
$P_L(\varphi)+Q_L(\varphi)$ (respectively, 
$P_R(\varphi)+Q_R(\varphi)$) exists if and only if
either $\varphi$ is left-recognizing 
(right-recognizing) or $X$ is finite.
\item
$P_L(\varphi)+Q_L(\varphi)$ and $P_R(\varphi)+Q_R(\varphi)$ 
are shift-invariant, if they exist. 
\end{enumerate}
\end{remark}
\begin{proof} Statements (1), (2), (3) are valid by definition. 
Since $P_L(\varphi)+Q_R(\varphi)$ does not exist if and only if 
$P_L(\varphi)=\infty$ and $Q_L(\varphi)=-\infty$, (4) follows 
from (1),(2),(3).  By Propositions 3.4 and 6.4, (5) follows.
\end{proof} 

The following proposition is a counterpart of Proposition 5.4.
\begin{proposition} 
If right-recognizing (respectively, left-recognizing) 
endomorphisms $\varphi$ and $\psi$ of 
infinite subshifts are topologically-conjugate, then 
there exists an integer 
$s\geq 0$ such that for all $i\geq 0$, 
$|P_L(\varphi^i)-P_L(\psi^i)|\leq s$ (respectively,
$|P_R(\varphi^i)-P_R(\psi^i)|\leq s$). 
\end{proposition}
\begin{proof} 
The proposition is proved in a similar manner to
the proof of Proposition 5.4 
by using Proposition 6.7 and (6.1) and by using the fact that 
topological conjugacy $\theta$ between subshifts are 
biclosing, hence bi-recognizing 
and hence 
$P_L(\theta),P_L(\theta^{-1}), P_R(\theta),P_R(\theta^{-1})\in\Z$.
\end{proof}

\begin{proposition} Let $\varphi$ be an endomorphism 
of a subshift. 
\begin{enumerate} 
\item 
If $\varphi$ essentially has nonnegative $p$-L degree  
(respectively, nonnegative $p$-R degree ), 
then there exists an integer $s\geq 0$ 
such that for all $i\geq 0$, $P_L(\varphi^i)\geq -s$ 
(respectively, $P_R(\varphi^i)\geq -s$). 
\item
If $\varphi$ essentially has memory (respectively, anticipation)
zero, then there exists an integer $s\geq 0$ 
such that for all $i\geq 0$, 
$\varphi^i$ has memory (respectively, anticipation) $s$.  
\item 
If $\varphi$ is essentially weakly $p$-L (respectively, 
essentially weakly $p$-R), 
then there exists an integer $s\geq 0$ such that for all $i\geq 0$,
$\varphi^i\sigma^s$ is weakly $p$-L (respectively, 
$\varphi^i\sigma^{-s}$ is weakly $p$-R). 
\item 
If $\varphi$ is essentially weakly $p$-L (respectively, 
essentially weakly $p$-R), then there exists $s\geq 0$ such that
for any $x,y\in X$ it holds that 
if $x$ and $y$ coincide in the $k$-th coordinate 
(i.e., $x=\aseq$ and $y=\bseq$ with $a_k=b_k$) 
for all $k\geq 0$ (respectively, for all $k\leq 0$), 
then for all $i\geq 1$, $\varphi^i(x)$ and $\varphi^i(y)$ coincide
in the $k$-th coordinate for all $k\geq s$ (respectively, $k\leq -s$). 
\item
If in addition,
$(X,\sigma)$ is a topological Markov shift, then  
(3) and (4) with all ``weakly'' in them deleted hold. 
\end{enumerate}
\end{proposition}
\begin{proof} The proof of (1) 
is given in a similar manner to the 
proof of Proposition 5.5(1) 
by using Propositions 6.15 and 6.7. 

We see that (2), (3) and (4) are proved by (1) and Theorem 6.11.
By (1) and Proposition 6.9, (5) is proved. 
\end{proof} 

\begin{lemma} Let $m,m', n,n'\geq 0$. 
If a homomorphism  $\phi$ of a topological Markov shift 
$(X,\sigma)$ into a subshift $(Y,\sigma_Y)$ is 
of $(m,n')$-type and also of $(m',n)$-type, then 
$\phi$ is of $(m,n)$-type. 
\end{lemma}
\begin{proof} 
Suppose that 
$\phi$ is of $(m,n')$-type with 
a local rule $f_1:L_{m+n'+1}(X)\to L_1(Y)$ and that
$\phi$ is of $(m',n)$-type 
with a local rule $f_2:L_{m'+n+1}(X)\to L_1(Y)$. 
We may assume that $m'\geq m$ and $n'\geq n$ 
(by adding right redundancy by $n-n'$ to $f_1$ if $n> n'$, 
and by adding left redundancy by $m-m'$ to $f_2$ if $m>m'$).  
Let $g_1: L_{m'+n'+1}(X)\to L_1(Y)$ and 
$g_2: L_{m'+n'+1}(X)\to L_1(Y)$ be the local rules defined by
$g_1(w'ww'')=f_1(ww'')$ and $g_2(w'ww'')=f_2(w'w)$, 
where $w'ww''\in L_{m'+n'+1}(X)$ with $w'\in L_{m'-m}(X)$ 
and $w''\in L_{n'-n}(X)$.  
Then we see that $g_1=g_2$ as a mapping, 
because $g_1$ and $g_2$
give the same block-map $\phi$ of $(m',n')$-type. 
Let $w\in L_{m+n+1}(X)$ and let $w'_1,w'_2\in L_{m'-m}(X)$ 
with $w'_1w, w'_2w \in L_{m'+n+1}(X)$. Since 
$(X,\sigma)$ is a topological Markov shift, there exists 
$w''\in L_{n'-n}(X)$ such that 
$w'_1ww'',w'_2ww''\in L_{m'+n'+1}(X)$. We see that 
$f_2(w'_1w)=g_2(w'_1ww'')=g_1(w'_1ww'')=f_1(ww'')
=g_1(w'_2ww'')=g_2(w'_2ww'')=f_2(w'_2w)$.
Therefore $f_2$ is $m'-m$ left-redundant. Let
$f:L_{m+n}(X)\to L_1(Y)$ be the local rule 
obtained from $f_2$ by deleting left redundancy by $m'-m$. 
Then  $\phi$ is of $(m,n)$-type and given by $f$. 
\end{proof} 

\begin{proposition}
Let $\varphi$ be an 
endomorphism of a subshift $(X,\sigma)$. 
\begin{enumerate} 
\item 
The following conditions are pairwise equivalent 
(the equivalence of (a),(b) and (c) 
is due to K\uo rka \cite{Kur1}): 
\begin{enumerate}
\item[(a)] $\varphi$ is uniformly equicontinuous;
\item[(b)] there exist $m,n\geq 0$ such that for all $i\geq 0$ 
$\varphi^i$ is of $(m,n)$-type;
\item[(c)] there exist $r,s\geq 0$ such that $\varphi^{r+s}=\varphi^r$; 
\item[(d)] $\varphi$ is an essentially symbolic 
endomorphism of $(X,\sigma)$.
\end{enumerate}
\item 
If $\varphi$ is onto,  
then (b) (as a representative of (a),(b),(c) and (d)) 
is equivalent to each of: 
\begin{enumerate} 
\item[(e)] 
there exists $s\geq 0$ such that $\varphi^s=i_X$; 
\item[(f)]
$\varphi$ is an essentially symbolic automorphism of $(X,\sigma)$;
\item[(g)]
$\varphi$ is an essentially weakly $p$-biresolving endomorphism    
of $(X,\sigma)$. 
\end{enumerate} 
\item 
If $(X,\sigma)$ is an SFT, then (b) 
is equivalent to: 
\begin{enumerate}
\item[(h)] $\varphi$ essentially has memory zero and essentially 
has anticipation zero. 
\end{enumerate}
\item 
If $(X,\sigma)$ is an SFT and $\varphi$ is onto, then (b) 
is equivalent to each of: 
\begin{enumerate}
\item[(i)] $\varphi$ is essentially $p$-L and essentially $p$-R; 
\item[(j)] $\varphi$ is essentially $p$-biresolving. 
\end{enumerate}
\end{enumerate}
\end{proposition}
\begin{proof} 
(1) It is direct that (a) is equivalent to (b).  
If $\varphi^i$ is of $(m,n)$-type 
and given by a local rule $f_i:L_{m+n+1}(X)\to L_1(X)$ for $i\geq 0$, 
then there exist $r,s\geq 0$ such that $f_{r+s}=f_r$ and 
hence $\varphi^{r+s}=\varphi^r$. Hence (b) implies (c). 
We see that (c) implies (d) because 
$\varphi^{[^*r+s]}$ is a symbolic endomorphism of 
$(X^{[^*_\varphi r+s]},\sigma^{[^*_\varphi r+s]})$ with the notation 
defined in Subsection 8.1.  (A similar proof 
for proving that (e) implies (f) appears in \cite[Proof of 
Proposition 8.1]{Nasu-t}.) 
If $\varphi=\theta^{-1}\varphi_0\theta$, where $\varphi_0$ 
is a symbolic endomorphism of 
a subshift $(X_0,\sigma_0)$ and $\theta:(X,\sigma)\to(X_0,\sigma_0)$ 
is a conjugacy, then $\varphi$ clearly satisfies (b). 
Hence (d) implies (b). 

(2) Since $\varphi$ is onto, it follows that if (c) holds 
then all $\varphi$-orbits have period $s$, and hence (e) is 
equivalent to (c). Since $\varphi$ is onto, (d) is equivalent to (f). 
The equivalence of (f) and (g) follows from Proposition 
7.1 (the equivalence of (c) and (d) in it). 
Therefore (2) follows from (1). 

(3) Suppose (h). Since $(X,\sigma)$ is an SFT, there 
exists a topological Markov shift $(\bar{X},\bar{\sigma})$
and a conjugacy $\theta:(X,\sigma)\to (\bar{X},\bar{\sigma})$.
Since $\theta\varphi\theta^{-1}$ essentially has memory zero 
and essentially has anticipation zero, 
it follows from Proposition 6.16(2) 
that there exist $\bar{m},\bar{n}\geq 0$ 
such that for all $i\geq 0$, 
$(\theta\varphi\theta^{-1})^i$ is of $(\bar{m},n_i)$-type
and also $(m_i,\bar{n})$-type for some $m_i,n_i\geq 0$.  
Therefore, by Lemma 6.17 $\theta\varphi^i\theta^{-1}$ 
is of $(\bar{m},\bar{n})$-type for all $i\geq 1$, and 
hence (b) follows.

Therefore, since (d) implies (h), (3) follows from (1). 

(4) Since $(X,\sigma)$ is an SFT and $\varphi$ is onto, 
it follows from Proposition 6.9 that (h) and (i) 
are equivalent. 
By Remark 2.10(1), (g) and (j) are equivalent. Hence 
(4) follows from (2) and (3). 
\end{proof}

By Proposition 6.18(4), we have covered 
\cite[Proposition 8.6(1)]{Nasu-te} 
(i.e. the equivalence of (i) and (j)
in Proposition 6.18(4)) 
by a much easier 
proof than the proof given in \cite{Nasu-te}  
by using results \cite[Theorems 7.19 and 7,20]{Nasu-t} on 
resolvable textile systems, which need very long proofs  
in \cite{Nasu-t}. However, we cannot give an easier proof 
to the following, which was proved 
in \cite[Proposition 6.19(2)(3)]{Nasu-te}
by using \cite[Theorems 7.19 and 7,20]{Nasu-t}. 
\begin{proposition} [\cite{Nasu-te}]
Let $\varphi$ be an onto endomorphism of an SFT $(X,\sigma)$. 
\begin{enumerate} 
\item 
If $\varphi$ is essentially $p$-L and essentially $q$-R, 
then $\varphi$ is essentially LR. 
\item 
If $\varphi$ is essentially $q$-L and essentially $q$-R, 
then $\varphi$ is essentially $q$-biresolving. 
\end{enumerate} 
\end{proposition} 

Moreover, we cannot answer the question whether or not 
an essentially weakly $p$-L, 
essentially weakly $p$-R onto endomorphism of a subshift
is a essentially weakly $p$-biresolving endomorphism 
of the subshift, 
the question whether or not 
an essentially weakly $p$-L, 
essentially weakly $q$-R onto endomorphism of a subshift
is a essentially weakly LR endomorphism 
of the subshift, and so on.  

\section{Endomorphisms with two nonnegative degrees} 

By Theorems 4.6 and 5.3, 
we know what are the endomorphisms of the shift with 
both $q$-R and $q$-L degrees nonnegative. 
In this section, we study further the endomorphisms of the 
shift with two nonnegative degrees of onesided resolvingness. 

\begin{proposition}
Let $\varphi$ be an endomorphism of 
a subshift $(X,\sigma)$. 
\begin{enumerate}
\item The following conditions are equivalent:
  \begin{enumerate}
  \item[(a)]
  $P_L(\varphi)\geq 0$ and $P_R(\varphi)\geq 0$; 
  \item[(b)]
  $\varphi$ is weakly $p$-L and weakly $p$-R;
  \item[(c)] 
  $\varphi^{[t]}$ is a symbolic endomorphism  
  of $(X^{[t]},\sigma^{[t]})$ with some $t\geq 0$; 
  \item[(d)] 
  $\varphi^{[t]}$ is a weakly $p$-biresolving endomorphism 
  of $(X^{[t]},\sigma^{[t]})$ with some $t\geq 0$. 
\end{enumerate}
\noindent If in addition, $\varphi^i(X)$ is infinite for 
all $i\geq 0$, then (a) implies 
that $P_L(\varphi)= P_R(\varphi)= 0$. 

\item If $\varphi$ is onto, then all ``endomorphism'' can be 
changed to ``automorphism'' in (c) and (d), and 
if in addition, $(X,\sigma)$ is a topological Markov shift, 
then in addition, all ``weakly'' 
can be deleted from (b) and (d). 
\end{enumerate}
\end{proposition}
\begin{proof}
(1) By Theorem 6.11, (a) implies (b). 
If (b) holds, then 
by Theorem 6.11 
$\varphi$ is of $(0,n)$-type and of 
$(m,0)$ type with some $m,n\geq 0$ and 
therefore $\varphi^{[m+n]}$ 
is a 1-block endomorphism 
of $(X^{[m+n]},\sigma^{[m+n]})$, which proves (c). 
By Remark 2.14(2) (c) implies (d). 
By Theorem 6.11 and Proposition 6.5, (d) implies (a). 
The remainder of (1) follows from Proposition 6.13(1). 

(2) Since an onto symbolic endomorphism of a subshift is an symbolic 
automorphism of the subshift, the first part of (2) holds. To prove 
the remainder, use Proposition 6.9 instead of Theorem 6.11 and 
use Remark 2.14(1) instead of Remark 2.14(2) in the proof of (1). 
\end{proof} 

\begin{proposition} Suppose that 
$\phi:(X,\sigma_X)\to (Y,\sigma_Y)$ is a factor map of an SFT onto 
a sofic system. 
\begin{enumerate}
\item 
If $\phi^{[t]}$ is weakly LR 
(respectively, weakly RL) with some $t\geq 1$, then 
$P_L(\phi)\geq 0$ and $Q_R(\phi)\geq 0$ (respectively, 
$P_R(\phi)\geq 0$ and $Q_L(\phi)\geq 0$). 
\item 
If $P_L(\phi)\geq 0$ and $Q_R(\phi)\geq 0$ (respectively, 
$P_R(\phi)\geq 0$ and $Q_L(\phi)\geq 0$), 
then for some $t\geq 1$ (which we can find),  
we can construct 
a onesided 1-1, weakly LR 
(respectively, weakly RL) textile relation system 
$T=(p:\G\to G_A^{[t]},q:\G\to G_B^{[t]})$  
with $\phi_T=\phi^{[t]}$, where $A=L_1(X)$ and $B=L_1(Y)$, 
and moreover, we have a onesided 1-1, weakly LR 
(respectively, weakly RL) textile relation system 
$T_s=(p_s:\G\to G_A^{[t-s]},q_s:\G\to G_B^{[t-s]})$  
with $\phi_{T_s}=(\phi\sigma_X^s)^{[t-s]}$ for s=1,\dots,t-1. 
\end{enumerate}
\end{proposition}
\begin{proof} 
(1) By straightforward modification in 
\cite[Proof of Proposition 5.1(1)]{Nasu-te}, it is 
proved that
any weakly $p$-L factor map $\psi$ of 
an SFT onto a sofic system has memory zero and hence 
$P_L(\psi)\geq 0$. Therefore, (1) follows from 
Propositions 4.1, 3.5 and 6.5. 

(2) Let $J\geq 1$ such that $(X^{[J]},\sigma_{X^{[J]}})$ is a topological 
Markov shift and let $G_J$ be its defining graph. 
Suppose that $\phi$ is of $(m,n)$-type and given by a 
strictly $k$ right-mergible local rule $f$ of neighborhood-size $m+n$. 
Since by Proposition 6.5 $P_L(\phi^{[J]})=P_L(\phi)\geq 0$, 
it follows that $\phi^{[J]}$ is $m$ left-redundant.
Therefore by Lemma 6.8 
we can decrease the left redundancy of $f^{[J]}$ by $m$ to obtain 
a local rule $g: L_{n+1}(G_J)\to L_1(G_B^{[J]})$ 
such that $\phi^{[J]}$ is of $(0,n)$-type and given by $g$. 
Since by Proposition 3.5 $f^{[J]}$ is strictly $k$ right-mergible, 
so is $g$ (by Lemma 3.2 when $X$ is infinite). 
We have $n-k\geq 0$, because $n-k=Q_R(\phi)\geq 0$ 
if $X$ is infinite, and otherwise, $k=0$. 

Let $G=G_A^{[J]}$ and let $H=G_B^{[J]}$. 
Then $\phi^{[J]}$ is of $(0,n)$-type given by 
the $k$ right-mergible local rule 
$g:L_{n+1}(G)\to L_1(H)$ with $n\geq k$. 
We use the graph-homomorphism 
\[q^+_{g;k}:G^+_{g;k}\to H^{[k+1]}\]
(see Subsection 2.2).  Each arc in $G^+_{g;k}$ 
is written in the form
\[D^+_{g;k}(\alpha_1\dots \alpha_{n+k+1}),\] 
where $\alpha_1\dots \alpha_{n+k+1}\in L_{k+n+1}(G)$ 
with $\alpha_j\in A_G$.  
Arc $D^+_{g;k}(\alpha_1\dots \alpha_{n+k+1})$ goes from 
vertex $D^+_{g;k}(\alpha_1\dots \alpha_{n+k})$ to vertex 
$D^+_{g;k}(\alpha_2\dots \alpha_{n+k+1})$ with 
\[q^+_{g;k}(D^+_{g;k}(\alpha_1\dots \alpha_{n+k+1}))
=g(\alpha_1\dots \alpha_{n+k+1}).\]
By Proposition 2.4(1), $q^+_{g;k}$ is weakly right-resolving. 
Since $\alpha_1\dots \alpha_{n+1}$ is uniquely determined by 
arc $D^+_{g;k}(\alpha_1\dots \alpha_{n+k+1})$ and $n\geq k$, 
we can define another graph-homomorphism 
$p:G^+_{g;k}\to G^{[k+1]}$ by
\[p(D^+_{g;k}(\alpha_1\dots \alpha_{n+k+1}))=\alpha_1\dots \alpha_{k+1}.\] 
Since $\alpha_1\dots \alpha_{k+1}$ 
and $D^+_{g;k}(\alpha_2\dots \alpha_{n+k+1})$ uniquely 
determine $D^+_{g;k}(\alpha_1\dots \alpha_{n+k+1})$, we see that 
$p$ is weakly left-resolving. 
Thus we have a weakly LR textile-relation system 
\[T=(p:G^+_{g;k}\to G^{{[k+1]}}, 
q^+_{g;k}:G^+_{g;k}\to H^{[k+1]}).\]
Clearly $T$ is onesided 1-1 
and $\phi_T=(\phi^{[J]})^{[k+1]}=\phi^{[J+k]}$. 

Let $t=J+k$ and let $0\leq s\leq t-1$. From $T$, we have a 
textile relation system 
$T_s=(p_s:\G\to G_A^{[t-s]},q_s:\G\to G_B^{[t-s]})$ 
such that  $\G=G^+_{g;k}$ and for 
$W\in L_{n+k+1}(G)= L_{n+k+1}(G_A^{[J]})$, 
if 
$p(D^+_{g;k}(W)) 
=\alpha_1\dots \alpha_{k+1}$ 
and $\alpha_j=a_j\dots a_{j+J-1}$ 
for $j=1,\dots,k+1$ with $a_1,\dots,a_t\in A$, 
then $p(D^+_{g;k}(W))=a_1\dots a_{t-s}$, and if 
$q^+_{g;k}(D^+_{g;k}(W)=\beta_1\dots \beta_{k+1}$ 
and $\beta_j=b_j\dots b_{j+J-1}$ 
for $j=1,\dots,k+1$ with $b_1,\dots,b_t\in B$, 
then $p(D^+_{g;k}(W))=b_{s+1}\dots b_{t}$. 
We easily see that $T_s$ is onesided 1-1 and weakly LR 
and $\phi_{T_s}=(\phi\sigma_X^s)^{[t-s]}$.
Therefore (2) is proved. 
\end{proof} 

\begin{theorem} 
Let $\varphi$ be an onto endomorphism of a 
topological Markov shift $(X_G,\sigma_G)$ with 
$\varphi$ one-to-one or $G$ irreducible. Suppose that
$P_L(\varphi)\geq 0$ and  $Q_R(\varphi)\geq 0$  
(respectively, $P_R(\varphi)\geq 0$ and $Q_L(\varphi)\geq 0$). 
\begin{enumerate}
\item If $\varphi$ is of $(m,n)$-type 
(respectively, of $(n,m)$ type) and given by 
a local rule $f:L_{m+n+1}(G)\to A_G$ which is 
strictly $k$ right-mergible (respectively, strictly 
k left-mergible), 
then $k\leq n$ and 
$\varphi$ is of $(0,n)$-type (respectively, $(n,0)$-type) 
and given by the strictly $k$ right-mergible 
(respectively, strictly $k$ left-mergible) local rule 
$g:L_{n+1}(G)\to A_G$ obtained from $f$ by deleting 
left-redundancy (respectively, right-redundancy) by $m$. 
\item
Suppose that $\varphi$ is 
of $(0,n)$-type (respectively, $(n,0)$-type) 
and given by a
$k$ right-mergible 
(respectively, k left-mergible) local rule 
$g:L_{n+1}(G)\to A_G$ with $k\leq n$. 
Let
$T=(p,q^+_{g;k}:G^+_{g;k}\to G^{[k+1]})$ (respectively, 
$T=(p,q^-_{g;k}:G^-_{g;k}\to G^{[k+1]})$, where $p$ is defined 
as follows: for $w\in L_{n+k+1}(G)$ 
$p(D^+_{g;k}(w))$ (respectively, $p(D^-_{g;k}(w)$) 
is the initial (respectively, terminal) subpath of length $k+1$ of $w$. 
Then $T$ is a onesided 1-1 and  LR (respectively, RL) textile system 
with $\varphi_T=\varphi^{[k+1]}$. 
Moreover $(\varphi\sigma_G^s)^{[k+1-s]}$ 
(respectively, $(\varphi\sigma_G^{-s})^{[k+1-s]}$) is an LR 
(respectively, RL) endomorphism 
of $(X_G,\sigma_G)$ for $s=1,\dots,k$. 
\item
We can decide the expansiveness situation of $\varphi$. 
\end{enumerate}
\end{theorem} 
\begin{proof} 
For the proof Proposition 7.2(2), 
consider the case 
that $(X,\sigma_X)=(Y,\sigma_Y)=(X_G,\sigma_G)$ 
and $\phi=\varphi$. Then we can take 
$J=1$ and $H=G$ and the proof thus modified shows that 
(1) is valid and also shows, 
under the hypothesis of (2), that for $s=0,\dots, k$ we can construct  
a onesided 1-1, weakly LR textile system $T_s$ with $T_0=T$ 
such that $\phi_{T_s}=(\varphi\sigma_G^s)^{[t-s]}$. 
Using Lemma 2.1, we see that $T_s$ is LR for $s=0,\dots,k$. 
Hence (2) is proved. 

Since $T^*$ is LR, we can decide whether $\xi_{T^*}$ is 
one-to-one or not, and  whether $\eta_{T^*}$ is 
one-to-one or not. 
(This follows from \cite[Lemma 6.25]{Nasu-t}.)  Hence we can 
decide the expansiveness situation of $\varphi^{[k+1]}$ and 
hence of $\varphi$, by 
\cite[Theorem 2.5]{Nasu-t} and \cite[Proposition 6.2]{Nasu-te}. 
Hence (3) is proved. 
\end{proof} 

For an onto endomorphism $\varphi$ of a subshift , define
\begin{align*}
C_R(\varphi)&=\max\{-P_L(\varphi), -Q_R(\varphi)\},\q\q
C_L(\varphi)=\min\{P_R(\varphi), Q_L(\varphi)\},\\ 
D_R(\varphi)&=\min\{-P_L(\varphi), -Q_R(\varphi)\},\q\q
D_L(\varphi)=\max\{P_R(\varphi), Q_L(\varphi)\}. 
\end{align*} 

\begin{corollary} 
Let $\varphi$ be an onto endomorphism of 
an infinite topological Markov shift $(X,\sigma)$ 
such that $\varphi$ is one-to-one or 
$\sigma$ is topologically transitive. Let $s\in\Z$. 
\begin{enumerate} 
\item  Suppose that $\varphi$ is 
given by a strictly $k$ right-mergible local rule with $k\neq\infty$. 
\begin{enumerate}
\item If $s<C_R(\varphi)$, then there exists no $t\geq 1$ with  
$(\varphi\sigma^s)^{[t]}$ is LR.
\item If $C_R(\varphi)\leq s\leq C_R(\varphi)+k$, then 
$(\varphi\sigma^s)^{[k+1+C_R(\varphi)-s]}$ is LR;
\item If $s\geq k+C_R(\varphi)$ then $\varphi\sigma^s$ is LR. 
\item We can decide the expansiveness situation of 
$\varphi\sigma^{C_R(\varphi)}$. 
\item If $s>C_R(\varphi)$, then $\varphi\sigma^s$ is expansive. 
\end{enumerate}
\item Suppose that $\varphi$ is 
given by a strictly $l$ left-mergible local rule with $l\neq\infty$. 
\begin{enumerate}
\item If $s>C_L(\varphi)$, then there exists no $t\geq 1$ with  
$(\varphi\sigma^s)^{[t]}$ is RL.
\item If $-l+C_L(\varphi)\leq s\leq C_L(\varphi)$, then 
$(\varphi\sigma^s)^{[l+1-C_L(\varphi)+s]}$ is $RL$. 
\item If $s\leq -l+C_L(\varphi)$, then $\varphi\sigma^s$ is RL. 
\item We can decide the expansiveness situation of 
$\varphi\sigma^{C_L(\varphi)}$. 
\item If $s<C_L(\varphi)$, then $\varphi\sigma^s$ is expansive.
\end{enumerate}
\end{enumerate} 
\end{corollary}
\begin{proof} 
(1) Statement (a) follows 
from Propositions 4.1(1), 3.5(3), 6.12(3) and 6.5(2). 

To prove (b), Let $\psi=\varphi\sigma^{C_R(\varphi)}$. 
Then $\psi$ is given by a strictly $k$ right mergible local 
rule. (For suppose that $\varphi$ is of $(m_0,n_0)$-type and 
given by a strictly $k$ right-mergible local rule 
$f:L_{m_0+n_0+1}(X)\to L_1(X)$. Then, since $n_0=k+Q_R(\varphi)$
(when $X$ is infinite; recall Standing convention II for proofs), 
we have
$C_R(\varphi)=-n_0+k+Q_R(\varphi)+C_R(\varphi)\geq -n_0+k\geq -n_0$.
If $-n_0\leq C_R(\varphi)\leq m_0$, 
then $\psi$ is given by $f$. If $C_R(\varphi)\geq m_0$, then 
$\psi$ is given by the 
local rule obtained from $f$ 
by adding left-redundancy by $C_R(\varphi)-m_0$, which is still
strictly $k$ right-mergible by Lemma 3.2.) 
Therefore, since $P_R(\psi)=P_L(\varphi)+C_R(\varphi)\geq 0$ 
and $Q_R(\psi)=Q_R(\varphi)+C_R(\varphi)\geq 0$, it follows 
by Theorem 7.3(1) that for $n=k+Q_R(\psi)$, 
$\psi$ is of $(0,n)$-type and there exists
a strictly $k$ right-mergible local rule $g:L_{n+1}(X)\to L_1(X)$ 
which gives $\psi$. 
Applying Theorem 7.3(2) to $\psi$, we see that
$(\psi\sigma^s)^{[k+1-s]}$ is LR 
for $s=0,\dots,k$. Thus (b) is proved. 

Since $\varphi\sigma^{k+C_R(\varphi)}$ is LR, 
$\varphi\sigma^{k+C_R(\varphi)+s}$ is LR for all $s\geq 0$, 
by \cite[Corollary 3.18(1)]{Nasu-t} (because $\sigma$ is LR), so 
that (c) is proved. By Theorem 7.3, (d) follows. 

By Proposition 6.12(3) and Theorem 4.2, 
it follows that if  $s>C_R(\varphi)$, 
then $\varphi\sigma^s$ is 
right $\sigma$-expansive and left $\sigma$-expansive, and hence 
expansive (see Subsection 2.1). Therefore (e) is proved. 

(2) By symmetric arguments to the above.
\end{proof}

NB \cite[Proposition 6.30]{Nasu-t} erroneously asserts 
that an endomorphism $\varphi$ 
of a topological Markov shift 
$(X,\sigma)$ is right-closing 
(respectively, left-closing) if and only if 
$\varphi\sigma^l$ is  LR 
($\varphi\sigma^{-l}$ is  RL) 
for some (all sufficiently
large) $l$. Corollary 7.4  give 
a correction to it with refinements.   
If ``LR'' (respectively,``RL'') in it is replaced by 
``weakly LR'' (respectively, ``weakly RL''), 
we have another corrected version of it, 
which \cite[Proof of Proposition 6.30]{Nasu-t} with obvious corrections 
proves and is generalized to factor maps between subshifts 
in \cite[Proposition 7.10]{Nasu-te}. 

Here we remark the following. Suppose that $\varphi$ is 
an onto endomorphism of a topological Markov shift $(X,\sigma)$. 
Then $\varphi$ is LR up to higher block conjugacy if and only 
If $\varphi$ is $p$-L and $q$-R. (This is proved by 
Proposition 6.9, Theorem 4.2(1) 
and Propositions 6.5(2) and 3.5(3) and Theorem 7.3.) 
If $\varphi$ is $p$-L and $q$-R, then $\varphi$ is 
not necessarily LR. (See 
\cite{Nasu-note} together with 
\cite[Section 10, Example 2]{Nasu-t}). 

Let $G$ be a nondegenerate graph. 
Let $M_G=(m_{u,v})_{u,v\in V_G}$ be the 
adjacency matrix of $G$ (defined in Section 2). 
Let $L_G=(l_{u,a})_{u\in V_G,a\in A_G}$
denote the rectangular matrix 
with indexing set $V_G\times A_G$ such that 
$l_{u,a}=1$ if $i_G(a)=u$ and otherwise $l_{u,a}=0$. 
Let $R_G=(r_{a,v})_{a\in A_G,v\in V_G,}$
denote the rectangular matrix with indexing set 
$A_G\times V_G$ such that 
$r_{a,v}=1$ if $t_G(a)=v$ and
otherwise $r_{a,v}=0$. Then $M_G=L_GR_G$ 
and $M_{G^{[2]}}=R_GL_G$. 

Recall the definition of the product of graphs 
$G$ and $H$ with $V_G=V_H$ (given before Remark 2.15).
For $k\geq 0$, let $G^k$ denote the graph defined as follows: 
$G^0$ is the graph such that $V_{G^0}=V_G$
and $M_{G^0}$ is the identity matrix, and for $k\geq 1$, $G^k=GG^{k-1}$. 
 
\begin{lemma}
Let $G$ be a nondegenerate graph. Let $t\geq 1$. 
Suppose that $H$ is a
nondegenerate graph with $V_H=V_{G^{[t]}}$
such that \[M_HM_{G^{[t]}}=M_{G^{[t]}}M_H.\] 
Then there exists a nondegenerate graph $K$ with $V_K=V_G$ 
such that 
\[L_GM_H=M_KL_G,\q\q\q\q M_HR_G=R_GM_K\] 
and hence
\[M_KM_G=M_GM_K \]
and such that 
the onesided topological Markov shifts 
$(\tilde{X}_{H^i(G^{[t]})^{j+1}},\tilde{\sigma}_{H^i(G^{[t]})^{j+1}})$ 
and 
$(\tilde{X}_{K^iG^{j+1}},\tilde{\sigma}_{K^iG^{j+1}})$ 
are topologically conjugate for all $i,j\geq 0$. 
\end{lemma} 
\begin{proof} Since $G^{[t]}=(G^{[t-1]})^{[2]}$, if we show 
the lemma in case $t=2$, then the lemma follows by induction. 

Let $M_{HG^{[2]}}=(m_{a,b})_{a,b\in A_G}$. 
Since $M_{HG^{[2]}}=M_HM_{G^{[2]}}=M_{G^{[2]}}M_H$, we have 
\[M_{HG^{[2]}}=(m_{a,b})_{a,b\in A_G}=M_HR_GL_G=R_GL_GM_H.\] 
Since $M_{HG^{[2]}}=R_G(L_GM_H)$, it follows that for each $u\in V_G$, 
all $a$-rows of $M_{HG^{[2]}}$ with $a\in t_G^{-1}(u)$ are the same 
and equal to the $u$-row of $L_GM_H$. 
Since $M_{HG^{[2]}}=(M_HR_G)L_G$, it follows that for each $v\in V_G$, 
all $a$-columns of $M_{HG^{[2]}}$ with $a\in i_G^{-1}(v)$ are the same  
and equal to the $v$-column of $M_HR_G$. 
Hence we can define a graph $K$ by $M_K=(k_{u,v})_{u,v\in V_G}$ 
with $k_{u,v}=m_{a,b}$ if $a\in t_G^{-1}(u)$ and $b\in i_G^{-1}(v)$. 
We have 
\[R_GM_KL_G=M_{HG^{[2]}}.\] 
Therefore,
since $M_{HG^{[2]}}=R_G(L_GM_H)$, we have $M_KL_G=L_GM_H$, and
since $M_{HG^{[2]}}=(M_HR_G)L_G$, we have $R_GM_K=M_HR_G$. 
Therefore,
we have \[M_KM_G=M_KL_GR_G=L_GM_HR_G=L_GR_GM_K=M_GM_K.\] 
Since 
\begin{align*}
M_{H^i(G^{[2]})^{j+1}}&=M_H^iM_{G^{[2]}}^{j+1}
=M_H^i(R_GL_G)^jR_GL_G=M_H^iR_G(L_GR_G)^jL_G,\\
M_{K^iG^{j+1}}&=M_K^iM_G^{j+1}=M_K^iL_GR_G(L_GR_G)^j=
L_GM_H^iR_G(L_GR_G)^j
\end{align*}
and $L_G$ is a zero-one matrix with exactly one $1$ in each column 
and no rows all zero,
by well-known Williams' conjugacy result
for onesided topological Markov shifts \cite{Wil} 
we see that 
the onesided topological Markov shifts 
$(\tilde{X}_{H^i(G^{[2]})^{j+1}},\tilde{\sigma}_{H^i(G^{[2]})^{j+1}})$
and $(\tilde{X}_{K^iG^{j+1}},\tilde{\sigma}_{K^iG^{j+1}})$ 
are topologically conjugate. 
\end{proof} 

For an onto endomorphism $\varphi$ of a transitive 
topological Markov shift, 
let $r_\varphi$ and $l_\varphi$ denote 
the right and left \itl{multipliers} 
of Boyle \cite{Boy-c}. (An explanation of them 
is found in \cite[pp.52--56]{Nasu-t}. 
Note that $r_\varphi=1/R(\varphi)$ and $l_\varphi=1/L(\varphi)$ 
for $R(\varphi)$ and $L(\varphi)$ 
appearing in \cite{Boy-c} and \cite{Nasu-t}.)
 
Let $h(\varphi)$ denote the topological entropy of 
an endomorphism $\varphi$ of a subshift. 

For a subshift $(X,\sigma)$ and 
its endomorphism with memory zero, let 
$(X^\sim,\sigma^\sim)$ and $\varphi^\sim$
denote $(\tilde{X},\tilde{\sigma})$ and $\tilde{\varphi}$, 
respectively. 

\begin{theorem} 
\begin{enumerate}
\item 
If $\varphi$ is an onto endomorphism of 
a topological Markov shift 
$(X_G,\sigma_G)$ such that 
for some $t\geq 1$, $\varphi^{[t]}$ 
is an LR endomorphism of 
$(X_{G^{[t]}},\sigma_{G^{[t]}})$, 
then $\varphi$ has memory zero and 
there exists a nondegenerate graph 
$K$ with $M_GM_K=M_KM_G$ 
such that for all $i\geq 0, j\geq 1$, 
$(\tilde{X}_G, \tilde{\varphi}^i\tilde{\sigma}_G^j)$ 
is topologically conjugate to the onesided topological 
Markov shift
$(\tilde{X}_{K^iG^j},\tilde{\sigma}_{K^iG^j})$ 
and
the inverse limit system of $\varphi^i\sigma_G^j$
is topologically conjugate to 
$(X_{K^iG^j},\sigma_{K^iG^j})$. 
\item 
If $\varphi$ is an onto endomorphism of the full $n$-shift 
$(X,\sigma)$ such that for some
$t\geq 1$ $\varphi^{[t]}$ is an 
LR endomorphism of $(X^{[t]},\sigma^{[t]})$, then 
for some integer $k\geq 1$,  
\begin{equation}
h(\varphi)=h(\tilde{\varphi})=\log r_\varphi=\log k 
\end{equation}
and for all $i\geq 0, j\geq 1$, 
$(\tilde{X}, \tilde{\varphi}^i\tilde{\sigma}^j)$ is 
topologically conjugate to 
the onesided full $k^in^j$-shift and 
the inverse limit system of 
$\varphi^i\sigma^j$ is topologically conjugate 
to the full $k^in^j$-shift.
\end{enumerate}
\end{theorem}  
\begin{proof} (1) 
Let $i\geq 0$ and $j\geq 1$. 
Suppose that $\varphi^{[t]}$ 
is an LR endomorphism 
of $(X_{G^{[t]}},\sigma_{G^{[t]}})$. 
Since $\varphi^{[t]}$ 
is a $p$-L endomorphism 
of $(X_{G^{[t]}},\sigma_{G^{[t]}})$, 
$\varphi$ has memory zero 
(by Propositions 6.9 and 6.5(2)), and hence 
the induced endomorphism 
$\tilde{\varphi}$ of the induced onesided topological 
Markov shift $(\tilde{X},\tilde{\sigma})$ can be defined. 
Since $\varphi^{[t]}$ 
is an LR endomorphism 
of $(X_{G^{[t]}},\sigma_{G^{[t]}})$, 
there exists 
a onesided 1-1, LR textile system 
$T$ over $G^{[t]}$ such that $\varphi_T=\varphi^{[t]}$. 
Let $G^*$ be the nondegenerate graph over which the dual $T^*$ is 
defined. Then by \cite[Proposition 6.1]{Nasu-t} that 
$M_{G^{[t]}}M_{G^*}=M_{G^*}M_{G^{[t]}}$. 
It follows from 
\cite[Corollary 6.7(2)]{Nasu-t} that  
$(({X}_G^{[t]})^\sim, ((\varphi^{[t]})^i(\sigma_G^{[t]})^j)^\sim)$
is topologically-conjugate to
$(\tilde{X}_{(G^*)^i(G^{[t]})^j},\tilde{\sigma}_{(G^*)^i(G^{[t]})^j})$. 
By Lemma 7.5, there exists a nondegenerate graph $K$ with 
$M_KM_G=M_GM_K$ such that 
$(\tilde{X}_{K^iG^j}, \tilde{\sigma}_{K^iG^j})$ 
is topologically-conjugate to
$(\tilde{X}_{(G^*)^i(G^{[t]})^j},\tilde{\sigma}_{(G^*)^i(G^{[t]})^j})$. 
Therefore 
$(\tilde{X}_G, \tilde{\varphi}^i\tilde{\sigma}_G^j)$ 
is topologically conjugate to 
$(\tilde{X}_{K^iG^j}, \tilde{\sigma}_{K^iG^j})$.

The proof of the remainder of (1) is similar, but use 
\cite[Theorem 6.3(3)]{Nasu-t} instead of 
\cite[Corollary 6.7(2)]{Nasu-t}. 

(2)  We consider the proof of (1) for the case that 
$(X_G,\sigma_G)=(X,\sigma)$ 
(with notation being the same as in it). 

Since the defining graph $G$ of $(X,\sigma)$ is 
a one-vertex graph, so is $K$, because $M_KM_G=M_GM_K$.  
If $K$ has exactly $k$ arcs, 
then the conclusions of (2) follow 
from (1) except for (7.1). 

By the proof of (1)
$(({X}_G^{[t]})^\sim, (\varphi^{[t]}\sigma_G^{[t]})^\sim)$ 
is topologically-conjugate to 
$(\tilde{X}_{KG},\tilde{\sigma}_{KG})$. 
Therefore, 
\[h(\varphi^{[t]}\sigma^{[t]})
=h((\varphi^{[t]}\sigma_G^{[t]})^\sim)
=h(\tilde{\sigma}_{KG})=\log kn=\log k +\log n.\] 
By Boyle \cite{Boy-c} we know that if $\psi$ and $\psi'$
are onto endomorphisms of an irreducible
topological Markov shift, then $r_{\psi\psi'}=r_\psi r_{\psi'}$.  
Using this and 
\cite[Theorem 6.31(2)]{Nasu-t}
(noting that
 $\varphi^{[t]}$ and $\sigma^{[t]}$ are LR endomorphisms 
of $(X^{[t]},\sigma^{[t]})$ and hence so is 
$\varphi^{[t]}\sigma^{[t]}$ by \cite[Corollary 3.18(1)]{Nasu-t}), 
we see that  
\begin{align*} 
h(\varphi^{[t]}\sigma^{[t]})
=&\log r_{\varphi^{[t]}\sigma^{[t]}}
=\log r_{\varphi^{[t]}}+\log r_{\sigma^{[t]}}\\
=&h(\varphi^{[t]})+h(\sigma^{[t]})
=h(\varphi)+ h(\sigma)=h(\varphi) + \log n.
\end{align*}
Therefore we have $h(\varphi)=\log k.$ 

By \cite[Theorem 6.31(2)]{Nasu-t} 
$h(\varphi^{[t]})=h((\varphi^{[t]})^\sim)=r_{\varphi^{[t]}}$, 
and hence $h(\varphi)=h(\tilde{\varphi})=r_{\varphi}$, 
because by Boyle \cite{Boy-c}, $r_\varphi$ is an invariant 
of topological conjugacy between endomorphisms of 
irreducible topological Markov shifts. 

Therefore (7.1) is proved. 
\end{proof} 

\begin{corollary} Let $\varphi$ be an onto endomorphism of 
the full $n$-shift $(X,\sigma)$. 
\begin{enumerate} 
\item 
If $P_L(\varphi)\geq 0$ and $Q_R(\varphi)\geq 0$, then for 
some integer $k\geq 1$,
\[h(\varphi)=h(\tilde{\varphi})=\log r_\varphi=\log k\] and
for all $i\geq 0$ and $j\geq 1$, 
the dynamical system 
$(\tilde{X},\tilde{\varphi}^i\tilde{\sigma}^j)$ is 
topologically conjugate to 
the onesided full $k^in^j$-shift and 
the inverse limit system of 
$\varphi^i\sigma^j$ is topologically conjugate 
to the full $k^in^j$-shift. 
\item
If $P_R(\varphi)\geq 0$ and $Q_L(\varphi)\geq 0$, then 
for some integer $k\geq 1$,
\[h(\varphi)=\log l_\varphi=\log k\] 
and for all $i\geq 0$ and $j\geq 1$, 
the inverse limit system of 
$\varphi^i\sigma^{-j}$ is topologically conjugate 
to the full $k^in^j$-shift. 
\end{enumerate}
\end{corollary}
\begin{proof}
By Theorems 7.3(1),(2) and 7.6(2), (1) is proved. 
By symmetry, (2) follows from (1). 
\end{proof}

\begin{theorem} 
Let $\varphi$ be an endomorphism of a subshift $(X,\sigma)$. 
Let $s\in\Z$. 
\begin{enumerate} 
\item  
\begin{enumerate} 
\item $\varphi^{[t]}$ is weakly LR with some $t\geq 1$ if and only if 
$P_L(\varphi)\geq 0$ and $Q_R(\varphi)\geq 0$. 
\item If $(\varphi\sigma^s)^{[t]}$ is weakly LR with some $t\geq 1$, 
then $s\geq C_R(\varphi)$;  
there exists $t\geq 1$ such that for 
$C_R(\varphi)\leq s\leq C_R(\varphi)+t-1$,
$(\varphi\sigma^s)^{[t-s+C_R(\varphi)]}$ is weakly LR, 
and so is
$\varphi\sigma^s$ for $s\geq C_R(\varphi)+t$. 
\item If $s>C_R(\varphi)$ and $\varphi$ is onto, 
then $\varphi\sigma^s$ is expansive. 
\end{enumerate}
\item 
\begin{enumerate} 
\item $\varphi^{[t]}$ is weakly RL with some $t\geq 1$ 
if and only if 
$P_R(\varphi)\geq 0$ and $Q_L(\varphi)\geq 0$. 
\item If 
$(\varphi\sigma^s)^{[t]}$ is weakly RL with some $t\geq 1$, 
then $s\leq C_L(\varphi)$; there exists $t\geq 1$ such that 
for $C_L(\varphi)-t+1\leq s\leq C_L(\varphi)$,
$(\varphi\sigma^s)^{[t+s-C_L(\varphi)]}$ is weakly RL, 
and so is $\varphi\sigma^s$
for all $s\leq C_L(\varphi)-t$. 
\item If $s<C_L(\varphi)$ and $\varphi$ is 
onto, then $\varphi\sigma^s$ is expansive. 
\end{enumerate}
\end{enumerate}
\end{theorem}
\begin{proof} 
(1)(a) The ``only-if'' part follows from Theorems 5.2 and 6.11 
and Propositions 3.5 and 6.5. .  

To prove the ``if'' part 
together with more for proving (b), 
assume that 
$P_L(\varphi)\geq 0$ and $Q_R(\varphi)\geq 0$. It follows from 
Lemmas 5.1 and 6.10 that there exists $J\in\N$ such that 
$Q_R(\phi^{\langle J,\varphi\rangle})=Q_R(\varphi)$ and 
$P_L(\phi^{\langle J,\varphi\rangle})=P_L(\varphi)$. 
Since 
$P_L(\phi^{\langle J,\varphi\rangle})\geq 0$ and 
$Q_R(\phi^{\langle J,\varphi\rangle})\geq 0$,  
it follows from Proposition 7.2 that there exists $t\geq 1$ 
such that for $s=0,\dots, t-1$ 
we have a onesided 1-1, weakly LR textile-relation system 
\[T_s=(p_s:\G\to G_A^{[t-s]}, q_s:\G\to G_A^{[t-s]}) 
\q \text{with}\;\, A=L_1(X)\]
such that $\phi_{T_s}=
(\phi^{\langle J,\varphi\rangle}\sigma_{X^{\langle J\rangle}}^s)^{[t-s]}$. 
We can regard $T_s$ as a textile system 
$T_s=(p_s,q_s:\G\to G_A^{[t-s]})$, 
which is onesided 1-1 and weakly LR. 
Let $Z_s=\hf{\xi}_{T_s}^{-1}(X^{[t-s]})$ 
for $s=0,\dots,t-1$. Then, since 
$Z_s\subset\hf{Z}_{T_s}$ with 
$\hf{\xi}_{T_s}(Z_s)\supset\hf{\eta}_{T_s}(Z_s)$ and 
$\hf{\xi}_{T_s}|Z_s$ one-to-one, 
there exists a unique onesided 1-1 
half-textile-subsystem $\hf{U}_s$ of $T_s$ 
with $Z_{\hf{U}_s}=Z_s$. 
Since $(\varphi\sigma^s)^{[t-s]}=\varphi_{\hf{U}_s}$, 
we conclude that $(\varphi\sigma^s)^{[t-s]}$ 
is weakly LR for $s=0,\dots, t-1$, and hence, 
in particular the ``if'' part is proved.  Further, since 
$\varphi\sigma^{t-1}$ is weakly LR, it follows from 
Proposition 2.8(3) that 
$\varphi\sigma^s$ is weakly LR for all $s\geq t$. Therefore
we have proved the following $(*)$. 

$(*)$  If $P_L(\varphi)\geq 0$ and $Q_R(\varphi)\geq 0$, 
then there exists $t\geq 1$ such that 
$(\varphi\sigma^s)^{[t-s]}$ is weakly LR 
for $s=0,\dots,t-1$ and $\varphi\sigma^s$ is weakly LR for all 
$s\geq t$. 

(b) For the first assertion,   
substitute $\varphi\sigma^s$ for $\varphi$ in (a); for the second,  
apply $(*)$ to $\varphi\sigma^{C_R(\varphi)}$. 
Then by using Propositions 3.4 and 6.4, (b) is proved.

(c) By \cite[Proposition 7.6(3)]{Nasu-te}. 

(2) By symmetric arguments to the above, (2) is proved. 
\end{proof}

\begin{proposition}
Let $\varphi$ be an endomorphism of 
an infinite subshift $(X,\sigma)$. 
\begin{enumerate}
\item 
$P_L(\varphi)+Q_L(\varphi)$ and $P_R(\varphi)+Q_R(\varphi)$ 
are nonpositive and hence 
\[Q_L(\varphi)\leq -P_L(\varphi),\q
P_R(\varphi)\leq -Q_R(\varphi)\q\text{and}\q
C_L(\varphi)\leq C_R(\varphi).\] 
\item 
If $\varphi$ is $(m,n)$-type
then $-n\leq P_R(\varphi)\leq -Q_R(\varphi)$, 
$Q_L(\varphi)\leq -P_L(\varphi)\leq m$, and hence 
\begin{enumerate}
\item[(a)] if $Q_R(\varphi)+Q_L(\varphi)\geq 0$ then 
\[-n\leq P_R(\varphi)\leq -Q_R(\varphi)\leq 
Q_L(\varphi)\leq -P_L(\varphi)\leq m,\]
\item[(b)] if $Q_R(\varphi)+Q_L(\varphi)\leq 0$ then 
\[-n\leq D_L(\varphi)\leq D_R(\varphi)\leq m.\]
\end{enumerate} 
\end{enumerate}
\end{proposition}
\begin{proof} (1) Assume that 
$P_L(\varphi)+Q_L(\varphi)>0$. 
Then $Q_L(\varphi)\in\Z$, because
$X$ is infinite and hence $Q_L(\varphi)<\infty$. 
By Proposition 5.2(2), $\varphi\sigma^{Q_L(\varphi)-1}$ 
is (weakly $q$-L and) positively right $\sigma$-expansive. 
Since $Q_L(\varphi)-1\geq -P_L(\varphi)$ (by assumption), 
it follows from Theorem 6.11 that 
$\varphi\sigma^{Q_L(\varphi)-1}$ is weakly $p$-L. 
Therefore we are led to a contradiction, because
there exists no endomorphism of $(X,\sigma)$
which is weakly $p$-L 
and positively right $\sigma$-expansive. 
We shall prove this in the following.

Assume that there exists an endomorphism $\psi$ 
which is weakly $p$-L
and positively right $\sigma$-expansive. 
Let $x,y\in X$ and put $\psi^i(x)=(a_{i,j})_{i\geq 0,j\in\Z}$
with $a_{i,j},b_{i,j}\in L_1(X)$. Since $\psi$ is 
positively right $\sigma$-expansive, there exists $\delta\geq 0$
such that if $d_X(\sigma^j\psi^j(x),\sigma^j\psi^j(y))\leq \delta$ 
for all $i,j\geq 0$ then $x=y$. Hence there exists $m>0$ such that 
if $a_{i,j}=b_{i,j}$ for all $i\geq 0, j\geq -m$, 
then $x=y$. Since $\psi$ is 
weakly $p$-L, by Theorem 6.11 $\psi$ has memory zero. 
Therefore, if $a_{0,j}=b_{0,j}$ for all $j\geq -m$, then 
$a_{i,j}=b_{i,j}$ for all $i\geq 0, j\geq -m$ and hence $x=y$. 
This implies that for any $x=(a_{0,j})_{j\in\Z}\in X$, 
$(a_{0,j})_{j\geq -m}$ uniquely determines $x$, which 
cannot be the case because $(X,\sigma)$ is an infinite subshift. 

We have proved that
$P_L(\varphi)+Q_L(\varphi)\leq 0$. By symmetry we have 
$P_R(\varphi)+Q_R(\varphi)\leq 0$. 

(2) By (6.3) and (1).
\end{proof}

A direct proof of the first inequality of Proposition 7.9(1) 
for the case that $\varphi$ is an automorphism of a 
topological Markov shift or $\varphi$ is an onto endomorphism of a 
transitive topological Markov shift is given 
by replacing Theorem 6.11 by Proposition 6.9 and 
Theorem 5.2(2) by Theorem 4.2(2) and 
by deleting ``weakly'' in the general proof above;
the second one follows from the first by symmetry. 
\begin{proposition}
Let $\varphi$ be an endomorphism of 
a subshift $(X,\sigma)$. 
\begin{enumerate} 
\item 
The following conditions are equivalent:
  \begin{enumerate}
  \item
  $\varphi$ is weakly LL (respectively, weakly RR);
  \item
  $P_L(\varphi)\geq 0$ and $Q_L(\varphi)\geq 0$ 
(respectively, $P_R(\varphi)\geq 0$ and $Q_R(\varphi)\geq0$);
  \item 
  For some $n\geq 0$, 
  $\varphi$ is of $(0,n)$-type (respectively, $(n,0)$-type) 
and given by a zero left-mergible 
(respectively, zero right-mergible) local rule 
  $f:L_{n+1}(X)\to L_1(X)$. 
  \end{enumerate}
\noindent If in addition $X$ is infinite, then (b) is equivalent to
the condition that $P_L(\varphi)= Q_L(\varphi)= 0$ (respectively, 
$P_R(\varphi)= Q_R(\varphi)= 0$). 
\item 
The following conditions are equivalent if $X$ is infinite:
  \begin{enumerate}
  \item[(a)]
There exists $s\in\Z$ with $\varphi\sigma^s$ weakly LL 
(respectively, weakly RR);
  \item[(b)]
  $P_L(\varphi)+Q_L(\varphi)=0$ 
(respectively, $P_R(\varphi)+Q_R(\varphi)=0$);
  \item[(c)]  
  $\varphi\sigma^{-P_L(\varphi)}$ is weakly LL 
  (respectively, $\varphi\sigma^{-Q_R(\varphi)}$ is weakly RR).  
\end{enumerate}
\item 
When $(X,\sigma)$ is a topological Markov shift, if 
$\varphi$ is bijective or if $\varphi$ is onto with 
$\sigma$ transitive,  
then all ``weakly'' can be deleted in (1) and (2). 
\end{enumerate} 
\end{proposition}
\begin{proof}
(1) Assume (a). Then $P_L(\varphi)\geq 0$ by
Theorem 6.11, and 
$Q_L(\varphi)\geq 0$ by Theorem 5.2(2). 
Hence (b) follows by Proposition 7.9(1). 

Assume (b). Since $P_L(\varphi)\geq 0$, it follows from  
Theorem 6.11 that
$\varphi$ is of $(0,n)$-type and given by a local rule 
$f: L_{n+1}(X)\to L_1(X)$ with some $n\geq 0$. If 
$f$ is strictly $l$ left mergible, 
then $0-l=Q_L(\varphi)\geq 0$. Hence $l=0$ and (c) is proved. 

Assume (c). Let $T=(p,q:G[X,n+1]\to G_A)$ with $A=L_1(X)$ be the 
textile system defined as follows: for 
each arc $\alpha=a_0\dots a_n$ in $G[X,n+1]$ with $a_j\in A$, 
$p(\alpha)=a_0$ and $q(\alpha)=f(a_0\dots a_n)$. (Recall that 
$G[X,s]$ was defined at the beginning of Section 5.) Then 
$\hf{\xi}_T$ is one-to-one. Clearly $p$ is weakly left-resolving. 
Since $f$ is zero left-mergible, $q$ is weakly left-resolving. 
Let $Z=\hf{\xi}_T^{-1}(X)$. Then, since $Z\subset\hf{Z}_T$ with 
$\hf{\xi}_T(Z)\supset\hf{\eta}_T(Z)$ and 
$\hf{\xi}_T|Z$ is one-to-one, 
there exists a unique onesided 1-1 
half-textile-subsystem $\hf{U}$ of $T$ with $Z_{\hf{U}}=Z$. 
Since $\varphi=\varphi_{\hf{U}}$, 
we see that $\varphi$ is weakly LL and (a) follows. 

The remainder follows from Proposition 7.9(1).  

(2) By (1) and Propositions 6.4 and 3.4. 

(3) The proof of (1) with ``weakly'' deleted 
is the same as that of the original (1) except for the following. 
To prove that (a)(without ``weakly'') implies (b) 
and that (b) implies (c), 
use Proposition 6.9 and Theorem 4.2(2)
instead of Theorem 6.11 and Theorem 5.2(2), respectively.  
To prove that (c) implies (a) with ``weakly'' deleted, 
a onesided 1-1, LL textile system 
$T$ with $\varphi_T=\varphi$ 
is directly given as 
follows. If $(X,\sigma)=(X_G,\sigma_G)$, then 
define $T=(p,q: G^{[n+1]}\to G)$, where for 
$\alpha=a_0\dots a_n\in A_{G^{[n+1]}}$ with $a_j\in A_G$, 
$p(\alpha)=a_0$ and $q(\alpha)=f(a_0\dots a_n)$. 
Clearly, $T$ is weakly LL. Hence, 
by Lemma 2.1, $T$ is LL.

The proof of the claim (2) with ``weakly'' deleted 
is given by (1) with ``weakly'' deleted  
and Propositions 6.4 and 3.4.
\end{proof} 

For more information on onto LL endomorphisms and LL automorphisms
of topological Markov shifts, see \cite[Section 5]{Nasu-d}. 

\begin{proposition} 
\begin{enumerate} 
\item
If $\varphi$ is an automorphism of a subshift, 
then 
\begin{gather*}
Q_L(\varphi)=P_L(\varphi^{-1}),\q\q\q\q
Q_R(\varphi)=P_R(\varphi^{-1}),\\ 
C_L(\varphi^{-1})=-C_R(\varphi),\q\q\q\q 
D_L(\varphi^{-1})=-D_R(\varphi).
\end{gather*} 
\item
If $\varphi$ is an automorphism of an infinite subshift, 
then 
\[Q_R(\varphi)+Q_L(\varphi)\leq 0\q\text{and}\q
D_L(\varphi)\leq D_R(\varphi).\]
\item
If $\varphi$ is a symbolic automorphism of an 
infinite subshift, then
$P_L(\varphi)=Q_R(\varphi)=P_R(\varphi)=Q_L(\varphi)=0$. 
\item 
If $\varphi$ is a symbolic automorphism of 
an infinite subshift $(X,\sigma)$, then 
for any endomorphism $\psi$ of 
$(X,\sigma)$,  $\varphi\psi$ and $\psi\varphi$ have
the same onesided resolving degrees as $\psi$.
\end{enumerate}
\end{proposition}
\begin{proof} 
(1) Let $s\in\Z$. 
By Proposition 6.12(2), $s\geq -P_L(\varphi^{-1})$ 
if and only $\varphi^{-1}\sigma^s$ is weakly $p$-L. 
The latter holds if and only if
$\varphi\sigma^{-s}$ is weakly $q$-L. 
The latter holds if and only if 
$-s\leq Q_L(\varphi)$, by Theorem 5.2(2). 
Therefore, $s\geq -P_L(\varphi^{-1})$ 
if and only if $s\geq -Q_L(\varphi)$, 
and hence $P_L(\varphi^{-1})=Q_L(\varphi)$. 
The proof that $Q_R(\varphi)=P_R(\varphi^{-1})$ is similar. 
The remainder follow from the first two equations. 

(2) By Proposition 6.13(1), 
$P_R(\varphi^{-1})+P_L(\varphi^{-1})\leq 0$. Hence by (1), 
$Q_R(\varphi)+Q_L(\varphi)\leq 0$. Therefore 
by Proposition 7.9(2)(b), $D_L(\varphi)\leq D_R(\varphi)$.

(3) By Proposition 7.1, 
$P_L(\varphi)=P_R(\varphi)=0$. 
Since $\varphi^{-1}$ is a symbolic automorphism, 
we have $P_L(\varphi^{-1})=P_R(\varphi^{-1})=0$.  
Hence by (1), $Q_L(\varphi)=Q_R(\varphi)=0$. 

(4) Since $\varphi^{-1}$ is also a symbolic automorphism, 
the results follow from Propositions 3.8 and 6.7 and (3).
\end{proof}

\section{Endomorphisms having resolving powers} 

\subsection{Endomorphisms having resolving powers}
Let $\varphi$ be an endomorphism of a subshift $(X,\sigma)$. Let 
$s\in\N$. For $x\in X$, let 
\[\rho^*_{\varphi,s}(x)=((a_{i,j})_{0\leq i\leq {s}-1})_{j\in\Z},\]
where $(a_{i,j})_{j\in\Z}=\varphi^i(x)$ for $i=0,\dots, {s}-1$ 
with $a_{i,j}\in L_1(X)$. We often write 
\[\rho^*_{\varphi,s}(x)=(a_{i,j})_{0\leq i\leq {s}-1,j\in\Z}.\]  
Let 
\[X^{[^*_\varphi s]}=\{\rho^*_{\varphi,s}(x) \bigm| x\in X\}.\]
Then we have a subshift 
$(X^{[^*_\varphi s]},\sigma^{[^*_\varphi s]})$, which is 
called the \itl{$\varphi$-dual higher-block system of order $s$ of   
$(X,\sigma)$}, such that 
\[L_1(X^{[^*_\varphi s]})\subset\{(a_i)_{0\leq i\leq {s}-1}
\bigm| a_i\in L_1(X)\}\]
and have a topological conjugacy     
$\rho^*_{\varphi,s}:(X,\sigma)\to 
(X^{[^*_\varphi s]},\sigma^{[^*_\varphi s]})$, 
which is called the 
\itl{$\varphi$-dual higher-block conjugacy of order $s$} 
or a \itl{dual higher-block conjugacy}. 
We also have an endomorphism $\varphi^{[^*s]}$ of 
$(X^{[^*_\varphi s]},\sigma^{[^*_\varphi s]})$ by 
\[\varphi^{[^*s]}
=\rho^*_{\varphi,s}\varphi(\rho^*_{\varphi,s})^{-1},\] 
which is called the \itl{dual higher-block presentation 
of order $s$} of $\varphi$. 
Clearly, the endomorphism $\varphi$ of $(X,\sigma)$ is topologically
conjugate to the endomorphism $\varphi^{[^*s]}$ of 
$(X^{[^*_\varphi s]},\sigma^{[^*_\varphi s]})$ 
through $\rho^*_{\varphi,s}:(X,\sigma,\varphi)\to 
(X^{[^*_\varphi s]},\sigma^{[^*_\varphi s]},\varphi^{[^*s]})$. 

Let $N\geq 1$.  For convenience, 
each element $\alpha_{1}\dots\alpha_{N}$ 
of $L_{N}(X^{[^*_\varphi s]})$ such that
$\alpha_j=(a_{i,j})_{0\leq i\leq s-1} \in L_1(X^{[^*_\varphi s]})$ 
with $a_{i,j}\in L_1(X)$
will often be written 
\[(a_{i,j})_{0\leq i\leq s-1, 1\leq j\leq N}.\]
Clearly, if 
$(a_{i,j})_{0\leq i\leq s-1, 1\leq j\leq N}\in L_{N}(X^{[^*_\varphi s]})$, 
then
$a_{i,1}\dots a_{i,N}\in L_{N}(X)$ for $i=0,\dots,s-1$. 

In 
\cite[Proof of Proposition 8.2]{Nasu-t}, it was shown that 
for an automorphism $\varphi$ of a subshift $(X,\sigma)$, if
$\varphi^s$ is the identity map with $s\geq 1$, 
then $\varphi^{[^*s]}$ is a symbolic 
automorphism of 
$(X^{[^*_\varphi s]},\sigma^{[^*_\varphi s]})$, and hence $\varphi$ is 
an essentially symbolic automorphism of $(X,\sigma)$. 

Mike Boyle \cite{Boy-p} proved  
the theorem  below 
(for onto endomorphisms of subshifts) 
and suggested the possibility of 
proving the main results of \cite[Section 7]{Nasu-t} on 
resolving endomorphisms of topological Markov shifts
without using the long theory of 
``resolvable textile systems'' developed there. 

\begin{theorem}[Boyle \cite{Boy-p}]
Let $\varphi$ be an endomorphism of a subshift $(X,\sigma)$. 
Let $s\geq 1$. 
\begin{enumerate} 
\item
If $\varphi^s$ is weakly $p$-L 
(respectively, weakly $p$-R), then so is $\varphi^{[^*s]}$, and hence
$\varphi$ is essentially weakly $p$-L 
(respectively, essentially weakly $p$-R). 
\item
If $\varphi^s$ is a weakly LR automorphism of $(X,\sigma)$, then 
$\varphi$ is essentially weakly LR. 
\item 
When $(X,\sigma)$ is a topological Markov shift, 
if $\varphi^s$ is $p$-L 
(respectively, $p$-R), then so is $(\varphi^{[^*s]})^{[t]}$ for 
some $t\geq 1$.
\end{enumerate}
\end{theorem} 
\begin{proof} 
(1) Since $\varphi^s$ is weakly $p$-L, it follows from 
Theorem 6.11 
(the equivalence of (2) and (3)) 
$\varphi^s$ is of $(0,n)$-type and given by a 
local-rule $f:L_{n+1}(X)\to L_1(X)$ with 
some $n\geq 0$.   We define a local-rule 
$F:L_{n+1}(X^{[^*_\varphi s]})\to L_1(X^{[^*_\varphi s]})$ 
as follows: $F=f$ if $s=1$;  
if $s\geq 2$, then 
for each 
\[W=(a_{i,j})_{0\leq i\leq s-1, 1\leq j\leq n+1}\in 
L_{n+1}(X^{[^*_\varphi s]}), 
\q\q a_{i,j}\in L_1(X),\]   
$F(W)=(a_i)_{1\leq i\leq s}$,
where 
\[a_i=a_{i,1}\;\text{for}\; i=1,\dots,s-1,  
\q\text{and}\q a_s=f(a_{0,1}\dots a_{0,n+1}).\]
Then it is clear that 
$\varphi^{[^*s]}$ is a block-map of $(0,n)$-type given by 
the local-rule $F$. 
Therefore, again by Theorem 6.11 
(the equivalence of (2) and (3)), 
$\varphi^{[^*s]}$ is weakly $p$-L, and hence 
$\varphi$ is essentially weakly $p$-L. 

(2) Since $\varphi^s$ is weakly $p$-L 
and $\varphi^{-s}$ is weakly $p$-R 
(because $\varphi^s$ is weakly $q$-R), 
it follows from (1) that $\varphi^{[^*s]}$ is weakly $p$-L 
and $(\varphi^{-1})^{[^*s]}=(\varphi^{[^*s]})^{-1}$ is 
weakly $p$-R. Hence it follows from 
\cite[Proposition 7.5]{Nasu-te} that for some $t\geq 1$
$(\varphi^{[^*s]})^{[t]}$ is weakly LR, and hence 
$\varphi$ is essentially weakly LR. 

(3) Since $(X^{[^*_\varphi s]},\sigma^{[^*_\varphi s]})$ is an SFT, 
$(\varphi^{[^*s]})^{[t]}$ is an endomorphism of 
a topological Markov shift for some $t\geq 1$. 
Using this and \cite[Proposition 5.1(1)]{Nasu-te} 
(instead of \cite[Proposition 7.5(1)]{Nasu-te}), (3) is 
proved in a similar way to the proof of (1).  
\end{proof} 

Moreover, Theorem 8.1(3) gives a very short proof for 
\cite[Theorem 7.22(1)]{Nasu-t} that if $\varphi^s$ is a $p$-L 
endomorphism of a topological Markov shift, then 
$\varphi$ is an essentially $p$-L endomorphism of the shift, and 
also a similar proof to that of 
(2) using it together with \cite[Proposition 5.1]{Nasu-te} 
(instead of 
\cite[Proposition 7.5]{Nasu-te}) proves 
that if $\varphi^s$ is an LR automorphism of a topological Markov shift, 
then $\varphi$ is an essentially LR automorphism of the shift, 
which is the automorphism case of \cite[Theorem 7.22(2)]{Nasu-t}. 

\begin{proposition} 
Let $\varphi$ be an endomorphism of a subshift $(X,\sigma)$. 
Let $s\geq 1$. 
\begin{enumerate}
\item If $P_L(\varphi^s)\geq 0$, then $P_L(\varphi^{[^*s]})\geq 0$. 
\item If $P_R(\varphi^s)\geq 0$, then $P_R(\varphi^{[^*s]})\geq 0$.  
\end{enumerate}
\end{proposition}
\begin{proof} 
By Theorems 8.1(1) and 6.11. 
\end{proof}

\begin{remark} Suppose that
$\varphi$ is an endomorphism of a subshift $(X,\sigma)$ 
and $s\geq 1$. Then 
$P_L(\varphi^{[^*s]})\leq P_L(\varphi^{[^*s+1]})$ 
and $P_R(\varphi^{[^*s]})\leq P_R(\varphi^{[^*s+1]})$.
\end{remark}
\begin{proof} 
Since $\varphi^{[^*s+1]}=(\varphi^{[^*s]})^{[^*2]}$, 
it suffices to prove the remark for $s=1$. 

Suppose that $\varphi$ is 
of $(m,n)$-type 
and given by 
a local rule $f:L_{N+1}(X)\to L_1(X)$ with $N=m+n$. 
Define $g:L_{N+1}(X^{[^*_\varphi 2]})\to L_1(X^{[^*_\varphi 2]})$ 
to be the local-rule such that  
for each  
$W=(a_{i,j})_{0\leq i\leq 1, 1\leq j\leq N+1}\in L_{N+1}(X^{[^*_\varphi 2]})$
with $a_{i,j}\in L_1(X)$,   
$g(W)=(a_i)_{1\leq i\leq 2}$, 
where $a_1=f(a_{0,1}\dots a_{0,N+1})(=a_{1,m+1})$ 
and $a_2=f(a_{1,1}\dots a_{1,N+1})$.
Then $\varphi^{[^*2]}$ is of $(m,n)$-type and 
given by $g$. 
We see that if $f$ is $I$ left-redundant 
then so is $g$ and that
if $f$ is $J$ right-redundant then so is $g$. Therefore, 
since both $\varphi$ and $\varphi^{[^*2]}$ are of $(m,n)$-type 
and
given by $f$ and $g$, respectively,  
we have  
$P_L(\varphi)\leq P_L(\varphi^{[^*2]})$ and 
$P_R(\varphi)\leq P_R(\varphi^{[^*2]})$. 
\end{proof} 

\begin{proposition} 
Let $\varphi$ be an endomorphism 
of a subshift $(X,\sigma)$. Let 
$s\geq 1$. 
\begin{enumerate}
\item If $Q_R(\varphi^s)\geq 0$, then $Q_R(\varphi^{[^*s]})\geq 0$, 
and if $Q_R(\varphi^s)\leq 0$, 
then $Q_R(\varphi^{[^*s]})\geq Q_R(\varphi^s)$.
\item If $Q_L(\varphi^s)\geq 0$, then $Q_L(\varphi^{[^*s]})\geq 0$, 
and if $Q_L(\varphi^s)\leq 0$, 
then $Q_L(\varphi^{[^*s]})\geq Q_L(\varphi^s)$. 
\end{enumerate}
\end{proposition}
\begin{proof}  
(1) Suppose that $\varphi^s$ is of $(m_s,n_s)$-type and given by a 
local-rule $f_s:L_{N_s+1}(X)\to L_1(X)$ with 
$N_s=m_s+n_s$. Suppose that 
$f_s$ is strictly $k_s$ right-mergible. 

We define a local-rule 
$F_s:L_{N_s+1}(X^{[^*_\varphi s]})\to L_1(X^{[^*_\varphi s]})$ 
as follows: $F_1=f_1$;  
if $s\geq 2$, then 
for each 
\[W=(a_{i,j})_{0\leq i\leq s-1, 1\leq j\leq N_s+1}\in 
L_{N_s+1}(X^{[^*_\varphi s]}), 
\q\q a_{i,j}\in L_1(X),\]   
$F_s(W)=(a_i)_{1\leq i\leq s}$,
where 
\[a_i=a_{i,m_s+1}\;\text{for}\; i=1,\dots,s-1,  
\q\text{and}\q a_s=f_s(a_{0,1}\dots a_{0,N_s+1}).\]
Then it is clear that 
$\varphi^{[^*s]}$ is a block-map of $(m_s,n_s)$-type given by 
$F_s$. 

We may assume that $s\geq 2$. 
Suppose that $Q_R(\varphi^s)\geq 0$.  Then $k_s\leq n_s$, and hence 
$f_s$ is $n_s$ right-mergible. 
To prove that $F_s$ is $n_s$ right-mergible, 
let
\[(a_{i,j})_{0\leq i\leq s-1,j\in\Z}\q\text{and}\q
(b_{i,j})_{0\leq i\leq s-1,j\in\Z}, 
\q\q a_{i,j},b_{i,j}\in L_1(X)\]
be two points in $X^{[^*_\varphi s]}$
such that 
$(a_{i,j})_{0\leq i\leq s-1}=(b_{i,j})_{0\leq i\leq s-1}$
for all $j\leq 0$ and 
\[
F_s((a_{i,J})_{0\leq i\leq s-1,j-N_s\leq J\leq j})
=F_s((b_{i,J})_{0\leq i\leq s-1,j-N_s\leq J\leq j})\q\text{for}\;
j=1,\dots, n_s+1.
\] 
Then, by the definition of $F_s$, we have
\begin{gather}
(a_{i,j-n_s})_{1\leq i\leq s-1}=(b_{i,j-n_s})_{1\leq i\leq s-1} 
\q\text{for}\;j=1,\dots, n_s+1,\\
f_s((a_{0,J})_{j-N_s\leq J\leq j})
=f_s((b_{0,J})_{j-N_s\leq J\leq j})\q\text{for}\;
j=1,\dots, n_s+1.
\end{gather}
Since $(a_{0,j})_{j\in\Z}, (b_{0,j})_{j\in\Z}\in X$ with 
$(a_{0,j})_{j\leq 0}=(b_{0,j})_{j\leq 0}$ and
$f_s$ is $n_s$ right-mergible, it follows from (8.2) that 
$a_{0,1}=b_{0,1}$. Combining this with (8.1) for the case that $j=n_s+1$, 
we have $(a_{i,1})_{0\leq i\leq s-1}=(b_{i,1})_{0\leq i\leq s-1}$.
Therefore $F_s$ is $n_s$ right-mergible, and hence 
$Q_R(\varphi^{[^*s]})\geq 0$. 

Assume that  $Q_R(\varphi^s) \leq 0$. Then $k_s\geq n_s$. 
To prove that $F_s$ is $k_s$ right-mergible, 
let
\[(a_{i,j})_{0\leq i\leq s-1,j\in\Z}\q\text{and}\q
(b_{i,j})_{0\leq i\leq s-1,j\in\Z}, 
\q\q a_{i,j},b_{i,j}\in L_1(X)\]
be two points in $X^{[^*_\varphi s]}$
such that 
$(a_{i,j})_{0\leq i\leq s-1}=(b_{i,j})_{0\leq i\leq s-1}$
for all $j\leq 0$ and 
\[
F_s((a_{i,J})_{0\leq i\leq s-1,j-N_s\leq J\leq j})
=F_s((b_{i,J})_{0\leq i\leq s-1,j-N_s\leq J\leq j})\q\text{for}\;
j=1,\dots, k_s+1.
\] 
Then, by the definition of $F_s$, we have
\begin{gather}
(a_{i,j-n_s})_{1\leq i\leq s-1}=(b_{i,j-n_s})_{1\leq i\leq s-1} 
\q\text{for}\;j=1,\dots, k_s+1,\\
f_s((a_{i,J})_{0,j-N_s\leq J\leq j})
=f_s((b_{i,J})_{0,j-N_s\leq J\leq j})\q\text{for}\;
j=1,\dots, k_s+1.
\end{gather}
Since $(a_{0,j})_{j\in\Z}, (b_{0,j})_{j\in\Z}\in X$ with 
$(a_{0,j})_{j\leq 0}=(b_{0,j})_{j\leq 0}$ and
$f_s$ is $k_s$ right-mergible, it follows from (8.4) that 
$a_{0,1}=b_{0,1}$. Combining this with (8.3) for the case that $j=n_s+1$, 
we have $(a_{i,1})_{0\leq i\leq s-1}=(b_{i,1})_{1\leq i\leq s-1}$.
Therefore $F_s$ is $k_s$ right-mergible, and hence 
$Q_R(\varphi^{[^*s]})\geq n_s-k_s=Q_R(\varphi^s)$. 

(2) By symmetry, (2) follows from (1). 
\end{proof} 

\begin{remark} Suppose that 
$\varphi$ is an endomorphism of a subshift $(X,\sigma)$ 
and $s\geq 1$. Then 
$Q_R(\varphi^{[^*s]})\leq Q_R(\varphi^{[^*s+1]})$
and $Q_L(\varphi^{[^*s]})\leq Q_L(\varphi^{[^*s+1]})$. 
\end{remark}
\begin{proof} The proof is similar to that of Remark 8.3. 
For the same $f$ and $g$ there, 
we see that if $f$ is $k$ right-mergible 
then so is $g$ and if $f$ is $l$ left-mergible then so is $g$, and
hence we have  
$Q_R(\varphi)\leq Q_R(\varphi^{[^*2]})$ and 
$Q_L(\varphi)\leq Q_L(\varphi^{[^*2]})$. 
\end{proof}

\begin{theorem}
Let $\varphi$ be 
an endomorphism of a subshift $(X,\sigma)$. 
Let $s\geq 1$. 
\begin{enumerate} 
\item
If $\varphi^s$ is weakly $q$-R
(respectively, weakly $q$-L), then so is $\varphi^{[^*s]}$, 
and hence $\varphi$ is essentially weakly $q$-R 
(respectively, essentially weakly $q$-L). 
\item When $(X,\sigma)$ is a topological 
Markov shift, if 
$\varphi$ is bijective or if $\varphi$ is onto with $\sigma$ transitive, 
then it holds that
if $\varphi^s$ is $q$-R
(respectively, $q$-L), then so is $(\varphi^{[^*s]})^{[t]}$ for 
some  $t\geq 1$. 
\end{enumerate}
\end{theorem} 
\begin{proof} 
(1) By Proposition 8.4 and Theorem 5.2. 

(2) Since $\varphi^s$ is $q$-R, $Q_R(\varphi^s)\geq 0$, by Theorem 4.2. 
Hence, it follows from 
Propositions 8.4 and 3.5(3) that 
$Q_R((\varphi^{[^*s]})^{[t]})\geq 0$ for all $t\geq 1$. 
Since $(X^{[^*_\varphi s]},\sigma^{[^*_\varphi s]})$ is an SFT, 
there exists $t\geq 1$ such that 
$((X^{[^*_\varphi s]})^{[t]},(\sigma^{[^*_\varphi s]})^{[t]})$  
is a topological Markov shift 
and hence, by 
Theorem 4.2, $(\varphi^{[^*s]})^{[t]}$ 
is a $q$-R endomorphism of the shift. 
\end{proof} 

One can prove using the method of \cite[Section 7]{Nasu-t}
that if $\varphi^s$ is a $q$-R (respectively, $q$-L) 
endomorphism of a topological Markov shift, then 
$\varphi$ is an essentially $q$-R 
(respectively, essentially $q$-L) endomorphism of the shift.
Theorem 8.6(2) provides a very short proof for an important 
special case of the result. 

\begin{theorem} 
Let $\varphi$ be an endomorphism 
of a subshift $(X,\sigma)$. 
Let $s\geq 1$.
\begin{enumerate} 
\item 
    \begin{enumerate}
    \item[(a)] If  
$\varphi^s$ is weakly $p$-L and 
weakly $p$-R, then $(\varphi^{[^*s]})^{[t]}$ is 
weakly $p$-biresolving for some $t\geq 1$, and hence 
$\varphi$ is essentially weakly $p$-biresolving. 
    \item[(b)]If $\varphi^s$ is weakly $p$-L and weakly $q$-R 
(respectively, weakly $p$-R and weakly $q$-L), then 
$(\varphi^{[^*s]})^{[t]}$ is weakly LR (respectively, weakly RL) 
with some $t\geq 1$, and hence 
$\varphi$ is essentially weakly LR 
(respectively, essentially weakly RL).
    \item[(c)] If $\varphi^s$ is weakly $q$-L and weakly $q$-R, then 
$(\varphi^{[^*s]})^{[t]}$ is 
weakly $q$-biresolving with some $t\geq 1$, and hence 
$\varphi$ is essentially weakly $q$-biresolving. 
    \item[(d)] If $\varphi^s$ is weakly $p$-L and weakly $q$-L 
(respectively, weakly $p$-R and weakly $q$-R), then 
$\varphi^{[^*s]}$ is weakly LL (respectively, weakly RR), and hence 
$\varphi$ is essentially weakly LL 
(respectively, essentially weakly RR). 
     \end{enumerate} 
\item Suppose that $(X,\sigma)$ is a topological Markov shift 
and $\varphi$ is onto.   
     \begin{enumerate}
     \item[(a)]
If $\varphi^s$ is $p$-L and $p$-R, then 
$(\varphi^{[^*s]})^{[t]}$ is $p$-biresolving  
with some $t\geq 1$. 
     \end{enumerate}
If in addition, $\varphi$ is one-to-one or $\sigma$ is transitive, 
then the following hold.
     \begin{enumerate}
     \item[(b)] If $\varphi^s$ is $p$-L and $q$-R 
(respectively, $p$-R and $q$-L), then 
$(\varphi^{[^*s]})^{[t]}$ is LR (respectively, RL) 
with some $t\geq 1$. 
     \item[(c)]  
If $\varphi^s$ is $q$-L and $q$-R, 
then 
$(\varphi^{[^*s]})^{[t]}$ is $q$-biresolving  
with some $t\geq 1$.
     \item[(d)] 
If $\varphi^s$ is $p$-L and $q$-L 
(respectively, $p$-R and $q$-R), then 
$(\varphi^{[^*s]})^{[t]}$ is LL (respectively, RR) 
for some $t\geq 1$, and hence 
$\varphi$ is essentially LL (respectively, essentially RR). 
     \end{enumerate}
\end{enumerate}
\end{theorem}
\begin{proof}  
(1)(a)  By Theorem 8.1(1) $\varphi^{[^*s]}$ is 
weakly $p$-L and weakly $p$-R, and hence 
 by Proposition 7.1(1) $(\varphi^{[^*s]})^{[t]}$ is 
weakly $p$-biresolving with some $t\geq 1$. 

(1)(b)  By Theorems 6.11 and 5.2(1), $P_L(\varphi^s)\geq 0$ and 
$Q_R(\varphi^s)\geq 0$. Hence by Propositions 8.2(1) and 8.4(1), 
$P_L(\varphi^{[^*s]})\geq 0$ and $Q_R(\varphi^{[^*s]})\geq 0$. 
Therefore, by Theorem 7.8(1) $(\varphi^{[^*s]})^{[t]}$ is 
weakly LR with some $t\geq 1$. 

(1)(c) By Theorem 5.2, 
$Q_R(\varphi^s)\geq 0$ and $Q_L(\varphi^s)\geq 0$ and hence, 
by Proposition 8.4,     
$Q_R(\varphi^{[^*s]})\geq 0$ and $Q_L(\varphi^{[^*s]})\geq 0$. 
Therefore, by Theorem 5.3 
$(\varphi^{[^*s]})^{[t]}$ is 
weakly $q$-biresolving with some $t\geq 1$. 

(1)(d) 
By Theorems 6.11 and 5.2(2), $P_L(\varphi^s)\geq 0$ and 
$Q_L(\varphi^s)\geq 0$. Hence by Propositions 8.2(1) and 8.4(2), 
$P_L(\varphi^{[^*s]})\geq 0$ and $Q_L(\varphi^{[^*s]})\geq 0$. 
Therefore, by Propositions 7.10(1) $\varphi^{[^*s]}$ is weakly LL. 

(2)(a)  
By Proposition 6.9, 
$P_L(\varphi^s)\geq 0$ and $P_R(\varphi^s)\geq 0$. 
Hence $P_L(\varphi^{[^*s]})\geq 0$ and $P_R(\varphi^{[^*s]})\geq 0$ 
by Proposition 8.2. 
Since $(X^{[^*_\varphi s]},\sigma^{[^*_\varphi s]})$
is an SFT, for sufficiently large $t'$
$((X^{[^*_\varphi s]})^{[t']},(\sigma^{[^*_\varphi s]})^{[t']})$
is a topological Markov shift . Hence,  
using Proposition 6.5(2) 
we can apply Proposition 7.1 to $(\varphi^{[^*s]})^{[t']}$ 
to prove (2)(a). 

(2)(b) By Proposition 6.9 and Theorem 4.2(1), 
$P_L(\varphi^s)\geq 0$ and  
$Q_R(\varphi^s)\geq 0$. Hence by Propositions 8.2(1) and 8.4(1), 
$P_L(\varphi^{[^*s]})\geq 0$ and  
$Q_R(\varphi^{[^*s]})\geq 0$. 
For the same reason as above, 
using Propositions 6.5(2) and 3.5(3) we can apply 
Theorem 7.3 to $(\varphi^{[^*s]})^{[t']}$ 
for sufficiently large $t'$ to prove (2)(b).

(2)(c) By Theorem 4.2, 
$Q_L(\varphi^s)\geq 0$ and $Q_R(\varphi^s)\geq 0$. Hence 
by Proposition 8.4, 
$Q_L(\varphi^{[^*s]})\geq 0$ and $Q_R(\varphi^{[^*s]})\geq 0$. 
For the same reason as above, using Proposition 3.5(3) we can apply 
Theorem 4.4 to $(\varphi^{[^*s]})^{[t']}$ 
for sufficiently large $t'$ to prove (2)(c).

(2)(d) By Proposition 6.9 and Theorem 4.2(1), 
$P_L(\varphi^s)\geq 0$ and $Q_L(\varphi^s)\geq 0$, and hence
by Propositions 8.2(1) and 8.4(2), 
$P_L(\varphi^{[^*s]})\geq 0$ and $Q_L(\varphi^{[^*s]})\geq 0$.  
For the same reason as above, 
using Propositions 6.5(2) and 3.5(3) 
we can apply 
Proposition 7.10(3) to  
$(\varphi^{[^*s]})^{[t]}$ for sufficiently large $t$ to 
prove (2)(d). 
\end{proof} 

Theorem 8.7(2)(b) provides 
a short proof for an important special case of 
\cite[Theorem 7.22(2)]{Nasu-t} that if $\varphi^s$ is an LR 
endomorphism of a topological Markov shift, then 
$\varphi$ is an essentially LR endomorphism of the shift. 
Theorem 8.7(2)(c) provides 
a short proof for an important special case of 
\cite[Theorem 7.22(3)]{Nasu-t} that if $\varphi^s$ is a $q$-biresolving  
endomorphism of a topological Markov shift, then 
$\varphi$ is an essentially $q$-biresolving endomorphism of the shift.
By Theorem 8.7(2)(d), 
we obtain the result that if  
$\varphi^s$ is an LL endomorphism 
of a topological Markov shift $(X,\sigma)$ 
with $\varphi$ one-to-one or $\sigma$ transitive, then $\varphi$ is 
essentially LL. This is a new result 
which was not given by \cite{Nasu-t}. 

\begin{corollary} 
Let $\varphi$ be an endomorphism of a subshift $(X,\sigma)$. 
Let \rt be any resolving term. Let $i,i'\in\N$ and $j,j'\in\Z$ 
with $j/i=j'/i'$. 
\begin{enumerate} 
\item 
$\varphi^i\sigma^j$ essentially weakly \rt if and only if 
so is $\varphi^{i'}\sigma^{j'}$. 
\item 
If $(X,\sigma)$ is an SFT and $\varphi$ is onto, 
then $\varphi^i\sigma^j$ is essentially $p$-L (respectively, 
essentially $p$-R) if and only if so is 
$\varphi^{i'}\sigma^{j'}$, 
and $\varphi^i\sigma^j$ is essentially $p$-biresolving 
if and only if so is $\varphi^{i'}\sigma^{j'}$. 
If in addition, $\varphi$ is 
one-to-one or $\sigma$ is transitive, then  
$\varphi^i\sigma^j$ is essentially \rt 
if and only if so is $\varphi^{i'}\sigma^{j'}$. 
\end{enumerate}
\end{corollary} 
\begin{proof} 
(1) We describe a proof for the case that \rt is ``LR''. 
Suppose that $\varphi^i\sigma^j$ is essentially 
weakly LR. Then by Remark 2.15(2), so is
$(\varphi^i\sigma^j)^{i'}$. Hence so is 
$(\varphi^{i'}\sigma^{j'})^i$ because $j/i=j'/i'$. 
Therefore, by Theorem 8.7(1)(b) $\varphi^{i'}\sigma^{j'}$ is 
essentially weakly LR. Hence we see that 
$\varphi^i\sigma^j$ is essentially 
weakly LR if and only if so is $\varphi^{i'}\sigma^{j'}$. 

Similarly (1) is proved for each \rt of the other resolving terms 
using Theorems 8.1(1), 8.6(1) and 8.7(1).

(2) The proof is similar to the above, 
but use Remark 2.15(1) instead of Remark 2.15(2) and use 
Theorems 8.1(3), 8.6(2) and 8.7(2) instead of 
Theorems 8.1(1), 8.6(1) and 8.7(1).
\end{proof}

\subsection{Endomorphisms having some power with positive degree}

The following two theorems are  key results 
for understanding the relation between ``resolvingness'' and 
``expansiveness''.

\begin{theorem} 
Let $\varphi$ be an endomorphism of 
a subshift $(X,\sigma)$. 
\begin{enumerate}
\item
If $\varphi$ is positively left $\sigma$-expansive 
(respectively, positively right $\sigma$-expansive), then
$\varphi$ is essentially weakly $q$-R 
(respectively, essentially weakly $q$-L), 
\item 
When $\varphi$ is onto, $\varphi$ is 
essentially weakly $q$-R and left $\sigma$-expansive
(respectively, essentially weakly $q$-L and right $\sigma$-expansive) 
if and only if $\varphi$ is positively left $\sigma$-expansive 
(respectively, positively right $\sigma$-expansive). 
\item 
If $(X,\sigma)$ is an SFT and 
$\varphi$ is one-to-one or $\sigma$ is topologically transitive, 
then (2) with all ``weakly'' in them deleted hold. 
\end{enumerate}
\end{theorem}
\begin{proof} 
(1) Let $A=L_1(X)$. Let $\N_0=\N\cup\{0\}$. We
define $\hf{O}_{\varphi,\sigma}$ to be 
the set of all two-dimensional 
configurations $(a_{i,j})_{i\in\N_0,j\in\Z}$ with $a_{i,j}\in A$ 
such that there exists $x\in X$ with 
$\varphi^i(x)=(a_{i,j})_{j\in\Z}$ for 
all $i\in \N_0$. Since $\varphi$ is a positively left 
$\sigma$-expansive, it follows that for any 
$(a_{i,j})_{i\in\N_0,j\in\Z}\in\hf{O}_{\varphi,\sigma}$, 
$(a_{i,j})_{i\in\N_0,j\leq 0}$ 
uniquely determines $x=(a_{0,j})_{j\in\Z}$
(for a similar reason to that in the proof of the statement (1) 
given after the statements (1),(2),(3) 
at the beginning of Subsection 11.2)
and hence in particular $a_{0,1}$. 
(See (3) at the beginning of Subsection 11.2.)
Therefore, since
$\hf{O}_{\varphi,\sigma}$ is a compact 
subspace of the product topological space $A^{\N\times\Z}$ 
with a compatible metric, a standard compactness argument 
shows that there exist $k,l\geq 1$ such that 
for any $(a_{i,j})_{i\in\N_0,j\in\Z}\in\hf{O}_{\varphi,\sigma}$, 
the subconfiguration $(a_{i,j})_{0\leq i\leq k, -l+1\leq j\leq 0}$
uniquely determines $a_{0,1}$. 

Recalling the definitions of 
$(X^{[^*_\varphi s]},\sigma^{[^*_\varphi s]})$ (see the 
beginning of Section 8)  
$\varphi^{[^*s]}$ for $s\geq 1$ and the graph $G[X,s]$ 
(see the beginning of Section 5) we define a textile system 
$T=(p,q:\G\to G)$ as follows: 
$G=G_{l+1}[X^{[^*_\varphi k]}]$, that is, 
\begin{align*}
A_G&=\{(a_{i,j})_{0\leq i\leq k-1, 0\leq j\leq l}\,|\,
(a_{i,j})_{i\in\N_0,j\in\Z}\in\hf{O}_{\varphi,\sigma}\},\\ 
V_G&=\{(a_{i,j})_{0\leq i\leq k-1, 0\leq j\leq l-1}\,|\,
(a_{i,j})_{i\in\N_0, j\in\Z}\in\hf{O}_{\varphi,\sigma}\}
\end{align*} 
and for each arc 
$\alpha=(a_{i,j})_{0\leq i\leq k-1, 0\leq j\leq l}$,  
\[i_G(\alpha)=(a_{i,j})_{0\leq i\leq k-1, 0\leq j\leq l-1},\q 
t_G(\alpha)=(a_{i,j})_{0\leq i\leq k-1, 1\leq j\leq l};\]
$\G=G_{l+1}[X^{[^*_\varphi k+1]}]$, that is, 
\begin{align*}
A_\G&=\{(a_{i,j})_{0\leq i\leq k, 0\leq j\leq l}\,|\,
(a_{i,j})_{i\in\N_0,j\in\Z}\in\hf{O}_{\varphi,\sigma}\},\\ 
V_\G&=\{(a_{i,j})_{0\leq i\leq k, 0\leq j\leq l-1}\,|\,
(a_{i,j})_{i\in\N_0,j\in\Z}\in\hf{O}_{\varphi,\sigma}\}
\end{align*} 
and for each arc $\gamma=(a_{i,j})_{0\leq i\leq k, 0\leq j\leq l}$, 
\[
i_\G(\gamma)=(a_{i,j})_{0\leq i\leq k, 0\leq j\leq l-1},\q 
t_\G(\gamma)=(a_{i,j})_{0\leq i\leq k, 1\leq j\leq 1};\\
\]
$p$ and $q$ are such that for each arc $\gamma$ as above 
\[p(\gamma)=(a_{i,j})_{0\leq i\leq k-1, 0\leq j\leq l},\q 
q(\gamma)=(a_{i,j})_{1\leq i\leq k, 0\leq j\leq l}.\]
Since the subconfiguration 
$(a_{i,j})_{0\leq i\leq k, 0\leq j\leq l-1}$
determines $a_{0,l}$
in every configuration 
$(a_{i,j})_{i\in\N_0,j\in\Z}\in \hf{O}_{\varphi,\sigma}$, 
we see that for each arc 
$\gamma=(a_{i,j})_{0\leq i\leq k, 0\leq j\leq l}$, 
$i_\G(\gamma)=(a_{i,j})_{0\leq i\leq k,0\leq j\leq l-1}$
and $q(\gamma)=(a_{i,j})_{1\leq i\leq k, 0\leq j\leq l}$
determines $\gamma$. Hence $T$ is weakly $q$-R. 

Let $Z=(X^{[^*_\varphi k+1]})^{[l+1]}$. 
Then, since $Z\subset\hf{Z}_T$ with 
$\hf{\xi}_T(Z)\supset\hf{\eta}_T(Z)$ and 
$\hf{\xi}_T|Z$ one-to-one, there exists a unique 
onesided 1-1, half-textile-subsystem 
$\hf{U}$ of  $T$ 
with $Z_{\hf{U}}=Z$. Since $T$ is weakly $q$-R and 
$(X_{\hf{U}},\sigma_{\hf{U}},\varphi_{\hf{U}})
=((X^{[^*_\varphi k]})^{[l+1]},
(\sigma^{[^*_\varphi k]})^{[l+1]},(\varphi^{[^*k]})^{[l+1]})$ 
is topologically 
conjugate to $(X,\sigma,\varphi)$, 
we conclude that 
$\varphi$ is an essentially weakly 
$q$-R endomorphism of $(X,\sigma)$. 

(2) Since a positively left $\sigma$-expansive, onto 
endomorphism of $(X,\sigma)$ is left $\sigma$-expansive, 
the ``if'' part follows from (1). 

To prove the ``only-if'' part, suppose that 
$\varphi$ is 
essentially weakly $q$-R and left $\sigma$-expansive. 
There exists a onesided 1-1 textile-subsystem $U$ of 
a weakly $q$-R textile system $T$ such that $\varphi_U$ 
is left $\sigma_U$-expansive and 
$(X_U,\sigma_U,\varphi_U)$ and $(X,\sigma,\varphi)$ are
topologically conjugate. Since $\varphi$ is weakly 
left $\sigma_U$-expansive, 
by \cite[Proposition 7.2(1)]{Nasu-te} there exists 
$k\geq 1$ such that $(U^{[k]})^*$ is onesided 1-1. 
Since $(U^{[k]})^*$ is a textile-subsystem of the textile 
system $(T^{[k]})^*$ which is weakly $p$-L because $T$ is 
weakly $q$-R. Therefore by \cite[Proposition 7.5(1)]{Nasu-te}
$\varphi_{(U^{[k]})^*}$ has memory zero. 

Let $x\in X_U$ and let $\varphi_U^i(x)=(a_{i,j})_{j\in\Z}$ 
for $i\geq 0$. Suppose that $\xi_U^{-1}$ is of $(m,n)$-type. 
Then $(a_{i,j})_{i\in\N_0,j\leq 0}$ determines 
$(\alpha_{i,j})_{i\in\N,j\leq -m+1}$, 
where $u=(\alpha_{i,j})_{i,j\in\Z}$ is any textile (woven by $T$) 
in $U$ with $\alpha_{i,j}\in Z_U$ such that 
$\xi_U((\alpha_{0,j})_{j\in\Z})=x$, and hence it 
determines $(\alpha_{i,j})_{i\in\N,j\in\Z}$ because 
$\varphi_{(U^[k])^*}$ has memory zero. Therefore
$(a_{i,j})_{i\in\N_0,j\leq 0}$ determines
$\xi_U((\alpha_{1,j})_{j\in\Z})=x$. This implies that 
$\varphi_U$ is positively left $\sigma_U$-expansive. 
Hence $\varphi$ is positively left $\sigma$-expansive. 

(3) By (2) and Remark 2.10(2). 
\end{proof} 

\begin{theorem} Let $\varphi$ be an onto endomorphism 
of a subshift $(X,\sigma)$. 
\begin{enumerate}
\item 
$\varphi$ is 
essentially weakly $p$-L and right $\sigma$-expansive
(respectively, essentially weakly $p$-R 
and left $\sigma$-expansive) 
if and only if $\varphi$ is right $\sigma$-expansive 
(respectively, left $\sigma$-expansive) on the upper side. 
\item
If $(X,\sigma)$ is an SFT, then (1) with 
all ``weakly'' in it deleted holds. 
\end{enumerate}
\end{theorem}
\begin{proof} 
(1)
Let $A=L_1(X)$. We
define $O_{\varphi,\sigma}$ to be 
the set of all two-dimensional 
configurations $(a_{i,j})_{i,j\in\Z}$ with $a_{i,j}\in A$ 
such that there exists $x\in X$ with 
$\varphi^i(x)=(a_{i,j})_{j\in\Z}$ for 
all $i\in\Z$. Since $\varphi$ is right $\sigma$-expansive 
on the upper side, it follows that for any 
$(a_{i,j})_{i,j\in\Z}\in O_{\varphi,\sigma}$, 
$(a_{i,j})_{i\leq 0,j\geq0}$ uniquely determines 
$(a_{i,j})_{i\leq 0,j\in\Z}$ 
(for a similar reason to that in the proof of the statement (1) 
given after the statements (1),(2),(3) at  
the beginning of Subsection 11.2)
and hence in particular $a_{0,-1}$. 
(See (2) at the beginning of Subsection 11.2.)
Therefore, since $O_{\varphi,\sigma}$ 
is a compact subspace of the product space 
$A^{\Z^2}$ with a compatible metric, a standard 
compactness argument shows that 
there exist $k,l\in\N$ such that 
for every configuration 
$(a_{i,j})_{i,j\in\Z}\in O_{\varphi,\sigma}$, the 
subconfiguration 
$(a_{i,j})_{-k\leq i\leq 0, 0\leq j\leq l-1}$ 
determines $a_{0,-1}$. 

Let $T=(p,q:\G\to G)$ be the textile system 
(which is the same as defined 
in the proof of Theorem 8.9(1)) 
such that $G=G_{l+1}[X^{[^*_\varphi k]}]$ 
and $\G=G_{l+1}[X^{[^*_\varphi k+1]}]$ and 
for $\gamma=(a_{i,j})_{0\leq i\leq k, 0\leq j\leq l}\in A_\G$ 
with $a_{i,j}\in A$, that is, for 
$\gamma=(a_{i,j})_{0\leq i\leq k, 0\leq j\leq l}$ 
with $(a_{i,j})_{i,j\in\Z}\in O_{\varphi,\sigma}$, 
$p(\gamma)=(a_{i,j})_{0\leq i\leq k-1,0\leq j\leq l}$ and 
$q(\gamma)=(a_{i,j})_{1\leq i\leq k, 0\leq j\leq l}$.
Since the subconfiguration 
$(a_{i,j})_{0\leq i\leq k,1\leq j\leq l}$
determines $a_{k,0}$ in every configuration 
$(a_{i,j})_{i,j\in\Z}\in O_{\varphi,\sigma}$,
it follows that for each arc 
$\gamma=(a_{i,j})_{0\leq i\leq k, 0\leq j\leq l}$, 
$p(\gamma)=(a_{i,j})_{0\leq i\leq k-1,0\leq j\leq l}$
and $t_\G(\gamma)=(a_{i,j})_{0\leq i\leq k, 1\leq j\leq l}$
determines $\gamma$. Hence  $T$ is weakly $p$-L. 

Let $Z=(X^{[^*_\varphi k+1]})^{[l+1]}$. 
Since $Z\subset Z_T$ with
$\xi_T|Z$ one-to-one and $\xi_T(Z)=\eta_T(Z)$, 
there exists a unique 
onesided 1-1, textile-subsystem 
$U$ of $T$ with $Z_U=Z$. Since $T$ is weakly $p$-L and 
$(X_U,\sigma_U,\varphi_U)
=((X^{[^*_\varphi k]})^{[l+1]},
(\sigma^{[^*_\varphi k]})^{[l+1]},(\varphi^{[^*k]})^{[l+1]})$ 
is topologically 
conjugate to $(X,\sigma,\varphi)$, 
we conclude that 
$\varphi$ is an essentially weakly 
$p$-L endomorphism of $(X,\sigma)$. 
Since an endomorphism which is right $\sigma$-expansive 
on the upper side is right $\sigma$-expansive, 
``if'' part is proved. 

To prove ``only-if'' part, suppose that $\varphi$ is 
essentially weakly $p$-L and right $\sigma$-expansive. 
By \cite[Proposition 7.11(1)]{Nasu-te}, there exists 
$l\in N$ such that $\varphi^l$ is 
weakly $p$-L and hence has memory $0$ (by Theorem 6.11). 
Using this we easily observe that 
for any $(a_{i,j})_{i,j\in\Z}\in O_{\varphi,\sigma}$, 
$(a_{i,j})_{i\leq 0,j\geq 0}$
determines $(a_{i,j})_{i\in\Z,j\geq 0}$, which 
determines $(a_{i,j})_{i,j\in\Z}$ because 
$\varphi$ is right $\sigma$-expansive. 
Therefore, 
for any $(a_{i,j})_{i,j\in\Z}\in O_{\varphi,\sigma}$ 
it holds that $(a_{i,j})_{i\leq 0,j\geq 0}$
determines $(a_{i,j})_{i,j\in\Z}$. 
Hence $\varphi$ is 
right $\sigma$-expansive on the upper side. 

(2) By(1) and Remark 2.10(1). 
\end{proof} 

The following proposition is  
``not-necessarily-onto'' versions of 
\cite[Propositions 7.7 and 7.10]{Nasu-te}. 
A proof of the part (1) is given by 
using Theorem 8.9(1) and modifying straightforwardly  
\cite[Proof of Proposition 7.7(2)]{Nasu-te}(by 
using a half-textile-subsystem $\hf{U}$ instead of the 
textile-subsystem $U$ in it). However in the proof below
we describe another proof using Proposition 11.3  
without losing consistency in the subsequent use of 
the part (1). (The result \cite[Propositions 7.7]{Nasu-te} 
could have been shown by using Proposition 11.2.) 

\begin{proposition} 
Let $\varphi$ be a (not necessarily onto) endomorphism of 
a subshift. If $\varphi$ is positively left $\sigma$-expansive 
(respectively, positively right $\sigma$-expansive), then
\begin{enumerate}
\item
there exists $k\geq 1$ such that $\varphi^k\sigma^{-1}$
(respectively, $\varphi^k\sigma$) is  
essentially weakly $q$-R 
(respectively, essentially weakly $q$-L), and
\item 
there exists $m\geq 1$ such that 
$\varphi^n$ is weakly $q$-R (respectively, weakly $q$-L) for 
all $n\geq m$. 
\end{enumerate}
\end{proposition}
\begin{proof} 

(1) Since $\varphi$ is positively left $\sigma$ expansive, 
by Theorem 8.9(1) $\varphi$ is essentially weakly 
$q$-R, and by
Proposition 11.3(2) there exists $k\geq 1$ such that $-1/k$ is 
a positively left $\sigma$-expansive direction for $\varphi$, that is, 
$\varphi^k\sigma^{-1}$ is positively left 
$\sigma$-expansive and 
hence again by  Theorem 8.9(1) $\varphi^k\sigma^{-1}$ is 
essentially weakly $q$-R.  

(2) By (1), $\varphi$ is essentially weakly 
$q$-R. There exist a subshift $(X_1,\sigma_1)$ and a topological 
conjugacy $\theta:(X,\sigma)\to (X_1,\sigma_1)$ such that 
$\varphi_1=\theta\varphi\theta^{-1}$ is a weakly $q$-R endomorphisms 
of $(X_1,\sigma_1)$. Since $\theta$ and $\theta^{-1}$ are 
right closing, by \cite[Proposition 7.10]{Nasu-te} we may assume that 
$\theta$ and $\theta^{-1}\sigma_1^l$ 
is weakly LR for some $l\geq 0$. 

Since $\varphi_1$ is positively left $\sigma$ expansive, by (1)
there exists $k\geq 1$ such that $\varphi_1^k\sigma_1^{-1}$  
essentially weakly $q$-R.  
The remainder of the proof is similar to the corresponding part 
of \cite[Proof of Proposition 7.11]{Nasu-te}.
\end{proof}

\begin{theorem}
If $\varphi$ is an onto endomorphism of a subshift 
$(X,\sigma)$ then the following (1), (2), (3), (4) hold. 
\begin{enumerate}
\item[(1)]
$\varphi$ is essentially weakly $p$-L and 
right $\sigma$-expansive (respectively, essentially weakly $p$-R 
and left $\sigma$-expansive) 
if and only if there exists $s\geq 1$ such that 
$P_L(\varphi^s)>0$ (respectively, $P_R(\varphi^s)>0$). 
\item[(2)] 
$\varphi$ is essentially weakly $q$-R and left $\sigma$-expansive
(respectively, essentially weakly $q$-L and right $\sigma$-expansive) 
if and only if there exists $s\geq 1$ such that 
$Q_R(\varphi^s)>0$ (respectively, $Q_L(\varphi^s)>0$).
\item[(3)] 
The following three statements are equivalent:
\begin{enumerate}
\item[(a)] $\varphi$ is essentially weakly LR  
(respectively, essentially weakly RL) and expansive; 
\item[(b)] there exist $s, t\geq 1$ such that 
$P_L(\varphi^s)>0$ and $Q_R(\varphi^t)>0$ 
(respectively, $P_R(\varphi^s)>0$ and $Q_L(\varphi^t)>0$).
\item[(c)] $\varphi$ is 
right $\sigma$-expansive on the upper side and 
positively left $\sigma$-expansive (respectively, 
left $\sigma$-expansive on the upper side and 
positively right $\sigma$-expansive); 
\end{enumerate}
\item[(4)] 
If $(X,\sigma)$ is an SFT, then (1) with 
with all ``weakly'' in it deleted holds; if in addition, 
$\varphi$ is one-to-one or $\sigma$ is topologically transitive, 
then (2) with all ``weakly'' in it deleted holds 
and (3) with all ``weakly'' in (a) deleted holds. 
\end{enumerate}
If $\varphi$ 
is a (not necessarily onto) endomorphism of a subshift then 
the following hold. 
\begin{enumerate}
\item[(5)]
$\varphi$ is positively left $\sigma$-expansive
(respectively, positively right $\sigma$-expansive) 
if and only if there exists $s\geq 1$ such that 
$Q_R(\varphi^s)>0$ (respectively, $Q_L(\varphi^s)>0$). 
\item[(6)] 
$\varphi$ is positively expansive if and only
if there exist $s, t\geq 1$ such that 
$Q_R(\varphi^s)>0$ and $Q_L(\varphi^t)>0$.
\end{enumerate}
\end{theorem}
\begin{proof} 
(1) Assume that $\varphi$ is essentially weakly $p$-L
and right $\sigma$-expansive. Then it follows from 
\cite[Proposition 7.7]{Nasu-te} that 
there exists $s_1\geq 1$ such that $\varphi^{s_1}\sigma^{-1}$ 
is essentially weakly $p$-L endomorphism of $(X,\sigma)$. 
Hence so is
$(\varphi^{s_1}\sigma^{-1})^2$ (by Remark 2.15), 
and hence 
$(\varphi^{s_1}\sigma^{-1})^2\sigma$ is 
essentially weakly $p$-L and right $\sigma$-expansive 
by \cite[Proposition 7.6]{Nasu-te}.
Therefore,
it follows from 
\cite[Proposition 7.11]{Nasu-te} that 
there exists $s_2\geq 1$ such that 
$(\varphi^{2s_1}\sigma^{-1})^{s_2}$ is 
a weakly $p$-L endomorphism of $(X,\sigma)$. 
Hence, by Proposition 6.12(2) 
$P_L(\varphi^{2s_1s_2})\geq s_2>0$. 

Conversely assume that $P_L(\varphi^s)>0$.  
Then  it follows from 
Proposition 6.12(2) that $\varphi^s$ is weakly $p$-L and 
right $\sigma$-expansive. 
Hence, $\varphi$ is essentially weakly $p$-L, by Theorem 8.1(1), 
and right $\sigma$-expansive. 

(2) By straightforward modifications in the proof above using 
Theorem 5.2(1) (instead of Proposition 6.12(2)) and  
Theorem 8.6(1) (instead of Theorem 8.1(1)). 

(3) We give two different proofs. 

By \cite[Proposition 8.1]{Nasu-te}, $\varphi$ is essentially weakly 
LR and expansive if and only if $\varphi$ is essentially weakly 
$p$-L and right $\sigma$-expansive and 
essentially weakly 
$q$-R and left $\sigma$-expansive. Using this, (1),(2) and Theorems 
8.10 and 8,9 we see that (3) is valid. 

Another proof of (3) is given by showing the equivalence of (a) and (b) 
first and then by using this, (1),(2) and Theorems 
8.10 and 8,9. By (1) and (2), 
(a) implies (b). To prove that (b) implies (a) assume 
that $P_L(\varphi^s)>0$ and $Q_R(\varphi^t)>0$. 
Then by (1) and (2), 
$\varphi$ is right $\sigma$-expansive and left $\sigma$-expansive, 
and hence expansive
(by \cite[Corollary 7.3]{Nasu-te}). If we 
modify the proof of Theorem 8.7(1)(b) 
using Remarks 8.3 and 8.5, then we 
see that $\varphi$ is essentially weakly LR. 

(4) Let (1'), (2') and (a')  be 
the statements (1), (2) and (a)  
with all ``weakly'' in these deleted, respectively. We may assume 
that $(X,\sigma)$ is a topological Markov shift (see Propositions 
6.5(2) and 3.5(3)). 

The proof of (1') under the assumption that $(X,\sigma)$ is 
a topological Markov shift 
is similar to that of (1), but use 
Propositions 2.12, 2.11 and 2.13 (instead of 
\cite[Proposition 7.7, 7.6 and 7.11]{Nasu-te}), 
\cite[Fact 3.16]{Nasu-t} (instead of Remark 2.15), 
Proposition 6.12(3) (instead of 6.12(2)) 
and Theorem 8.1(3) (instead of Theorem 8.1(1)). 

The proof of (2') under the assumptions that $(X,\sigma)$ is 
a topological Markov shift and that 
$\varphi$ is one-to-one or $\sigma$ is topologically transitive 
is  similar to that of (1'), but 
use Theorem 4.2 (instead of Propositions 6.12(3)) and 
Theorem 8.6(2) (instead of Theorem 8.1(3)). 

To prove (4) it remains to show the equivalence of 
(a'), (b) and (c) under the 
assumptions that $(X,\sigma)$ is a topological Markov shift 
and that $\varphi$ is one-to-one or 
$\sigma$ is topologically transitive. 
By (1') and (2'), 
(a') implies (b). To prove that (b) implies (a') assume 
that $P_L(\varphi^s)>0$ and $Q_R(\varphi^t)>0$. 
Then by (1') and (2'), 
$\varphi$ is right $\sigma$-expansive and left $\sigma$-expansive, 
and hence expansive
(by \cite[Proposition 6.2]{Nasu-te} and \cite[Theorem 2.5]{Nasu-t}). 
Using Lemma 9.1 (appearing 
in the next section), 
we see that 
$P_L(\varphi^{st})>0$ and $Q_R(\varphi^{st})>0$. 
It follows from Theorem 7.3 that 
$\varphi^{st}$ is essentially LR. Hence, $\varphi$ is 
essentially LR, by Theorem 8.7(2)(b). 

The equivalence of (a') and (c) follows from  (3) (the equivalence of
(a) and (c)) and Remark 2.10(2). 

(5) Suppose that $Q_R(\varphi^s)>0$ with some $s\geq 1$.
Then by Theorem 5.2 
$\varphi^s$ is positively left $\sigma$-expansive, and hence 
so is $\varphi$. 

Conversely suppose that $\varphi$ is 
positively left $\sigma$-expansive. 
Then it follows from 
Proposition 8.11(1) that 
there exists $s_1\geq 1$ such that $\varphi^{s_1}\sigma^{-1}$ 
is essentially weakly $q$-R endomorphism of $(X,\sigma)$. 
Hence so is
$(\varphi^{s_1}\sigma^{-1})^2$ (by Remark 2.15), 
and hence by Theorem 5.2(1)
$(\varphi^{s_1}\sigma^{-1})^2\sigma$ is 
positively left $\sigma$-expansive .
Therefore,
it follows from Proposition 8.11(2) that 
there exists $s_2\geq 1$ such that 
$(\varphi^{2s_1}\sigma^{-1})^{s_2}$ is 
a weakly $q$-R endomorphism of $(X,\sigma)$. 
Hence, by Proposition  5.2(1) 
$Q_R(\varphi^{2s_1s_2})\geq s_2>0$. 

(6) By (5) and Proposition 11.3(3). 
(By (2) and Theorem 4.5 for an onto endomorphism $\varphi$.)
\end{proof} 

Though the result Theorem 8.12(5) is given 
for any (not necessarily onto) endomorphism 
of any subshift, a positively 
right $\sigma$-expansive or positively 
left $\sigma$-expansive endomorphism of 
a transitive SFT $(X,\sigma)$ 
is necessarily onto (\cite{Kur2},\cite{Sab}). 

A direct proof of Theorem 8.12(6)
for the case when $\varphi$ is an onto endomorphism of 
a transitive SFT $(X,\sigma)$, is given as follows. 
The ``only-if'' part follows from (2') in the proof above and 
the fact recalled after Theorem 4.5. If
$Q_R(\varphi^s)>0$ and $Q_L(\varphi^t)>0$, then 
Using Lemma 9.1 (appearing 
in the next section), 
we see that 
$Q_R(\varphi^{st})>0$ and $Q_L(\varphi^{st})>0$. 
Therefore, by Theorem 4.6(1), 
$\varphi^{st}$ is positively expansive, and hence so is $\varphi$.

Here we present some results on endomorphisms of onesided 
subshifts. 
\begin{theorem} Let $\tilde{\varphi}$ be an endomorphism of a 
onesided-subshift $(\tilde{X},\tilde{\sigma})$. 
Let $\varphi$ be its induced endomorphism
of the induced subshift $(X,\sigma)$ of $(\tilde{X},\tilde{\sigma})$. 
Then the following conditions are equivalent: 
\begin{enumerate}
\item $\tilde{\varphi}$ is positively expansive;
\item $\varphi$ is positively left $\sigma$-expansive;
\item $Q_R(\varphi^s)>0$ for some $s\geq 1$;
\item 
There exists a subshift $(X_0,\sigma_0)$, 
a positively left $\sigma_0$-expansive, 
weakly LR endomorphism $\varphi_0$ of $(X_0,\sigma_0)$ such that 
$(\tilde{X},\tilde{\sigma},\tilde{\varphi})$ and  
$(\tilde{X}_0,\tilde{\sigma}_0,\tilde{\varphi}_0)$ 
are topologically conjugate 
(i.e., there exists a conjugacy
$\tilde{\theta}:
(\tilde{X},\tilde{\sigma})\to 
(\tilde{X}_0,\tilde{\sigma}_0)$ 
(between onesided-subshifts) such that 
$\tilde{\varphi}_0=\tilde{\theta}\tilde{\varphi}\tilde{\theta}^{-1}$).  
\end{enumerate} 
If in addition, $(\tilde{X}, \tilde{\sigma})$ 
is a transitive onesided SFT, then 
(1) is equivalent to each of 
the following statements:
\begin{enumerate} 
\item[(5)] $\varphi$ is onto, essentially $q$-R and 
left $\sigma$-expansive;
\item[(6)] There exists a topological Markov shift $(X_0,\sigma_0)$, 
a left $\sigma_0$-expansive, LR (onto) endomorphism $\varphi_0$ 
of $(X_0,\sigma_0)$ such that 
$(\tilde{X},\tilde{\sigma},\tilde{\varphi})$ and 
$(\tilde{X}_0,\tilde{\sigma}_0,\tilde{\varphi}_0)$ are 
topologically conjugate. 
\end{enumerate}
\end{theorem} 
\begin{proof} First we prove the equivalence of (1) and (2). 
Assume (1). Then there exists $k\geq 1$ such that 
for any $x,y\in X$ the following holds: if
$\varphi^i(x)=(a_{i,j})_{j\in\Z}$ and 
$\varphi^i(y)=(b_{i,j})_{j\in\Z}$ for all $i\in\N_0$ with 
$a_{i,j}, b_{i,j}\in L_1(X)$, then 
it holds that if $(a_{i,j})_{i\in\N_0, 1\leq j\leq k}
=(b_{i,j})_{i\in\N_0, 1\leq j\leq k }$ then 
$(a_{0,j})_{j\geq 1}
=(b_{0,j})_{j\geq 1}$ and hence
it holds that if
$(a_{i,j})_{i\in\N_0,j\leq 0}
=(b_{i,j})_{i\in\N_0,j\leq 0}$ then
$(a_{0,j})_{j\geq -k+1}
=(b_{0,j})_{j\geq -k+1}$ 
and hence $x=y$. Therefore (2) follows.

Assume (2). Let
$\hf{O}_{\varphi,\sigma}$ be 
the set of all 
configurations $(a_{i,j})_{i\in\N_0,j\in\Z}$ with $a_{i,j}\in L_1(x)$ 
such that there exists $x\in X$ with 
$\varphi^i(x)=(a_{i,j})_{j\in\Z}$ for 
all $i\in \N_0$.
As we saw in the proof of Theorem 8.9(1), 
there exist $k,l\geq 1$ such that 
for any $(a_{i,j})_{i\in\N_0,j\in\Z}\in\hf{O}_{\varphi,\sigma}$, 
the subconfiguration $(a_{i,j})_{0\leq i\leq k,-l+1\leq j\leq 0}$
uniquely determines $a_{0,1}$. Therefore,  
for any $(a_{i,j})_{i\in\N_0,j\in\Z}\in\hf{O}_{\varphi,\sigma}$, 
$(a_{i,j})_{i\in\N_0, -l+1\leq j\leq 0}$
uniquely determines $(a_{i,j})_{i\in\N_0,-l+1\leq j\leq 1}$  
and hence $(a_{i,j})_{i\in\N_0,j\geq -l+1}$, 
so that it uniquely determines 
$(a_{0,j})_{j\geq -l+1}$. Consequently 
for any $(a_{i,j})_{i\in\N_0,j\in\Z}\in\hf{O}_{\varphi,\sigma}$,
$(a_{i,j})_{i\in\N_0, 1\leq j\leq l}$ determines
$(a_{0,j})_{j\geq 1}$, which implies (1). 

By Theorem 8.12(5), (2) is equivalent to (3). 

Assume (3). Then, since $\varphi$ has memory zero, so does 
$\varphi^s$ and hence by Theorem 6.11 $\varphi^s$ is weakly $p$-L.
Therefore, since $Q_R(\varphi^s)>0$, 
it follows from Theorem 5.2(1) that $\varphi^s$ is weakly $q$-R and 
positively left $\sigma$-expansive. 
Therefore, by Theorem 8.7(1)(b), there exists 
$t\geq 1$ such that $(\varphi^{[^* s]})^{[t]}$ is weakly LR and 
positively left $\sigma_0$-expansive, where 
$(X_0,\sigma_0)=(Y^{[t]}, \sigma_{Y^{[t]}})$ 
with $Y=X^{[^*_\varphi s]}$.
Since $\varphi$ has memory zero, so does $\rho^*_{\varphi,s}$, and 
clearly so does
$(\rho^*_{\varphi,s})^{-1}$.   
Let $\theta=\rho_{Y,0,t}\rho^*_{\varphi,s}$ 
(recalling that $\rho_{Y,0,t}$
is the higher-block conjugacy of $(0,t)$-type on $Y$)
and let $\varphi_0=\theta\varphi\theta^{-1}$. Then, since 
$\rho_{Y,0,t}$ and $\rho_{Y,0,t}^{-1}$  have memory zero, 
(4) is proved. 

Assume (4). Then since $(X,\sigma,\varphi)$ and 
$(X_0,\sigma_0,\varphi_0)$ are topologically conjugate, 
$\varphi$ is positively left $\sigma$-expansive and hence 
(2) follows.

To prove the remainder, 
suppose that and $(\tilde{X}, \tilde{\sigma})$ 
is a transitive onesided SFT. 
Passing through a higher-block conjugacy 
between endomorphisms of 
onesided subshifts, 
we may assume that $(\tilde{X},\tilde{\sigma})$ is 
a transitive onesided topological Markov shift. 

Since any positively expansive endomorphism of 
a transitive onesided SFT is onto and conjugate 
to a onesided topological Markov shift \cite{Kur1}, it follows 
from \cite[Theorem 3.13]{Nasu-t} and \cite[Proposition 6.2]{Nasu-te}
that (1) implies (6). 

Assume (6). 
Then there exists a onesided 1-1, nondegenerate, LR 
textile system T such that 
$(X_T,\sigma_T,\varphi_T)=(X_0,\sigma_0,\varphi_0)$. 
Since $\varphi_T$ is left $\sigma_T$-expansive, 
by \cite[Proposition 6.2]{Nasu-te} $T^*$ is onesided 1-1. 
Therefore,
since $T$ and $T^*$ are $p$-L (because $T$ is $q$-R), 
by \cite[Lemma 3.6(1)]{Nasu-t} 
$\tilde{\xi}_T$ and $\tilde{\xi}_{T^*}$ are 
one-to-one. Therefore, by \cite[Theorem 2.12]{Nasu-t} 
$\tilde{\varphi}_T$ is positively expansive, and 
hence (1) follows. Also (5) follows, since $(X,\sigma,\varphi)$ and 
$(X_0,\sigma_0,\varphi_0)$ are topologically conjugate.

By Theorem 8.12(4), (5) implies (3). 

Assume (3). Then, since $\varphi$ has memory zero, so is 
$\varphi^s$ and hence by Proposition 6.9 $\varphi^s$ is $p$-L. 
Since $Q_R(\varphi^s)>0$, by Theorem 4.2 $\varphi^s$ is $q$-R 
and left $\sigma$-expansive. 
Therefore,  using Theorem 8.7(2)(b) we see that 
$(\varphi^{[^* s]})^{[t]}$ is LR for some $t\geq 1$ 
and left $\sigma_0$-expansive, where 
$(X_0,\sigma_0)=(Y^{[t]}, \sigma_{Y^{[t]}})$ 
with $Y=X^{[^*_\varphi s]}$. 
Therefore (6) is proved in a similar way to to that in the above 
proof that (3) implies (4). 
\end{proof} 

\section{Limits of onesided resolving directions} 
\subsection{Limits of onesided resolving directions}

The following lemma directly follows from Propositions 6.7 and 3.8. 
\begin{lemma} 
Let $\varphi$ be an endomorphism of 
a subshift. Let $s,t\geq 0$.
\begin{align*}
&P_L(\varphi^{s+t})\geq P_L(\varphi^s)+P_L(\varphi^t);&
P_R(\varphi^{s+t})\geq P_R(\varphi^s)+P_R(\varphi^t);\\
&Q_R(\varphi^{s+t})\geq Q_R(\varphi^s)+Q_R(\varphi^t);&
Q_L(\varphi^{s+t})\geq Q_L(\varphi^s)+Q_L(\varphi^t).
\end{align*}
\end{lemma} 

In what follows, we follow the convention that 
if $\alpha_s=\pm\infty$ for all $s\in\N$ then 
$\lim_{s\to\infty}\alpha_s=\pm\infty$ and 
$\sup_s\alpha_s=\pm\infty$. 
\begin{theorem}
Let $\varphi$ be an endomorphism of a subshift. 
\begin{enumerate}
\item 
$\lim_{s\to\infty} P_L(\varphi^s)/s$ exists and equals
$\sup_s P_L(\varphi^s)/s$; 
\item
$\lim_{s\to\infty} P_R(\varphi^s)/s$ exists 
and equals $\sup_s P_R(\varphi^s)/s$; 
\item $\lim_{s\to\infty} Q_R(\varphi^s)/s$ exists and equals
$\sup_s Q_R(\varphi^s)/s$; 
\item $\lim_{s\to\infty} Q_L(\varphi^s)/s$ exists 
and equals $\sup_s Q_L(\varphi^s)/s$.
\end{enumerate}
\end{theorem} 
\begin{proof} 
(1) If $P_L(\varphi^s)\in\Z$ for all $s\geq 1$, then 
it follows from Lemma 9.1 and
\cite[Theorem 4.9]{Wal-2} that (1) is valid. If 
$P_L(\varphi^s)=\infty$ for some $s\geq 1$, then 
by Lemma 9.1 $P_L(\varphi^t)=\infty$ for all $t\geq s$ and 
hence (1) follows. 

(2) The proof of (2) is similar. 

(3) For the cases that  
$Q_R(\varphi^s)\in\Z$ for all $s\geq 1$ and that 
$Q_R(\varphi^s)=\infty$ for some $s\geq 1$, 
(3) is valid for the same reasons as in (1). 
If $Q_R(\varphi^s)=-\infty$ for some $s\geq 1$, then (3) is 
valid, because $\varphi$ is not right-closing and hence 
$Q_R(\varphi^s)=-\infty$ for all $s\geq 1$.  

(4) The proof of (4) is similar. 
\end{proof} 

\begin{definition}
Let $\varphi$ be an endomorphism of a subshift. Define 
\begin{align*}
p_L(\varphi)&=\lim_{s\to\infty} P_L(\varphi^s)/s,&
p_R(\varphi)=\lim_{s\to\infty} P_R(\varphi^s)/s,\\
q_L(\varphi)&=\lim_{s\to\infty} Q_L(\varphi^s)/s,&
q_R(\varphi)=\lim_{s\to\infty} Q_R(\varphi^s)/s. 
\end{align*}
We call $-p_L(\varphi)$ (respectively, $p_R(\varphi)$, 
$-q_R(\varphi)$, $q_L(\varphi)$)   
the \itl{limit of} \itl{$p$-L} (respectively, 
\itl{$p$-R} \itl{$q$-R}, 
\itl{$q$-L}) \itl{directions of $\varphi$};  
or the \itl{$p$-L} (respectively, 
\itl{$p$-R} \itl{$q$-R}, 
\itl{$q$-L}) \itl{limit of $\varphi$}, for short; 
$-p_L(\varphi)$, $p_R(\varphi)$, $-q_R(\varphi)$ and $q_L(\varphi)$
are generically 
called the 
\itl{limits of onesided resolving directions of $\varphi$}. 
\end{definition} 

We remark that the following hold 
for an endomorphism $\varphi$  of an infinite subshift 
$(X,\sigma)$: 
if $\varphi^i(X)$ is infinite for all $i\geq 1$ then 
$p_L(\varphi)\in\R, p_R(\varphi)\in\R$, and 
otherwise, $p_L(\varphi)=\infty$ and $p_R(\varphi)=\infty$ 
(by Proposition 9.9(2) which will appear shortly), and  
hence $p_L(\varphi)\in\R$ if and only if $p_R(\varphi)\in\R$;
if $\varphi$ is right-closing, then $q_R(\varphi)\in\R$, 
and otherwise, $q_R(\varphi)=-\infty$, and  
if $\varphi$ is left-closing, then $q_L(\varphi)\in\R$, 
and otherwise, $q_L(\varphi)=-\infty$ (by Theorem 9.2(3),(4),  
Proposition 9.9(1) and Definitions 3.1 and 9.3). We also remark that
if $\varphi$ is an endomorphism of 
a finite subshift, then all $p_L(\varphi)$, 
$p_R(\varphi)$, $q_L(\varphi)$, 
$q_R(\varphi)$ are $\infty$ (by Definitions 6.1, 3.1 and 9.3). 

For any infinite subshift $(X,\sigma)$, 
$p_L(i_X)=p_R(i_X)=q_R(i_X)=q_L(i_X)=0$ and 
$p_L(\sigma^s)=q_R(\sigma^s)=s$ and $p_R(\sigma^s)=q_L(\sigma^s)=-s$ 
for $s\in\Z$ (by (3.1),(3.2),(6.1),(6.2)). 

The following theorem shows that 
the limits of onesided resolving directions 
is an invariant of topological conjugacy between 
endomorphisms of subshifts. 
\begin{theorem} 
If endomorphisms $\varphi$ 
and $\psi$ of subshifts are topologically conjugate, then 
\[-p_L(\varphi)=-p_L(\psi),\q p_R(\varphi)=p_R(\psi),\q 
-q_R(\varphi)=-q_R(\psi),\q q_L(\varphi)=q_L(\psi).\]
\end{theorem} 
\begin{proof}  We may assume that $\varphi$ is an 
endomorphism of an infinite subshift
(otherwise, the theorem is clear). 
It follows from Propositions 6.13 and 6.15 and Theorem 9.2 
that
$-p_L(\varphi)=-p_L(\psi)$ and $p_R(\varphi)=p_R(\psi)$. 
If $\varphi$ is right-closing, then 
it follows from Proposition 5.4 and Theorem 9.2 that 
$q_R(\varphi)=q_R(\psi)$; if $\varphi$ is not right-closing then
$q_R(\varphi)=-\infty=q_R(\psi)$. Therefore 
we have $-q_R(\varphi)=-q_R(\psi)$. Similarly we see that
$q_L(\varphi)=q_L(\psi)$. 
\end{proof} 

\begin{theorem}
Let $\varphi$ and $\psi$ be endomorphisms of 
a subshift $(X,\sigma)$. 
\begin{enumerate}
\item 
If $\varphi^i(X)$ and $\psi^i(X)$ are infinite for 
all $i\geq 0$, 
and hence in particular, 
if $X$ is infinite and $\varphi$ and $\psi$ 
are onto, then 
$-p_L(\varphi\psi)=-p_L(\psi\varphi)$ and
$p_R(\varphi\psi)=p_R(\psi\varphi)$. 
\item If $\varphi$ and $\psi$ are right-closing then
$-q_R(\varphi\psi)=-q_R(\psi\varphi)$,  
and if $\varphi$ and $\psi$ are left-closing
then $q_L(\varphi\psi)=q_L(\psi\varphi)$. 
\end{enumerate}
\end{theorem} 
\begin{proof} 
(1) It follows from Proposition 6.7 that 
$P_L((\varphi\psi)^s)\geq 
P_L((\psi\varphi)^{s-1}) +P_L(\varphi)+P_L(\psi)$ 
for $s\geq 1$. Since 
$\varphi^i(X)$ and $\psi^i(X)$ are infinite for 
all $i\geq 0$, it follows from Proposition 6.13(1) that  
$P_L(\varphi),P_L(\psi)\in\Z$. Therefore by Theorem 9.2 
we have  $p_L(\varphi\psi)\geq p_L(\psi\varphi)$. 
Similarly we have 
$p_L(\psi\varphi)\geq p_L(\varphi\psi)$. Therefore
$-p_L(\varphi\psi)=-p_L(\psi\varphi)$. 

(2) We may assume that $X$ is infinite (because if 
$X$ is finite then (2) directly follows 
from Definitions 3.1 and 9.3).
It follows from Proposition 3.8 that 
$Q_R((\varphi\psi)^s)\geq 
Q_R((\psi\varphi)^{s-1}) +Q_R(\varphi)+Q_R(\psi)$ 
for $s\geq 1$. If $\varphi$ and $\psi$ are 
right-closing then  
$Q_R(\varphi),Q_R(\psi)\in\Z$ and hence by Theorem 9.2 we have  
$q_R(\varphi\psi)\geq q_R(\psi\varphi)$. 
Similarly we have 
$q_R(\psi\varphi)\geq q_R(\varphi\psi)$. Therefore
$-q_R(\varphi\psi)=-q_R(\psi\varphi)$. 

The proof of the remainder is similar.  
\end{proof} 

\begin{theorem}
Let $\varphi$ be an endomorphism of a subshift $(X,\sigma)$. 
Let $i\geq 0$ and $j\in \Z$. Then 
\begin{align*}
&p_L(\varphi^i\sigma^j)=ip_L(\varphi)+j,& 
p_R(\varphi^i\sigma^j)=ip_R(\varphi)-j,\\
&q_R(\varphi^i\sigma^j)=iq_R(\varphi)+j,& 
q_L(\varphi^i\sigma^j)=iq_L(\varphi)-j;
\end{align*} 
each of 
$p_L(\varphi)+p_R(\varphi)$, $p_L(\varphi)+q_L(\varphi)$, 
$p_R(\varphi)+q_R(\varphi)$ and $q_L(\varphi)+q_R(\varphi)$ 
exists and is shift-invariant. 
\end{theorem} 
\begin{proof}  Using Proposition 6.4,  we have
\begin{align*}
p_L(\varphi^i\sigma^j)&=\lim_{s\to\infty} P_L(\varphi^{is}\sigma^{js})/s 
=\lim_{s\to\infty} P_L(\varphi^{is})/s  + j\\
&=\lim_{s\to\infty}iP_L(\varphi^{is})/(is) +j =ip_L(\varphi)+j.
\end{align*} 
The other equations are similarly proved. Hence  
the sums $p_L(\varphi)+p_R(\varphi)$, $p_L(\varphi)+q_L(\varphi)$, 
$p_R(\varphi)+q_R(\varphi)$ and $q_L(\varphi)+q_R(\varphi)$ are 
shift-invariant if they exist. 
By the facts remarked after Definition 9.3, 
the existence of the sums follows. 
\end{proof} 

\begin{theorem}
Let $\varphi$ be an endomorphism of a subshift $(X,\sigma)$. 
Let $i\in\N$ and $j\in\Z$. If $\varphi$ is onto, then 
the following statements (1),(2) and (3) hold: 
\begin{enumerate} 
\item[(1)] 
$\varphi$ is essentially weakly $p$-L 
and right $\sigma$-expansive (respectively, essentially 
weakly $p$-R and left $\sigma$-expansive)
if and only if $p_L(\varphi)>0$ (respectively $p_R(\varphi)>0$), 
and moreover,  
$\varphi^i\sigma^j$ is essentially weakly $p$-L 
and right $\sigma$-expansive 
(respectively, essentially 
weakly $p$-R and left $\sigma$-expansive)
 if and only if $j/i> -p_L(\varphi)$
(respectively, $j/i< p_R(\varphi)$); 
\item[(2)]
$\varphi$ is essentially weakly $q$-R 
and left $\sigma$-expansive (respectively, essentially 
weakly $q$-L and right $\sigma$-expansive)
if and only if $q_R(\varphi)>0$ (respectively, $q_L(\varphi)>0$), 
and moreover,  
$\varphi^i\sigma^j$ is essentially weakly $q$-R 
and left $\sigma$-expansive 
(respectively, essentially 
weakly $q$-L and right $\sigma$-expansive) 
if and only if $j/i> -q_R(\varphi)$
(respectively, $j/i< q_L(\varphi)$);
\item[(3)] 
If $(X,\sigma)$ is an SFT, then (1) with 
all ``weakly'' in them deleted hold; if in addition, 
$\varphi$ is one-to-one or $\sigma$ is topologically transitive, 
then (2) with all ``weakly'' in them deleted hold. 
\end{enumerate}
If $\varphi$ is not necessarily onto, then 
the following statement hold: 
\begin{enumerate}
\item[(4)] 
$\varphi$ is positively left $\sigma$-expansive 
(respectively, 
positively right $\sigma$-expansive)
if and only if $q_R(\varphi)>0$ (respectively, $q_L(\varphi)>0$), 
and moreover,  
$\varphi^i\sigma^j$ is positively left $\sigma$-expansive 
(respectively, positively right $\sigma$-expansive) 
if and only if $j/i> -q_R(\varphi)$
(respectively, $j/i< q_L(\varphi)$)
\end{enumerate}
\end{theorem} 
\begin{proof} 
The theorem follows from Theorems 8.12, 9.2 and 9.6.
\end{proof}

\begin{corollary} Let $\tilde{\varphi}$ be an endomorphism 
of a onesided subshift $(\tilde{X},\tilde{\sigma})$. Let $\varphi$ be its 
induced endomorphism of the induced subshift $(X,\sigma)$ of 
$(\tilde{X},\tilde{\sigma})$. Then
$\tilde{\varphi}$ is positively expansive if and only if $q_R(\varphi)>0$.
Moreover, for $i\geq 1, j\geq 0$, 
$\tilde{\varphi}^i\tilde{\sigma}^j$ is positively expansive 
if and only if $j/i>-q_R(\varphi)$. 
\end{corollary} 
\begin{proof}
By Theorems 8.13 and 9.7(4). 
\end{proof}

For an onto endomorphism $\varphi$ of a subshift, define
\begin{align*}
c_R(\varphi)&=\max\{-p_L(\varphi), -q_R(\varphi)\},\q\q
c_L(\varphi)=\min\{p_R(\varphi), q_L(\varphi)\},\\ 
d_R(\varphi)&=\min\{-p_L(\varphi), -q_R(\varphi)\},\q\q
d_L(\varphi)=\max\{p_R(\varphi), q_L(\varphi)\}. 
\end{align*} 
Clearly, $c_L(\varphi)=c_R(\varphi)$ if and only if 
all limits of onesided resolving directions of $\varphi$ coincide. 

\begin{proposition} 
Let $\varphi$ be an endomorphism of an infinite subshift 
$(X,\sigma)$. 
\begin{enumerate} 
\item[(1)] 
If $\varphi$ is of $(m,n)$-type, then 
\begin{align*}
&-n\leq p_R(\varphi),\q\q\q\q -p_L(\varphi)\leq m,\\
&-n\leq -q_R(\varphi),\q\q\q\q q_L(\varphi)\leq m. 
\end{align*}
\item[(2)] 
If $\varphi^i(X)$ is infinite for all $i\geq 1$, then 
$p_L(\varphi)+p_R(\varphi)$ is nonpositive, and hence if in addition 
$\varphi$ is of $(m,n)$-type then 
\[-n\leq p_R(\varphi)\leq -p_L(\varphi)\leq m;\] 
if $\varphi^i(X)$ is finite for some $i\geq 1$, then 
\[p_L(\varphi)=\infty ,\q\q p_R(\varphi)=\infty.\]
\item[(3)] 
$p_L(\varphi)+q_L(\varphi)$ and 
$p_R(\varphi)+q_R(\varphi)$ are nonpositive and hence 
\[p_R(\varphi)\leq -q_R(\varphi), \q
q_L(\varphi)\leq -p_L(\varphi) \q\text{and}\q
c_L(\varphi)\leq c_R(\varphi).\]
\item [(4)]
If $\varphi$ is of $(m,n)$-type, then 
\begin{enumerate}
\item[(a)] 
 if $q_R(\varphi)+q_L(\varphi)\geq 0$ then
\[-n\leq p_R(\varphi)\leq -q_R(\varphi)\leq 
q_L(\varphi)\leq -p_L(\varphi)\leq m;\] 
\item [(b)]
if $q_R(\varphi)+q_L(\varphi)\leq 0$ then
\[-n\leq d_L(\varphi)\leq d_R(\varphi)\leq m.\] 
\end{enumerate}
\end{enumerate}
\end{proposition} 
\begin{proof} 
If $\varphi$ is of $(m,n)$-type, then $\varphi^s$ is of 
$(ms,ns)$-type for $s\geq 1$. Hence (1) follows 
by (6.3),(3.3) and Theorem 9.2, 
and (2) follows by Proposition 6.13, Theorem 9.2 and (1). 
By Proposition 7.9(1) and Theorem 9.2, (3) follows. 
By (1) and (3), (4) follows. 
\end{proof} 

\begin{theorem}
Let $\varphi$ be an onto endomorphism of a subshift $(X,\sigma)$.  
Let $i\in\N$ and $j\in\Z$. 
 \begin{enumerate}
 \item 
$\varphi^i\sigma^j$ is essentially weakly LR 
(respectively, essentially weakly RL) and  
 expansive   if and only if $j/i> c_R(\varphi)$ 
(respectively, $j/i< c_L(\varphi)$). 
\item 
When $(X,\sigma)$ is an SFT, if $\varphi$ is one-to-one 
or $\sigma$ is transitive, 
then all ``weakly'' can be deleted in the statements (1). 
\end{enumerate} 
\end{theorem} 
\begin{proof} (1) If $j/i> c_R(\varphi)$, then, 
by Theorem 9.7, $\varphi^i\sigma^j$ is essentially weakly $p$-L, 
right $\sigma$-expansive, essentially weakly $q$-R 
and left $\sigma$-expansive. Therefore, using Theorem 8.12(1),(2),(3), 
we see that $\varphi^i\sigma^j$ essentially weakly LR and expansive. 

Conversely assume that $\varphi^i\sigma^j$ is essentially weakly LR and  
expansive. Then it follows from Theorem 9.7 that
$j/i> c_R(\varphi)$. 

(2)  By using Theorem 8.12(4) and 9.7(3), (2) is proved. 
\end{proof} 

\begin{theorem} 
Let $\varphi$ be an endomorphism of a subshift $(X,\sigma)$. 
There exist $i\in\N$ and $j\in\Z$ 
such that $\varphi^i\sigma^j$ is positively expansive
if and only if 
\[q_R(\varphi)+q_L(\varphi)>0.\]
Moreover, for $i\in\N$ and $j\in\Z$, $\varphi^i\sigma^j$ is 
positively expansive if and only if 
\[-q_R(\varphi)<j/i<q_L(\varphi).\] 
\end{theorem}
\begin{proof}  By Theorem 9.7(4) and Proposition 11.3(3). 
(By Theorem 9.7(2) and Theorem 4.5 for an onto endomorphism $\varphi$.)
\end{proof} 
A direct proof of Theorem 9.11 
for the case that $\varphi$ is an onto endomorphism of 
a transitive topological Markov shift $(X,\sigma)$ 
is similarly given by Theorem 8.12(4)(with the comment after 
the proof of Theorem 8.12) and Theorem 9.7(3). 

\begin{proposition} \begin{enumerate}
\item If $\varphi$ is an automorphism of a subshift, then
\begin{gather*}
q_L(\varphi)=p_L(\varphi^{-1}),\q\q\q\q
q_R(\varphi)=p_R(\varphi^{-1}),\\ 
c_L(\varphi^{-1})=-c_R(\varphi),\q\q\q\q 
d_L(\varphi^{-1})=-d_R(\varphi). 
\end{gather*} 
\item If $\varphi$ is an automorphism of an infinite subshift, then
\[q_R(\varphi)+q_L(\varphi)\leq 0\q\q\text{and}\q\q 
d_L(\varphi)\leq d_R(\varphi).\]
\item
If $\varphi$ is an essentially symbolic automorphism 
of an infinite subshift, then 
\[p_L(\varphi)=q_R(\varphi)=p_R(\varphi)=q_L(\varphi)=0.\]
\end{enumerate}
\end{proposition}
\begin{proof}
By Proposition 7.11(1),(2),(3),Theorem 9.2 and Theorem 9.4.
\end{proof} 

Here we explain about our terminology of Definition 9.3. 
Let $\varphi$ be an endomorphism of 
a subshift $(X,\sigma)$.  
In view of Proposition 6.9, Theorem 6.11 and 
Theorems 4.2 and 5.2, let us call 
a rational number $r$ 
a \itl{$p$-L direction} (respectively, 
\itl{$p$-R direction}, \itl{$q$-R direction}, 
\itl{$q$-L direction}) 
of $\varphi$
if  there exist $i\in\N, j\in\Z$ such that 
$r=j/i$ and $P_L(\varphi^i\sigma^j)\geq 0$
(respectively, $P_R(\varphi^i\sigma^j)\geq 0$, 
$Q_R(\varphi^i\sigma^j)\geq 0$, 
$Q_L(\varphi^i\sigma^j)\geq 0$). Then 
it follows from Propositions 6.4 and 3.4 that 
$r\in\Q$ is 
a $p$-L direction (respectively, 
$p$-R direction, $q$-R direction, $q$-L direction)  
of $\varphi$ if and only if for some 
$i\in\N$, $r\geq -P_L(\varphi^i)/i$
(respectively, $r\leq P_R(\varphi^i)/i$, 
$r\geq -Q_R(\varphi^i)/i$, $r\leq Q_L(\varphi^i)/i$). 
Therefore it follows from Theorem 9.2 
that the limit $-p_L(\varphi)$ of $p$-L directions (respectively, 
the limit $-q_R(\varphi)$ of $q$-R directions) of $\varphi$ 
is the infimum of the set of 
$p$-L (respectively, $q$-R) directions of $\varphi$
and the limit $p_R(\varphi)$ of $p$-R directions (respectively, 
the limit $q_L(\varphi)$ of $q$-L directions) of $\varphi$ 
is the supremum of  $p$-R (respectively, $q$-L) 
directions of $\varphi$. 

Finally we discuss here about the relation 
between the limits $-p_L(\varphi)$ and $p_R(\varphi)$
of an endomorphism $\varphi$ of 
an infinite 
subshift $(X,\sigma)$  and 
the right and left Lyapunov 
exponents for a cellular automaton $\varphi$
which were defined and treated 
by M. A. Shereshevsky \cite{Sher} and 
further studied  by P. Tisseur \cite{Tiss} and  
M. Hochman \cite{Hoch}. 

We follow the notation which appeared in the  paragraph 
after the proof of Proposition 6.2. Let $\varphi$ be 
an endomorphism of an infinite subshift $(X,\sigma)$. 
For $s\in\Z$ and $x\in X$, define
\begin{align*}
\tilde{\Lambda}_k^+(x)&=\min\{s\geq 0\,|\,\varphi^k(W_0^+(x))\subset 
W_s^+(\varphi^k(x))\} \q \text{and}\\
\tilde{\Lambda}_k^-(x)&=\min\{s\geq 0\,|\,\varphi^k(W_0^-(x))\subset 
W_s^-(\varphi^k(x))\}.  
\end{align*} 
(The value $\tilde{\Lambda}^+_k(x)$ 
(respectively, $\tilde{\Lambda}^-_k(x)$) can be 
considered to show 
how far a perturbation front can exist on the right 
(respectively, on the left) at time $k$ 
if the front is initially
located at coordinate $0$ for the configuration $x$ 
in the cellular automaton $\varphi$. See \cite{Sher}.)
Define  
$\Lambda_k^+(x)=\max_{j\in\Z}\tilde{\Lambda}_k^+(\sigma^j(x))$ and
$\Lambda_k^-(x)=\max_{j\in\Z}\tilde{\Lambda}_k^-(\sigma^j(x))$. 
Shereshevsky \cite{Sher} proved that if a
probability measure $\mu$ on $X$ is $\sigma$-invariant and 
$\varphi$-invariant, then for $\mu$-almost all $x\in X$
the limits \[\lambda^+(x)\stackrel{\mathrm{def}}{=}
\lim_{k\to\infty}(1/k)\Lambda_k^+(x)\q\text{and}\q 
\lambda^-(x)\stackrel{\mathrm{def}}{=}  
\lim_{k\to\infty}(1/k)\Lambda_k^-(x)\] exist with the functions 
$\lambda^+$ and $\lambda^-$ $\varphi$-invariant,  
and called 
them the \itl{right and left Lyapunov exponents of $\varphi$ at $x$}, 
remarking that $\lambda^+$ and $\lambda^-$ are 
constant for $\mu$-almost all $x\in X$ if $\mu$ is $\varphi$-ergodic. 
Let 
\[\hat{\Lambda}_k^+=\max_{x\in X}\tilde{\Lambda}_k^+(x)\q\text{and}\q
\hat{\Lambda}_k^-=\max_{x\in X}\tilde{\Lambda}_k^-(x).\]
Then Tisseur \cite{Tiss} proved 
that $\lim_{k\to\infty}(1/k)\hat{\Lambda}_k^+$
and $\lim_{k\to\infty}(1/k)\hat{\Lambda}_k^-$ exist and 
that if $(X,\sigma)$ is an irreducible topological Markov shift 
and $\varphi$ is onto, then 
for $\mu$-almost all $x\in X$, where $\mu$ is the 
Parry measure for $(X,\sigma)$, and for all $k\geq 1$, 
$\Lambda^+_k(x)=\hat{\Lambda}_k^+$ 
and $\Lambda^-_k(x)=\hat{\Lambda}_k^-$ and hence 
$\lambda^+(x)=\lim_{k\to\infty}(1/k)\hat{\Lambda}_k^+$ 
and $\lambda^-(x)= \lim_{k\to\infty}(1/k)\hat{\Lambda}_k^-$ for 
$\mu$-almost all $x\in X$.
The \itl{right and left Lyapunov exponents of $\varphi$} 
which Hochman \cite{Hoch} defined are equal to 
$\lim_{k\to\infty}(1/k)\hat{\Lambda}_k^+$ 
and $\lim_{k\to\infty}(1/k)\hat{\Lambda}_k^-$, respectively. 

If $\varphi$ is of $(m,n)$-type given by a local rule 
$f: L_{m+n+1}(X)\to L_1(X)$, then we see that for $x\in X$
\[\tilde{\Lambda}^+_k(x)=\max\{km-\rho_{f^k,L}(x),\,0\}\q\text{and} 
\q \tilde{\Lambda}^-_k(x)=\max\{kn-\rho_{f^k,R}(x),\,0\}. \]
Hence, since $P_L(\varphi^k)=\min_{x\in X}\rho_{f^k,L}(x)-km$ 
and $P_R(\varphi^k)=\min_{x\in X}\rho_{f^k,R}(x)-kn$, 
it easily follows that 
\[
\hat{\Lambda}_k^+=\begin{cases}
0 & \mathrm{if}\; P_L(\varphi^k)\geq 0 \\
-P_L(\varphi^k) & \mathrm{if}\; P_L(\varphi^k)\leq 0
\end{cases}\q\text{and}\q
\hat{\Lambda}_k^-=\begin{cases}
0 & \mathrm{if}\; P_R(\varphi^k)\geq 0 \\
-P_R(\varphi^k) & \mathrm{if}\;  P_R(\varphi^k)\leq 0 
\end{cases}. 
\]
If there exists $k\geq 1$ such that $P_L(\varphi^k)\geq 0$, then 
for all $i\geq 1$, $P_L(\varphi^{ik})\geq 0$ (by Lemma 9.1), 
and hence $\lim_{k\to\infty}(1/k)\hat{\Lambda}_k^+=0$, 
because $((1/k)\hat{\Lambda}_k^+)_{k\in\N}$ is a 
convergent sequence. We similarly see that
if there exists $k\geq 1$ such that $P_R(\varphi^k)\geq 0$, 
then $\lim_{k\to\infty}(1/k)\hat{\Lambda}_k^-=0$. Therefore, 
using Theorems 6.11 and 9.2, we see: 
\begin{proposition}
Let $\varphi$ be an endomorphism of an infinite subshift $(X,\sigma)$.
If there exits $k\geq 1$ such that  $\varphi^k$ has memory zero 
then $\lim_{k\to\infty}(1/k)\hat{\Lambda}_k^+=0$, and otherwise 
$\lim_{k\to\infty}(1/k)\hat{\Lambda}_k^+=-p_L(\varphi)$. 
If there exits $k\geq 1$ such that  $\varphi^k$ has anticipation zero
then $\lim_{k\to\infty}(1/k)\hat{\Lambda}_k^-=0$, and otherwise 
$\lim_{k\to\infty}(1/k)\hat{\Lambda}_k^-=-p_R(\varphi)$.  
\end{proposition}

\subsection{Types of limits}

Let $\varphi$ be an endomorphism of an 
infinite subshift $(X,\sigma)$. 
By Proposition 9.9(2) 
we know that if $\varphi^i(X)$ is infinite 
for all $i\geq 1$ then $-p_L(\varphi), p_R(\varphi)\in\R$. 
By Proposition 9.9(1) we know that if $\varphi$ is 
right-closing then $-q_R(\varphi)\in\R$, and if $\varphi$ is 
left-closing then $q_L(\varphi)\in\R$.  
Each of the real limits of the onesided resolving directions 
of $\varphi$ is classified 
into the following three types: 
\begin{enumerate} 
\item[I.]  
The limit $-p_L(\varphi)$ 
(of $p$-L directions of $\varphi$) is 
said to be of \itl{type I} 
if there exist $i\in\N, j\in\Z$ such that  
$j/i=-p_L(\varphi)$ and $\varphi^i\sigma^j$ is an
essentially weakly $p$-L endomorphism of $(X,\sigma)$.
Similarly, the definition 
of being of \itl{type I} is given for each of 
the limits $-q_R(\varphi)$, $p_R(\varphi)$ and $q_L(\varphi)$. 
\item[II.] 
The limit $-p_L(\varphi)$ is 
said to be of \itl{type II} 
if there exist $i\in\N, j\in\Z$ such that  
$j/i=-p_L(\varphi)$ and $\varphi^i\sigma^j$ is not an 
essentially weakly $p$-L endomorphism of $(X,\sigma)$.  
Similarly, the definition  
of being of \itl{type II} is given for each of 
the limits $-q_R(\varphi)$, $p_R(\varphi)$ and $q_L(\varphi)$. 
\item[III.]  
The limit $-p_L(\varphi)$ 
is said to be of \itl{type III}
if $-p_L(\varphi)$ is an irrational number. 
Similarly, the definition  
of being of \itl{type III} is given for each of 
the limits $-q_R(\varphi)$, $p_R(\varphi)$ and $q_L(\varphi)$.  
\end{enumerate} 

By Corollary 8.8 we know that 
The definitions of type I limits and type II ones above 
do not contradict each other. 

For an endomorphism $\varphi$ of 
an infinite subshift $(X,\sigma)$, we see that 
the limit $-p_L(\varphi)$ 
is of type I if and only if some 
$p$-L direction of $\varphi$ essentially attains the limit, 
i.e., there exists 
a topological conjugacy $\theta:(X,\sigma)\to (X',\sigma')$
such that for some $i\in\N$ 
$-P_L((\theta\varphi\theta^{-1})^i)/i=-p_L(\varphi)$. 
Similar statements are also valid for the limits $p_R(\varphi)$, 
$-q_R(\varphi)$ and $q_L(\varphi)$. 

Here we define, in case $d=2$, the notion of an expansive line
for a $\Z^d$-action on an infinite compact metric space 
which was introduced by Boyle and Lind \cite{BoyLin}. 
Let $\alpha:(i,j)\mapsto \alpha^{(i,j)}$ 
be a $\Z^2$-action on an infinite compact metric space $X$. 
Let $\ell$ be a line on the plane $\R^2$. 
For $t\geq 0$ let $\ell^t$ denote the set of all points  
having distance not greater than $t$ from $\ell$. Then $\ell$ is 
\itl{expansive for $\alpha$} if there exist $\epsilon>0$ and $t\geq 0$
such that for any $x,y\in X$ 
it holds that if $d_X(\alpha^{(i,j)}(x),\alpha^{(i,j)}(y))\leq\epsilon$
for all $(i,j)\in \ell^t\cap\Z^2$ then $x=y$. 

We define the \itl{direction} 
of the line passing through 
points $(a,b),(c,d)\in\R^{2}$ with $a\neq c$ to be $(b-d)/(a-c)$. 
A number $r\in \R$ is called an \itl{expansive direction for $\alpha$} 
if there exists an expansive line of direction $r$ for $\alpha$. 

\begin{example} (1) Let $r_1,r_2\in\R$ with $r_1<0<r_2$. 
The result of Boyle and Lind \cite[Proof of Proposition 4.1]{BoyLin} 
shows that there exists an automorphism $\varphi$ 
of an infinite subshift $(X,\sigma)$ 
such that $r_1$ and $r_2$ are non-expansive, in fact, more precisely, 
neither left $\sigma$-expansive nor
right $\sigma$-expansive directions 
for the $\Z^2$-action $\alpha:(i,j)\mapsto \varphi^i\sigma^j$ 
and there are no more non-expansive directions for $\alpha$.  
The set-union of all expansive lines passing through $(0,0)$, with 
$\{(0,0)\}$ subtracted 
is divided into four connected sets which are open cones in $\R^2$.
These are the ``expansive components of 1-frames'' for $\alpha$ 
of Boyle and Lind \cite{BoyLin}. Let $\C$ be the component 
containing point $(0,1)$. 
Then $\C=\{(a,b)\in\R^2\,|\, b>0,\,b/r_1<a<b/r_2\}$. 
Let $C^\circ_K(\sigma)$ 
be the set of all essentially weakly LR, 
expansive automorphism of $(X,\sigma)$ 
in $K=\{\varphi^i\sigma^j\,|\, i,j\in\Z\}$. 
As is implied by \cite[Remark 9.6]{Nasu-te} 
together with Theorem 12.2 or by Theorem 11.12, 
$\{\varphi^i\sigma^j\,|\, (i,j)\in\C\cap\Z^2\}=C_K^\circ(\sigma)$. 
Therefore, it follows from Theorems 9.10 and 9.7  
that $c_L(\varphi)=p_R(\varphi)=q_L(\varphi)=r_1$ and 
$c_R(\varphi)=-p_L(\varphi)=-q_R(\varphi)=r_2$. 

(2) Let $r_0\in\R$. According to Mike Hochman \cite{Hoch}, 
there exists an automorphism $\pi$ of an infinite subshift 
$(X,\sigma)$ such that  $r_0$ is 
a unique non-expansive direction 
for the $\Z^2$-action $\beta:(i,j)\mapsto\pi^i\sigma^j$  
and such that if $r_0$ is rational then $\pi^i\sigma^j$ with $j/i=r_0$ is
not a periodic automorphism (or essentially-symbolic automorphism). 
The expansive components of 
1-frames for $\beta$ are the two open half-planes separated by 
the line passing through $(0,0)$ with direction $r_0$. 
Let $C^\circ_K(\sigma)$ 
be the set of all essentially weakly LR, 
expansive automorphism of $(X,\sigma)$  
in $K=\{\pi^i\sigma^j\,|\, i,j\in\Z\}$. 
As is implied by \cite[Remark 9.6]{Nasu-te} 
with Theorem 12.2 or by Theorem 11.12, 
$C_K^\circ(\sigma)=\{\varphi^i\sigma^j\,|\,j>r_0i \}$, and hence 
by Theorem 9.10 
we see that $c_L(\varphi)=c_R(\varphi)=r_0$. Hence 
we have $p_R(\varphi)=q_L(\varphi)
=-p_L(\varphi)=-q_R(\varphi)=r_0$. 

The automorphism $\pi$ constructed by
Hochman in \cite[Section 7]{Hoch} 
for $r_0=0$ is 
not essentially weakly \rt 
automorphisms of $(X,\sigma)$, 
where \rt is any one of the terms ``$p$-L'', 
``$p$-R'', ``$q$-R'' and ``$q$-L'', and hence all of the limits of 
onesided resolving directions of $\pi$ are of type II. 
This follows because neither $\pi$ nor $\pi^{-1}$ 
satisfies either one of the necessary conditions
given by Proposition 6.16(4) for essentially weakly 
$p$-L endomorphisms and for essentially weakly $p$-R ones, 
for the same reason as that presented 
by Hochman in \cite[Section 7]{Hoch} 
to show that the action of $\pi$ is not periodic.  
\end{example} 

By the above examples, 
we know the existence of automorphisms of subshifts with 
type III limits of onesided resolving directions 
and those with type II ones. 
However we cannot answer the following question. Are  
the (real) limits of onesided 
resolving directions of onto endomorphisms of an SFT always type I?
Related questions are found in \cite[Section 9]{BoyLin} and 
\cite[Section 11]{Boy-o}. 

For every endomorphism of every SFT, we can calculate its degrees of 
onesided resolvingness. However we know no general 
method for calculating the limits of onesided resolving directions 
even for onto endomorphisms of SFTs. 
The following proposition is useful to calculate 
type I limits of onesided resolving directions of onto 
endomorphisms, especially of SFTs. 

\begin{proposition}
Let $\varphi$ be an onto endomorphism of 
an infinite subshift $(X,\sigma)$. Let $k\geq 1$.
\begin{enumerate}
\item $-P_L(\varphi^k)/k=-p_L(\varphi)$ if and only if 
$\varphi^k\sigma^{-P_L(\varphi^k)}$ is not right $\sigma$-expansive. 
\item  If $\varphi$ is right-closing, then 
$-Q_R(\varphi^k)/k=-q_R(\varphi)$ if and only if 
$\varphi^k\sigma^{-Q_R(\varphi^k)}$ is not left $\sigma$-expansive. 
\item 
$P_R(\varphi^k)/k=p_R(\varphi)$ if and only if 
$\varphi^k\sigma^{P_R(\varphi^k)}$ is not left $\sigma$-expansive. 
\item  If $\varphi$ is left-closing, then 
$Q_L(\varphi^k)/k=q_L(\varphi)$ if and only if 
$\varphi^k\sigma^{Q_L(\varphi^k)}$ is not right $\sigma$-expansive.  
\end{enumerate}
\end{proposition} 
\begin{proof} By Theorem 9.6 it suffices to show the proposition 
for $k=1$. 

To prove (1) for $k=1$, let $\psi=\varphi\sigma^{-P_L(\varphi)}$. 
Assume that $-P_L(\varphi)\neq -p_L(\varphi)$. Since 
$p_L(\varphi)=\sup_{s\in\N} P_L(\varphi^s)/s$, there exists 
$s\geq 1$ such that $P_L(\varphi^s)/s>P_L(\varphi)$. 
It follows that 
\[P_L(\psi^s)=P_L(\varphi^s\sigma^{-sP_L(\varphi)})
=P_L(\varphi^s)-sP_L(\varphi)>0.\]
Therefore, since $P_L(\psi^s\sigma^{-1})=P_L(\psi^s)-1\geq 0$, 
$\psi^s\sigma^{-1}$ is a weakly $p$-L endomorphism of $(X,\sigma)$ 
(by Theorem 6.11). 
Hence it follows from \cite[Proposition 7.6]{Nasu-te} that 
$\psi^s$ is right $\sigma$-expansive, and hence so is $\psi$. 

Assume that $\psi$ is right $\sigma$-expansive. 
Since $P_L(\psi)=0$, $\psi$ is a 
weakly $p$-L endomorphism of $(X,\sigma)$, by Theorem 6.11. 
Therefore, it follows from 
Theorem 8.12 that there exists $s\geq 1$ such that $P_L(\psi^s)>0.$
Therefore $P_L(\psi^s\sigma^{-1})\geq 0$.   
Since 
\[P_L(\varphi^s)-sP_L(\varphi)-1=P_L(\varphi^s\sigma^{-sP_L(\varphi)-1})
=P_L(\psi^s\sigma^{-1})\geq 0,\]  
we have $P_L(\varphi^s)/s\geq P_L(\varphi)+1/s$. Therefore, 
$P_L(\varphi)\neq \sup_s P_L(\varphi^s)/s=p_L(\varphi)$. 

The proofs of (2),(3) and (4) for $s=1$ are similarly given 
by using Theorems 6.11, 5.2 and 8.12 and \cite[Proposition 7.6]{Nasu-te}. 
\end{proof} 

Direct proofs of (1) and (3) of the theorem above for the case that 
$\varphi$ is an onto endomorphism of 
an SFT $(X,\sigma)$  and those of (2) and (4) of it for the case that 
in addition, 
$\varphi$ is one-to-one or $\sigma$ transitive, 
are given similarly 
by using Proposition 6.9, Theorems 4.2 and 8.12(4) and 
Proposition 2.11. 

\begin{proposition}
Let $\varphi$ be an onto endomorphism 
of an infinite topological Markov shift 
$(X,\sigma)$ such that $\varphi$ is one-to-one 
or $\sigma$ is transitive. 
\begin{enumerate}
\item If $\varphi$ is right-closing and 
$-P_L(\varphi)\geq -Q_R(\varphi)$ 
(respectively, $-P_L(\varphi)\leq -Q_R(\varphi)$), then  
we can decide whether $\varphi\sigma^{-P_L(\varphi)}$ 
(respectively, $\varphi\sigma^{-Q_R(\varphi)}$) is right 
$\sigma$-expansive (respectively, left $\sigma$-expansive) or not. 
\item If $\varphi$ is left-closing and 
$P_R(\varphi)\leq Q_L(\varphi)$ 
(respectively, $P_R(\varphi)\geq Q_L(\varphi)$), then  
we can decide whether $\varphi\sigma^{P_R(\varphi)}$  
(respectively, $\varphi\sigma^{Q_L(\varphi)}$) is left 
$\sigma$-expansive (respectively, right $\sigma$-expansive) or not. 
\item If $-Q_R(\varphi)\leq Q_L(\varphi)$, then 
we can decide whether
$\varphi\sigma^{-Q_R(\varphi)}$ 
(respectively, $\varphi\sigma^{Q_L(\varphi)}$) is left $\sigma$-expansive
(respectively, right $\sigma$-expansive) or not. 
\end{enumerate}
\end{proposition}
\begin{proof} 
By Corollary 7.4(1)(d),(2)(d) and Theorem 4.6(2). 
\end{proof} 

By the proof of Proposition 9.16 and \cite[Lemma 6.25]{Nasu-t}, 
we see that 
the decision procedures whose existence is asserted by Proposition 9.16 
are intended to use a decision procedure for ``definiteness'' of 
right-resolving and left-resolving graph-homomorphisms.  
Definite right resolving and left 
resolving graph-homomorphisms 
appear in \cite[pp.96,97]{Nasu-t} as an extension of the definite 
transition diagram of a finite automaton, which was introduced 
Perles, Rabin and Shamir \cite{PerRabSha}. 
The properties of the definite transition diagrams of finite automata
and a practical decision procedure for them presented by 
\cite{PerRabSha} can be straightforwardly extended to 
right-resolving and left-resolving graph-homomorphisms.
 
\begin{example} Let $A=\{0,1\}$. 
Let $(X,\sigma)$ be the full 2-shift $(A^\Z,\sigma_A)$. 
Let $e:A\to A$ be 
the local rule such that $e(0)=1$ and $e(1)=0$. 
Let $f_0:A^4\to A$ be the local rule such that 
$f_0(1001)=1, f_0(1101)=0$ and $f_0(abcd)=b$  for $a,b,c,d\in A$
with $acd\neq 101$. 
Let $f:A^4\to A$ be the local rule such that $f=ef_0$. 

Let $\epsilon$, $\varphi_0$ and $\varphi$
be the automorphism of $(X,\sigma)$ as follows:
$\epsilon$ is the symbolic automorphism given by $e$;
$\varphi_0$ and $\varphi$ are 
of $(1,2)$-type and given by $f_0$ and $f$, respectively. 
These are the well-known automorphisms 
found in \cite[Section 20]{Hedlund}. 

Since $f_0$ and $f$ are zero left-redundant and zero 
right-redundant, 
$P_L(\varphi_0)=P_L(\varphi)=-1$ and 
$P_R(\varphi_0)=P_R(\varphi)=-2$. Since $\varphi_0^2$ is 
the identity $i_X$, 
we have, by Proposition 7.11(1), 
$Q_L(\varphi_0)=-1$ and $Q_R(\varphi_0)=-2$.  Since 
$\varphi=\epsilon\varphi_0$, we have, by Proposition 7.11(4), 
$Q_L(\varphi)=-1$ and $Q_R(\varphi)=-2$.
Hence, 
we know that $f_0$ and $f$ are 
strictly $2$ left-mergible and strictly $4$ right-mergible. 

Since $\varphi_0^2$ is the identity,  
$\varphi_0$ is 
an essentially symbolic automorphism of $(X,\sigma)$ 
(by \cite[Proposition 2.9]{BoyLinRud}). 
Hence 
there exists an SFT $(\bar{X},\bar{\sigma})$ and 
a conjugacy $\theta:(X,\sigma) \to (\bar{X},\bar{\sigma})$ such that 
$\bar{\varphi}_0=\theta^{-1}\varphi_0\theta$ 
is a symbolic automorphism of $(\bar{X},\bar{\sigma})$. 
By Proposition 7.11(3), 
$P_L(\bar{\varphi}_0)=Q_R(\bar{\varphi}_0)
=P_R(\bar{\varphi}_0)=Q_L(\bar{\varphi}_0)=0$. 
Thus we see that each of the degrees 
of onesided resolvingness is not 
an invariant of topological conjugacy 
(between endomorphisms of subshifts),  
and we also see, by Proposition 9.12, that   
$p_L(\varphi_0)=q_R(\varphi_0)=p_R(\varphi_0)=q_L(\varphi_0)=0$. 

Consider the textile system 
$T_2=(p_2,q_2:\G_2\to G)$ in \cite[page 201]{Nasu-t}. 
Then $T_2$ is 1-1 and LR and $\varphi_{T_2}=(\varphi\sigma^2)^{[2]}$. 
Since $\xi_{T_2^*}$ is not one-to-one, 
$(\varphi\sigma^2)^{[2]}$ is not left $\sigma^{[2]}$-expansive. 
Hence $\varphi\sigma^{-Q_R(\varphi)}$ is not left $\sigma$-expansive. 
Therefore, it follows from Proposition 9.15(2) that 
$-q_R(\varphi)=-Q_R(\varphi)=2$. 

Similar arguments show that $\varphi\sigma^{-2}$ is not 
left $\sigma$-expansive, and hence we know by 
Proposition 9.15(3) that $p_R(\varphi)=P_R(\varphi)=-2$. 

As is found in \cite[Example 8.12(2)]{Nasu-te}, 
$\varphi\sigma$ is not left $\sigma$-expansive 
(it is not right 
$\sigma$-expansive either) and 
$\varphi\sigma^{-1}$ is not right $\sigma$-expansive 
(it is not left $\sigma$-expansive either). 
Therefore, it follows from Proposition 9.15(1),(4) that 
$-p_L(\varphi)=-P_L(\varphi)=1$ and $q_L(\varphi)=Q_L(\varphi)=-1$. 
\end{example} 

\begin{example} 
The reader is referred to \cite[Section 10, Example 1]{Nasu-t}.
In it, we find an LR automorphism  
$\varphi$ of a topological Markov shift $(X_M,\sigma_M)$ 
and 1-1, LR textile systems $T$ and $T'$ such 
that $(X_M,\sigma_M)=(X_T,\sigma_T)=(X_{T'},\sigma_{T'})$,  
$\varphi_T=\varphi$, and $(\varphi\sigma_M^{-1})^{-1}=\varphi_{T'}$. 
The topological Markov shift $(X_M,\sigma_M)$ is 
defined by the graph $G$ with the 
representation matrix 
$\begin{pmatrix} a+b & c\\ d & e\end{pmatrix}$. 
Let $f:L_2(G)\to A_G$ be the local rule defined by 
$aa\mapsto b, ab\mapsto c, ac\mapsto b, ba\mapsto d, bb\mapsto e, 
bc\mapsto d, ce\mapsto a, cd\mapsto a, ee\mapsto a, ed\mapsto a,
da\mapsto b, db\mapsto c, dc\mapsto b$. 
Then $\varphi$ is of $(0,1)$-type given by $f$. 
It is easily seen that $f$ is strictly $0$ left-redundant 
and strictly 0 right-redundant. It is also seen that 
$f$ is strictly $1$ right-mergible and 
strictly $1$ left-mergible. Hence 
$-P_L(\varphi)=0$, $-Q_R(\varphi)=0$, $P_R(\varphi)=-1$ and 
$Q_L(\varphi)=-1$. Therefore we have
$\varphi\sigma^{-P_L(\varphi)}=\varphi\sigma^{-Q_R(\varphi)}=\varphi_T$ 
and 
$\varphi\sigma^{P_R(\varphi)}
=\varphi\sigma^{Q_L(\varphi)}=\varphi_{T'}^{-1}$.
Since $T$ and $T'$ are LR, using 
the decision procedure for definiteness of right-resolving 
and left-resolving graph-homomorphisms we easily see that
$\varphi_T$ and $\varphi_{T'}^{-1}$ are neither right $\sigma_M$-expansive
nor left $\sigma_M$-expansive  
(as seen in \cite[Example 8.12(1)]{Nasu-te}). 
Hence, by Proposition 9.15 we have 
$-p_L(\varphi)=-P_L(\varphi)=0$, $-q_R(\varphi)=-Q_R(\varphi)=0$, 
$p_R(\varphi)=P_R(\varphi)=-1$ and 
$q_L(\varphi)=Q_L(\varphi)=-1$. 
\end{example}

\begin{example} The main purpose here is to show the existence 
of a non-invertible onto endomorphism $\varphi_1$ of an SFT 
such that $q_L(\varphi_1)+q_R(\varphi_1)<0$. 

Let $A=\{0,1\}$. Let $f_0:A^4\to A$ be the same as in 
Example 9.17 and $g:A^2\to A$ be the local rule such that 
$g(ab)=a+b \pmod{2}$ for $a,b\in A$. Let $f_1:A^5\to A$ be the 
local rule defined by $f_1=gf_0$. 
Let $\varphi_1$ be the endomorphism of $(2,2)$-type given by $f_1$. 
Since $f_1$ is zero left-redundant and zero right-redundant, 
$P_L(\varphi_1)=0-2=-2$ and $P_R(\varphi_1)=0-2=-2$. 

Since $f_0$ is strictly $2$ left-mergible 
and strictly $4$ right-mergible 
(as seen in Example 9.17), so is $f_1=gf_0$ by Remark 3.7, 
because $g$ is zero right-mergible and zero left-mergible.  
Therefore, $Q_L(\varphi_1)=0$ and $Q_R(\varphi_1)=-2$. 

Since $P_L(\varphi_1\sigma^2)=0$ and $Q_R(\varphi_1\sigma^2)=0$, 
we know, by Theorem 7.3, that 
$\varphi_1\sigma^2$ is LR up to higher block 
conjugacy. Since $Q_L(\varphi_1)=0$, we know, by Theorem 4.2, that 
$\varphi_1$ is $q$-L. Since $P_R(\varphi_1\sigma^{-2})=0$ 
and $Q_L(\varphi_1\sigma^{-2})=2$, we know , by Theorem 7.3, that 
$\varphi_1\sigma^{-2}$ is $RL$ up to higher block conjugacy.

Let $T_{0,k}=(p_{0,k},q_{0,k}:G_A^{[5]}\to G_A)$ be 
the onesided 1-1, nondegenerate textile 
system such that for each arc $w=a_{-2}\dots a_2\in A^5$ in $G_A^{[5]}$ 
with $a_j\in A$, $p_{0,k}(w)=a_k$ and $q_{0,k}(w)=f_1(w)$, for 
$ -2 \leq k\leq 2$. Then $\varphi_{T_{0,k}}=\varphi_1\sigma^{-k}$ 
for $ -2 \leq k\leq 2$.

We can have a onesided 1-1, LR textile system 
$T_{-2}$ such that $\varphi_{T_{-2}}=(\varphi\sigma^2)^{[2]}$ 
(such a textile system is unique, by \cite[Corollary 7.25]{Nasu-t}). 
In fact, if $T'_{-2}$ is the textile system such that 
$\M_{T'_{-2}}$ is the representation matrix of order 6 obtained 
by applying ``same-column reductions'' to $\M_{T_{0,-2}}$ 
(see \cite[p.575]{Nasu-t-ex} for ``same-row and 
same-column reductions'' of representation matrices),
then $\M_{T_{-2}}$ is the 
representation matrix of order 8 obtained 
by applying same-column reductions to 
$\M_{(T'_{-2})^{[2]}}$ (which is of order 12). 
Since neither $\xi_{T_{-2}^*}$ nor $\eta_{T_{-2}^*}$
is one-to-one, $\varphi_1\sigma^2$ is neither 
left $\sigma$-expansive nor right $\sigma$-expansive. 
Therefore, 
since $-P_L(\varphi_1)=-Q_R(\varphi_1)=2$, it follows by 
Proposition 9.15 that 
$-p_L(\varphi_1)=-P_L(\varphi_1)=2$ and $-q_R(\varphi_1)=-Q_R(\varphi_1)=2$.
 
Let $T_0$ be 
the textile system such that $\M_{T_0}$ is 
obtained by applying same-column reductions to 
the representation matrix of order 8 obtained by applying 
same-row reductions to $\M_{T_{0,0}}$. 
Then $T_0$ is onesided 1-1, nondegenerate and $q$-L, 
$\varphi_{T_0}=\varphi_1$, 
\[\M_{T_0}=
\bordermatrix{
        & a& b& c& d& e\cr 
       a& \frac{0}{0}&\frac{0}{0}& & & \cr
       b& & &\frac{0}{0}& & \cr
       c& & & &\frac{1}{1}&\frac{1}{1}\cr
       d&\frac{0}{1}&\frac{1}{1}& &\frac{1}{0}&\frac{0}{0}\cr 
       e& & &\frac{0}{1}& &} 
\quad\text{and}\quad  \M_{T_0^*}=
\bordermatrix{
        & 0& 1\cr 
       0& \frac{a}{a}+\frac{a}{b}+\frac{b}{c}+\frac{d}{e}
        &\frac{d}{a}+\frac{e}{c}\cr
       1& \frac{d}{d}&\frac{c}{d}+\frac{c}{e}+\frac{d}{b} }.
\] 
Let $\dot{a}$ be the fixed point 
$(a_j)_{j\in\Z}$ in $(X_{T_0^*},\sigma_{T_0^*})$ with 
$a_j=a$ for all $j\in\Z$ and let $\dot{b}$, 
$\dot{c}$, $\dot{d}$ and $\dot{e}$ be the fixed points 
similarly defined. 
Observing the trace of the representation matrix $\M_{T_0^*}$, 
we know that $T_0^*$ weaves 
a textile on which 
$\dots,\dot{a},\dot{a},\dot{a},
\dot{b},\dot{c},\dot{d},\dot{b},\dot{c},\dot{d},\dots$
appear in order vertically, 
the periodic textile  
$\dots,\dot{a},\dot{a},\dot{a},\dots$ , and 
a textile on which 
$\dots,\dot{b},\dot{c},\dot{d},
\dot{b},\dot{c},\dot{d},\dot{b},\dot{c},\dot{d},\dots$. 
appear in order vertically. Therefore, we see that 
neither $\xi_{T_0^*}$ nor $\eta_{T_0^*}$ is one-to-one,  
and hence $\varphi_1$ is not onesided $\sigma$-expansive. 
Since $\varphi_1\sigma^{Q_L(\varphi_1)}$ is not right $\sigma$-expansive, 
we have, by Proposition 9.15, $q_L(\varphi_1)=Q_L(\varphi_1)=0$. 

Hence we have $q_L(\varphi_1)+q_R(\varphi_1)=-2$. 

If we apply same-column reductions to 
the representation matrix of order 10 
obtained by applying same-row reductions 
to $\M_{T_{0,2}}$, then we have the representation matrix 
$\M_{T_2}$ (of order 8) of a onesided 1-1, RL
textile system $T_2$ with 
$\varphi_{T_2}=\varphi_1\sigma^{-2}$ (which is unique by 
\cite[Corollary 7.25]{Nasu-t}). We easily check that 
neither $\xi_{T_2^*}$ nor $\eta_{T_2^*}$ is one-to-one. 
Since $\varphi_1\sigma^{P_R(\varphi_1)}$ is not 
left $\sigma$-expansive, we see, by Proposition 9.15, that 
$p_R(\varphi_1)=P_R(\varphi_1)=-2$. 
\end{example}

\begin{example} (1) Let $A$ be an alphabet. 
Let $f:A^{N+1}\to A$ be a \itl{bipermutive} local rule 
(i.e. rightmost permutive (or zero right-mergible) and 
leftmost permutive (or zero left-mergible)) local rule 
(on the full shift over $A$) introduced in \cite{Hedlund}. 
Let $i\geq 0$. Then $f^i:A^{Ni+1}\to A$ is a bipermutive 
local rule. Let $m,n\geq 0$ with $m+n=N$ and let 
$\varphi$ be the endomorphism of $(m,n)$-type of 
the full-shift $(X,\sigma)=(X_A,\sigma_A)$ given by $f$. 
Then $\varphi^i$ is of $(mi,ni)$-type and given by $f^i$. 
Hence $Q_L(\varphi^i)=mi$ and $Q_R(\varphi^i)=ni$,
and hence, by definition $q_L(\varphi)=m$ and $-q_R(\varphi)=-n$. 
Since $f^i$ is zero left-redundant and zero right-redundant, 
$P_L(\varphi^i)=-mi$ and $P_R(\varphi^i)=-ni$ and hence
by definition, $-p_L(\varphi)=m$ and $p_R(\varphi)=-n$. 

(2) Let $A,B$ be alphabets. Let $r,s\geq 0$. 
Let $f_1:A^{r+1}\to A$ and $f_2:B^{s+1}\to B$ be 
bipermutive local rules. 
Let $\hat{f}_1:A^{r+1+s}\to A$ and 
$\hat{f}_2:B^{r+1+s}\to B$ be
the local rules 
such that $\hat{f}_1(w_1aw_2)=f_1(w_1a)$ and 
$\hat{f}_2(w_1aw_2)=f_2(aw_2)$ 
for $w_1\in A^r$, $a\in A$, $w_2\in A^s$. Then 
$\hat{f}_1^t$ and $\hat{f}_2^t$
are a zero left-mergible, strictly $st$ right-mergible local rule 
and a strictly $st$ left-mergible, zero right-mergible local rule, 
respectively, for all $t\geq 1$ (by Lemma 3.2). 
Let $f:(A\times B)^{r+1+s}\to A\times B$ be a local rule 
such that for $a_1,\dots, a_{r+1+s}\in A$ and $b_1,\dots,b_{r+1+s}\in B$, 
$f((a_1,b_1)\dots (a_{r+1+s},b_{r+1+s}))
=(\hat{f}_1(a_1\dots a_{r+1+s}),\hat{f}_2(b_1\dots b_{r+1+s}))$.
Then $f^{t}$ is strictly $rt$ left-mergible and 
strictly $st$ right-mergible for all $t\geq 1$. 
Let $\varphi$ be the endomorphism of the full-shift of $(r,s)$-type of
$(X,\sigma)=(X_{A\times B},\sigma_{A\times B})$ 
given by $f$. 
Then for all $t\geq 1$, $\varphi^{t}$ is of $(rt, st)$-type and 
given by $f^t$.  Therefore $Q_L(\varphi^t)=0$ and $Q_R(\varphi^t)=0$ 
for all $t\geq 1$, and hence $q_L(\varphi)=-q_R(\varphi)=0$.
Though $\varphi$ is essentially $q$-biresolving (by Theorem 4.6(1)), 
there exist no $i\in\N, j\in\Z$ such that $\varphi^i\sigma^j$ 
is positively expansive (by Theorem 9.11). Since $f^t$ is 
zero left-redundant and zero right-redundant, we have 
$P_L(\varphi^t)=rt$ and $P_R(\varphi^t)=st$, and hence 
$p_R(\varphi)=-s$ and $-p_L(\varphi)=r$. 
\end{example}

A nontrivial example  
of an essentially $q$-biresolving 
endomorphism $\varphi$ of a full-shift $(X,\sigma)$ such that  
there exist no $i\in\N,j\in\Z$ with $\varphi^i\sigma^j$ positively
expansive was
given in \cite[p.890]{Nasu-n} by using an arithmetic constraint 
\cite[Theorem 4.2]{Nasu-n} for positively expansive 
endomorphisms. Example 20(2) generalizes the example and provides 
examples for which the constraint does not work. 

More examples are found in \cite{Jad-Nasu-Yaz} together with 
applications of results obtained in this 
paper to ``permutation cellular automata''. 

\section{Essentially-weakly-LR cones and extended districts}  

This section treats ``resolving'' endomorphisms of a 
subshift $(X,\sigma)$ from 
a much more general viewpoint 
than that of the overall dynamics 
of a single endomorphism of $(X,\sigma)$. 
We consider all ``resolving'' onto endomorphisms and 
all ``resolving'' automorphisms of $(X,\sigma)$. 
The notion of 
``resolving'' automorphism is 
used for classifying the 
expansive homeomorphisms of $X$. 

As will be shown by Theorem 12.2 
in the final section,  a ``directionally essentially weakly
\st'' endomorphism of a subshift is  
an ``essentially weakly \st'' endomorphism of the subshift, 
where \rt is any one of the resolving terms 
``$p$-L''  
``$p$-R'', ``$q$-R'', ``$q$-L'', 
``LL'', ``RR'', 
``LR'', ``RL'' , ``$p$-biresolving''  
and ``$q$-biresolving''. Therefore the definitions 
and the  statements of results 
appearing in the following are compatible with those in 
\cite[pp.193--208]{Nasu-te}.  

Let $X$ be a zero-dimensional compact metric space. 
Let $S(X)$ denote the monoid of  surjective continuous maps 
of $X$ (onto itself) and  
$H(X)$ denote the group of homeomorphisms of $X$. 
Let $J$ be a submonoid of $S(X)$. 
Suppose that $\tau$ is an expansive element in $J\cap H(X)$  
with $\tau^{-1}\in J$. 
Let $PL_J(\tau)$ denote the set of all essentially weakly
$p$-L endomorphisms of $(X,\tau)$ in $J$. 
Replacing ``$p$-L'' in this by ``$p$-R'', 
``$q$-R'' and ``$q$-L'', respectively, 
we define $PR_J(\tau)$, $QR_J(\tau)$ 
and $QL_J(\tau)$. Clearly 
$PR_J(\tau)=PL_J(\tau^{-1})$ and 
$QL_J(\tau)=QR_J(\tau^{-1})$.
Let $PL_J^\circ(\tau)$, $PR_J^\circ(\tau)$, 
$QR_J^\circ(\tau)$,  
and $QL_J^\circ(\tau)$ be
the set of all right $\tau$-expansive elements of 
$PL_J(\tau)$, the set of all left $\tau$-expansive 
elements of $PR_J(\tau)$, 
the set of all left $\tau$-expansive elements of 
$QR_J(\tau)$ and the set of all right $\tau$-expansive 
elements of $QL_J(\tau)$, respectively. 
By Theorem 8.10, $PL_J^\circ(\tau)$ is the set of 
all right $\tau$-expansive elements on the upper side in $J$ 
and $PR_J^\circ(\tau)$ is the set of all left 
$\tau$-expansive elements on the upper side in $J$. 
By Theorem 8.9, 
$QR_J^\circ(\tau)$ is the set of 
all positively left $\tau$-expansive elements in $J$ 
and $QL_J^\circ(\tau)$ is the set of all positively right 
$\tau$-expansive elements in $J$. 

Let $C_J(\tau)$ denote 
the set of all  essentially 
weakly LR 
endomorphisms of $(X,\tau)$ in $J$. 
Then $C_J(\tau^{-1})$ is the set of all  essentially 
weakly RL endomorphisms of $(X,\tau)$ in $J$. 
Let $C_J^\circ(\tau)$ denote the set of all expansive 
elements of $C_J(\tau)$. By \cite[Proposition 8.1]{Nasu-te}, 
$C_J^\circ(\tau)=PL_J^\circ(\tau)\cap QR_J^\circ(\tau)$. 
(NB \cite[Propositions 8.1, 8.2 and 8.3]{Nasu-te} are given
with the condition that $J$ in them 
is a commutative submonoid of $S(X)$, but they are valid 
even if $J$ is any subset of $S(X)$.)
Define 
\[D_J(\tau)=PL_J(\tau)\cup QR_J(\tau)\q \text{and}\q 
D_J^\circ(\tau)=PL_J^\circ(\tau)\cup QR_J^\circ(\tau).\] 
We call $D_J(\tau)$ the 
\textit{district of $\tau$ in $J$}, whereas we call $C_J(\tau)$ the 
\itl{essentially-weakly-LR cone}, or  \itl{cone}, 
\itl{of $\tau$ in $J$}. 
Let $Q_J(\tau)$ denote the set of all essentially weakly 
$q$-biresolving endomorphisms of $(X,\tau)$ in $J$. We call 
$Q_J(\tau)$ the 
\itl{essentially-weakly-$q$-biresolving cone of $\tau$ in $J$}. 
Let $Q_J^\circ(\tau)$ denote the set of all expansive 
elements of $Q_J(\tau)$. 
Then $QR_J^\circ(\tau)\cap QL_J^\circ(\tau)=Q_J^\circ(\tau)$ 
and this is the set of all 
positively expansive endomorphisms of $(X,\tau)$
in $J$ \cite[Proposition 8.3(2)]{Nasu-te} (see Theorem 4.5).
We call $\alpha_J^\circ(\tau)$ the \itl{interior} 
of $\alpha_J(\tau)$
when ``$\alpha$'' represents one of the symbols 
``$PL$'', ``$QR$'', ``$PR$'', ``$QL$'', ``$C$'',  ``$D$'', ``$Q$''. 

If $J$ is commutative, then 
for any $\tau'\in H(X)\cap C^\circ_J(\tau)$, it holds that 
$C_J(\tau)=C_J(\tau')$, $PL_J(\tau)=PL_J(\tau')$ and 
$QR_J(\tau)=QR_J(\tau')$ (by \cite[Propositions 9.2]{Nasu-te}) and 
hence $D_J(\tau)=D_J(\tau')$. 

When $(X,\tau)$ is 
conjugate to a topological Markov shift, then 
we say that $\varphi\in S(X)$ is an \itl{essentially LR endomorphism of 
$(X,\tau)$} if there exists a onesided 1-1, nondegenerate LR 
textile system $T$ such that 
$(X_T,\sigma_T,\varphi_T)$ and
$(X,\tau,\varphi)$ are topologically conjugate commuting systems;  
the set of all 
essentially LR endomorphisms of $(X,\tau)$ will be called 
the \itl{essentially-LR cone}, or 
\itl{ELR} cone, \itl{of $\tau$} and denoted by 
$(C_0)_J(\tau)$; the set of all expansive elements in 
$(C_0)_J(\tau)$ will be  
called the \itl{interior} of 
the ELR cone and  
denoted by $(C_0)_J^{\circ}(\tau)$. 

Similarly, when $(X,\tau)$ is 
conjugate to a topological Markov shift, then 
we say that $\varphi\in S(X)$ is an 
\itl{essentially $p$-L} 
(respectively, \itl{essentially $p$-R, essentially $q$-R, 
essentially $q$-L}) \itl{endomorphism  of 
$(X,\tau)$} if there exists a onesided 1-1, nondegenerate 
$p$-L (respectively, $p$-R, $q$-R, $q$-L)
textile system $T$ such that 
$(X_T,\sigma_T,\varphi_T)$ and
$(X,\tau,\varphi)$ are topologically conjugate commuting systems;  
the set of all essentially $p$-L (respectively, 
essentially $p$-L, essentially $q$-R, essentially $q$-L) 
endomorphisms of $(X,\tau)$ will be denoted by 
$(PL_0)_J(\tau)$ (respectively, $(PR_0)_J(\tau)$, 
$(QR_0)_J(\tau)$, $(QL_0)_J(\tau)$); 
the set of all right $\tau$-expansive and essentially $p$-L 
(respectively, 
left $\tau$-expansive and essentially $p$-L, 
left $\tau$-expansive and essentially $q$-R, 
right $\tau$-expansive and essentially $q$-L) 
endomorphisms of $(X,\tau)$ will be denoted by 
$(PL_0)_J^\circ(\tau)$ (respectively, $(PR_0)_J^\circ(\tau)$, 
$(QR_0)_J^\circ(\tau)$, $(QL_0)_J^\circ(\tau)$).

We emphasize that 
for any $\varphi$ in the ELR cone $(C_0)_J(\tau)$ of $tau$, 
the dynamics of $\varphi^i\tau^j, i,j\geq 0$ is lucid and regular. 
The statement (1) of the following theorem is a 
version of \cite[Theorem 6.3]{Nasu-t}; 
a weaker version of (2) was given as \cite[Theorem 6.11]{Nasu-m}. 
Note that if $\psi\in H(X)$ then the inverse limit system of 
$\psi$ is conjugate to the dynamical system $(X,\psi)$. 

\begin{theorem} 
Let $X$ be a zero-dimensional compact metric 
space. 
Let $J=S(X)$ and suppose that $\tau\in H(X)$ with 
$(X,\tau)$ conjugate to a topological Markov shift. 
\begin{enumerate}
\item 
If $\varphi\in (C_0)_J(\tau)$, then
there exist nonnegative integral matrices $M$ 
and $N$ such that $MN=NM$ 
and for all $i\geq 0, j\geq 1$ 
the inverse limit system of $\varphi^i\sigma^j$ 
is conjugate to the topological Markov shift 
whose defining matrix is $N^iM^j$ 
and such that if $\varphi$ is expansive 
then for all $i\geq 1$ the inverse limit system 
of $\varphi^i$ is conjugate to 
the topological Markov shift whose defining matrix is $N^i$. 
\item
If $\varphi_1,\dots,\varphi_k$ are in $(C_0)^\circ_J(\tau)$ and 
$(X,\tau)$ is conjugate to 
a topological Markov shift whose defining matrix is $M$, 
then there exist $m\geq 0$ and 
nonnegative integral matrices $N_1,\dots,N_k$
with $MN_i=N_iM$ for $i=1,\dots,k$ 
such that 
for all $i_1,\dots,i_k\geq 0$ and all $j\geq 1$, 
the inverse limit system of 
$(\varphi_1^{m})^{i_1}\dots(\varphi_k^{m})^{i_k}\tau^j$
is conjugate to the topological Markov shift whose defining matrix is 
$N_1^{i_1}\dots N_k^{i_k}M^j$, and such that 
if in addition, $\varphi_1,\dots,\varphi_k$ 
are pairwise commuting, then 
$N_1,\dots,N_k$ are pairwise commuting 
and the inverse limit system of
$(\varphi_1^{m})^{i_1}\dots(\varphi_k^{m})^{i_k}$
is conjugate to the topological Markov shift whose defining matrix is 
$N_1^{i_1}\dots N_k^{i_k}$ for all $i_1,\dots,i_k\geq 0$ 
with $(i_1,\dots,i_k)\neq (0,\dots,0)$. 
\end{enumerate}
\end{theorem}
\begin{proof} (1) Since $\varphi\in(C_0)_J(\tau)$, there exists 
a onesided 1-1, nondegenerate LR textile system 
$T$ such that 
$(X_T,\sigma_T,\varphi_T)$ and $(X,\tau,\varphi)$ are conjugate.  
Let $G$ and $G^*$ be the graphs such that $T$ and $T^*$ are defined 
over $G$ and $G^*$, respectively. Let $M=M_G$ and $N=M_{G^*}$. 
Then (1) follows from \cite[Theorems 6.3 and 2.5]{Nasu-t}. 

(2) Since each $\varphi_i$ is an expansive, essentially LR endomorphism 
of the topological Markov shift $(X_G,\sigma_G)$, where $G$ is 
the graph with $M_G=M$, it follows from 
Proposition 2.13(3) that there exists $m\geq 1$ such that 
$\varphi_i^m$ is an expansive, LR endomorphism of $(X_G,\sigma_G)$ 
for all $i=1,\dots,k$. There exist onesided 1-1, 
nondegenerate LR textile systems 
$T_i, i=1,\dots, k$, over $G$ 
such that $\varphi_{T_i}=\varphi_i^m$. 
Since $T_i$ is LR, if the dual $T_i^*$ is defined over a graph 
$G_i^*$ with $M_{G_i^*}=N_i$, then by \cite[Proposition 6.1]{Nasu-t} 
$MN_i=N_iM$ for $i=1,\dots,k$. 
For any $i_1,\dots,i_k\geq 0$, let $T_{i_1,\dots,i_k}$
be the composition 
$T_1^{i_1}\circ\dots\circ T_k^{i_k}$ 
of textile systems 
(see the final remarks in Subsection 2.3). Then  $T_{i_1,\dots,i_k}$ is 
a onesided 1-1, nondegenerate LR textile system such that 
$T_{i_1,\dots,i_k}^*$ is defined over the graph whose adjacency matrix 
is $N_1^{i_1}\dots N_k^{i_k}$ and  
$\varphi_{T_{i_1,\dots,i_k}}=(\varphi_k^m)^{i_k}\dots(\varphi_1^m)^{i_1}$ 
(by \cite[Corollaries 3.17, 3.18]{Nasu-t} and Remark 2.16(1)). 
Therefore, the first part of (2) follows from 
\cite[Theorem 6.3(3)]{Nasu-t}. 

Suppose in addition, that $\varphi_1,\dots,\varphi_k$ 
are pairwise commuting. Then it follows from Remark 2.16(2) that 
$N_1,\dots,N_k$ are pairwise commuting. 
Since $\varphi_1,\dots,\varphi_k$ 
are essentially LR and expansive and $i_1,\dots,i_k\geq 0$ 
with $(i_1,\dots,i_k)\neq (0,\dots,0)$, 
it follows from \cite[Proposition 7.12(3)]{Nasu-te} that 
$(\varphi_k^m)^{i_k}\dots(\varphi_1^m)^{i_1}$ is expansive 
and hence $T_{i_1,\dots,i_k}^*$ is 1-1 
(by \cite[Theorem 2.5]{Nasu-t}). 
Therefore, by \cite[Theorem 6.3(2)]{Nasu-t} the second part of (2) 
is proved. 
\end{proof} 

Let $\varphi:X\to X$ be a continuous map 
of a compact metric 
space. A bisequence $(x_i)_{i\in\Z}$ of points $x_i\in X$ is 
called a \itl{$\delta$-pseudo orbit} of $\varphi$ 
if $d_X(\varphi(x_i),x_{i+1})\leq\delta$ for all $i\in\Z$. 
For $\epsilon>0$, 
an orbit $(x_i)_{i\in\Z}$ of $\varphi$ is said to 
\itl{$\epsilon$-trace} 
a $\delta$-pseudo orbit $(x'_i)_{i\in\Z}$ of $\varphi$ 
if $d_X(x_i,x'_i) <\epsilon$ for all $i\in\Z$. 
We say that $\varphi$ has 
\itl{ the pseudo orbit tracing property}, 
or \itl{POTP} for short, if for any $\epsilon >0$, 
there exists $\delta >0$
such that any $\delta$-pseudo orbit of $\varphi$ 
is $\epsilon$-traced by some orbit of $\varphi$. 
We recall that for $\varphi\in H(X)$, where
$X$ is a zero-dimensional, compact metric space, $(X,\varphi)$ 
is conjugate to a topological Markov shift 
if and only if $\varphi$ is expansive 
and has POTP (Walters \cite{Wal-1}). 

\begin{proposition} 
Let $X$ be a zero-dimensional compact metric space and  
$J$ a submonoid of $S(X)$. Suppose that  
$\tau \in J\cap H(X)$ with $\tau^{-1}\in J$ and 
$(X,\tau)$ is conjugate to an SFT.  
\begin{enumerate} 
\item 
If $J=K$ with a subgroup $K$ of $H(X)$ or  
$\tau$ is topologically transitive, then $C_J(\tau)=(C_0)_J(\tau)$,  
$PL_J(\tau)=(PL_0)_J(\tau)$, $PR_J(\tau)=(PR_0)_J(\tau)$, 
$QR_J(\tau)=(QR_0)_J(\tau)$, $QL_J(\tau)=(QL_0)_J(\tau)$.

\item  All elements 
in $(C_0)_J(\tau)$ have POTP and 
the inverse limit system of every element of $(C_0)_J^\circ(\tau)$ 
is conjugate to an SFT; if $\tau$ is topologically mixing 
(or alternatively, $\tau$ commutes with some 
topologically mixing element in $S(X)$), 
then the inverse limit system of every element of 
$(C_0)_J^\circ(\tau)$ is conjugate to a mixing SFT.   
\item 
If $\tau$ is topologically transitive, then 
the inverse limit system of every expansive, 
right closing or left closing element of $D_J(\tau)$ 
is conjugate to an SFT; 
if in addition 
$\tau$ is topologically mixing 
(or alternatively, $\tau$ commutes with some 
topologically mixing element in $S(X)$), then 
the inverse limit system of every 
expansive, right closing or left closing element 
in $D_J(\tau)$ is conjugate to a mixing SFT. 
\end{enumerate} 
\end{proposition} 
\begin{proof}
(1) By Remark 2.10(2). 

(2) By (1) we see that $C_J(\tau)=(C_0)_J(\tau)$. 
(NB \cite[Proof of Proposition 9.9(2)]{Nasu-te}
has been revised to be much simpler.)
Hence, since any essentially LR endomorphism of an SFT
has POTP (\cite[Proposition 9.9(1)]{Nasu-te}), (2) is proved 
by using \cite[Theorem 1.6]{Nasu-e} (For let $\varphi$ be an onto 
endomorphism of $(X,\tau)$. If $\varphi$ has POTP then so does 
$\sigma_\varphi$, and if $\varphi$ is topologically mixing then so is
$\sigma_\varphi$ (see e.g. \cite{AokHir}). )

(3) It follows from \cite[Corollary 6.6]{Nasu-te} that 
if $(X,\tau)$ is conjugate to a transitive SFT, 
then the inverse limit system of every 
expansive, right closing or left closing element in 
$(PL_0)(\tau)\cup(QR_0)(\tau)$, is conjugate to a 
topological Markov shift. Since 
$D_J(\tau)=(PL_0)(\tau)\cup(QR_0)(\tau)$ by (1),
(3) follows by \cite[Theorem 1.6]{Nasu-e}.
\end{proof}

Let $X$ be a zero-dimensional compact metric space. 
Let $\varphi,\tau\in H(X)$. 
We write 
$\varphi\pl\tau$ if 
$\tau$ is expansive and $\varphi$ is 
an essentially weakly $p$-L automorphism of $(X,\tau)$, and
if in addition, $\varphi$ is right $\tau$-expansive, then 
we write $\varphi\pl^\circ\tau$, which is 
equivalent to the condition that 
$\tau$ is expansive and $\varphi$ 
is right $\tau$-expansive on 
the upper side (by Theorem 8.10(1)); 
we write  
$\varphi\qr\tau$ if 
$\tau$ is expansive and $\varphi$ is 
an essentially weakly $q$-R automorphism of $(X,\tau)$, and
if in addition, $\varphi$ is left $\tau$-expansive, then 
we write $\varphi\qr^\circ\tau$, which is equivalent to 
the condition that $\tau$ is expansive and 
$\varphi$ is positively left $\tau$-expansive 
(by Theorem 8.9(2)). 
We write $\varphi\lr\tau$ to mean that 
$\tau$ is expansive and $\varphi$ is an essentially weakly LR 
automorphism of $(X,\tau)$, and 
$\varphi\lr^\circ\tau$ to mean that
$\tau$ is expansive and $\varphi$ is an 
expansive, essentially weakly LR
automorphism of $(X,\tau)$, which is equivalent to
the condition that 
$\varphi\pl^\circ\tau$ and $\varphi\qr^\circ\tau$ 
(by Theorem 8.12(3)). 

Let $K$ be any subgroup of $H(X)$. (We do not 
assume that $K$ is finitely generated.)
Let $E(K)$ denote the 
set of all expansive elements in $K$. 
Suppose that 
$\tau\in E(K)$. Then 
$C_K(\tau)=\{\varphi\in K\,|\,\varphi\lr\tau\}$ and 
$C_K^\circ(\tau)=
\{\varphi\in K\,|\,\varphi\lr^\circ\tau\}=
\{\varphi\in E(K)\,|\,\varphi\lr\tau\}$, 
which is the interior of $C_K(\tau)$. 
Further, 
$PL_K(\tau)=\{\varphi\in K\,|\,\varphi\pl\tau\}$, 
$QR_K(\tau)=\{\varphi\in K\,|\,\varphi\qr\tau\}$, 
$C_K(\tau)=\{\varphi\in K\,|\,\varphi\lr\tau\}$ and 
$PL_K^\circ(\tau)=\{\varphi\in K\,|\,\varphi\pl^\circ\tau\}$, 
$QR_K^\circ(\tau)=\{\varphi\in K\,|\,\varphi\qr^\circ\tau\}$. 

For any subset $F$ of $H(X)$, let $F^{-1}$ denote
$\{\varphi^{-1}\,|\,\varphi\in F\}$. Then, as is easily seen, 
$PR_K(\tau)=PL_K(\tau^{-1})=(QR_K(\tau))^{-1}$, 
$QL_K(\tau)=QR_K(\tau^{-1})=(PL_K(\tau))^{-1}$, 
$C_K(\tau^{-1})=(C_K(\tau))^{-1}$ and 
$D_K(\tau^{-1})=(D_K(\tau))^{-1}$. 
It follows from \cite[Proposition 8.2]{Nasu-te} 
that if $X$ is infinite, then 
$PL_K(\tau)\cap PL^\circ_K(\tau^{-1})=\emptyset$, 
$QR_K(\tau)\cap QR^\circ_K(\tau^{-1})=\emptyset$ and 
$C_K(\tau)\cap C^\circ_K(\tau^{-1})=\emptyset$. 
More properties of the relations $\pl$, $\qr$ and 
$\lr$ on $H(X)$ were presented 
in \cite[Section 9]{Nasu-te} (for summaries of them, 
see \cite[Propositions 9.4 and 9.5]{Nasu-te}). 

\begin{proposition} Let $X$ be 
a zero-dimensional compact metric space. Let 
$K$ be a commutative subgroup of $H(X)$. Then 
the relation $\lr^\circ$ on $E(K)$
is an equivalence relation. 
For $\tau\in E(K)$, $C_K^\circ(\tau)$ is 
the equivalence class containing $\tau$ 
with respect to $\lr^\circ$ and $C_K(\tau)=C_K(\tau')$ 
for every $\tau'\in C^\circ_K(\tau)$. 
\end{proposition} 
\begin{proof} 
The proof is given by the explanation in \cite[p. 206]{Nasu-te}
(which uses \cite[Proposition 9.5(1)(2)(3)]{Nasu-te})
and Theorem 12.2.
\end{proof}

Let $K$ be a commutative subgroup of $H(X)$ and $\tau\in E(K)$, 
where $X$ is a zero-dimensional compact metric space. For the 
reason shown in Proposition 10.3, 
we call $C_K(\tau)$ 
\itl{an essentially-weakly-LR cone} or \itl{a cone} 
\itl{in $K$}. For a cone $C$ in $K$, $C^\circ=C\cap E(K)$ is 
the interior of $C$. If a cone $C$ in $K$
contains an element having POTP in  $C^\circ$, 
then $C=(C_0)_K(\tau)$ for all elements 
$\tau\in C^\circ$ (by Propositions 10.3 and 10.2) 
and all elements 
of the cone $C$ have POTP (by Proposition 10.2), and hence 
we call $C$ \itl{an ELR} (\itl{essentially LR}) \itl{cone in $K$}. 

Let $X$ be a zero-dimensional compact metric space. 
Let $\varphi,\tau\in H(X)$.  
Let us write $\varphi\mm\tau$ to mean 
that $\varphi\pl\tau$ or $\varphi\qr\tau$. 
Let $K$ be any (not necessarily commutative) subgroup of $H(X)$. 
Let us write $\varphi\mm^\ast_K\tau$ to mean 
that there exist $r\geq 1$ and (expansive) elements 
$\tau_i, i=1,\dots, r$, in $K$ such that 
$\varphi=\tau_0\mm\tau_1, \tau_1\mm\tau_2, \dots, 
\tau_{r-1}\mm\tau_r=\tau$. 
Let us write $\varphi\mm^\circ\tau$ to mean that 
$\varphi\pl^\circ\tau$ or $\varphi\qr^\circ\tau$, and let us write  
$\varphi\mm_K^{\circ,\ast}\tau$ to mean that 
there exist $r\geq 1$ and (expansive) elements 
$\tau_i, i=1,\dots, r$, in $K$ such that 
$\varphi=\tau_0\mm^\circ\tau_1, \tau_1\mm^\circ\tau_2, \dots, 
\tau_{r-1}\mm^\circ\tau_r=\tau$
(or equivalently, to mean that
there exists $\tau'\in E(K)$ such that 
$\varphi\mm^\circ\tau'$ and $\tau'\mm^\ast_K\tau$). 

Let $\tau\in E(K)$. Then  
$D_K(\tau)=\{\varphi\in K\,|\,\varphi\mm\tau\}$ and 
$D_K^\circ(\tau)=\{\varphi\in K\,|\,\varphi\mm^\circ\tau\}$. 
Define 
\[
D_K^\ast(\tau)
=\{\varphi\in K\,|\,\varphi\mm^\ast_K\tau\}\q 
\text{and}\q 
(D_K^\ast)^\circ(\tau)=
\{\varphi\in K\,|\, \varphi\mm_K^{\circ,\ast}\tau\}.
\]
We call $D_K^\ast(\tau)$  
the \itl{extended district of $\tau$ in $K$},
and also call \itl{an extended district in $K$} 
for the reason shown in Proposition 10.4 below. 
We call the set $D_K^\ast(\tau)\cap E(K)= 
\{\varphi\in E(K)\,|\,\varphi\mm^\ast_K\tau\}$
of all expansive elements in $D_K^\ast(\tau)$ 
the \itl{core} of $D_K^\ast(\tau)$ 
in $K$. We call $(D_K^\ast)^\circ(\tau)$
the \itl{interior} of $D_K^\ast(\tau)$ in $K$. 
Clearly we have  \[D_K^\ast(\tau)=
\bigcup_{\tau'\in D_K^\ast(\tau)\cap E(K)} D_K(\tau'),\q\q 
(D_K^\ast)^\circ(\tau)=
\bigcup_{\tau'\in D_K^\ast(\tau)\cap E(K)}D_K^\circ(\tau').\] 

\begin{proposition} Let $X$ be a zero-dimensional compact metric 
space, 
$K$ a (not necessarily commutative) subgroup of $H(X)$. Then 
the relation $\mm^\ast_K$ 
($\mm^{\circ,\ast}_K$)   
restricted on $E(K)$
is an equivalence relation on $E(K)$. For $\tau\in E(K)$, 
the core of $D_K^\ast(\tau)$ (i.e. $D_K^\ast(\tau)\cap E(K)$)
is the equivalence class containing $\tau$ with respect 
to $\mm^\ast_K$ on $E(K)$. 
For every $\tau'$ in the core of $D_K^\ast(\tau)$, 
$D^\ast_K(\tau)=D^\ast_K(\tau')$ and  
$(D^\ast_K)^\circ(\tau)=(D^\ast_K)^\circ(\tau')$.
\end{proposition}
\begin{proof}
If $\tau\in H(X)$ is expansive, then $\tau\mm\tau$ 
(because $\tau\lr\tau$). If $\tau,\tau'\in H(X)$ are expansive, 
then $\tau\mm\tau'$ implies $\tau'\mm\tau$, 
because $\tau\pl\tau'$ if and only if $\tau'\qr\tau$ 
(\cite[Proposition 9.4(1)]{Nasu-te}). By definition 
$\mm^\ast_K$ is transitive on $E(K)$. Therefore 
$\mm^\ast_K$ restricted on $E(K)$ is an equivalence relation on $E(K)$, 
and hence the remainder follow. 
\end{proof} 

The following theorem is 
an extension of the main part of \cite[Proposition 9.11]{Nasu-te}. 

\begin{theorem} 
Let $X$ be a zero-dimensional compact metric space and $K$ be 
$H(X)$  or any subgroup of $H(X)$. Let 
$D^*$ be any extended district in $K$. If
$D^*$ contains an element $\tau$ with $(X,\tau)$ conjugate to 
a mixing SFT, then for every expansive element 
$\varphi$ in $D^*$, the following hold: 
\begin{enumerate} 
\item
$(X,\varphi)$ is conjugate to a mixing SFT; 
\item
all left $\varphi$-expansive elements in $PL_K(\varphi)$ and all 
right $\varphi$-expansive elements in $QR_K(\varphi)$ 
have POTP and are topologically mixing;
\item 
all elements in $C_K(\varphi)$ have POTP.
\end{enumerate}
\end{theorem}
\begin{proof} By Proposition 10.4, $D^*=D_K^*(\tau)$. 
We know, by \cite[Corollary 6.6]{Nasu-te}, 
that for any $\psi_1,\psi_2\in H(X)$ with 
$\psi_1\pl\psi_2$ (respectively, $\psi_1\qr\psi_2$) and 
$\psi_1$ left $\psi_2$-expansive (respectively, 
right $\psi_2$-expansive), it holds that if $\psi_2$ has POTP and is 
topologically-mixing, then so does and so is $\psi_1$.  Using this, 
(1) and (2) are proved. 

By (1) and Proposition 10.2, (3) is proved.
\end{proof}

Let $K$ be a commutative subgroup of $H(X)$.  
Let $\varphi,\tau\in E(K)$.   
If $\varphi\pl\tau$ (respectively, $\varphi\qr\tau$), 
$\varphi'\in C_K^\circ(\varphi)$ and $\tau'\in C_K^\circ(\tau)$, 
then 
$\varphi'\pl\tau'$ (respectively, $\varphi'\qr\tau'$) 
(\cite[Proposition 9.7]{Nasu-te}).
Hence the relations ``$\pl$'', ``$\qr$'', ``$\mm$'', 
and ``$\mm_K^\ast$'' are 
defined on the family  of interiors $C^\circ$'s of cones $C$'s
in $K$ as the quotient 
relations of them with respect to 
the equivalence relation $\lr$ on $E(K)$, 
where $C^\circ$ denotes the interior of a cone $C$ in $K$. 
By Theorem 10.5 we see that if 
$C,C'$ are cones in $K$ with $C^\circ{\mm_K^\ast}C'^\circ$
and if one of $C, C'$ is an ELR cone 
containing a topologically mixing element in its interior,  
then the other is also an ELR cone with 
all elements in its interior 
topologically mixing. 

Here, in connection with Theorem 10.5, we note that 
the question of whether or not there exists an expansive automorphism 
$\varphi$ of a transitive SFT $(X,\sigma)$ with $(X,\varphi)$ not 
conjugate to any SFT is still open. See \cite[Section 12]{Boy-o}. 

If we restrict our attention to the overall dynamics of 
an onto endomorphism $\varphi$ of a subshift $(X,\sigma)$, then letting
$J$ be the submonoid of $S(X)$ generated by 
$\{\sigma,\varphi,\sigma^{-1}\}$ we know, 
by Theorems 9.7, 9.10 and 9.11, the following:  
\begin{align*}
PL_J^\circ(\sigma)&=\{\varphi^i\sigma^j\,\bigm|\, 
j>-p_L(\varphi)i, i\geq 0\},\q
PR_J^\circ(\sigma)=\{\varphi^i\sigma^j\,\bigm|\, 
j<p_R(\varphi)i, i\geq 0\},\\
QR_J^\circ(\sigma)&=\{\varphi^i\sigma^j\,\bigm|\, 
j>-q_R(\varphi)i, i\geq 0\},\q
QL_J^\circ(\sigma)=\{\varphi^i\sigma^j\,\bigm|\, 
j<q_L(\varphi)i, i\geq 0\},\\
C_J^\circ(\sigma)&=\{\varphi^i\sigma^j\,\bigm|\, 
j>c_R(\varphi)i, i\geq 0\},\q\:\;
C_J^\circ(\sigma^{-1})=\{\varphi^i\sigma^j\,\bigm|\, 
j<c_L(\varphi)i, i\geq 0\},\\
D_J^\circ(\sigma)&=\{\varphi^i\sigma^j\,\bigm|\, 
j>d_R(\varphi)i, i\geq 0\},\q\:
D_J^\circ(\sigma^{-1})=\{\varphi^i\sigma^j\,\bigm|\, 
j<d_L(\varphi)i, i\geq 0\},\\ 
Q_J^\circ(\sigma)&=
\{\varphi^i\sigma^j\bigm| q_L(\varphi)i <j< -q_R(\varphi)i, i> 0\}. 
\end{align*} 

\begin{proposition} 
Let $\varphi$ be an onto endomorphism of a 
subshift $(X,\sigma)$. Let $J$ be the submonoid of $S(X)$ 
generated by $\{\varphi,\sigma,\sigma^{-1}\}$. 
If $-p_L(\varphi)$ is of type I, then 
$PL_J(\sigma)=\{\varphi^i\sigma^j\,|\, 
j\geq-p_L(\varphi)i, i\geq 0\}$ and 
$PL_J(\sigma)\supsetneq PL_J^\circ(\sigma)$; 
if $-p_L(\varphi)$ is of type II or of type III, 
then $PL_J(\sigma)=PL_J^\circ(\sigma)$; similar facts 
for $QR_J(\sigma)$, $PR_J(\sigma)$ and $QL_J(\sigma)$ 
also hold. 
\end{proposition} 
\begin{proof} We note that $\sigma^j$ is an LR automorphism 
of $(X,\sigma)$ for all $j\geq 0$. 
Let $i,j\in\Z$ with $i>0$. 
If $j/i<-p_L(\varphi)$, then $\varphi^i\sigma^j$ 
is not essentially weakly $p$-L. For otherwise, there would 
be an integer $m>0$ with $(mj+1)/(mi)<-p_L(\varphi)$
but $(\varphi^i\sigma^j)^m\sigma$ would be 
essentially weakly $p$-L and right $\sigma$-expansive 
(by \cite[Remark 7.9, Proposition 7.6]{Nasu-te}), contrary 
to the fact that $j/i>-p_L(\varphi)$ if and only if 
$\varphi^i\sigma^j$ is essentially weakly $p$-L 
and right $\sigma$-expansive (by Theorem 9.7). 
If $j/i=-p_L(\varphi)$, then $\varphi^i\sigma^j$ is 
essentially weakly $p$-L if and only if $-p_L(\varphi)$ is of type I.  
Hence we have proved the statements for $PL_J(\sigma)$.

The remainder are proved similarly. 
\end{proof}

We cannot answer the question whether or not 
$C_J(\sigma)= PL_J(\sigma)\cap QR_J(\sigma)$ generally holds, 
because we cannot answer the question whether or not, 
for the case that $-p_L(\varphi)$ and $-q_(\varphi)$ 
are the same value and both of type I, 
$\varphi^i\sigma^{c_R(\varphi)i}$
with $i>0$ is essentially weakly LR endomorphism of $(X,\sigma)$. 
We cannot answer the question whether or not 
$Q_J(\sigma)= QL_J(\sigma)\cap QR_J(\sigma)$ generally holds,  
because we cannot answer the question whether or not, 
for the case that $q_L(\varphi)$ and $-q_R(\varphi)$ 
are the same value and both of type I, 
$\varphi^i\sigma^{q_L(\varphi)i}$
with $i>0$ is an essentially weakly 
$q$-biresolving endomorphism 
of $(X,\sigma)$.
The following proposition together with Proposition 10.6 
is useful for describing $C_J(\sigma)$ and $Q_J(\sigma)$
in terms of the limits of onesided resolving directions 
of $\varphi$. 

\begin{proposition} 
Let $\varphi$ be an onto endomorphism of a 
subshift $(X,\sigma)$. Let $J$ be the submonoid of $S(X)$ 
generated by $\{\varphi,\sigma,\sigma^{-1}\}$.
\begin{enumerate} 
\item
If either $-p_L(\varphi)<-q_R(\varphi)$ with $-q_R(\varphi)$ of type I
(respectively, $q_L(\varphi)<p_R(\varphi)$ with $q_L(\varphi)$ of  
type I) or $-p_L(\varphi)>-q_R(\varphi)$ with  
$-p_L(\varphi)$ of type I (respectively, 
$q_L(\varphi)>p_R(\varphi)$ with $p_R(\varphi)$ of type I), then  
$\varphi^i\sigma^j$ with $j=c_R(\varphi)i, i\geq 0$ is contained 
in $C_J(\sigma)$
(respectively, 
$\varphi^i\sigma^j$ with $j=c_L(\varphi)i, i\geq 0$ is contained 
in $C_J(\sigma^{-1})$). 
\item 
If $-p_L(\varphi)=-q_R(\varphi)$ 
with both $-p_L(\varphi)$ and $-q_R(\varphi)$ of type I
(respectively, $p_R(\varphi)=q_L(\varphi)$ with both 
$p_R(\varphi)$ and $q_L(\varphi)$ of type I), then 
with the additional condition that $(X,\sigma)$ is an SFT and that  
$\varphi$ is invertible or $\sigma$ is topologically transitive, 
$\varphi^i\sigma^j$ with $j=c_R(\varphi)i, i\geq 0$ is contained 
in $C_J(\sigma)$
(respectively, 
$\varphi^i\sigma^j$ with $j=c_L(\varphi)i, i\geq 0$ is contained 
in $C_J(\sigma^{-1})$). 
\item 
If $-q_R(\varphi)<q_L(\varphi)$ with $-q_R(\varphi)$ of type I
(respectively, with $q_L(\varphi)$ of  
type I), then  
$\varphi^i\sigma^j$ with $j=-q_R(\varphi)i, i\geq 0$ is contained 
in $Q_J(\sigma)$
(respectively, 
$\varphi^i\sigma^j$ with $j=q_L(\varphi)i, i\geq 0$ is contained 
in $Q_J(\sigma)$). 
\item 
If $-q_R(\varphi)=q_L(\varphi)$ 
with both $-q_R(\varphi)$ and $q_L(\varphi)$ of type I, then 
with the additional condition that $(X,\sigma)$ is a 
topologically transitive SFT, it holds that 
$\{\varphi^i\sigma^j\,|\, j=-q_R(\varphi)i, i\geq 0\}=Q_J(\sigma)$. 
\end{enumerate}
\end{proposition} 
\begin{proof} 
(1) Suppose that $j=c_R(\varphi)i, i\geq 0$.
Suppose that $-p_L(\varphi)<-q_R(\varphi)$ and 
$-q_R(\varphi)$ is type I. Then 
$j=-q_R(\varphi)i$ and $\varphi^i\sigma^j$ is 
essentially weakly $q$-R. 
There exists an onto 
endomorphism $\varphi_1$ of a subshift $(X_1,\sigma_1)$ 
such that $(X_1,\sigma_1,\varphi_1)$ 
is conjugate to $(X,\sigma,\varphi)$ and 
$\varphi_1^i\sigma_1^j$ is a weakly $q$-R endomorphism of 
$(X_1,\sigma_1)$. Since $i_X\in C_J(\sigma)$, 
we may assume that $i>0$. 
Since $j/i>-p_L(\varphi)=-p_L(\varphi_1)$ (by Theorem 9.4), 
it follows from Theorem 9.7 that $\varphi_1^i\sigma_1^j$ 
is essentially weakly $p$-L and right $\sigma_1$-expansive. 
Therefore, by \cite[Proposition 7.11]{Nasu-te} 
there exists $k>0$ such that 
$(\varphi_1^i\sigma_1^j)^k$ is weakly $p$-L. 
Since $\varphi_1^i\sigma_1^j$ is weakly $q$-R, 
so is $(\varphi_1^i\sigma_1^j)^k$ 
(by Remark 2.15). Therefore by 
Theorem 8.7(1)(b) $\varphi_1^i\sigma_1^j$ is 
an essentially weakly LR endomorphism of $(X_1,\sigma_1)$, 
and hence $\varphi^i\sigma^j\in C_J(\sigma)$. 
The proof for the case that 
$-p_L(\varphi)>-q_R(\varphi)$ and 
$-p_L(\varphi)$ is of type I is similar. 

(2) By \cite[Theorem 8.8(4)]{Nasu-te}. 

(3) The proof is similar to that of (1); use 
Theorem 8.7(1)(c) instead of Theorem 8.7(1)(b).

(4) By \cite[Theorem 8.8(5)]{Nasu-te}. 
\end{proof} 

Let $\varphi$ be an automorphism of $(X,\sigma)$. 
In order to discuss the overall dynamics of the 
automorphism $\varphi$ of $(X,\sigma)$, 
it is appropriate to let $K$ be the subgroup 
of $H(X)$ generated by $\{\sigma,\varphi\}$.  
For any subset $F$ of $H(X)$, 
let $F^{-1}$ denote the set $\{\psi^{-1}\,|\,\psi\in F\}$. 
Let $J$ be the submonoid of $S(X)$ generated by 
$\{\sigma, \sigma^{-1}, \varphi\}$. Then  
$PL_K(\sigma)=PL_J(\sigma)\cup PL_{J^{-1}}(\sigma)$. 
$PL_{J^{-1}}(\sigma)$ is the set of all 
$\varphi^{-i}\sigma^{-j}$ such that $\varphi^{-i}\sigma^{-j}\pl\sigma$ 
with $i\geq 0,j\in\Z$, this set equals the set
of all $\varphi^{-i}\sigma^{-j}$ 
such that $\varphi^i\sigma^j$ is essentially 
weakly $q$-L automorphism of $(X,\sigma)$ with $i\geq 0,j\in\Z$, and 
this set is $(QL_J(\sigma))^{-1}$. Therefore,  
$PL_K(\sigma)=PL_J(\sigma)\cup(QL_J(\sigma))^{-1}$. 
Similarly   
$QR_K(\sigma)=QR_J(\sigma)\cup (PR_{J}(\sigma))^{-1}$ and 
$C_K(\sigma)=C_J(\sigma)\cup (C_J(\sigma^{-1}))^{-1}$. We also have 
$PL_K^\circ(\sigma)=PL_J^\circ(\sigma)\cup (QL_J^\circ(\sigma))^{-1}$, 
$QR_K^\circ(\sigma)=QR_J^\circ(\sigma)\cup (PR_J^\circ(\sigma))^{-1}$ 
and 
$C_K^\circ(\sigma)=C_J^\circ(\sigma)\cup(C_J^\circ(\sigma^{-1}))^{-1}$.
Hence, by the equations before Proposition 10.6, we see: 
\begin{align*}
PL_K^\circ(\sigma)&=\{\varphi^i\sigma^j \,\bigm|\, 
j>q_L(\varphi)i\;\,\text{if}\;\,i\leq 0,\q 
j>-p_L(\varphi)i\;\,\text{if}\;\,i\geq 0\},\\
QR_K^\circ(\sigma)&=\{\varphi^i\sigma^j \,\bigm|\, 
j>p_R(\varphi)i\;\,\text{if}\;\,i\leq 0,\q 
j>-q_R(\varphi)i\;\,\text{if}\;\,i\geq 0\},\\
C_K^\circ(\sigma)&=\{\varphi^i\sigma^j \,\bigm|\, 
j>c_L(\varphi)i\;\,\text{if}\;\,i\leq 0,\q 
j>c_R(\varphi)i\;\,\text{if}\;\,i\geq 0\},\\ 
D_K^\circ(\sigma)&=\{\varphi^i\sigma^j \,\bigm|\, 
j>d_L(\varphi)i\;\,\text{if}\;\,i\leq 0,\q 
j>d_R(\varphi)i\;\,\text{if}\;\,i\geq 0\}. 
\end{align*} 
By using $PL_K(\sigma)=PL_J(\sigma)\cup (QL_J(\sigma))^{-1}$, 
$QR_K(\sigma)=QR_J(\sigma)\cup (PR_J(\sigma))^{-1}$ 
and Proposition 10.6, 
$PL_K(\sigma)$ and $QR_K(\sigma)$ can be described in terms 
of the limits of onesided resolving directions of $\varphi$, 
and so can be done 
$C_K(\sigma)$ in many cases by using
$C_K(\sigma)=C_J(\sigma)\cup (C_J(\sigma^{-1}))^{-1}$ and 
Proposition 10.7. 

For an automorphism $\varphi$
of a subshift $(X,\sigma)$ and 
for $(k,l)\in\Z^2$ with 
$\varphi^k\sigma^l$ expansive, 
Proposition 10.8(2) below 
is often useful to determine 
$PL_K(\varphi^k\sigma^l)$, $QR_K(\varphi^k\sigma^l)$ and 
$C_K(\varphi^k\sigma^l)$,  
and its continuation Corollary 11.14(2) in the next 
section is often useful 
to determine $D_K^\ast(\varphi^k\sigma^l)$. 

For $(a_1,b_1),(a_2,b_2)\in\R^2$ with $a_1b_2\neq a_1b_2$, 
let us define  
$\cH[(a_1,b_1),(a_2,b_2)]$, $\cH[(a_1,b_1),(a_2,b_2))$, 
$\cH((a_1,b_1),(a_2,b_2)]$, and $\cH((a_1,b_1),(a_2,b_2))$ 
to be the convex cones in $\R^2$ 
(with apex $(0,0)$) consisting 
of all points $(ca_1+da_2,cb_1+db_2)\in\R^2$
with $c,d\geq 0$, with $c>0,d\geq 0$, with $c\geq 0,d>0$, and 
with $c>0,d>0$, respectively. 

\begin{proposition} Let $\varphi$ be an automorphism of 
an infinite subshift $(X,\sigma)$. Let $K$ be the subgroup of 
$H(X)$ generated by $\{\sigma,\varphi\}$. Let $(k,l)\in\Z^2$ 
with $\varphi^k\sigma^l$ expansive. 
\begin{enumerate} 
\item If 
$\varphi^i\sigma^j\pl\varphi^k\sigma^l$ 
(respectively, $\varphi^i\sigma^j\qr\varphi^k\sigma^l$, 
$\varphi^i\sigma^j\lr\varphi^k\sigma^l$,)
with $il\neq jk$, then for all   
$(i',j')\in\cH[(i,j),(k,l)]\cap\Z^2$, 
$\varphi^{i'}\sigma^{j'}\in PL_K(\varphi^k\sigma^l)$ 
(respectively, $\varphi^{i'}\sigma^{j'}\in QR_K(\varphi^k\sigma^l)$, 
$\varphi^{i'}\sigma^{j'}\in C_K(\varphi^k\sigma^l)$), 
and for all $(i',j')\in\cH((i,j),(k,l)]\cap\Z^2$, 
$\varphi^{i'}\sigma^{j'}\in PL_K^\circ(\varphi^k\sigma^l)$
(respectively, $\varphi^{i'}\sigma^{j'}\in QR_K^\circ(\varphi^k\sigma^l)$, 
$\varphi^{i'}\sigma^{j'}\in C_K^\circ(\varphi^k\sigma^l)$). 
\item 
Let $(i_1,j_1), (i_2,j_2)\in\Z^2$ with $i_1j_2\neq i_2j_1$. 
If $\varphi^{i_t}\sigma^{j_t}\pl\varphi^k\sigma^l$ 
(respectively, $\varphi^{i_t}\sigma^{j_t}\qr\varphi^k\sigma^l$
$\varphi^{i_t}\sigma^{j_t}\lr\varphi^k\sigma^l$)
with 
$\varphi^{i_t}\sigma^{j_t}$ not right $\varphi^k\sigma^l$-expansive
(respectively, not left $\varphi^k\sigma^l$-expansive, 
nonexpansive) for $t=1,2$,  
then  
$(i,j)\in\cH[(i_1,j_1),(i_2,j_2)]\cap\Z^2$ if and only if 
$\varphi^i\sigma^j\in PL_K(\varphi^k\sigma^l)$ 
(respectively, $\varphi^i\sigma^j\in QR_K(\varphi^k\sigma^l)$, 
$\varphi^i\sigma^j\in C_K(\varphi^k\sigma^l)$), and 
$(i,j)\in\cH((i_1,j_1),(i_2,j_2))\cap\Z^2$ if and only if  
$\varphi^i\sigma^j\in PL_K^\circ(\varphi^k\sigma^l)$ 
(respectively, $\varphi^i\sigma^j\in QR_K^\circ(\varphi^k\sigma^l)$, 
$\varphi^i\sigma^j\in C_K^\circ(\varphi^k\sigma^l)$). 
\end{enumerate}
\end{proposition} 
\begin{proof}  We describe proofs 
only for the first versions of 
(1) and (2). The other versions are 
proved by modifying them straightforwardly 
(by using Theorems 8.6 and 8.7 
instead of Theorem 8.1). 

(1) Suppose that $(i',j')\in\cH[(i,j),(k,l)]$.  Then 
there exist integers $m>0$, $r\geq 0, s\geq 0$ such that 
$(mi',mj')=(ri+sk,rj+sl)$. 
Noting that $i_X\in C_K(\varphi^k\sigma^l)$,  
we know, by Remark 2.15, that 
$(\varphi^i\sigma^j)^r\in PL_K(\varphi^k\sigma^l)$
and hence it follows from \cite[Proposition 7.6]{Nasu-te} that 
$(\varphi^{i'}\sigma^{j'})^m=
(\varphi^i\sigma^j)^r(\varphi^k\sigma^l)^s$ is in 
$PL_K(\varphi^k\sigma^l)$. 
Hence by Theorem 8.1, 
$\varphi^{i'}\sigma^{j'}\in PL_K(\varphi^k\sigma^l)$. 

If $(i',j')\in\cH((i,j),(k,l)]$, 
then a similar argument to the above (replaced by ``$s\geq 0$'' by 
``$s>0$'') shows that 
$\varphi^{i'}\sigma^{j'}\in PL_K^\circ(\varphi^k\sigma^l)$. 

(2) Suppose that $\varphi^{i_t}\sigma^{j_t}\pl\varphi^k\sigma^l$ 
with $\varphi^{i_t}\sigma^{j_t}$ 
not right $\varphi^k\sigma^l$-expansive for $t=1,2$. 
Since $\varphi^k\sigma^l$ is expansive and 
$\varphi^{i_t}\sigma^{j_t}$ is nonexpansive 
for $t=1,2$, it follows that $i_tl\neq j_tk$ for $t=1,2$. 
It follows from (1) that if 
$(i,j)\in\cH[(i_1,j_1),(k,l)]\cup\cH[(i_2,j_2),(k,l)]$ then 
$\varphi^i\sigma^j \in PL_K(\varphi^k\sigma^l)$.
Since $i_1j_2\neq i_2j_1$,  
neither one of $\cH[(i_1,j_1),(k,l)]$ and $\cH[(i_2,j_2),(k,l)]$ 
can include the other, because if 
$\cH[(i_1,j_1),(k,l)]\supsetneq\cH[(i_2,j_2),(k,l)]$ then by (1)
$\varphi^{i_2}\sigma^{j_2}\in PL_K^\circ(\varphi^k\sigma^l)$, 
which contradicts the hypothesis that $\varphi^{i_2}\sigma^{j_2}$  
is not right $\varphi^k\sigma^l$-expansive. Therefore,  
$\cH[(i_1,j_1),(i_2,j_2)]
=\cH[(i_1,j_1),(k,l)]\cup\cH[(i_2,j_2),(k,l)]$,  
and if $(i,j)\in \cH[(i_1,j_1),(i_2,j_2)]$ 
then $\varphi^i\sigma^j\in PL_K(\varphi^k\sigma^l)$. 

Assume that there were $(i',j')\in\Z^2$ 
such that $\varphi^{i'}\sigma^{j'}\in PL_K(\varphi^k\sigma^l)$ 
but $(i',j')\notin \cH[(i_1,j_1),(i_2,j_2)]$. Then, since 
$(i',j')\notin \cH[(i_1,j_1),(i_2,j_2)]\ni (k,l)$ and 
$\varphi^{i'}\sigma^{j'}\notin PL_K^\circ(\varphi^{-k}\sigma^{-l})$ 
(because 
$PL_K(\varphi^k\sigma^l)\cap PL_K^\circ(\varphi^{-k}\sigma^{-l})
=\emptyset$ as seen by using \cite[Proposition 8.2]{Nasu-te}), 
it follows that $i'l\neq j'k$.  
Therefore, it follows from (1) that for all 
$(i,j)\in\cH((i',j'),(k,l)]$,  
$\varphi^i\sigma^j\in PL_K^\circ(\varphi^k\sigma^l)$. 
Since either $(i_1,j_1)$ or $(i_2,j_2)$ is in  
$\cH((i',j'),(k,l)]$, 
either $\varphi^{i_1}\sigma^{j_1}$ or $\varphi^{i_2}\sigma^{j_2}$ 
would be right $\varphi^k\sigma^l$-expansive, 
contrary to the hypothesis. 

Hence we have proved that
$(i,j)\in\cH[(i_1,j_1),(i_2,j_2)]\cap\Z^2$ if and 
only if $\varphi^i\sigma^j\in PL_K(\varphi^k\sigma^l)$. 
The proof that 
$(i,j)\in\cH((i_1,j_1),(i_2,j_2))\cap\Z^2$ if and 
only if $\varphi^i\sigma^j\in PL_K^\circ(\varphi^k\sigma^l)$ is 
given by modifying the proof above straightforwardly. 
\end{proof} 

\begin{example} Consider the automorphism $\varphi$ of $(X,\sigma)$ 
in Example 9.17. Let $K$ be the subgroup of $H(X)$ generated by 
$\{\sigma,\varphi\}$. Applying Propositions 10.6 and 10.7 to
the limits of onesided directions of $\varphi$ 
obtained in Example 9.17, we see: 
\begin{align*}
D_K(\sigma)=PL_K(\sigma)=\{\varphi^i\sigma^j\,|\, j\geq |i|\},\q
C_K(\sigma)=QR_K(\sigma)=\{\varphi^i\sigma^j\,|\, j\geq 2|i|\}.
\end{align*} 
These are also seen by using Proposition 10.8(2) as follows:
since $\varphi\sigma^{-P_L(\varphi)}\pl\sigma$ (by Proposition 6.12) 
and $(\varphi\sigma^{Q_L(\varphi)})^{-1}\pl\sigma$ 
(by Theorem 4.2(2)) 
and $\varphi\sigma^{-P_L(\varphi)}$ and 
$(\varphi\sigma^{Q_L(\varphi)})^{-1}$ 
are not right $\sigma$-expansive (by Example 9.17), 
by Proposition 10.8(2)  
$PL_K(\sigma)=\{\varphi^i\sigma^j\,|\,(i,j)\in
\cH[(1,-P_L(\varphi)),(-1,-Q_L(\varphi))]\}$; 
since $\varphi\sigma^{-Q_R(\varphi)}\qr\sigma$ 
(by Theorem 4.2(1)) and 
$(\varphi\sigma^{P_R(\varphi)})^{-1}\qr\sigma$ 
(by Proposition 6.12) and 
$\varphi\sigma^{-Q_R(\varphi)}$ and 
$(\varphi\sigma^{P_R(\varphi)})^{-1}$
are not left $\sigma$-expansive (by Example 9.17), 
by Proposition 10.8(2)  
$QR_K(\sigma)=\{\varphi^i\sigma^j\,|\,(i,j)\in
\cH[(1,Q_R(\varphi)),(-1,-P_R(\varphi))]\}$; 
since $\varphi\sigma^{C_R(\varphi)}\lr\sigma$ 
(by Corollary 7.4(1)(b)) and 
$(\varphi\sigma^{C_L(\varphi)})^{-1}\lr\sigma$ 
(by Corollary 7.4(2)(b)) and 
$\varphi\sigma^{C_R(\varphi)}$ and 
$(\varphi\sigma^{C_L(\varphi)})^{-1}$ are not 
expansive, by 
Proposition 10.8(2)  
$C_K(\sigma)=\{\varphi^i\sigma^j\,|\,(i,j)\in
\cH[(1,C_R(\varphi)),(-1,-C_L(\varphi))]\}$. 
Moreover, by Corollary 11.14(2) 
\[D^\ast_K(\sigma)=
\{\varphi^i\sigma^j\,|\,(i,j)\in
\cH[(1,-P_L(\varphi)),(-1,-Q_L(\varphi))]\}=PL_K(\sigma)\]
because $\varphi\sigma^{-P_L(\varphi)}\mm\sigma$  
and $(\varphi\sigma^{Q_L(\varphi)})^{-1}\mm\sigma$ 
and $\varphi\sigma^{-P_L(\varphi)}$ and 
$(\varphi\sigma^{Q_L(\varphi)})^{-1}$ are 
neither right $\sigma$-expansive 
nor left $\sigma$-expansive (by Example 9.17).  

Since $(X,\sigma)$ is an SFT, 
$C_K(\sigma)$ and $C_K(\sigma^{-1})$ are 
ELR cones in $K$. 
It is known (\cite[p.202]{Nasu-t}) that no element 
in $D_K(\sigma)\setminus C^\circ_K(\sigma)$ 
is left $\sigma$-expansive
and no element in $D_K(\sigma^{-1})\setminus C_K^\circ(\sigma^{-1})$ 
is left $\sigma$-expansive.  
However no more is known about 
the dynamics of the elements in 
$K\setminus(D_K(\sigma^{-1})\cup D_K(\sigma))$ except that 
$\varphi$ is neither left $\sigma$-expansive 
nor right $\sigma$-expansive.  
\end{example} 

\begin{example}
Consider the LR automorphism $\varphi$ of $(X_M,\sigma_M)$ 
in Example 9.18. Let $K$ be the subgroup of $H(X)$ generated by 
$\{\sigma_M,\varphi\}$. As was seen in Example 9. 18, 
$\varphi\lr\sigma_M$ and $\varphi^{-1}\sigma_M\lr\sigma_M$ and 
$\varphi$ and $\varphi^{-1}\sigma_M$ are neither left $\sigma$-expansive 
nor right $\sigma$-expansive. Therefore by Proposition 10.8(2) and 
Corollary 11.14(2) 
\begin{gather*}
D_K^\ast(\sigma_M)=C_K(\sigma_M)=\cH[(1,0),(-1,1)]
=\{\varphi^i\sigma_M^j \,|\, 
j\geq -i\;\,\text{if}\;\,i<0,
j\geq 0\;\,\text{if}\;\,i\geq 0\}. 
\end{gather*}
Let us refer to the 1-1 textile system $\bar{T}$ in 
\cite[Section 10, Example 1]{Nasu-t}. As is 
shown there, $(X_{\bar{T}},\sigma_{\bar{T}},\varphi_{\bar{T}})$
is conjugate to $(X_M. \varphi^2\sigma_M^{-1}, \sigma_M)$.
One can show (more strongly than stated there) that 
$\varphi_{\bar{T}}\sigma_{\bar{T}}^3$ and
$\varphi_{\bar{T}}\sigma_{\bar{T}}^{-3}$ are LR automorphisms 
of $(X_{\bar{T}},\sigma_{\bar{T}})$ which 
are neither left $\sigma_{\bar{T}}$-expansive nor 
right $\sigma_{\bar{T}}$-expansive. This 
implies that $\varphi^6 \sigma_M^{-2}\lr\varphi^2\sigma_M^{-1}$ and 
$\varphi^6\sigma_M^{-4}\lr\varphi^2\sigma_M^{-1}$
and $\varphi^6 \sigma_M^{-2}$ and $\varphi^6 \sigma_M^{-4}$
are neither left $\varphi^2\sigma_M^{-1}$-expansive nor 
right $\varphi^2\sigma_M^{-1}$-expansive. Hence 
by Proposition 10.8(2) and Corollary 11.14(2) 
\begin{gather*}
D_K^\ast(\varphi^2\sigma_M^{-1})=C_K(\varphi^2\sigma_M^{-1})=
\cH[(6,-2),(6,-4)]\\
=\{\varphi^i\sigma_M^j\,|\, 
-(2/3)i\leq j\leq -(1/3)i, i\geq 0\}. 
\end{gather*}

It is very important to recall 
that an expansive element
in the ELR cone $C_K(\sigma_M)$ and that  
in the ELR cone $C_K(\varphi^2\sigma^{-1})$ define
drastically different topological Markov shifts 
(see \cite[Section 10, Example 1]{Nasu-t}; 
see \cite[Example 5.1]{BoyLin}). 

Nothing is known about the dynamics of elements 
outside the ELR cones 
$C_K(\sigma_M)$, $C_K(\sigma_M^{-1})$, 
$C_K(\varphi^2\sigma_M^{-1})$ and $C_K(\varphi^{-2}\sigma_M)$.
\end{example}
\begin{example} 
Consider the onto endomorphism $\varphi_1$ of $(X,\sigma)$ in 
Example 9.19. Let $J$ be the submonoid of $S(X)$ generated by 
$\{\sigma,\sigma^{-1},\varphi_1\}$. Let $i\geq 0,j\in\Z$.  
Applying Propositions 10.6 and 10.7 to
the results of Example 9.19, we have 
\begin{gather*}
C_J(\sigma^{-1})=PR_J(\sigma)=\{\varphi_1^i\sigma^j\,|\, j\leq -2i\},\q
D_J(\sigma^{-1})=QL_J(\sigma)=\{\varphi_1^i\sigma^j\,|\, j\leq 0\},\\
D_J(\sigma)=PL_J(\sigma)=QR_J(\sigma)=C_J(\sigma) 
=\{\varphi_1^i\sigma^j\,|\, j\geq 2i\}. 
\end{gather*}
One can see that no element in 
$QL_J(\sigma)\setminus PR_J(\sigma)$ is left $\sigma$-expansive
by a similar argument to that in 
\cite[Section 10, Example 2 (p.202)]{Nasu-t}. 
However nothing is known about the dynamics of elements 
in $J\setminus(D_J(\sigma^{-1})\cup D_J(\sigma))$. 
\end{example} 

\begin{example} 
(1) Consider the onto endomorphism $\varphi$ of $(X,\sigma)$ in 
Example 9.20(1). Let $J$ be the submonoid of $S(X)$ 
generated by $\{\sigma,\sigma^{-1},\varphi\}$. Let $i\geq 0,
j\in\Z$. 
Then applying Propositions 10.6 and 10.7 to the results of 
Example 9.20(1), we see the following:
\begin{align*}
C_J(\sigma^{-1})&=PR_J(\sigma)=\{\varphi^i\sigma^j\,|\, j\leq -ni\},&
D_J(\sigma^{-1})&=QL_J(\sigma)=\{\varphi^i\sigma^j\,|\, j\leq mi\},\\
D_J(\sigma)&=QR_J(\sigma)=\{\varphi^i\sigma^j\,|\, j\geq -ni\},&
C_J(\sigma)&=PL_J(\sigma)=\{\varphi^i\sigma^j\,|\, j\geq mi\},\\  
Q_J(\sigma)&=\{\varphi^i\sigma^j\,|\, -ni\leq j\leq mi\}.&
\end{align*} 

(2) Consider the onto endomorphism $\varphi$ of $(X,\sigma)$ in 
Example 9.20(2). Let $J$ be the submonoid of $S(X)$ 
generated by $\{\sigma,\sigma^{-1},\varphi\}$.  Let $i\geq 0,
j\in\Z$. Then applying Propositions 10.6 and 10.7 to the results of 
Example 9.20(2), we see the following:
\begin{align*}
C_J(\sigma^{-1})&=PR_J(\sigma)=\{\varphi^i\sigma^j\,|\, j\leq -si\},&
D_J(\sigma^{-1})&=QL_J(\sigma)=\{\varphi^i\sigma^j\,|\, j\leq 0\},\\
D_J(\sigma)&=QR_J(\sigma)=\{\varphi^i\sigma^j\,|\, j\geq 0\},&
C_J(\sigma)&=PL_J(\sigma)=\{\varphi^i\sigma^j\,|\, j\geq ri\},\\ 
Q_J(\sigma)&=\{\varphi^i\}.&
\end{align*} 
\end{example}

\section{Onesided expansive lines and directions} 

In this section, we define ``onesided expansive'' analogues 
of the notions of expansive lines and expansive components 
for $\Z^2$-actions due to Boyle and Lind \cite{BoyLin} 
(who treated more generally corresponding notions
for $\Z^d$-actions in a general framework), for 
endomorphisms of compact, invertible dynamical systems. 
Then using them we understand 
what the limits of onesided resolving directions of 
endomorphisms of subshifts are and understand the relation between 
``expansiveness'' and ``resolvingness'' for
onto endomorphisms and automorphisms of subshifts 
and then generally for $\Z^d$ actions 
on zero-dimensional compact metric spaces. 

\subsection{The set of onesided expansive directions}

We consider the plane $\R^2$ and the coordinate system 
for it such that a point $(a,b)$ in $\R^2$ goes down 
if $a$ increases, and goes right 
if $b$ increases (like the indexing for matrices). 
The \itl{direction} of 
a non-horizontal line passing through 
points $(a,b)$ and $(c,d)$ with $a\neq c$ is $(b-d)/(a-c)$. 
(For a horizontal line, we do not define its direction.) 
Hereafter a ``line'' means a non-horizontal line unless 
otherwise stated.

For any line $\ell$ in the plane $\R^2$, 
let $\R^2(\cdot,\ell]$ (respectively, $\R^2[\ell,\cdot))$
denote the set of points in $\R^2$ that lie 
on the left (respectively, right) of $\ell$ or on $\ell$. 
For parallel lines $\ell,\ell'=\ell+(0,t)$ with $t\geq 0$ , let 
$\R^2[\ell,\ell']=\R^2[\ell,\cdot)\cap \R^2(\cdot,\ell']$. 
For a horizontal line $\hbar$, let $\R^2(\cdot,\hbar]$ (respectively, 
$\R^2[\hbar,\cdot)$) denote the set of points in $\R^2$ that lie either on 
or above (respectively, below) $\hbar$.  (Here and in what follows, 
for any subset $R$ of $\R^2$ and a point $(a,b)\in \R^2$, 
$R+(a,b)$ denotes $\{(c+a,d+b)\,|\,(c,d)\in R\}$, and hence
for a line $\ell$ passing through a point $(c,d)\in \R^2$, 
$\ell+(a,b)$ denotes the line which 
is parallel with $\ell$ and passes 
through the point $(a+c,b+d)$.)

Let $\varphi$ be an onto endomorphism of 
an invertible dynamical system $(X,\tau)$. 
A line $\ell$ in $\R^2$ is said to be 
\itl{left $\tau$-expansive for $\varphi$} 
(respectively, \itl{right $\tau$-expansive for $\varphi$}), 
if there exists $\delta>0$ such that 
for any $\varphi$-orbits $(x_i)_{i\in\Z},(y_i)_{i\in\Z}$, 
it holds that if 
$d_X(\tau^j(x_i),\tau^j(y_i))\leq\delta$ 
for all $(i,j)\in\R^2(\cdot,\ell]\cap\Z^2$ 
(respectively, for all $(i,j)\in\R^2[\ell,\cdot)\cap\Z^2$) then 
$(x_i)_{i\in\Z}=(y_i)_{i\in\Z}$; 
we call such $\delta$ a 
\itl{left $\tau$-expansive constant} 
(respectively, \itl{right $\tau$-expansive constant}) 
\itl{for $\ell$ for $\varphi$}. 
A line $\ell$ is said to be 
\itl{expansive for $\varphi$}, if there exist $\delta>0$ and 
$t\geq 0$ such that 
for any $\varphi$-orbits $(x_i)_{i\in\Z},(y_i)_{i\in\Z}$, 
it holds that if 
$d_X(\tau^j(x_i),\tau^j(y_i))\leq\delta$ 
for all $(i,j)\in\R^2[\ell-(0,t),\ell+(0,t)]\cap\Z^2$ then 
$(x_i)_{i\in\Z}=(y_i)_{i\in\Z}$; we call such $\delta$ an 
\itl{expansive constant for $\ell$ for $\varphi$}.  
When $\varphi$ is invertible, 
an expansive line for $\varphi$ is  
an expansive line 
for the $\Z^2$-action $(i,j)\mapsto\varphi^i\tau^j$ in the sense of
Boyle and Lind \cite[Definition 2.2]{BoyLin}. 

For any line $\ell$ and $(I,J)\in\Z^2$, 
the line 
$\ell+(I,J)$ is left $\tau$-expansive (respectively, 
right $\tau$-expansive, expansive) for 
$\varphi$ with a left $\tau$-expansive (respectively, 
right $\tau$-expansive, expansive) constant $\delta$  
if and only if so is $\ell$ with the same left 
$\tau$-expansive (respectively, 
right $\tau$-expansive, expansive) 
constant $\delta$ 
(because if $(x_i)_{i\in\Z}$ 
is a $\varphi$-orbit then so is $(\tau^J(x_{i+I}))_{i\in\Z}$). 
For any parallel lines $\ell,\ell'$, 
there exists an integer $J\geq 0$ 
such that $\ell'+(0,J) \subset\R^2[\ell,\cdot)$ and hence, 
if $\ell$ is left $\tau$-expansive for 
$\varphi$ with a left $\tau$-expansive constant $\delta$ 
then so is $\ell'+(0,J)$ with the same 
left $\tau$-expansive constant $\delta$. Therefore,  
if $\ell$ is left $\tau$-expansive for 
$\varphi$ with a left $\tau$-expansive constant $\delta$, 
then so is every parallel line $\ell'$ with $\ell$ 
with the same left $\tau$-expansive constant $\delta$. 
Similarly we know that if 
one of two parallel lines is right $\tau$-expansive 
for $\varphi$ with 
a right $\tau$-expansive constant $\delta$, then so is the other 
with the same right $\tau$-expansive constant $\delta$. 
If $\ell,\ell'$ are parallel lines 
and $t\geq 0$, 
then there exists $J\in\Z$ such that 
$\R^2[\ell-(0,t),\ell+(0,t)]+(0,J)\subset \R^2[\ell'-(0,t+1),\ell'+(0,t+1)]$. 
Hence if $\ell$ is expansive for $\varphi$ with 
an expansive constant $\delta$ , then 
so are all lines $\ell'$ parallel with $\ell$ 
with the same expansive constant $\delta$. 

Therefore we can say that a real number $r$ is 
a \itl{left $\tau$-expansive} 
(respectively, \itl{right $\tau$-expansive, 
expansive}) \itl{direction for $\varphi$} if there exists 
a line $\ell$ of direction $r$ which is left $\tau$-expansive 
(respectively, right $\tau$-expansive, expansive) for $\varphi$; 
we can define similarly a 
\itl{left $\tau$-expansive} 
(respectively, \itl{right $\tau$-expansive}, \itl{expansive}) 
\itl{constant}  for a 
left $\tau$-expansive 
(respectively, right $\tau$-expansive, expansive) 
direction $r$ for $\varphi$. 

We note that if $r$ is a rational number and $r=j/i$ 
with $i\in\N,j\in\Z$, then $r$ is a left $\tau$-expansive 
(respectively, right $\tau$-expansive, expansive) 
direction for $\varphi$ 
if and only if $\varphi^i\tau^j$ is 
left $\tau$-expansive (respectively, right $\tau$-expansive, expansive). 

In \cite{BoyLin}, Boyle and Lind used a geometric coding method 
(as an analogue of sliding block coding of symbolic dynamics), which 
can be called ``sliding pattern coding''. We shall use 
their method taking parallelograms as patterns. 
Let $r\in\R$. Let $m, n$ be nonnegative integers and $l,e\geq 0$. 
We define the \itl{parallelogram $ABCD$ with the left head $EF$}
(respectively, \itl{right head $GH$}) 
\itl{determined by} $(r;m,n,l;e)$ as follows:
the sides $AD$ and $BC$ are horizontal and of length $l$,    
the sides $AB$ and $DC$ are of direction $r$  
and of length $(m+n)\sqrt{1+r^2}$, hence the \itl{height} of the 
parallelogram $ABCD$ is $m+n$, $EF$ is the horizontal edge of 
length $e$ such that
the node $F$ is on the 
side $AB$ and the length of $AF$ is 
$m\sqrt{1+r^2}$ (hence the length of $FB$ is 
$n\sqrt{1+r^2}$) 
(respectively, $GH$ is the horizontal edge of length $e$
such that the node $G$ is on the 
side $DC$ and the length of $DG$ is 
$m\sqrt{1+r^2}$ (hence the length of $GC$ is equal to  
$n\sqrt{1+r^2}$)).  

\begin{proposition} 
Let $\varphi$ be an onto endomorphism of 
an invertible dynamical system $(X,\tau)$. 
Let $r\in\R$. Then $r$ is an expansive direction for 
$\varphi$ if and only if 
$r$ is both a left $\tau$-expansive direction for $\varphi$  
and a right $\tau$-expansive direction for $\varphi$.  
\end{proposition} 
\begin{proof} The ``only-if'' part is clear.  
To prove ``if'' part, 
suppose that $r$ is both 
a left $\tau$-expansive direction for $\varphi$ 
and a right $\tau$-expansive direction for $\varphi$. 

For $\delta>0$ and subsets $R,S$ of $\R^2$, 
let us say, following the idea of ``coding'' of Boyle and Lind 
\cite{BoyLin}, that \itl{$R$ $\delta$-codes $S$ 
for $\varphi$} if for any $\varphi$-orbits 
$(x_i)_{i\in\Z}, (y_i)_{i\in\Z}$ it holds that 
if $d_X(\tau^j(x_i),\tau^j(y_i))\leq\delta$ for all $(i,j)\in R\cap\Z^2$
then $d_X(\tau^j(x_i),\tau^j(y_i))\leq\delta$ for all $(i,j)\in S\cap\Z^2$.  

Let $\ell$ be the line of direction $r$ in $\R^2$
passing through the origin $(0,0)$ . Since $\ell$ is 
left $\tau$-expansive and right $\tau$-expansive, 
there exist a left $\tau$-expansive constant 
$\delta_L>0$ for $\ell$ for $\varphi$ and 
a right $\tau$-expansive constant 
$\delta_R>0$ for $\ell$ for $\varphi$. 
Let $\delta=\min\{\delta_L,\delta_R\}$. Since $\R^2(\cdot,\ell]$ 
$\delta$-codes $\R^2$ for $\varphi$, $\R^2(\cdot,\ell]$ 
 $\delta$-codes in particular the horizontal line segment 
$\{(0,b)\,|\, 0\leq b\leq 1\}$ in $\R^2$ for $\varphi$. 

We consider the product topological space 
$X^{\Z^2}$, which is a compact metric space 
with a compatible metric. 
Let $\Oh_{\varphi,\tau}
=\{(x_{i,j})_{i,j\in\Z}\in X^{\Z^2}\,|\, 
x_{i,j}=\tau^j(x_i), (x_i)_{i\in\Z}\in \Oh_\varphi\}$, 
where $\Oh_\varphi$ denotes 
the set of all $\varphi$-orbits. 
Then $\Oh_{\varphi,\tau}$ is closed and hence compact. 
A standard compactness argument shows that there exist integers 
$m,n,l_0\geq 1$ such that if the parallelogram $A_0B_0C_0D_0$ 
with the right head $G_0H_0$ determined by $(r;m,n,l_0;1)$ is put on 
the plane $\R^2$ in such a way that $G_0$ is put on $(0,0)$ (hence, 
$H_0$ is put on (0,1), $A_0$ is put on $(-m,-l_0-mr)$, 
$B_0$ is put on $(n,-l_0+nr)$, $C_0$ is put on $(n,nr)$ and $D_0$ is 
put on $(-m,-mr)$), then the subset of $\R^2$ \itl{covered} by 
the parallelogram $A_0B_0C_0D_0$ (i.e. the set of points 
on the sides and interior of the parallelogram hulled by the points
$A_0(-m,-l_0-mr), B_0(n,-l_0+nr), C_0(n,nr), D_0(-m,-mr)$) 
$\delta$-codes the subset of $\R^2$ covered by 
the edge $G_0H_0$ (i.e. the line segment with the 
endpoints $G_0(0,0)$ and $H_0(0,1)$). Since 
if $(x_{i,j})_{i,j\in\Z}\in\Oh_{\varphi,\tau}$ then 
$(x_{i+I,j+J})_{i,j\in\Z}\in\Oh_{\varphi,\tau}$ for any 
$I,J\in\Z$, if $G_0$ is put on a point $(I,J)$ with $I,J\in\Z$, 
then the subset of $\R^2$ covered by 
the parallelogram $A_0B_0C_0D_0$  
(i.e. the set of points 
on the sides and interior of the parallelogram hulled by the points
$A_0(-m+I,-l_0-mr+J), B_0(n+I,-l_0+nr+J), C_0(n+I,nr+J), 
D_0(-m+I,-mr+J)$)
$\delta$-codes the set of points covered by the edge $G_0H_0$
(i.e. the line segment with the 
endpoints $G_0(I,J)$ and $H_0(I,1+J)$). 

Let $l=l_0+1$ and consider the parallelogram $ABCD$ with 
the right head $GH$ determined by $(r;m,n,l;1)$. Let $(a,b)\in\R^2$ and 
suppose that the parallelogram $ABCD$ with 
the right head $GH$ is put on the plane $\R^2$ in such a way that 
$G$ is put on the point $(a,b)$. 
If $a\notin\Z$, then the subset of $\R^2$ covered by 
the parallelogram $ABCD$  
$\delta$-codes the subset of $\R^2$ covered by the edge $GH$ 
(because this subset contains no point in $\Z^2$). If 
$a\in\Z$, then for the parallelogram $A_0B_0C_0D_0$ 
with the right head 
$G_0H_0$ such that $G_0$ is put on $(a,\lfloor b\rfloor)$, 
the subset of $\R^2$ covered by 
the parallelogram $A_0B_0C_0D_0$, say $R_0$, is 
included in the subset of $\R^2$ covered by 
the parallelogram  $ABCD$, say $R$. 
Since $R_0$ $\delta$-codes the subset of $\R^2$ 
covered by the edge $G_0H_0$, say $S_0$, 
$R$ $\delta$-codes $S_0$ and hence 
the subset of $\R^2$ covered by the edge $GH$, say $S$, 
because 
$S\setminus S_0$ contains no point in $\Z^2$. 

Therefore we have seen that 
the parallelogram $ABCD$ with the 
right head $GH$ determined by $(r;m,n,l;1)$ is a 
\itl{sliding 
parallelogram for $\delta$-coding for $\varphi$} 
i.e, if the parallelogram $ABCD$ with the 
right head $GH$ determined by $(r;m,n,l;1)$ is located 
on any place of the plane $\R^2$, then the subset of $\R^2$ 
covered by the parallelogram $ABCD$ 
$\delta$-codes the subset of $\R^2$ covered by the edge $GH$. 
Hence by sliding the parallelogram $ABCD$ with 
the edge $GH$ forward 
and backward in the direction $r$
in such a way that $G$ is on the line $\ell+(0,l)$, we know that 
$\R^2[\ell,\ell+(0,l)]$ $\delta$-codes $\R^2[\ell,\ell+(0,l+1)]$. 
By sliding the parallelogram $ABCD$ with the edge $GH$ forward 
and backward in the direction $r$
in such a way that the node $G$ is on the line $\ell+(0,l+1)$, 
we know that 
$\R^2[\ell,\ell+(0,l)]$ $\delta$-codes $\R^2[\ell,\ell+(0,l+2)]$. 
Continuing this process, we know that 
$\R^2[\ell,\ell+(0,l)]$ $\delta$-codes 
$\R^2[\ell,\ell+(0,l+k)]$ for all $k\geq 0$ 
and hence $\R^2[\ell,\cdot)$. 

Therefore, if 
$(x_i)_{i\in\Z},(y_i)_{i\in\Z}$ are  
any $\varphi$-orbits and  
$d_X(\tau^j(x_i),\tau^j(y_i))\leq\delta$ 
for all $(i,j)\in\R^2[\ell,\ell+(0,l)]\cap\Z^2$, 
then $d_X(\tau^j(x_i),\tau^j(y_i))\leq\delta$ 
for all $(i,j)\in\R^2[\ell,\cdot)\cap\Z^2$ and hence 
$(x_i)_{i\in\Z}=(y_i)_{i\in\Z}$, because 
$\delta$ is 
a right $\tau$-expansive constant for $\ell$ for $\varphi$. 
Therefore $r$ is 
an expansive direction for $\varphi$. 
\end{proof} 

Throughout the remainder of this section, 
the terminology given in this proof will be used. 

\begin{proposition} 
Let $\varphi$ be an onto endomorphism 
of an invertible dynamical system $(X,\tau)$. 
The set of all left $\tau$-expansive 
(respectively, right $\tau$-expansive) 
directions for $\varphi$ is an open subset of $\R$. 
\end{proposition} 
\begin{proof} 
Let $r$ be a left $\tau$-expansive direction for $\varphi$, and
$\delta$ a left 
$\tau$-expansive constant for $r$ for $\varphi$. 
Then the proof of Proposition 11.1 shows that 
there exist integers $m,n,l\geq 1$ such that 
the parallelogram $ABCD$ with the right head $GH$ determined 
by $(r; m,n,l;1)$ is
a sliding parallelogram for $\delta$-coding for $\varphi$. 

Put a node $G'$ on the edge $GH$ 
such that $G'\neq G$ and $G'\neq H$. Let $r'$ be the direction 
of the edge $G'C$. Then we readily see that 
the parallelogram $AB'CD'$ with the right head $G'H$   
determined by $(r';m,n,l';e')$ is 
a sliding parallelogram for $\delta$-coding for $\varphi$, 
where $l'$ is $l$ plus the $(m+n)/n$ times of the length of 
$GG'$, and $e'$ is the length of $G'H$. Using this 
sliding parallelogram for $\delta$-coding for $\varphi$
and using a similar argument to that in the proof of 
Proposition 11.1, we see that for a line $\ell'$ of direction 
$r'$ in $\R^2$, $\R^2(\cdot,\ell']$ $\delta$-codes $\R^2$, which 
obviously $\delta$-codes $\R^2(\cdot,\ell]$ for any line $\ell$
of direction $r$. Therefore
for any $\varphi$-orbits $(x_i)_{i\in\Z}$ and 
$(y_i)_{i\in\Z}$, it holds that 
if $d_X(\tau^j(x_i),\tau^j(y_i))\leq\delta$ for all 
$(i,j)\in\R^2(\cdot,\ell']\cap\Z^2$ then 
$(x_i)_{i\in\Z}=(y_i)_{i\in\Z}$. Therefore we know that 
$r'$ is a left $\tau$-expansive direction for $\varphi$. 

If $r''$ is the direction 
of the edge $DG'$, then  
the parallelogram $A''BC''D$ with the right head $G'H$   
determined by $(r'';m,n,l'';e')$ is 
a sliding parallelogram for $\delta$-coding for $\varphi$, 
where $l''$ is $l$ plus the $(m+n)/m$ times of the length of $GG'$. 
Hence we similarly know that $r''$ is 
a left $\tau$-expansive direction for $\varphi$. 
Since $r'$ can be any real number greater than 
the direction $r-1/n$
of the edge $CH$ and less than $r$ and since $r''$ 
can be any real number greater than $r$ and less than 
the direction $r+1/m$
of the edge $HD$, 
for all $r-1/n<s<r+1/m$, $s$ is 
an left $\tau$-expansive direction for $\varphi$. \end{proof} 

Let $\varphi$ be an onto endomorphism of 
an invertible dynamical system 
$(X,\tau)$. Let $\ell$ be a line in 
the plane $\R^2$. 
We say that $\ell$ is 
\itl{right $\tau$-expansive on the upper side for $\varphi$} 
(respectively, \itl{left $\tau$-expansive 
on the upper side for $\varphi$}), 
if there exists $\delta>0$ such that for 
a horizontal line $\hbar$ and 
for any $\varphi$-orbits $(x_i)_{i\in\Z},(y_i)_{i\in\Z}$, 
it holds that if 
$d_X(\tau^j(x_i),\tau^j(y_i))\leq\delta$ 
for all $(i,j)\in\R^2(\cdot,\hbar]\cap\R^2[\ell,\cdot)\cap\Z^2$ 
(respectively, for all 
$(i,j)\in\R^2(\cdot,\hbar]\cap\R^2(\cdot,\ell]\cap\Z^2$) then 
$x_i=y_i$ for all $(i,0)\in \R^2(\cdot,\hbar]$
(or equivalently, $(x_i)_{i\in\Z}=(y_i)_{i\in\Z}$); 
we call such $\delta$ a 
\itl{right $\tau$-expansive constant} 
(respectively, \itl{left $\tau$-expansive constant}) \itl{on the 
upper side for $\ell$ for $\varphi$}. 

Though the notions of 
``left $\tau$-expansive line on the lower side for $\varphi$''
and ``left $\tau$-expansive line on the lower side for $\varphi$''
can be defined similarly and  naturally for onto endomorphisms, 
they shall be included in more general notions for 
not necessarily onto 
endomorphisms of invertible dynamical systems.

Let $\varphi$ be a not necessarily onto endomorphism  
of an invertible dynamical system $(X,\tau)$. 
Let $\hbar(0,0)$ denote the horizontal line in $\R^2$ 
passing through the origin $(0,0)$. 
We say that a line $\ell$ in the plane $\R^2$ is 
\itl{positively left $\tau$-expansive} 
(respectively, 
\itl{positively right $\tau$-expansive}) \itl{for $\varphi$}  
if there exists $\delta>0$ such that 
for any $x,y\in X$, 
it holds that if 
$d_X(\tau^j\varphi^i(x),\tau^j\varphi^i(y))\leq\delta$ 
for all $(i,j)\in\R^2[\hbar(0,0),\cdot)\cap\R^2(\cdot,\ell]\cap\Z^2$ 
(respectively, for all 
$(i,j)\in\R^2[\hbar(0,0),\cdot)\cap\R^2[\ell,\cdot)\cap\Z^2$) then 
$x=y$; we call 
such $\delta$ a \itl{positively left $\tau$-expansive constant} 
(respectively, \itl{positively right $\tau$-expansive constant}) 
\itl{for $\ell$ for $\varphi$}. 
We also say that $\ell$ is   
\itl{positively expansive for $\varphi$}, if there exist  
$\delta>0$ and $t\geq 0$ such that 
for any $x,y\in X$, it holds that if 
$d_X(\tau^j\varphi^i(x),\tau^j\varphi^i(y))\leq\delta$ 
for all 
$(i,j)\in\R^2[\hbar(0,0),\cdot)\cap\R^2[\ell-(0,t),\ell+(0,t)]\cap\Z^2$ 
then $x=y$; 
we call such $\delta$ a \itl{positively expansive constant 
for $\ell$ for $\varphi$}. 

As is seen by a similar argument 
to that given before Proposition 11.1, 
the notions defined above for a line in $\R^2$ depend  
only on its direction. 
Hence we can define each of the terms for lines above 
for directions: for example, we say that a direction $r\in \R$ is 
\itl{left $\tau$-expansive on the upper side for $\varphi$}
if there exists a line of direction $r$ 
which is right $\tau$-expansive on the upper side for $\varphi$, 
and we also say that $\delta>0$ is a  a 
\itl{right $\tau$-expansive constant on the 
upper side for $r$ for $\varphi$} if there exists 
a line $\ell$ of direction $r$ such that $\delta$ 
is a right $\tau$-expansive constant 
on the upper side for $\ell$ for $\varphi$. 

We remark the following easily proved facts for an onto 
endomorphism of an invertible dynamical system $(X,\tau)$: 
if $\ell$ is right $\tau$-expansive 
(respectively, left $\tau$-expansive) on the upper side for $\varphi$
with a right $\tau$-expansive  
(respectively, left $\tau$-expansive) constant $\delta$ 
on the upper side for $\varphi$, then 
$\ell$ is right $\tau$-expansive 
(respectively, left $\tau$-expansive) for $\varphi$ 
with the same $\delta$ as a right $\tau$-expansive  
(respectively, left $\tau$-expansive) constant 
for $\varphi$; 
if $\ell$ is positively right $\tau$-expansive 
(respectively, positively left $\tau$-expansive) for $\varphi$
with a positively right $\tau$-expansive  
(respectively, positively left $\tau$-expansive) constant $\delta$ 
for $\varphi$, 
then $\ell$ is right $\tau$-expansive 
(respectively, left $\tau$-expansive) for $\varphi$ 
with the same $\delta$ as a right $\tau$-expansive 
(respectively, left $\tau$-expansive) constant for $\varphi$.

Here we add the definition of 
$\delta$-coding for a not onto endomorphism $\varphi$ of 
an invertible dynamical system $(X,\tau)$:
for any subsets $R,S$ of $\R^2[\hbar(0,0),\cdot)$, 
we say that \itl{$R$ $\delta$-codes $S$ 
for $\varphi$} if for any $x,y\in X$ it holds that 
if $d_X(\tau^j\varphi^i(x),\tau^j\varphi^j(y))\leq\delta$ 
for all $(i,j)\in R\cap\Z^2$
then $d_X(\tau^j\varphi^i(x),\tau^j\varphi^j(y))\leq\delta$ 
for all $(i,j)\in S\cap\Z^2$. (Recall that the definition of 
$\delta$-coding for an onto endomorphism $\varphi$ of 
an invertible dynamical system $(X,\tau)$ was given 
in the proof of Proposition 11.1.) 

For $s\in\R$ and $(a,b)\in\R^2$, 
let $\ell_s(a,b)$ denote the line of direction $s$  
passing through the point $(a,b)$ and let   
$\hbar(a,b)$ denote the horizontal line passing through $(a,b)$. 
Define 
$R_s(a,b)=
\R^2(\cdot,\hbar(a,b)]\cap\R^2[\ell_s(a,b),\cdot)$ and
$\bar{R}_s(a,b)=
\R^2[\hbar(a,b),\cdot)\cap\R^2(\cdot,\ell_s(a,b)]$. 

\begin{proposition} 
Let $\varphi$ be an endomorphism of 
an invertible dynamical system $(X,\tau)$. 
\begin{enumerate} 
\item When $\varphi$ is onto, if the set of right $\tau$-expansive 
(respectively, left $\tau$-expansive) directions on the upper side 
for $\varphi$ is nonempty, 
then it is an open interval $(\alpha,\infty)$ 
for some $\alpha\in{\R\cup\{-\infty\}}$ (respectively, 
$(-\infty,\alpha')$ for some $\alpha'\in{\R\cup\{\infty\}}$). 
\item 
If the set of positively left $\tau$-expansive 
(respectively, positively right $\tau$-expansive) directions   
for $\varphi$ is nonempty, 
then 
it is an open interval $(\beta,\infty)$ 
for some $\beta\in{\R\cup\{-\infty\}}$ (respectively, 
$(-\infty,\beta')$ for some $\beta'\in{\R\cup\{\infty\}})$. 
\item 
A real number $r$ is a positively expansive 
direction for $\varphi$ 
if and only if $r$ is both a positively 
left $\tau$-expansive direction 
for $\varphi$ and a positively right $\tau$-expansive direction 
for $\varphi$. 
\end{enumerate}
\end{proposition} 
\begin{proof}
(1) Suppose that $r$ is a right $\tau$-expansive direction on the 
upper side for $\varphi$ with a right $\tau$-expansive constant 
$\delta$ on the upper side for $\varphi$. 
Then, so is any $s\in [r,\infty)$, 
because $R_s(a,b)\supset R_r(a,b)$ for $(a,b)\in\R^2$ 
and hence $R_s(a,b)$ $\delta$-codes $R_r(a,b)$. 

For a similar argument to that in the proof of Proposition 11.1 
we see that there exist integers $m,l\geq 1$ 
such that the parallelogram 
$ABCD$ with the left head $EF$  determined by $(r;m,n,l;1)$ with 
$n=0$ (hence $F=B$) is 
a sliding parallelogram for $\delta$-coding for $\varphi$. 
Put a node $F'$ on the edge $EF=EB$ 
so that $F'\neq E$. Let $r'$ be the direction 
of the edge $AF'$. Then we readily see that 
the parallelogram $AB'CD'$ with the right head $EF'=EB'$  
determined by $(r';m,0,l';e')$ is 
a sliding parallelogram for $\delta$-coding for $\varphi$, 
where $l'$ is $l$ plus the length of 
$F'F$, and $e'$ is the length of $EF'$. 
Therefore by a similar argument to that in the proof of 
Proposition 11.1, we know that $r'$ is 
a right $\tau$-expansive direction on the upper side for $\varphi$. 

Since $r'$ can be any real number greater than 
the direction $r-1/m$ 
of the edge $AE$ and 
not greater than $r$, any real number in the 
interval $(r-1/m,r]$ is a right $\tau$-expansive direction 
on the upper side for $\varphi$. Hence, so is
any real number in the interval $(r-1/m,\infty)$ by the above. 
Therefore the set of right $\tau$-expansive directions 
on the upper side for $\varphi$ is a right-unbounded open 
interval.  

(2) 
Suppose that $r$ is a positively left $\tau$-expansive 
direction for $\varphi$ with 
a positively left $\tau$-expansive constant 
$\delta$  for $\varphi$. 
Then, so is any $s\in [r,\infty)$, 
because $\bar{R}_s(a,b)\supset \bar{R}_r(a,b)$ 
for $(a,b)\in\R^2[\hbar(0,0),\cdot)$ 
and hence $\bar{R}_s(a,b)$ $\delta$-codes $\bar{R}_r(a,b)$. 

Let $\delta$ be a positively left $\tau$-expansive constant 
for $r$ for $\varphi$. Let $\N_0=\N\cup\{0\}$.
We consider the product topological space $X^{\N_0\times\Z}$, 
which is a compact metric space with a compatible metric. 
Let $\hf{\Oh}_{\varphi,\tau}
=\{(x_{i,j})_{i\in\N,j\in\Z}\in X^{\N_0\times\Z}\,|\, 
x_{i,j}=\tau^j\varphi^i(x), x\in X$.  
Then $\hf{\Oh}_{\varphi,\tau}$ is closed and hence compact. 
Using a standard compactness argument and 
a similar argument to that in the proof of Proposition 11.1, 
we see that there exist integers $n,l\geq 1$ such that 
the parallelogram $ABCD$ with the left head $GH$ 
determined by $(r;m,n,l;1)$ with $m=0$ (hence $G=D$) is 
a sliding parallelogram for $\delta$-coding for $\varphi$ 
i.e., if the parallelogram $ABCD$ with the 
right head $GH$ determined by $(r;0,n,l;1)$ is located 
on any place of the lower half-plane $\R^2[\hbar(0,0),\cdot)$, 
then the subset of $\R^2$ covered by the parallelogram $ABCD$ 
$\delta$-codes the subset of 
$\R^2[\hbar(0,0),\cdot)$ covered by the edge $GH$. 
  
Put a node $G'$ on the edge $GH=DH$ 
such that $G'\neq H$. Let $r'$ be the direction 
of the edge $G'C$. Then we readily see that 
the parallelogram $AB'CD'$ with the right head $G'H=D'H$  
determined by $(r';0,n,l';e')$ is 
a sliding parallelogram for $\delta$-coding for $\varphi$, 
where $l'$ is $l$ plus  the length of 
$GG'$, and $e'$ is the length of $G'H$. 
Therefore by a similar argument to that in the proof of 
Proposition 11.1, we know that $r'$ is 
a positively left $\tau$-expansive direction  for $\varphi$. 

Since $r'$ can be any real number greater than the direction 
of the edge $HC$ (i.e. $r-1/n$) and not greater than $r$, 
any real number in the interval $(r-1/n,r]$ 
is a positively left $\tau$-expansive direction 
for $\varphi$. Therefore, so is
any real number in the interval $(r-1/n,\infty)$ by the above. 
Hence the set of positively left $\tau$-expansive directions 
for $\varphi$ is a right-unbounded open 
interval.  

(3) The proof is similar to that of Proposition 11.1.
\end{proof} 

Let $\varphi$ be an endomorphism of an invertible 
dynamical system $(X,\tau)$ such that $\tau$ is expansive.  
Let $\delta$ be an expansive constant for $\tau$. 
As is well known, there exist $m,n\geq 0$ 
such that \itl{$\varphi$ is of $(m,n)$-type 
with respect to $\delta$}, i.e., the subset  
$\{(0,j)\,|\, -m\leq j\leq n\}+(i',j')$ of lattice points 
in $\R^2$  $\delta$-codes the singleton
$\{(1,0)\}+(i',j')$ for $\varphi$ for all $(i',j')\in\Z^2$ 
if $\varphi$ is onto, 
and otherwise, for all $(i',j')\in\N_0\times\Z$. (In fact, 
a standard compactness argument shows that 
for any $\epsilon>0$ there exist $m,n\geq 0$ such that for any 
$x,y\in X$ it holds that
if $d_X(\tau^j(x),\tau^j(y))<\delta$ for $-m\leq j\leq n$ then 
$d_X(x,y)\leq \epsilon$.) 

An endomorphism of an invertible dynamical 
system $(X,\tau)$ is said to be 
\itl{right $\tau$-closing}  
(respectively, \itl{left $\tau$-closing}) if 
$\varphi$ never collapses distinct left 
$\tau$-asymptotic (respectively, right $\tau$-asymptotic) points.

For any $(a,b)\in\R^2$ and any $r\in\R$, 
let $\bar{S}_r(a,b)=
\R^2[\hbar(a,b),\cdot)\cap\R^2[\ell_r(a,b),\cdot)$, let 
$\hbar[(a,b),\cdot)=\{(a,c)\in\R^2\,|\,c\geq b\}$ (i.e., 
$\hbar[(a,b),\cdot)$ denotes 
the horizontal, closed right half-line with endpoint $(a,b)$), 
and let $\hbar(\cdot,(a,b)]=\{(a,c)\in\R^2\,|\,c\leq b\}$.

\begin{proposition} 
Let $\varphi$ be an endomorphism of an expansive, invertible 
dynamical system $(X,\tau)$. 
\begin{enumerate}
\item[(1)] 
When $\varphi$ is onto, 
there exist a right $\tau$-expansive direction 
on the upper side for $\varphi$ 
and a left $\tau$-expansive direction 
on the upper side for $\varphi$. 
\item[(2)]
When $\varphi$ is onto, if $\varphi$ is of $(m,n)$-type 
with respect to 
an expansive constant $\delta$ for $\tau$ then for
$\alpha,\alpha'$ in (1) of Proposition 11.3, 
$\alpha\leq m$ and $\alpha'\geq -n$.
\item[(3)] 
There exists a positively left $\tau$-expansive 
(positively right $\tau$-expansive) direction for 
$\varphi$ if and only if 
$\varphi$ is right $\tau$-closing (respectively, 
left $\tau$-closing). 
\item[(4)]
If $\varphi$ is of $(m,n)$-type with respect to 
an expansive constant $\delta$ for $\tau$ and 
there exist left $\tau$-asymptotic 
(respectively, right $\tau$-asymptotic) 
distinct points in $X$, then for 
$\beta,\beta'$ in (2) of Proposition 11.3, 
$\beta\geq -n$ (respectively, $\beta'\leq m$).
\end{enumerate}
\end{proposition}
\begin{proof} 
(1) Let $\delta$ be an expansive constant for $\tau$. 
Suppose that $\varphi$ is of $(m,n)$-type with respect to 
$\delta$. Then it is easily observed that
$R_{m+1}(0,0)$ $\delta$-codes $R_{m+1}(0,-1)$ 
for $\varphi$ and hence  
$R_{m+1}(0,-j)$ for all $j\geq 0$. 
Thus $R_{m+1}(0,0)$ $\delta$-codes 
$\R^2(\cdot,\hbar(0,0)]$ for $\varphi$. 
Hence, since for
any $\varphi$-orbits $\xseq,\yseq$ it holds that if
$d_X(\tau^j(x_i),\tau^j(y_i))\leq \delta$
for all $(i,j)\in \R^2(\cdot,\hbar(0,0)]$ then 
$(x_i)_{i\leq 0}=(y_i)_{i\leq 0}$ (because $\delta$ is an 
expansive constant for $\tau$), 
$m+1$ is a left $\tau$-expansive direction on the upper side 
for $\varphi$. 

We symmetrically see that $-n-1$ is 
a left $\tau$-expansive direction $\varphi$. 

(2) To see that $\alpha\leq m$, 
it suffices to show that any $s>m$ is 
a right $\tau$-expansive direction 
on the upper side for $\varphi$. Let $s>m$. 
Since $\varphi$ is of $(m,n)$-type 
with respect to 
an expansive constant $\delta$ for $\tau$, 
it easily follows that 
for any point $(a,b)\in\R^2$ with $a\in\Z$, $\hbar[(a,b),\cdot)$ 
$\delta$-codes $\bar{S}_m(a,b)$. 
Since $R_s(0,0)\supset\hbar[(i,si),\cdot)$ for all $i\leq 0$, 
$R_s(0,0)$ $\delta$-codes $\bar{S}_m(i,si)$ for all $i\leq 0$.
Since $s>m$, it follows that 
$\cup_{i\leq 0}\bar{S}_m(i,si)=\R^2$, and 
hence $R_s(0,0)$ $\delta$-codes $\R^2$. Therefore $s$ is 
a right $\tau$-expansive direction 
on the upper side for $\varphi$. 

A similar proof shows that $-n \leq \alpha'$.
 
(3) Suppose that $\varphi$ is not left $\tau$-closing. 
Then there exist distinct left 
$\tau$-asymptotic points $x,y\in X$ 
such that $\varphi(x)=\varphi(y)$.
Let $\epsilon>0$. Then there exists $t\in\Z$ such that 
$d_X(\tau^s(x),\tau^s(y))<\epsilon$ for all $s\leq t$. 
It is easily observed that for any $r\in\R$, 
for any lattice point $(i,j)$ in 
$\bar{R}_r(0,t)$, 
$d_X(\tau^j\varphi^i(x),\tau^j\varphi^i(y))<\epsilon$. 
Since $\epsilon$ is arbitrary, 
there exists no positively left $\tau$-expansive direction 
for $\varphi$. Hence the ``only-if'' part of (2) is proved. 
(Note that this is done without using the expansiveness of $\tau$). 

To prove the ``if'' part, suppose that 
$\varphi$ is left $\tau$-closing. Let $\delta$ be  
an expansive constant for $\tau$. 
Then the subset 
$(\{(0,j)\,|\,j\leq 0\}\cup\{(1,j)\,|\,j\in\Z\})+(i',j')$ 
of lattice points in $\R^2$ $\delta$-codes the singleton 
$\{(0,1)\}+(i',j')$ 
for all $(i',j')\in\N_0\times\Z$. (For otherwise, using 
the condition that $\delta$ is an expansive constant 
for $\tau$, 
we see that there would exist $x,y\in X$ such that 
for any $\epsilon>0$
there is $j_\epsilon\leq 0$ with 
$d_X(\tau^j(x),\tau^j(y))\leq \epsilon$ 
for all $j\leq j_\epsilon$, 
$d_X(\tau(x),\tau(y))> \delta$ 
and $\varphi(x)=\varphi(y)$, which cannot be the case 
because $\varphi$ is left $\tau$-closing.) 
Since $\delta$ is an expansive constant for $\tau$, 
a standard compactness argument shows that there exists 
$k\in\Z$ such that 
$(\{(0,j)\,|\,j\leq 0\}\cup\{(1,j)\,|\,j\leq k\})+(i',j')$ 
$\delta$-codes 
$\{(0,1)\}+(i',j')$ 
for all $(i',j')\in\N_0\times\Z$. 
Therefore, it is easily observed that 
$\bar{R}_k(0,0)$ $\delta$-codes $\bar{R}_k(0,1)$ 
and hence    
$\bar{R}_k(0,j)$ for all $j\geq 0$. Thus 
$\bar{R}_k(0,0)$ $\delta$-codes 
$\R^2[\hbar(0,0),\cdot)$, which implies that 
$k$ is a positively left $\tau$-expansive direction for $\varphi$. 

(4) It suffices to show that $-n$ is not a 
positively left $\tau$-expansive direction for $\varphi$. 
To do this, 
assume the contrary.  
Since $\varphi$ is of $(m,n)$-type 
with respect to 
an expansive constant $\delta$ for $\tau$, 
$\hbar(\cdot,(0,0)]$ 
$\delta$-codes $\bar{R}_{-n}(0,0)$, which would  
$\delta$-codes $\R^2[(\hbar(0,0),\cdot)$ for $\varphi$ 
by the assumption that 
$-n$ were a positively left $\tau$-expansive 
direction for $\varphi$. 
Therefore, $\hbar(\cdot,(0,0)]$ 
would $\delta$-codes 
$\hbar(0,0)$ for $\varphi$, which contradicts 
the hypothesis that there exist left $\tau$-asymptotic 
distinct points in $X$. 
\end{proof} 

\begin{remark} If
$\varphi$ is an onto endomorphism of an expansive, 
invertible dynamical system 
$(X,\tau)$ and $\delta$ is an expansive constant 
for $\tau$, then $\delta$ is 
a \itl{universal onesided expansive constant for $\varphi$}
i.e., 
$\delta$ has all of the following properties (1),(2) and (3) 
for every line $\ell$ in $\R^2$: 
(1) if $\ell$ is left $\tau$-expansive 
(respectively, right $\tau$-expansive) 
for $\varphi$, then 
$\delta$ is  a left $\tau$-expansive 
(respectively, right $\tau$-expansive) constant for $\ell$ 
for $\varphi$; (2) if $\ell$ is right $\tau$-expansive 
(respectively, left $\tau$-expansive)
on the upper side for $\varphi$, 
then $\delta$ is  a right $\tau$-expansive 
(respectively, left $\tau$-expansive) constant 
on the upper side for $\ell$ for $\varphi$; 
(3) if $\ell$ is positively left $\tau$-expansive 
(positively right $\tau$-expansive) 
for $\varphi$, then 
$\delta$ is  a positively left $\tau$-expansive 
(respectively, positively right $\tau$-expansive) 
constant for $\ell$ for $\varphi$. 
\end{remark} 
\begin{proof}
Let $\ell$ be 
any line in $\R^2$. To see (1) suppose that 
$\delta_0$ is a left $\tau$-expansive constant for $\ell$
for $\varphi$. A standard compactness argument shows that
there exists an integer $t\geq 0$ such that for any $x,y\in X$ 
it holds that 
if $d_X(\tau^j(x),\tau^j(y))\leq\delta$ for $-t\leq j\leq t$ 
then $d_X(x,y)\leq \delta_0$. Hence it follows that 
for any $\varphi$-orbits $\xseq,\yseq$ 
it holds that if 
$d_X(\tau^j(x_i),\tau^j(y_i))\leq\delta$ 
for all $(i,j)\in\R^2(\cdot,\ell+(0,t)]\cap\Z^2$ then 
$d_X(\tau^j(x_i),\tau^j(y_i))\leq\delta_0$ 
for all $(i,j)\in\R^2(\cdot,\ell]\cap\Z^2$ and hence 
$(x_i)_{i\in\Z}=(y_i)_{i\in\Z}$. Therefore $\delta$ is an 
left $\tau$-expansive constant for $\ell+(0,t)$ for $\varphi$ 
and hence for $\ell$ for $\varphi$. 

By similar arguments, (2) and (3) are proved.
\end{proof} 

\begin{proposition} 
Let $\varphi$ be an onto endomorphism of 
an expansive, invertible 
dynamical system $(X,\tau)$. 
\begin{enumerate} 
\item[(1)] 
The set of right $\tau$-expansive 
(respectively, left $\tau$-expansive) 
directions on the upper side for $\varphi$
is some open interval $(\alpha,\infty)$ 
(respectively, $(-\infty, \alpha')$)
such that  
$\alpha$ is not a right $\tau$-expansive direction
(respectively, $\alpha'$ is not a left $\tau$-expansive direction) 
for $\varphi$, 
and hence it is a (connected) component of the set of 
right $\tau$-expansive (respectively, left $\tau$-expansive) 
directions for $\varphi$ in $\R$. 
\item[(2)] 
If $\varphi$ is right $\tau$-closing (respectively, 
left $\tau$-closing), then
the set of positively left $\tau$-expansive  
(respectively, positively right $\tau$-expansive) directions  
for $\varphi$ is some open interval $(\beta,\infty)$ 
(respectively, $(-\infty, \beta')$)
such that $\beta$ is not a left $\tau$-expansive direction
(respectively, $\beta'$ is not a right $\tau$-expansive direction) 
for $\varphi$, 
and hence it is a component of the set of 
left $\tau$-expansive (respectively, right $\tau$-expansive) 
directions for $\varphi$ in $\R$ and 
 $\beta\geq\alpha'$ for  
$\alpha'$ in (1) unless $\beta=-\infty$ and $\alpha'=\infty$
(respectively, $\beta'\leq\alpha$ for $\alpha$ in (1) 
unless $\beta=\infty$ and $\alpha'=-\infty$ ). 
\item[(3)] 
If the set of positively expansive directions for $\varphi$  
is nonempty, then it is some open interval $(\beta,\beta')$
such that $\beta$ is not a left $\tau$-expansive 
direction for $\varphi$
and $\beta'$ is not a right $\tau$-expansive 
direction for $\varphi$, 
and hence it is a component of 
the set of expansive directions  
for $\varphi$ in $\R$. 
\end{enumerate} 
\end{proposition} 
\begin{proof} 
(1) Let $\delta$ be an expansive constant for $\tau$. 
By Propositions 11.3(1) and 11.4(1) there exists 
an open interval $(\alpha,\infty)$ which equals 
the set of all right $\tau$-expansive 
directions on the upper side 
for $\varphi$.  
We shall prove that 
$\alpha$ is not a right $\tau$-expansive direction for $\varphi$. 
This is true when $\alpha=-\infty$ 
because, by definition, a 
right $\tau$-expansive direction for $\varphi$ 
belongs to $\R$. 

Suppose that $\alpha\in\R$, and assume that $\alpha$
were a right $\tau$-expansive direction for $\varphi$. 
Then, 
since by Remark 11.5 $\delta$ is a universal onesided 
expansive constant for $\varphi$, 
by the proof of Proposition 11.1 
there exist positive integers $m,n,l$ such that the parallelogram 
$ABCD$ with the left head $EF$ determined by $(\alpha;m,n,l;1)$
is a sliding  parallelogram for $\delta$-coding for $\varphi$. 
Suppose that $\alpha- 1/(m+n)< s<\alpha$. 
Let $e=1-m(\alpha-s)$. Then $n/(m+n)<e<1$. Put a node $F'$ on the 
edge $EF$ in such a way that the edge $EF'$ is of length $e$. Then 
we have the parallelogram $AB'CD'$ with the left head $EF'$ 
determined by $(s;m,n,l+((m+n)/m)(1-e);e)$, which is a sliding 
parallelogram for $\delta$-coding for $\varphi$. 
If we put the parallelogram $AB'CD'$ with the head $EF'$ on the 
plane $\R^2$ in such a way that $B'$ is put on a point $(a,b)\in\R^2$ 
and then slide the parallelogram $AB'CD'$ with the edge $EF'$ 
upward in such a way that the side $AB'$ is always on $\ell_s(a,b)$, 
then we see that $R_s(a,b)$ $\delta$-codes $R_s(a-n,b-ns-e)$ 
for $\varphi$. 
Using this, we know that $R_s(0,0)$
$\delta$-codes $R_s(-kn,-k(ns+e))$ for $\varphi$ 
for all $k\geq 0$.  
Let $t=(ns+e)/n$. Then, since 
$R_t(0,1)\subset\cup_{k=0}^\infty R_s(-kn,-k(ns+e))$, 
$R_s(0,0)$ $\delta$-codes $R_t(0,1)$ for $\varphi$. 
Since $t> \alpha$ (because $t=s+e/n>(\alpha-1/(m+n))+1/(m+n)=\alpha$), 
$\ell_t(0,1)$ is right 
$\tau$-expansive on the upper side for $\varphi$. Moreover, 
since $\delta$ is a universal onesided 
expansive constant for $\varphi$ (by Remark 11.5), 
$\delta$ is a right $\tau$-expansive constant on the upper side 
for $\ell_t(0,1)$ for $\varphi$. 
Therefore, since $R_s(0,0)$ $\delta$-codes $R_t(0,1)$ 
for $\varphi$, 
$s$ would be a right $\tau$-expansive direction 
on the upper side for $\varphi$, which contradicts the fact that 
$s<\alpha$ and $(\alpha,\infty)$ is the set of all
right $\tau$-expansive directions 
on the upper side for $\varphi$. 

Therefore we have proved the first part of (1). By this 
and the fact that any right $\tau$-expansive direction on the 
upper side for $\varphi$ is 
a right $\tau$-expansive direction for $\varphi$
we see that $(\alpha,\infty)$ is a component of the set of 
right $\tau$-expansive directions for $\varphi$ in $\R$.

(2)  Let $\delta$ be an expansive constant for $\tau$.
Suppose that $\varphi$ is right $\tau$-closing. 
Then by Propositions 11.3(2) and 11.4(3),
the set of all 
positively left $\tau$-expansive directions for $\varphi$    
is some open interval $(\beta,\infty)$. 
We shall prove that $\beta$ is  
not a left $\tau$-expansive direction for $\varphi$. 
This is true when $\beta=-\infty$ (by definition). 

Suppose that $\beta\in\R$, and assume that $\beta$ 
were a left $\tau$-expansive direction for $\varphi$. 
Then, 
since by Remark 11.5 $\delta$ is a universal onesided 
expansive constant for $\varphi$,  
by the proof of Proposition 11.1  
there exist positive integers $m,n,l$ such that the parallelogram 
$ABCD$ with the right head $GH$ determined by $(\beta;m,n,l;1)$
is a sliding  parallelogram for $\delta$-coding for $\varphi$. 
Suppose that $\beta- 1/(m+n)< s<\beta$. 
Let 
$e=1-n(\beta-s)$. Then $m/(m+n)<e<1$. Put a node $G'$ on the 
edge $GH$ in such a way that the edge $G'H$ is of length $e$. Then 
we have the parallelogram $AB'CD'$ with the right head $G'H$ 
determined by $(s;m,n,l+((m+n)/n)(1-e);e)$, which is a sliding 
parallelogram for $\delta$-coding for $\varphi$. 
If we put the parallelogram $AB'CD'$ with the head $G'H$ on the 
plane $\R^2$ in such a way that $D'$ is put on a point $(a,b)\in\R^2$ 
and then slide the parallelogram $AB'CD'$ with the edge $G'H$ 
downward in such a way that the side $D'C$ is always on $\ell_s(a,b)$, 
then we see that $\bar{R}_s(a,b)$ $\delta$-codes $\bar{R}_s(a+m,b+ms+e)$ 
for $\varphi$. 
Using this, we know that $\bar{R}_s(0,0)$
$\delta$-codes $\bar{R}_s(km,k(ms+e))$ for $\varphi$ 
for all $k\geq 0$.  
Let $t=(ms+e)/m$. Then, since 
$\bar{R}_t(0,-1)\subset\cup_{k=0}^\infty \bar{R}_s(km,k(ms+e))$, 
$\bar{R}_s(0,0)$ $\delta$-codes $\bar{R}_t(0,-1)$ for $\varphi$. 
Since $t> \beta$ (because $t=s+e/m>\beta$), 
$\ell_t(0,-1)$ is positively left 
$\tau$-expansive  for $\varphi$. 
Moreover, 
since $\delta$ is a universal onesided 
expansive constant for $\varphi$ (by Remark 11.5), 
$\delta$ is a positively left $\tau$-expansive constant  
for $\ell_t(0,-1)$ for $\varphi$. 
Therefore, since $\bar{R}_s(0,0)$ $\delta$-codes $\bar{R}_t(0,-1)$ 
for $\varphi$, 
$s$ would be a positively left $\tau$-expansive direction 
for $\varphi$, 
which contradicts the fact that 
$s<\beta$ and $(\beta,\infty)$ is the set of all positively 
left $\tau$-expansive directions  for $\varphi$. 

Therefore we have proved the first part of (2). From
this and the 
fact that any positively left $\tau$-expansive direction 
for $\varphi$ is a left $\tau$-expansive direction 
for $\varphi$, it follows that $(\beta,\infty)$ is  a component 
of the set of left $\tau$-expansive directions for $\varphi$ in $\R$. 
Since so is $(-\infty,\alpha')$ by (1), it follows that 
$\alpha'\leq \beta$ unless $\beta=-\infty$ and $\alpha'=\infty$. 

(3) By (2) and Propositions 11.3(3) and 11.1. 
\end{proof} 

We remark that 
the proposition above includes a generalization of
Proposition 9.9(3) to onto endomorphisms of 
expansive  invertible dynamical systems, which is seen by 
comparing it with Theorem 11.8 in the next subsection.

\subsection{Onesided expansiveness and onesided resolvingness}

Now we return to symbolic dynamics. 
Let $(X,\sigma)$ be a subshift over an alphabet $A$. 
For an onto endomorphism $\varphi$ of $(X,\sigma)$, 
define $O_{\varphi,\sigma}$ to be the set of all two-dimensional 
configurations $(a_{i,j})_{i,j\in\Z}$ with $a_{i,j}\in A$ such that
if we put $(a_{i,j})_{j\in\Z}=x_i$ for $i\in\Z$, 
then $(x_i)_{i\in\Z}$ is a $\varphi$-orbit. 
For a not necessarily onto endomorphism $\varphi$ of $(X,\sigma)$
define $\hf{O}_{\varphi,\sigma}$ to be 
the set of all two-dimensional 
configurations $(a_{i,j})_{i\in\N_0,j\in\Z}$ with $a_{i,j}\in A$ 
such that there exists $x\in X$ with 
$\varphi^i(x)=(a_{i,j})_{j\in\Z}$ for 
all $i\in \N_0$, where $\N_0=\N\cup\{0\}$. 

Let $\varphi$ be an endomorphism of $(X,\sigma)$. 
Following the idea of ``coding'' of Boyle and Lind in 
\cite{BoyLin} we define the following: if $\varphi$ is onto, 
then for any subsets $R,S$ of the plane $\R^2$ we say that 
$R$ \itl{codes} $S$ \itl{for $\varphi$} if for any 
$(a_{i,j})_{i,j\in\Z}, (b_{i,j})_{i,j\in\Z}\in O_{\varphi,\sigma}$
it holds that if $a_{i,j}=b_{i,j}$ for all $(i,j)\in R\cap\Z^2$ then 
$a_{i,j}=b_{i,j}$ for all $(i,j)\in S\cap\Z^2$; if 
$\varphi$ is not onto, then 
for any subsets $R,S$ of the half-plane $\R^2[\hbar(0,0),\cdot)$ 
we say that $R$ \itl{codes} $S$ \itl{for $\varphi$} if for any 
$(a_{i,j})_{i\in\N_0,j\in\Z}, (b_{i,j})_{i\in\N_0,j\in\Z}
\in \hf{O}_{\varphi,\sigma}$
it holds that if $a_{i,j}=b_{i,j}$ for all $(i,j)\in R\cap\Z^2$ 
then $a_{i,j}=b_{i,j}$ for all $(i,j)\in S\cap\Z^2$. 

Let $\varphi$ be an endomorphism of a subshift $(X,\sigma)$. 
Then the following hold for any line $\ell$ on the 
plane $\R^2$.
\begin{enumerate}
\item[(1)]  When $\varphi$ is onto, 
$\ell$ is 
left $\sigma$-expansive (respectively, right 
$\sigma$-expansive) for $\varphi$
if and only if $\R^2(\cdot,\ell]$ 
(respectively, $\R^2[\ell,\cdot)$) codes $\R^2$ for $\varphi$.
\item[(2)] When $\varphi$ is onto, 
$\ell$ is 
left $\sigma$-expansive (respectively, right 
$\sigma$-expansive) on the upper side for $\varphi$
if and only if for a horizontal line $\hbar$, 
$\R^2(\cdot,\ell]\cap\R^2(\cdot,\hbar]$ 
(respectively, $\R^2[\ell,\cdot)\cap\R^2(\cdot,\hbar]$) 
codes $\R^2(\cdot,\hbar]$ (or equivalently, $\R^2$) 
for $\varphi$.  
\item[(3)] When $\varphi$ is not necessarily onto, 
$\ell$ is positively 
left $\sigma$-expansive (respectively, positively right 
$\sigma$-expansive) for $\varphi$ 
if and only if 
$\R^2(\cdot,\ell]\cap\R^2[\hbar(0,0),\cdot)$ 
(respectively, $\R^2[\ell,\cdot)\cap\R^2[\hbar(0,0),\cdot)$) 
codes $\hbar(0,0)$ for $\varphi$.  
\end{enumerate}
(We show (1). For $\delta>0$,
there exists an integer $s\geq 0$ such that 
for any $x=\aseq,y=\bseq\in X$, $d_X(x,y)\leq\delta$ if and 
only if $a_j=b_j$ for all $|j|\leq s$, and hence 
$\R^2(\cdot,\ell]$ (respectively, $\R^2[\ell,\cdot)$) 
$\delta$-codes $\R^2$ for $\varphi$ if and only if 
$\R^2(\cdot,\ell+(0,s)]$ 
(respectively, $\R^2[\ell-(0,s),\cdot)$) codes $\R^2$ for $\varphi$. 
Therefore, (1) follows  
because for any $t\in\Z$ it holds that $\R^2(\cdot,\ell+(0,t)]$
(respectively, $\R^2[\ell+(0,t),\cdot)$) codes $\R^2$ if and 
only if so does $\R^2(\cdot,\ell]$
(respectively, $\R^2[\ell,\cdot)$). Similarly 
(2) and (3) are shown.) 

\begin{proposition} 
Let $\varphi$ be an onto endomorphism of a subshift $(X,\sigma)$. 
Let $i\in\N, j\in\Z$. 
\begin{enumerate}
\item 
$\varphi^i\sigma^j$ is 
essentially weakly $p$-L (respectively, essentially weakly $p$-R) 
and right $\sigma$-expansive (respectively, left $\sigma$-expansive)
if and only if $j/i$ is a right $\sigma$-expansive (respectively, 
left $\sigma$-expansive) direction on the upper side for $\varphi$. 
\item
$\varphi^i\sigma^j$ is 
an essentially weakly $q$-R 
(respectively, essentially weakly $q$-L) 
and left $\sigma$-expansive (respectively, right $\sigma$-expansive)
if and only if $j/i$ is a positively left $\sigma$-expansive 
(respectively, positively 
right $\sigma$-expansive) direction for $\varphi$. 
\item
$\varphi^i\sigma^j$ is 
essentially weakly LR 
(respectively, essentially weakly RL) and expansive 
if and only if $j/i$ is both a right $\sigma$-expansive 
(respectively, left 
$\sigma$-expansive) direction on the upper side for $\varphi$ 
and a positively left  $\sigma$-expansive 
(respectively, positively right $\sigma$-expansive) direction 
for $\varphi$. 
\item 
If $(X,\sigma)$ is an SFT, then (1) with 
all ``weakly'' in it deleted holds; if in addition, 
$\varphi$ is one-to-one or $\sigma$ is topologically transitive, 
then (2) and (3) with all ``weakly'' in them deleted hold.  
\end{enumerate} 
\end{proposition} 
\begin{proof} 
(1) Since $\varphi^i\sigma^j$ is 
right $\sigma$-expansive on the upper side if and only if 
$j/i$ is a right $\sigma$-expansive 
direction on the upper side for $\varphi$, 
(1) follows Theorem 8.10.

(2)  The result follows from Theorem 8.9(2). 

(3) The ``only-if'' part follows from (1),(2) and Proposition 11.1. 
The ``if'' part follows from
(1),(2) and Theorems 9.7 and 9.10 and Proposition 11.1. 

(4) By (1),(2),(3) and Remark 2.10. 
\end{proof}

Let $\varphi$ be an onto endomorphism of a subshift $(X,\sigma)$. 
Define $E_L(\varphi)$ (respectively, $E_R(\varphi)$) to be the set of 
all left $\sigma$-expansive (respectively, right $\sigma$-expansive)
directions for $\varphi$. 
Define $E(\varphi)$ to be the set of 
all expansive directions for $\varphi$. Then by 
Proposition 11.1 $E(\varphi)=E_L(\varphi)\cap E_R(\varphi)$,
and by Proposition 11.2, $E_L(\varphi)$, $E_R(\varphi)$ 
and $E(\varphi)$ are open subsets of $\R$. We call $E_L(\varphi)$
(respectively, $E_R(\varphi)$, $E(\varphi)$) the 
\itl{left $\sigma$-expansive} 
(respectively, \itl{right $\sigma$-expansive}, \itl{expansive})
\itl{direction-set for $\varphi$}. 

The following theorem is described under the convention that  
for $\alpha,\beta\in\R\cup\{\infty,-\infty\}$, the interval
$(\alpha,\beta)$ with $\alpha\geq\beta$ means the empty set. 
Recall the remarks about the values 
of the limits of onesided resolving directions of endomorphisms 
of subshifts given immediately after Definition 9.3. 
\begin{theorem}
Let $\varphi$ be an endomorphism of a subshift $(X,\sigma)$. 
If $\varphi$ is onto, then the following (1), (2), (3) and (4) 
hold. 
\begin{enumerate}
\item 
The interval $(-p_L(\varphi),\infty)$  
(respectively, $(-\infty,p_R(\varphi))$) is the set of all
right $\sigma$-expansive (respectively, left $\sigma$-expansive) 
directions on the upper side for $\varphi$ and the right-unbounded  
(connected) component of $E_R(\varphi)$ (respectively, 
the left-unbounded component of $E_L(\varphi)$).  
\item 
The interval $(-q_R(\varphi),\infty)$  
(respectively, $(-\infty,q_L(\varphi))$) is 
the set of all positively
left $\sigma$-expansive 
(respectively, positively right $\sigma$-expansive) 
directions for $\varphi$, and it is the right-unbounded 
component of $E_L(\varphi)$ (respectively, 
the left-unbounded component of $E_R(\varphi)$) if it is nonempty. 
\item Particularly, 
neither $-p_L(\varphi)$ nor $q_L(\varphi)$ is a 
right $\sigma$-expansive direction for $\varphi$, and neither
$p_R(\varphi)$ nor $-q_R(\varphi)$ is a
left $\sigma$-expansive direction for $\varphi$. 
\item 
The interval  
$(c_R(\varphi),\infty)$ (respectively, 
$(-\infty, c_L(\varphi))$) is the right-unbounded 
(respectively, left-unbounded) component of $E(\varphi)$ 
if it is nonempty. Further, the interval 
$(-q_R(\varphi),q_L(\varphi))$ is a 
component of $E(\varphi)$ if it is nonempty. 
\end{enumerate} 
If $\varphi$ is not necessarily onto, then the following holds.
\begin{enumerate}
\item[(5)]
The interval $(-q_R(\varphi),\infty)$  
(respectively, $(-\infty,q_L(\varphi))$) is 
the set of all positively
left $\sigma$-expansive 
(respectively, positively right $\sigma$-expansive) 
directions for $\varphi$, and the interval 
$(-q_R(\varphi),q_L(\varphi))$ is the set of all positively 
expansive directions for $\varphi$. 
\end{enumerate} 
\end{theorem}
\begin{proof}
By Proposition 11.7(1) and Theorem 9.7(1) it holds 
for $i\in\N,j\in\Z$ that 
$j/i\in(-p_L(\varphi),\infty)$ if and only if 
$j/i$ is a right $\sigma$-expansive direction 
on the upper side for $\varphi$. By Proposition 11.6(1)
there exists $\alpha\in\R\cup\{-\infty\}$ such that 
$(\alpha,\infty)$ is 
the set of all right $\sigma$-expansive directions 
on the upper side for $\varphi$ and $\alpha$ is not a 
right $\sigma$-expansive direction for $\varphi$. 
Therefore, since 
$(-p_L(\varphi),\infty)\cap\Q=(\alpha,\infty)\cap\Q$, 
we have $\alpha=-p_L(\varphi)$ and (1) follows. 

The proof of (2) is given in the same way by using
Proposition 11.7(2) and Theorem 9.7(2) and Proposition 11.6(2). 

The statement (3) follows from  (1) and (2),  
and the statement (4) is proved by (1), (2) and Proposition 11.1. 

The proof of the first part of (5) is 
given by Theorem 9.7(4) and Proposition 11.3(2) in the same way 
as in the proof of (1), and the second part of (5) follows from the 
first one and Proposition 11.3(3). 
\end{proof} 

Here we explain the relation between some results of 
Sablik \cite{Sab} and  of K\uo rka \cite{Kur2} and those 
in this paper.

Let $\varphi$ be an endomorphism of a subshift $(X,\sigma)$. 
Using symbolic dynamics versions of sliding parallelograms 
in the proofs of Propositions 11.1 and 11.3 we easily see that 
the following statements (1'),(2') and (3') follow from 
the the statements (1),(2) and (3), respectively, 
in the  third paragraph of this subsection. 
\begin{enumerate}
\item[(1')] 
When $\varphi$ is onto, $\ell$ is 
left $\sigma$-expansive (respectively, right 
$\sigma$-expansive) for $\varphi$
if and only if there exists $t\geq 0$ 
such that $\R^2[\ell-(0,t), \ell+(0,t)]$ 
codes $\R^2[\ell,\cdot)$ (respectively, $\R^2(\cdot,\ell]$) 
for $\varphi$;
\item[(2')] 
When $\varphi$ is onto, $\ell$ is 
left $\sigma$-expansive (respectively, right 
$\sigma$-expansive) on the upper side for $\varphi$
if and only if there exists  $t\geq 0$ 
such that for a horizontal line $\hbar$, 
$\R^2[\ell-(0,t), \ell+(0,t)]\cap \R^2(\cdot,\hbar]$ 
codes $\R^2[\ell,\cdot)\cap \R^2(\cdot,\hbar]$ 
(respectively, $\R^2(\cdot,\ell]\cap \R^2(\cdot,\hbar]$) 
for $\varphi$;
\item[(3')] $\ell$  is positively 
left $\sigma$-expansive (respectively, positively right 
$\sigma$-expansive) for $\varphi$ 
if and only if there exists $t\geq 0$ 
such that  
$\R^2[\ell-(0,t), \ell+(0,t)]\cap\R^2[\hbar(0,0),\cdot)$ 
codes $\R^2[\ell,\cdot)\cap\hbar(0,0)$ 
(respectively, $\R^2(\cdot,\ell]\cap\hbar(0,0)$) for $\varphi$. 
\end{enumerate} 

Let $r\in\R$.  
Sablik's definitions in \cite{Sab} can be stated as follows 
(terminology is not as it is
in \cite{Sab}):
when $\varphi$ is bijective, 
$\varphi$ is said to be \itl{$\Z$-right-expansive} 
(respectively, \itl{$\Z$-left-expansive}) 
\itl{of slope $r$}
if there exist a line $\ell$ of direction $r$ 
and $t\geq 0$ 
such that $\R^2[\ell-(0,t), \ell+(0,t)]$ 
codes $\R^2[\ell,\cdot)\cap\hbar(0,0)$ 
(respectively, $\R^2(\cdot,\ell]\cap\hbar(0,0)$) for $\varphi$ 
(see (1')); (when $\varphi$ is not necessarily onto,)  
$\varphi$ is said to be \itl{$\N_0$-right-expansive} 
(respectively, \itl{$\N_0$-left-expansive}) 
\itl{of slope $r$}
if there exist a line $\ell$ of direction $r$  and 
$t\geq 0$ with the same condition as in (3') for $\ell$. 
Therefore, noting that 
if $\varphi$ is bijective with $\varphi$ of $(m,n)$-type 
and $\varphi^{-1}$ of $(m',n')$-type then  
$t$ in the definition of 
$\Z$-right-expansiveness (respectively, $\Z$-left-expansiveness) 
of slope $r$ of $\varphi$ can be chosen  
so that $t\geq \max\{m,m'\}$ (respectively, 
$t\geq \max\{n,n'\}$), we see that when 
$\varphi$ is bijective,  
$\varphi$ is $\Z$-right-expansive (respectively, 
$\Z$-left-expansive) of slope $r$ if and only if 
$r$ is a left $\sigma$-expansive (respectively, 
right $\sigma$-expansive) direction for $\varphi$.  
We also see similarly that 
$\varphi$ is $\N_0$-right-expansive (respectively, 
$\N_0$-left-expansive) of slope $r$ if and only if 
$r$ is a positively left $\sigma$-expansive (respectively, 
positively right $\sigma$-expansive) direction for $\varphi$. 
Note that ``left'' and ``right'' are used reversely in each pair of  
equivalent conditions in these statements. 

Sablik showed \cite[Theorem 5.2]{Sab} that 
for any endomorphism $\varphi$ of $(m,n)$-type  
of a subshift $(X,\sigma)$ it holds that if the set 
of real numbers $r$ such that $\varphi$ is 
$\N_0$-right-expansive (respectively, 
$\N_0$-left-expansive) of slope $r$ is nonempty, 
then it is an open 
interval $(\beta,\infty)$ with some $-n\leq \beta$
(respectively, $(-\infty,\beta')$ with some $\beta'\leq m$) 
and hence it holds that if the set of real numbers
$r$ such that $\varphi$ is 
positively expansive of slope $r$ is nonempty, then it is 
an open interval $(\beta,\beta')$ with 
some $-n\leq\beta<\beta'\leq m$. 
Propositions 11.3(2),(3) and 11,4(4) extend these.
K\uo rka proved \cite[Theorem 33]{Kur2} that 
the set of real numbers $r$ such that $\varphi$ is 
$\N_0$-right-expansive (respectively, 
$\N_0$-left-expansive) of slope $r$ is nonempty
if and only if $\varphi$ is right-closing 
(respectively, left-closing).  
Proposition 11.4(3) extends this. 
Theorem 11.8(5) refines the above results 
of Sablik and of K\uo rka (see  
Proposition 9.9 and the facts remarked after Definition 
9.3). By Proposition 11.2 the set of real numbers $r$ 
such that $\varphi$ is  
$\Z$-right-expansive (respectively, 
$\Z$-left-expansive) of slope $r$ is an open subset of $\R$, though 
the opposite claim is found in \cite[Section 5.3]{Sab}. 

Recall that $\ell_r(a,b)$ denote the 
non-horizontal line in $\R^2$ having direction $r\in\R$ and
passing through a point $(a,b)\in \R^2$. 

For a direction $r\in\R$, 
define $\ell^+_r(0,0)=
(\ell_r(0,0)\cap\R^2[\hbar(0,0),\cdot))\setminus\{(0,0)\}$ and 
$\ell^-_r(0,0)=(\ell_r(0,0)\cap\R^2(\cdot,\hbar(0,0)])\setminus\{(0,0)\}$. 
Then $\ell_r(0,0)$ equals the disjoint 
union $\ell^-_r(0,0)\cup\{(0,0)\}\cup\ell^+_r(0,0)$. 
Similarly, define $\hbar^+(0,0)=
(\hbar(0,0)\cap\R^2[\ell_0(0,0),\cdot))\setminus\{(0,0)\}$ and 
$\hbar^-(0,0)=
(\hbar(0,0)\cap\R^2(\cdot,\ell_0(0,0)])\setminus\{(0,0)\}$.
We call $\hbar^+(0,0)$ 
(respectively,  $\hbar^-(0,0)$) \itl{the positive} 
(respectively, \itl{negative}) \itl{horizontal half-line}.  
For any subset $R$ of $\R$, define 
\[\Ll^+(R)=\cup_{r\in R}\ell^+_r(0,0),\q 
\Ll^-(R)=\cup_{r\in R}\ell^-_r(0,0).\] 
For any subset $S$ of $\R^2$, we shall call the subset
$-S=\{(-a,-b)\,|\, (a,b)\in S\}$ of $\R^2$ the \itl{reverse of $S$}.  
Then $\Ll^+(R)$ and $\Ll^-(R)$ are the reverses of each other.

Let $\varphi$ be an onto endomorphism of a 
subshift $(X,\sigma)$. Define 
\begin{gather*}
\E_L^+(\varphi)=\Ll^+(E_L(\varphi)), \q
\E_R^+(\varphi)=\Ll^+(E_R(\varphi)),\q
\E^+(\varphi)=\Ll^+(E(\varphi)). 
\end{gather*}
Then it holds for $(i,j)\in\N\times\Z$ that $(i,j)\in \E_L^+(\varphi)$ 
(respectively, $(i,j)\in\E_R^+(\varphi)$, 
$(i,j)\in \E^+(\varphi)$) if and only if 
$\varphi^i\sigma^j$ is left $\sigma$-expansive (respectively, 
right $\sigma$-expansive, expansive). When $X$ is infinite, 
it holds for $(i,j)\in\N_0\times\Z$ that 
$(i,j)$ belongs to the union of $\hbar^-(0,0)\cup\hbar^+(0,0)$ 
and $\E_L^+(\varphi)$
(respectively, $\E_R^+(\varphi)$, 
$\E^+(\varphi)$)  if and only if 
$\varphi^i\sigma^j$ is left $\sigma$-expansive (respectively, 
right $\sigma$-expansive, expansive).
However we call 
$\E_L^+(\varphi)$ (respectively, $\E_R^+(\varphi)$, $\E^+(\varphi)$) 
the \itl{left $\sigma$-expansive}  
(\itl{right $\sigma$-expansive, expansive})
\itl{set of endomorphism type for $\varphi$}. 
$\E_L^+(\varphi)$ is the open subset of $\R^2$ such that 
each (connected) component of it is an open cone (in $\R^2$ 
with apex $(0,0)$) written in the form 
$\Ll^+(I^L)$, where $I^L$ is a 
component of $E_L(\varphi)$; 
$\E_R^+(\varphi)$ is the open subset of $\R^2$ such that 
each component of it is an open cone written in the form 
$\Ll^+(I^R)$, where $I^R$ is a component of $E_R(\varphi)$;  
$\E^+(\varphi)$ is the open subset of $\R^2$ such that 
each component of it is an open cone  written in the form 
$\Ll^+(I^L\cap I^R)$, where 
$I^L$ is a component of $E_L(\varphi)$ and $I^R$ is a 
component of $E_R(\varphi)$ with $I^L\cap I^R\neq\emptyset$. 
We call each component of $\E^+_L(\varphi)$
(respectively, $\E_R^+(\varphi)$, $\E^+(\varphi)$) 
a \itl{left $\sigma$-expansive} 
(respectively, \itl{right $\sigma$-expansive, 
expansive}) \itl{component of endomorphism-type for $\varphi$}. 

By Theorem 11.8 and Theorems 9.7, 9.10 and 9.11 we have:
\begin{corollary} Let $\varphi$ be an onto endomorphism of 
an infinite subshift $(X,\sigma)$. Let $J$ be the
submonoid of $S(X)$ generated by $\{\sigma, \sigma^{-1}, \varphi\}$. 
For any subset $F$ of $J$, let $\Lambda(F)$ denote 
the set $\{(i,j)\in\Z^2\,|\,\varphi^i\sigma^j\in F\}$. 
\begin{enumerate} 
\item If $S=\Ll^+((-p_L(\varphi),\infty))$
(respectively,  $S=\Ll^+((-\infty,p_R(\varphi)))$ then $S$ 
is a right $\sigma$-expansive
(respectively, left $\sigma$-expansive) component 
of endomorphism-type for $\varphi$ and  
$(S\cup\hbar^+(0,0))\cap\Z^2=\Lambda(PL_J^\circ(\sigma))$ 
(respectively,  
$(S\cup\hbar^-(0,0))\cap\Z^2=\Lambda(PR_J^\circ(\sigma))$).
\item 
If $\varphi$ is right-closing (respectively, left-closing) and 
$S=\Ll^+((-q_R(\varphi),\infty))$
(respectively, $S=\Ll^+((-\infty,q_L(\varphi)))$, 
then $S$ is a left $\sigma$-expansive 
(respectively, right $\sigma$-expansive) component 
of endomorphism-type for $\varphi$ and  
$(S\cup\hbar^+(0,0))\cap\Z^2= \Lambda(QR_J^\circ(\sigma))$ (respectively,  
$(S\cup\hbar^-(0,0))\cap\Z^2=\Lambda(QL_J^\circ(\sigma)$)). 
\item 
If $\varphi$ is right-closing (respectively, left-closing) and 
$S=\Ll^+((c_R(\varphi),\infty))$
(respectively, $S=\Ll^+((-\infty,c_L(\varphi)))$, 
then $S$ is an expansive component 
of endomorphism-type for $\varphi$ and  
$(S\cup\hbar^+(0,0))\cap\Z^2= \Lambda(C_J^\circ(\sigma))$ (respectively,  
$(S\cup\hbar^-(0,0))\cap\Z^2=\Lambda(C_J^\circ(\sigma^{-1})$)).
\item 
If $q_L(\varphi)+q_R(\varphi)>0$, then
$S=\Ll^+((-q_R(\varphi), q_L(\varphi)))$ is an 
expansive component of endomorphism-type and  
$S\cap\Z^2=\Lambda(Q_J^\circ(\sigma))$.  
\end{enumerate}
\end{corollary}

\subsection{Onesided-expansive components 
and extended districts}  

Let $\varphi$ be an automorphism of an infinite 
subshift $(X,\sigma)$. We define 
a horizontal line $\hbar$ in $\R^2$ to be 
right $\sigma$-expansive, left $\sigma$-expansive 
and expansive for $\varphi$. Define
\begin{gather*}
\E_L(\varphi)
=\Ll^+(E_L(\varphi))\cup\Ll^-(E_L(\varphi))
\cup\hbar^+(0,0)\cup\hbar^-(0,0),\\ 
\E_R(\varphi)
=\Ll^+(E_R(\varphi))\cup\Ll^-(E_R(\varphi))
\cup\hbar^+(0,0)\cup\hbar^-(0,0),\\ 
\E(\varphi)=\E_L(\varphi)\cap\E_R(\varphi).
\end{gather*}
Then $\E_L(\varphi)$ (respectively, $\E_R(\varphi)$, $\E(\varphi)$)
is the set-union of all left $\sigma$-expansive 
(respectively, right $\sigma$-expansive, expansive) lines 
passing through $(0,0)$ 
for $\varphi$ with $\{(0,0)\}$ subtracted. 
A lattice point $(i,j)$ in $\R^2$ is in $\E_L(\varphi)$ 
(respectively, in $\E_R(\varphi)$, in $\E(\varphi)$)
if and only if $\varphi^i\sigma^j$ is left $\sigma$-expansive 
(respectively, right $\sigma$-expansive, expansive). 
Each of $\E_L(\varphi)$, $\E_R(\varphi)$ and $\E(\varphi)$
is an open subset of $\R^2$. We call $\E_L(\varphi)$ 
(respectively, $\E_R(\varphi)$,  $\E(\varphi)$) the 
\itl{left $\sigma$-expansive} (respectively, 
\itl{right $\sigma$-expansive, expansive}) \itl{set 
for $\varphi$}. Further we call $\E_L(\varphi)\cup\E_R(\varphi)$ 
the \itl{onesided $\sigma$-expansive set} for $\varphi$. 
Define 
\begin{gather*}
\C^R_\varphi
=\Ll^+((-p_L(\varphi),\infty))\cup\Ll^-((-\infty,q_L(\varphi)))
\cup\hbar^+(0,0), \\
\C^L_\varphi 
=\Ll^+((-q_R(\varphi),\infty))\cup\Ll^-((-\infty,p_R(\varphi),))\cup
\hbar^+(0,0), \\
\C_\varphi=\C^R_\varphi\cap\C^L_\varphi, \\
\bar{\C}^L_\varphi
=\Ll^+((-\infty,p_R(\varphi)))\cup\Ll^-((-q_R(\varphi),\infty))\cup
\hbar^-(0,0), \\
\bar{\C}^R_\varphi 
=\Ll^+((-\infty,q_L(\varphi)))\cup\Ll^-((-p_L(\varphi),\infty))\cup
\hbar^-(0,0), \\
\bar{\C}_\varphi=\bar{\C}^L_\varphi\cap\bar{\C}^R_\varphi. 
\end{gather*}
Then it follows from Theorem 11.8(1),(2) that
$\C^R_\varphi$ and $\bar{\C}^R_\varphi$, 
(respectively, $\C^L_\varphi$ and 
$\bar{\C}^L_\varphi$) are the (connected) components of 
$\E_R(\varphi)$ (respectively, $\E_L(\varphi)$) that 
include the positive and negative horizontal 
half-lines, respectively. They are open cones in $\R^2$ 
and are the reverses of each other.
It also follows from Theorem 11.8(4)  
that $\C_\varphi$ and $\bar{\C}_\varphi$ are the 
components of $\E(\varphi)$ that 
include the positive and negative horizontal 
half-lines $\hbar^+(0,0)$ and $\hbar^-(0,0)$, respectively; 
they are open cones in $\R^2$ 
and are the reverses of each other.
The other components of 
$\E_R(\varphi)$ (respectively, $\E_L(\varphi)$)  
appear as pairs 
of open cones in $\R^2$ 
written in the form $\Ll^+(I^R)$ and $\Ll^-(I^R)$ 
(respectively, $\Ll^+(I^L)$ and $\Ll^-(I^L)$),
which are the reverses of each other,   
where $I^R$ (respectively, $I^L)$ is a component of $E_R(\varphi)$ 
(respectively, $E_L(\varphi)$) which is 
a bounded interval in $\R$. 
All components of 
$\E(\varphi)$ including $\C_\varphi$ and $\bar{\C}_\varphi$ above 
appear as pairs of open cones $\C$ and $\bar{\C}$ which are 
the reverses of each other and of the form 
$\C=\C^L\cap\C^R$ and $\bar{\C}=\bar{\C}^L\cap\bar{\C}^R$, 
where $\C^L$ and $\bar{\C}^L$ are the components of 
$\E_L(\varphi)$ with $\C^L\supset\C$ 
which are the reverses of each other, and 
$\C^R$ and $\bar{\C}^R$ are the components of 
$\E_R(\varphi)$ with $\C^R\supset\C$
which are the reverses of each other. 
(Since $\C^L\cap\C^R\supset \C$ and $\C^L\cap\C^R$ 
is a cone in $\R^2$ and hence connected, 
it follows from Proposition 11.1 that $\C^L\cap\C^R=\C$.)
We call each component of $\E_L(\varphi)$
(respectively, of $\E_R(\varphi)$, of $\E(\varphi)$, 
of $\E_R(\varphi)\cup\E_R(\varphi)$) 
a \itl{left $\sigma$-expansive component} 
(respectively, \itl{right $\sigma$-expansive component, 
expansive component, onesided $\sigma$-expansive component})
\itl{for $\varphi$}. 

An expansive component for an automorphism 
$\varphi$ of a subshift $(X,\sigma)$ (i.e. a component 
of $\E(\varphi)$)  
is an ``expansive component 
of $1$-frames'' for the 
$\Z^2$-action $\alpha:(i,j)\mapsto\varphi^i\sigma^j$ 
in the sense of 
Boyle and Lind \cite{BoyLin}. 
As is implied by \cite[Remark 9.6]{Nasu-t}  
in view of Theorem 12.2 (in the next section),      
for any lattice points $(i,j),(k,l)$ in $\R^2$, 
they belong to the same expansive component of 
1-frames for $\alpha$ if and only 
if $\varphi^i\sigma^j\lr^\circ\varphi^k\sigma^l$; or
alternatively, if $\C$ is any component of $\E(\varphi)$ 
and $(k,l)\in\C\cap\Z^2$ then $K$ being the subgroup of $H(X)$ 
generated by $\{\sigma,\varphi\}$,  
\[\C\cap\Z^2
=\{(i,j)\in \Z^2\,|\, 
\varphi^i\sigma^j\in C ^\circ_K(\varphi^k\sigma^l)\}.\]
A complete, self-contained proof of this will be given
in the proof of Theorem 11.12.

Let $K$ be the subgroup of $H(X)$ 
generated by $\{\sigma,\varphi\}$. Then 
by Theorems 9.7 and 9.10 and Proposition 9.12
(or by the equations appearing after Proof of 
Proposition 10.7) and the definitions of 
$\C_\varphi^R$, $\C_\varphi^L$ and $\C_\varphi$ above, 
we have 
\begin{gather} 
\C^R_\varphi\cap\Z^2= 
\{(i,j)\in\Z^2 \,|\, \varphi^i\sigma^j\in PL^\circ_K(\sigma)\},\\
\C^L_\varphi\cap\Z^2=
\{(i,j)\in\Z^2 \,|\, \varphi^i\sigma^j\in QR^\circ_K(\sigma)\},\\
\C_\varphi\cap\Z^2=
\{(i,j)\in\Z^2 \,|\, \varphi^i\sigma^j\in C^\circ_K(\sigma)\}. 
\end{gather}

Let $\bar{\ell}$ denote 
a line in the plane $\R^2$ passing through 
$(0,0)$ which may be horizontal, and for 
$(a,b)\in \R^2$ with $(a,b)\neq(0,0)$,  
let $\bar{\ell}(a,b)$ denote the line passing through $(0,0)$ and 
$(a,b)$ (which may be horizontal). 
For $(a,b)\in\R^2\setminus \bar{\ell}$, 
let $\R^2[\bar{\ell},(a,b))$ denote the 
closed half-plane whose boundary is $\bar{\ell}$ 
and whose interior contains $(a,b)$. 

Let $(k,l)\in\Z^2$ with $\varphi^k\sigma^l$ expansive. 
Let us define $\E_R(\varphi,(k,l))$ 
(respectively, $\E_L(\varphi,(k,l))$) 
to be the set-union of $\bar{\ell}(k,l)$ 
and all lines $\bar{\ell}$ (passing through $(0,0)$)
such that $\plane[\bar{\ell},(k,l))$ 
(respectively, $\plane[\bar{\ell},(-k,-l))$) 
codes $\plane$ for $\varphi$, 
with $\{(0,0)\}$ subtracted. 
We call $\E_L(\varphi,(k,l))$ 
(respectively, $\E_R(\varphi,(k,l))$, 
$\E_L(\varphi,(k,l))\cup\E_R(\varphi,(k,l))$) the 
\itl{left $\varphi^k\sigma^l$-expansive} (respectively, 
\itl{right $\varphi^k\sigma^l$-expansive, 
onesided $\varphi^k\sigma^l$-expansive}) \itl{set for $\varphi$}. 
This definition is consistent with the definition of the 
left $\sigma$-expansive (respectively,right $\sigma$-expansive, 
onesided $\sigma$-expansive) set for $\varphi$; in fact, 
\[\E_L(\varphi)=\E_L(\varphi,(0,1)),\q\q
\E_R(\varphi)=\E_R(\varphi,(0,1)).\] 

Let $E_L(\varphi,(k,l)$ and $E_R(\varphi,(k,l))$ 
be the subsets of $\R$ defined as follows:
if $l>0$ then $E_L(\varphi,(0,l))=E_L(\varphi)$ and
$E_R(\varphi,(0,l))=E_R(\varphi)$; if $l<0$ then
$E_L(\varphi,(0,l))=E_R(\varphi)$ and
$E_R(\varphi,(0,l))=E_L(\varphi)$; if $k\neq 0$ then 
\begin{align*}
E_L(\varphi,(k,l))=(E_L(\varphi)\cap(-\infty,l/k])\cup 
(E_R(\varphi)\cap[l/k,\infty)),\\
E_R(\varphi,(k,l))=(E_R(\varphi)\cap(-\infty,l/k])\cup 
(E_L(\varphi)\cap[l/k,\infty)). 
\end{align*} 
Then clearly we have
\begin{align*} 
E_L(\varphi,(k,l))\cap E_R(\varphi,(k,l))
&=E_L(\varphi)\cap E_R(\varphi)=E(\varphi),\\ 
E_L(\varphi,(k,l))\cup E_R(\varphi,(k,l))
&=E_L(\varphi)\cup E_R(\varphi). 
\end{align*}
Since $E_L(\varphi)$ and $E_R(\varphi)$ are open subsets of 
$\R$ and if $k\neq 0$ then 
$l/k\in E(\varphi)$, 
we see that 
$E_L(\varphi,(k,l))$ and $E_R(\varphi,(k,l))$ 
are open subsets of $\R$.

\begin{proposition}
Let $\varphi$ be an automorphism of 
an infinite subshift $(X,\sigma)$. 
Let $K$ be the subgroup of $H(X)$ 
generated by $\{\sigma,\varphi\}$. 
For any $(k,l)\in\Z^2$ 
with $\varphi^k\sigma^l$ expansive the following hold. 
\begin{enumerate}
\item
$\E_R(\varphi,(k,l))$ and $\E_L(\varphi,(k,l))$ are the open subsets  
of $\plane$ given by 
\begin{align}
\E_L(\varphi,(k,l))
=\Ll^+(E_L(\varphi,(k,l)))\cup\Ll^-(E_L(\varphi,(k,l)))
\cup\hbar^+(0,0)\cup\hbar^-(0,0),\\ 
\E_R(\varphi,(k,l))
=\Ll^+(E_R(\varphi,(k,l)))\cup\Ll^-(E_R(\varphi,(k,l)))
\cup\hbar^+(0,0)\cup\hbar^-(0,0),
\end{align}
and have the properties: 
\begin{align*} 
\E_L(\varphi,(k,l))\cap\E_R(\varphi,(k,l))
&=\E_L(\varphi)\cap\E_R(\varphi)=\E(\varphi),\\ 
\E_L(\varphi,(k,l))\cup\E_R(\varphi,(k,l))
&=\E_L(\varphi)\cup\E_R(\varphi). 
\end{align*}
\item 
$\{\varphi^i\sigma^j\,|\,(i,j)\in\E_L(\varphi,(k,l))\cap\Z^2\}$
is the set of all left $\varphi^k\sigma^l$-expansive elements in $K$, 
$\{\varphi^i\sigma^j\,|\,(i,j)\in\E_R(\varphi,(k,l))\cap\Z^2\}$ is 
the set of all right $\varphi^k\sigma^l$-expansive elements 
in $K$. 
\item 
$\{\varphi^i\sigma^j\,|
\,(i,j)\in(\E_L(\varphi)\cup\E_R(\varphi))\cap\Z^2\}$, 
i.e. the set of all
onesided $\sigma$-expansive elements in $K$, is the 
set of all
onesided $\varphi^k\sigma^l$-expansive elements in $K$. 
\end{enumerate}
\end{proposition} 
\begin{proof} 
(1) Let
$\A$ be the subset of $\R^2$ such that
if $k\neq 0$ or $l>0$ then 
$\A=\cH[(0,1),(k,l)]\cup(-\cH[(0,1),(k,l)])$ 
and if $k=0$ and $l<0$ then $\A=\R^2$. 
For any line $\bar{\ell}$ (passing through $(0,0)$) in 
$\R^2$ with $\bar{\ell}\neq\hbar(0,0)$ 
and $\bar{\ell}\neq\bar{\ell}(k,l)$, 
the following hold: 
if $\R^2[\bar{\ell},(k,l))=\R^2[\bar{\ell},(0,1))$, then 
$\plane[\bar{\ell},(k,l))$ (respectively, $\plane[\bar{\ell},(-k,-l))$  
codes $\plane$ for $\varphi$ if and only if 
$\bar{\ell}$ is right $\sigma$-expansive 
(respectively, left $\sigma$-expansive) for $\varphi$; 
if $\R^2[\bar{\ell},(k,l))\neq\R^2[\bar{\ell},(0,1))$ 
(i.e,, $\R^2[\bar{\ell},(k,l))=\R^2[\bar{\ell},(0,-1))$), then 
$\plane[\bar{\ell},(k,l))$ (respectively, $\plane[\bar{\ell},(-k,-l))$)
codes $\plane$ for $\varphi$ if and only if 
$\bar{\ell}$ is left $\sigma$-expansive 
(respectively, right $\sigma$-expansive) for $\varphi$. 
Each of $\E_L(\varphi)$ and $\E_R(\varphi)$ 
includes $\hbar(0,0)\setminus\{(0,0)\}$, and it 
also includes $\bar{\ell}(k,l)\setminus\{(0,0)\}$, because 
$\bar{\ell}(k,l)$ is both right $\sigma$-expansive 
and left $\sigma$-expansive for $\varphi$ 
(because $\varphi^k\sigma^l$ is expansive). 
Each of $\E_L(\varphi,(k,l))$ and $\E_R(\varphi,(k,l))$ 
includes $\bar{\ell}(k,l)\setminus\{(0,0)\}$, and it 
also includes $\hbar(0,0)\setminus\{(0,0)\}$, 
for if $k=0$ then $\hbar(0,0)=\bar{\ell}(k,l)$, 
and otherwise, each of $\R^2[\hbar(0,0),(k,l))$ and 
$\R^2[\hbar(0,0),(-k,-l))$ codes $\R^2$ for $\varphi$. 
Therefore, 
since $\A$ equals the set-union of $\hbar(0,0)$, $\bar{\ell}(k,l)$
and all lines $\bar{\ell}$ (passing through $(0,0)$) such that 
$\R^2[\bar{\ell},(0,1))\neq\R^2[\bar{\ell},(k,l))$, we have 
\begin{align*}
\E_L(\varphi,(k,l))
&=(\E_L(\varphi)\setminus\A)\cup(\E_R(\varphi)\cap\A), \\
\E_R(\varphi,(k,l))
&=(\E_R(\varphi)\setminus\A)\cup(\E_L(\varphi)\cap\A). 
\end{align*} 
From these we have the properties of $\E_L(\varphi,(k,l))$ 
and $\E_L(\varphi,(k,l))$ in (1), 
and noting that
$\A=\Ll^+([l/k,\infty))\cup\Ll^-([l/k,\infty))\cup\hbar(0,0)$ 
if $k\neq 0$, 
we derive (11.4) and (11.5) using the definitons of 
$\E_L(\varphi)$ and $\E_R(\varphi)$. 
Since $E_L(\varphi,(k,l))$ and $E_R(\varphi,(k,l))$ are 
open subsets of $\R$, it follows 
from (11.4) and (11.5) 
that $\E_L(\varphi,(k,l))$ and $\E_R(\varphi,(k,l))$ are 
open subsets of $\R^2$.

(2) Let $(i,j)\in\Z^2/\{(0,0)\}$. 
By definition, $(i,j)$ satisfies the condition 
that $\varphi^i\sigma^j$ is 
left $\varphi^k\sigma^l$-expansive  
if and only if it satisfies 
the condition that 
there exists $\delta>0$  such that 
for any $x,y\in X$ it holds that 
if for all $m\in\Z, n\leq 0$, 
\[d_X((\varphi^i\sigma^j)^m(\varphi^k\sigma^l)^n(x), 
(\varphi^i\sigma^j)^m(\varphi^k\sigma^l)^n(y))\leq\delta,\] then $x=y$. 
This condition is equivalent to the condition that 
either $\bar{\ell}(i,j)=\bar{\ell}(k,l)$ or there exists $\delta>0$ 
such that 
$\R^2[\bar{\ell}(i,j),(-k,-l))$ $\delta$-codes $\R^2$ for 
$\varphi$. This condition  
is equivalent to the condition that 
either $\bar{\ell}(i,j)=\bar{\ell}(k,l)$ or there exists an integer  
$t\geq 0$ such that 
$\cup_{-t\leq s\leq t} \R^2[\bar{\ell}(i,j),(-k,-l))+(0,s)$ 
codes $\R^2$ for $\varphi$, 
which is equivalent to the condition that 
either $\bar{\ell}(i,j)=\bar{\ell}(k,l)$ or 
$\R^2[\bar{\ell}(i,j),(-k,-l))$ codes $\R^2$ for $\varphi$. 
Since this condition is equivalent to the conditin that 
$(i,j)\in\E_L(\varphi,(k,l))$, we see that 
$\{\varphi^i\sigma^j\,|\,(i,j)\in\E_L(\varphi,(k,l))\cap\Z^2\}$ is
the set of all left $\varphi^k\sigma^l$-expansive elements 
in $K$. 

By symmetry we see that
$\{\varphi^i\sigma^j\,|\,(i,j)\in\E_R(\varphi,(k,l))\cap\Z^2\}$ is 
the set of all right $\varphi^k\sigma^l$-expansive elements in $K$. 

(3) By (2) and the last equation in (1), we know that
the set of all onesided $\varphi^k\sigma^l$-expansive  
elements in $K$ equals 
$\{\varphi^i\sigma^j\,|\,
(i,j)\in(\E_L(\varphi)\cup\E_R(\varphi))\cap\Z^2\}$. 
\end{proof} 

By definition (see Section 2),
an automorphism $\varphi$ of a subshift $(X,\sigma)$ is 
onesided-expansive if $\varphi$ is 
onesided $\varphi^k\sigma^l$-expansive for some $(k,l)\in\Z^2$ with   
$\varphi^k\sigma^l$ expansive.  Proposition 11.10(3) shows that if 
an automorphism $\varphi$ of a subshift $(X,\sigma)$ is 
onesided-expansive, then  $\varphi$ is onesided 
$\varphi^k\sigma^l$-expansive for all $(k,l)\in \Z^2$ with 
$\varphi^k\sigma^l$ expansive. 
More generally, as stated in Subection 2.1, 
if an automorphism $\varphi$ of an invertible 
dynamical system $(X,\tau)$ is 
onesided-expansive, then $\varphi$ is onesided 
$\varphi^k\tau^l$-expansive for all $(k,l)\in \Z^2$ with 
$\varphi^k\tau^l$ expansive (see Remark 11.23(3)).

\begin{proposition} 
Let $\varphi$ be an automorphism of 
an infinite subshift $(X,\sigma)$.  
Let $(k,l)\in\Z^2$ 
with  $\varphi^k\sigma^l$ expansive. 
Let $\C^L$ and $\C^R$ be 
the components of $\E_L(\varphi)$ and 
$\E_R(\varphi)$, respectively, with  
$\C^L\cap\C^R\ni(k,l)$, 
let $\Hat{\C}^L$ and $\Hat{\C}^R$ be 
the components of $\E_L(\varphi,(k,l))$ and 
$\E_R(\varphi,(k,l))$, respectively, with  
$\Hat{\C}^L\cap\Hat{\C}^R\ni(k,l)$ 
and let $\C$ be the component of $\E(\varphi)$ 
with $\C\ni (k,l)$. 
Then 
\begin{align}\Hat{\C}^L\cap\Hat{\C}^R&=\C^L\cap\C^R=\C,\\
\Hat{\C}^L\cup\Hat{\C}^R&=\C^L\cup\C^R 
\end{align}
\end{proposition} 
\begin{proof} 
Since 
$\E_L(\varphi,(k,l))\cap\E_R(\varphi,(k,l))=\E(\varphi)$ 
(by Proposition 11.10(1)), we see that 
$\C$ is a connected subset of $\E_L(\varphi,(k,l))$ 
containing $(k,l)$. Hence $\C$ is included in 
the component $\hat{\C}^L$ of  
$\E_L(\varphi,(k,l))$ containing  $(k,l)$. Similarly we see that 
$\C\subset\hat{\C}^R$, and hence $\hat{\C}^L\cap\hat{\C}^R\supset\C$. 
On the other hand, 
$\hat{\C}^L\cap\hat{\C}^R\subset
\E_L(\varphi,(k,l))\cap\E_R(\varphi,(k,l))=\E(\varphi)$. 
Since $\Hat{\C}^L$ and $\Hat{\C}^R$ are cones  with 
apex $(0,0)$ in $\R^2$ (by (11.4) and (11.5)),
$\hat{\C}^L\cap\hat{\C}^R$ is a cone and hence connected. 
Since $\hat{\C}^L\cap\hat{\C}^R$ 
is a connected subset of 
$\E(\varphi)$ containing $(k,l)$, we have 
$\hat{\C}^L\cap\hat{\C}^R\subset\C$. Therefore   
$\hat{\C}^L\cap\hat{\C}^R=\C$, and (11.6) is proved. 

If $k=0$, then (11.7) easily follows. 
Assume that $k>0$. 
Let $I^L$ and $I^R$ be the components 
of $E_L(\varphi)$ and $E_R(\varphi)$, respectively, with 
$I^L\cap I^R\ni l/k$. Let 
$\Hat{I}^L=(I^L\cap(-\infty,l/k])\cup(I^R\cap[l/k,\infty))$ 
and 
$\Hat{I}^R=(I^R\cap(-\infty,l/k])\cup(I^L\cap[l/k,\infty))$. 
Then 
$\Hat{I}^L$ and $\Hat{I}^R$ are the componens of $E_L(\varphi,(k,l))$ 
and $E_R(\varphi,(k,l))$, respectively, 
with $\Hat{I}^L\cap\Hat{I}^R\ni (k,l)$, and 
$\Hat{I}^L\cup\Hat{I}^R=I^L\cup I^R$. Using these 
and notihg that $d_L(\varphi)\leq d_R(\varphi)$ (Proposition 9.12(2)) 
and that $d_R(\varphi)\neq l/k$ and $d_L(\varphi)\neq l/k$ (by 
Theorem 11.8(3)),  
we see the following by the observation of (11.4) and (11.5). 
If $d_L(\varphi)<l/k<d_R(\varphi)$ then $I^L\cup I^R$ is a 
bounded interval and hence
$\Hat{\C}^L\cup\Hat{\C}^R=\Ll^+(I^L\cup I^R)=\C^L\cup\C^R$.
If $l/k>d_R(\varphi)$ then $I^L\cup I^R$ is 
a left-bounded, right-unbounded interval and hence  
$\Hat{\C}^L\cup\Hat{\C}^R=
\Ll^+(I^L\cup I^R)\cup\Ll^-((-\infty,d_L(\varphi)))\cup\hbar^+(0,0)
=\C^L\cup\C^R$. 
If $l/k< d_L(\varphi)$ then 
$I^L\cup I^R$ is 
a left-unbounded, right-bounded interval and hence  
$\Hat{\C}^L\cup\Hat{\C}^R=
\Ll^+(I^L\cup I^R)\cup\Ll^-((d_R(\varphi),\infty)\cup\hbar^-(0,0)
=\C^L\cup\C^R$. Hence (11.7) is proved in case $k>0$. 
From this it 
follows that if $k<0$ then $(-\Hat{\C}^L)\cup(-\Hat{\C}^R)=
(-\C^L)\cup(-\C^R)$, because by (11.4) and (11.5) $-\Hat{\C}^L$ and 
$-\Hat{\C}^R$ are also components of $\E_L(\varphi,(k,l))$
and $\E_R(\varphi,(k,l))$, respectively. 
Hence $(11.7)$ is proved in case $k<0$. 
\end{proof}

\begin{theorem} 
Let $\varphi$ be an automorphism of 
an infinite subshift $(X,\sigma)$. 
Let $K$ be the subgroup of $H(X)$ 
generated by $\{\sigma,\varphi\}$. For any $F\subset K$, 
let $\Lambda(F)$ denote the set 
$\{(i,j)\in\Z^2 \,|\, \varphi^i\sigma^j\in F\}$. 
Let $(k,l)\in\Z^2$ 
with  $\varphi^k\sigma^l$ expansive. 
Let $\C^L$ and $\C^R$ be 
the components of $\E_L(\varphi)$ and 
$\E_R(\varphi)$, respectively, with  
$\C^L\cap\C^R\ni(k,l)$, 
let $\Hat{\C}^L$ and $\Hat{\C}^R$ be 
the components of $\E_L(\varphi,(k,l))$ and 
$\E_R(\varphi,(k,l))$, respectively, with  
$\Hat{\C}^L\cap\Hat{\C}^R\ni(k,l)$ 
and let $\C$ be the component of $\E(\varphi)$ 
with $\C\ni (k,l)$. 
Then the following hold:
\begin{gather}
\Hat{\C}^R\cap\Z^2=\Lambda(PL^\circ_K(\varphi^k\sigma^l)),\q\q 
\Hat{\C}^L\cap\Z^2=\Lambda(QR^\circ_K(\varphi^k\sigma^l))\\
\C\cap\Z^2=\Lambda(C^\circ_K(\varphi^k\sigma^l))\\
(\C^L\cup\C^R)\cap\Z^2=\Lambda(D^\circ_K(\varphi^k\sigma^l)). 
\end{gather} 
\end{theorem} 
\begin{proof} 
Since  
$D_K^\circ(\varphi^k\sigma^l)
=PL^\circ_K(\varphi^k\sigma^l)\cup QR^\circ_K(\varphi^k\sigma^l)$, 
(11.10) follows from (11.8) and Proposition 11.11((11.7)). 
Hence it suffices to prove (11.8) and (11.9). 
(Since $C^\circ(\varphi^k\sigma^l)
=PL^\circ_K(\varphi^k\sigma^l)\cap QR^\circ_K(\varphi^k\sigma^l)$
(by \cite[Proposition 8.1]{Nasu-te}), 
(11.9) follows from (11.8) and Proposition 11.11((11.6)).
However, we shall prove (11.9) directly together with (11.8) below.)

Suppose that $k=0$ and $l>0$. 
Then $\E_L(\varphi,(0,l))=\E_L(\varphi)$ 
and $\E_R(\varphi,(0,l))=\E_R(\varphi)$, and hence 
$\Hat{\C}^L=\C^L=\C_\varphi^L$, $\Hat{\C}^R=\C^R=\C_\varphi^R$ 
and $\C=\C_\varphi$. 
Using (11.1) and Proposition 12.1 (in the next section), we see that
$\Hat{\C}^R\cap\Z^2
=\C^R_\varphi\cap\Z^2=\Lambda(PL^\circ_K(\sigma))=
\Lambda(PL^\circ_K(\sigma^l))$. Similarly using (11.2), (11.3) 
and Proposition 12.1 we have $\Hat{\C}^L\cap\Z^2=
\Lambda(QR^\circ_K(\sigma^l))$ and (11.9). 

Suppose that $k=0$ and $l<0$. 
Then $\E_L(\varphi,(0,l))=\E_R(\varphi)$ 
and $\E_R(\varphi,(0,l))=\E_L(\varphi)$, and hence 
$\Hat{\C}^L=\C^R=\bar{\C}_\varphi^R$, 
$\Hat{\C}^R=\C^L=\bar{\C}_\varphi^L$, and $\C=\bar{\C}_\varphi$. 
Using (11.1) and Proposition 12.1, we see that
$\Hat{\C}^L\cap\Z^2
=\bar{\C}_\varphi^R\cap\Z^2=\Lambda(PL^\circ_K(\sigma)^{-1})
=\Lambda(QR_K^\circ(\sigma^{-1}))=\Lambda(QR_K^\circ(\sigma^l))$. 
The remainder of (11.8), and (11.9) are similarly proved. 

Hereafter in this proof we suppose that $k\neq 0$. 
Since $\varphi^k\sigma^l$ is expansive
and $\sigma$ is an automorphism of the dynamical 
system  $(X,\varphi^k\sigma^l)$, 
there exist a subshift $(X',\sigma')$, an automorphism $\varphi'$
of $(X',\sigma')$ and a conjugacy 
$\theta: (X,\varphi^k\sigma^l,\sigma) \to (X',\sigma',\varphi')$
between commuting systems. 
Let $\pplane$ be a copy  
of the plane $\plane$. Let $(\Z^2)'$ be the set of lattice points 
in $\pplane$. 
Let $\pi:\plane\to \pplane$ be the linear map 
such that $\pi((a,b))=(-(l/k)a+b,(1/k)a)$ for $(a,b)\in\plane$; 
it is invertible 
with $\pi^{-1}((a',b'))=(kb',a'+lb')$ for $(a',b')\in \pplane$. 

Let $\bar{\ell}$ be any line (passing through $(0,0)$) in $\plane$ 
with $\bar{\ell}\neq \bar{\ell}(k,l)$. 
Since $\pi((k,l))=(0,1)$ and 
$(\R^2)'[\pi(\bar{\ell}),(0,1))=(\R^2)'[\pi(\bar{\ell}),\cdot)$, 
the following conditions are seen, 
in the order of their appearance, to be
pairwise equivalent: 

(a) there exists $\delta>0$ such that for any $x,y\in X$ it holds that if
\[d_{X'}((\varphi')^{i'}(\sigma')^{j'}(\theta(x)),
(\varphi')^{i'}(\sigma')^{j'}(\theta(y)))\leq\delta 
\q\forall (i',j')\in (\R^2)'[\pi(\bar{\ell}),\cdot)\cap(\Z^2)'\] 
then $\theta(x)=\theta(y)$; 

(b) there exists $\delta>0$ such that for any $x,y\in X$ it holds that if
\[d_{X'}(\theta\varphi^i\sigma^j(x),\theta\varphi^i\sigma^j(y))\leq\delta 
\q\forall (i,j)\in \plane[\bar{\ell},(k,l))\cap\pi^{-1}((\Z^2)')\]
then $x=y$; 

(c) there exists $\delta>0$ such that for any $x,y\in X$ it holds that if
\[  
d_{X}(\varphi^i\sigma^j(x),\varphi^i\sigma^j(y))\leq\delta
\q\forall (i,j)\in \plane[\bar{\ell},(k,l))\cap\pi^{-1}((\Z^2)')\]
then $x=y$; 

(d) there exists $\delta>0$ such that for any $x,y\in X$ it holds that if
\[d_{X}(\varphi^i\sigma^j(x),\varphi^i\sigma^j(y))\leq\delta
\q\forall (i,j)\in \plane[\bar{\ell},(k,l))\cap\Z^2\]
then $x=y$. 

\noindent(For by putting $(i,j)=\pi^{-1}((i',j'))$
the equivalence of (a) and (b) is seen; 
the equivalence of (c) and (d) follows because 
$\varphi$ and $\sigma$ are homeomorphisms.) 
Using the equivalence of (a) and (d) we know that 
$\pi(\bar{\ell})$ is right 
$\sigma'$-expansive for the automorphism $\varphi'$ of $(X',\sigma')$ 
if and only if 
either $\bar{\ell}=\bar{\ell}(k,l)$ or $\plane[\bar{\ell},(k,l))$ 
codes $\plane$ for the automorphism $\varphi$ of $(X,\sigma)$. 
Therefore we have 
\[\E_R(\varphi')=\{\pi({\bar{\ell}})\,|\,\bar{\ell}
\subset\E_R(\varphi,(k,l))\cup\{(0,0)\}\}\setminus\{(0,0)\}
=\pi(\E_R(\varphi,(k,l))).\]
Since $\pi((-k,-l))=(0,-1)$ and 
$(\R^2)'[\pi(\bar{\ell}),(0,-1))=(\R^2)'(\cdot,\pi(\bar{\ell})]$, 
a similar argument shows that 
$\pi(\bar{\ell})$ is left
$\sigma'$-expansive for $\varphi'$
if and only if either $\bar{\ell}=\bar{\ell}(k,l)$ or 
$\plane[\bar{\ell},(-k,-l))$
codes $\plane$ for $\varphi$. Therefore we have
\[\E_L(\varphi')=\{\pi({\bar{\ell}})\,|\,\bar{\ell}\subset
\E_L(\varphi,(k,l))\cup\{(0,0)\}\}\setminus\{(0,0)\}
=\pi(\E_L(\varphi,(k,l))).\]
Since $\E(\varphi)=\E_L(\varphi,(k,l))\cap\E_R(\varphi,(k,l))$ 
(by Proposition 11.10), we have 
\[\pi(\E(\varphi))=\E_L(\varphi')\cap\E_R(\varphi')=\E(\varphi').\]
Since $\Hat{\C}^L$, $\Hat{\C}^R$ and $\C$ are components of 
$\E_L(\varphi,(k,l))$, $\E_R(\varphi,(k,l))$ and 
$\E(\varphi)$, respectively, with   
$\Hat{\C}^L\cap\Hat{\C}^R=\C\ni (k,l)$ (by Proposition 11.11) and 
$\pi$ is a homeomorphism of $\plane$ onto $\pplane$ with
$\pi(k,l))=(0,1)$, 
it follows that 
$\pi(\Hat{\C}^L)$, $\pi(\Hat{\C}^R)$ 
and $\pi(\C)$ 
are components of 
$\E_L(\varphi')$, $\E_R(\varphi')$ 
and $\E(\varphi')$, 
respectively, such that    
$\pi(\Hat{\C}^L)\cap\pi(\Hat{\C}^R)=\pi(\C)$ and 
this includes positive horizontal half line in $\pplane$. 
Therefore, we have
$\pi(\Hat{\C}^R)=\C^R_{\varphi'}$, $\pi(\Hat{\C}^L)=\C^L_{\varphi'}$
and $\pi(\C)=\C_{\varphi'}$. Hence, if we let $K'$ be the 
subgroup generated by $\{\sigma', \varphi'\}$, 
then by (11.1),(11.2),(11.3) we have 
\begin{gather} \pi(\Hat{\C}^R)\cap\Z^2=
\{(i',j')\in\Z^2 \,|\,(\varphi')^{i'}(\sigma')^{j'}\in
PL_{K'}^\circ(\sigma')\},\\
\pi(\Hat{\C}^L)\cap\Z^2=
\{(i',j')\in\Z^2 \,|\,(\varphi')^{i'}(\sigma')^{j'}\in
QR_{K'}^\circ(\sigma')\}, \\
\pi(\C)\cap\Z^2=
\{(i',j')\in\Z^2 \,|\,(\varphi')^{i'}(\sigma')^{j'}\in
C_{K'}^\circ(\sigma')\}. 
\end{gather}

Suppose that $\varphi^i\sigma^j$ is an 
essentially weakly $p$-L automorphism of $(X,\varphi^k\sigma^l)$ and 
right $\varphi^k\sigma^l$-expansive. 
Then so is $(\varphi^i\sigma^j)^{|k|}$. 
Hence $\theta\varphi^{|k|i}\sigma^{|k|j}$ is an essentially weakly 
$p$-L automorphism of $(\theta(X),\theta\varphi^k\sigma^l)
=(X',\sigma')$ and right $\sigma'$-expansive. 
If we put $\epsilon=|k|/k$ then  
\[
\theta\varphi^{|k|i}\sigma^{|k|j}=\theta(\varphi^{ki}\sigma^{kj})^\epsilon
=\theta((\varphi^k\sigma^l)^i\sigma^{-li+kj})^\epsilon
=(\sigma')^{\epsilon i}(\varphi')^{\epsilon(-li+kj)}, 
\] 
and hence $(\sigma')^{\epsilon i}(\varphi')^{\epsilon(-li+kj)}$ 
is an essentially weakly 
$p$-L automorphism of $(X',\sigma')$ and right $\sigma'$-expansive.
Therefore by (11.11) we see that 
$(\epsilon(-li+kj),\epsilon i)\in \pi(\Hat{\C}^R)$. From this 
it follows that $|k|(i,j)\in \Hat{\C}^R$, because 
\[|k|(i,j)=\pi^{-1}(|k|\pi((i,j)))
=\pi^{-1}(|k|(-(l/k)i+j, (1/k)i))
=\pi^{-1}((\epsilon(-li+kj),\epsilon i)).\]
Therefore 
$(i,j)\in \Hat{\C}^R$, because 
$\Hat{\C}^R$ is a cone in $\R^2$ with apex $(0,0)$. 

We have proved that 
$\Hat{\C}^R\cap\Z^2\supset \Lambda(PL^\circ_K(\varphi^k\sigma^l))$. 
To prove the converse inclusion, suppose that 
$(i,j)\in \Hat{\C}^R\cap\Z^2$. Then $(|k|i,|k|j)\in\Hat{\C}^R$. 
Let $(i',j')=\pi((|k|i,|k|j))$. Then by (11.11) 
$(\varphi')^{i'}(\sigma')^{j'}$ is an essentially 
weakly $p$-L automorphism of $(X',\sigma')$ and 
right $\sigma'$-expansive. Hence 
$\theta^{-1}(\varphi')^{i'}(\sigma')^{j'}$ is an essentially 
weakly $p$-L automorphism of $(\theta^{-1}(X'),\theta^{-1}\sigma')
=(X,\varphi^k\sigma^l)$ and right $\varphi^k\sigma^l$-expansive. 
Since $\theta^{-1}(\varphi')^{i'}(\sigma')^{j'}
=\sigma^{i'}(\varphi^k\sigma^l)^{j'}=(\varphi^i\sigma^j)^{|k|}$ 
(because $(kj',i'+lj')=\pi^{-1}((i',j'))=(|k|i,|k|j)$), 
$\varphi^i\sigma^j$ is essentially weakly $p$-L automorphism 
of $(X,\varphi^k\sigma^l)$ (by Theorem 8.1) and right 
$\varphi^k\sigma^l$-expansive. Hence we have proved that 
$\Hat{\C}^R\cap\Z^2= \Lambda(PL^\circ_K(\varphi^k\sigma^l))$. 

The proof that $\Hat{\C}^L\cap\Z^2= 
\Lambda(QR^\circ_K(\varphi^k\sigma^l))$ and that of (11.9) are
similarly given by using 
(11.12) and (11.13), respectively, and Theorem 8.1. 
\end{proof}

Let $\varphi$ be an automorphism of 
an infinite subshift $(X,\sigma)$. 
As is known by \cite[p. 71]{BoyLin}, 
if $\C$ is an expansive component for $\varphi$, i.e.  
a component of $\E(\varphi)$, then so is its reverse 
$-\C$ and furthermore, $\C\cap(-\C)=\emptyset$; 
hence there are  at least two 
expansive components for $\varphi$, which 
are open cones with apex $(0,0)$ 
in $\R^2$. 

In view of Proposition 11.10(3), 
$\E_L(\varphi)\cup\E_R(\varphi)$ can be called the 
\itl{onesided-expansive set for $\varphi$} and
a component of 
$\E_L(\varphi)\cup\E_R(\varphi)$ can be called 
a \itl{onesided-expansive component for $\varphi$}. 
Since $\E_L(\varphi)\cup\E_R(\varphi)=
\Ll^+(E_L(\varphi)\cup E_R(\varphi))\cup
\Ll^-(E_L(\varphi)\cup E_R(\varphi))\cup
\hbar^+(0,0)\cup\hbar^-(0,0)$, we see that
if $\cD$ is a onesided-expansive component for $\varphi$
then so is its reverse $-\cD$. 
Therefore, a onesided-expansive component for $\varphi$ is either
an open cone with apex $(0,0)$ 
in $\R^2$ with $\cD\cap(-\cD)=\emptyset$ or equals 
$\R^2/(0,0)$ (hence is unique).  
An example of an automorphism of a full-shift having a  
unique onesided-expansive component for it   
will be given in Example 11.17. 

We do not know whether or not there can be a onesided-expansive 
component for $\varphi$ containing no lattice point $(k,l)$ 
with $\varphi^k\sigma^l$ expansive. 

Let $\varphi$ be an automorphism of 
an infinite subshift $(X,\sigma)$. Let $K$ be the subgroup of 
$H(X)$ generated by $\{\sigma,\varphi\}$. 
The statement (1) of the following theorem is 
a restatement of Theorem 11.12((11.9)) 
(\cite[Remark 9.6]{Nasu-te} in view of Theorem 12.2), 
which asserts that the expansive components for 
$\varphi$ and 
the interiors of essentially-weakly-LR cones in $K$ 
are the same in a natural sense. We shall furthermore 
prove the statement (2) which asserts that 
the onesided-expansive components for $\varphi$ 
containing expansive elements and 
the interiors of extended districts in $K$ 
are the same in the same natural sense as above.  

\begin{theorem} Let $\varphi$ be an automorphism of 
an infinite subshift $(X,\sigma)$. Let 
$K$ be the subgroup of $H(X)$ generated 
by $\{\sigma,\varphi\}$. For any subset $F$ of $K$, 
let $\Lambda(F)$ denote the set 
$\{(i,j)\in\Z^2 \,|\, \varphi^i\sigma^j\in F\}$. Then 
the following are valid. 
\begin{enumerate}
\item 
A subset of $\Z^2$ is given as 
$\C\cap\Z^2$, where $\C$ is a (connected) component of   
$\E(\varphi)=\E_R(\varphi)\cap\E_L(\varphi)$,  
if and only if it is given as  
$\Lambda(C^\circ)$, where $C^\circ$ is the interior of  
an essentially-weakly-LR cone $C$ in $K$. 
\item 
A subset of $\Z^2$ is given as
$\cD\cap\Z^2$, where $\cD$ is a component  
of $\E_R(\varphi)\cup\E_L(\varphi)$ such that  
$\cD\cap\E(\varphi)\neq\emptyset$,  
if and only if it is given as $\Lambda((D^\ast)^\circ)$, where 
$(D^\ast)^\circ$ is the interior of an extended district 
$D^\ast$ in $K$. 
\end{enumerate}
\end{theorem}
\begin{proof}
(1) By Theorem 11.12((11.9)). 

(2)  First we shall prove the following.  

(i) For any $(k,l)\in \Z^2$ with $\varphi^k\sigma^l$ expansive, there 
exists a component $\cD$ of $\E_L(\varphi)\cup\E_R(\varphi)$ 
such that $\cD\supset\Lambda((D^\ast_K)^\circ(\varphi^k\sigma^l))$.

To prove this 
suppose that $(m,n)\in \Lambda((D^\ast_K)^\circ(\varphi^k\sigma^l))$. 
Then there exist $r\geq 1$ and $(k_i,l_i)\in\Z^2, i=0,\dots,r$, 
such that  $(k_0,l_0)=(m,n)$,
$(k_r,l_r)=(k,l)$, 
$\varphi^{k_i}\sigma^{l_i}$ is expansive for $i=1,\dots,r$
and 
\[\varphi^{k_0}\sigma^{l_0}\mm^\circ
\varphi^{k_1}\sigma^{l_1}\mm^\circ
\dots\mm^\circ\varphi^{k_r}\sigma^{l_r}.\] 
Since $\varphi^{k_{i-1}}\sigma^{l_{i-1}}
\in D_K^\circ(\varphi^{k_i}\sigma^{l_i})$ for $i=1,\dots,r$, 
it follows from Theorem 11.12((11.10)) that 
there exist open cones $\C_i, i=1,\dots,r$, in $\R^2$ 
each of which is 
a component of $\E_L(\varphi)$ or of $\E_R(\varphi)$ such that 
$(k_{i-1},l_{i-1}),(k_i,l_i)\in\C_i$ 
for $i=1,\dots,r$. Clearly $(m,n)=(k_0,l_0)\in\C_1$, 
$(k,l)=(k_r,l_r)\in\C_r$ and
$\C_{i-1}\cap C_i\neq\emptyset$ for $i=2,\dots, r$ if $r\geq 2$. 
Therefore, since all $\C_i$'s are 
connected subsets of $\R^2$ included in $\E_L(\varphi)\cup\E_R(\varphi)$, 
it follows that $(m,n)$ is contained in the component $\cD$
of $\E_L(\varphi)\cup\E_R(\varphi)$ 
with $(k,l)\in\cD$. Hence (i) is proved. 

Next we shall prove:

(ii) If $\cD$ is a compornent of $\E_L(\varphi)\cup\E_R(\varphi)$, 
then $\cD\cap\Z^2\subset\Lambda(D^\ast_K)^\circ(\varphi^k\sigma^l))$ 
for any $(k,l)\in \cD\cap\Z^2$ with $\varphi^k\sigma^l$ expansive. 

Let $s:\R^2\setminus\{(0,0)\}\to S^1$ be the mapping such that 
$s(a,b)=(a,b)/\sqrt{a^2+b^2}$ for $(a,b)\in\R^2\setminus\{(0,0)\}$. 
For any $R\subset\R^2\setminus\{(0,0)\}$, let $sR$ denote the subset 
$\{s(a,b)\,|\,(a,b)\in R\}$ of $S^1$. 

Let $\cD$ be a component of $\E_L(\varphi)\cup\E_R(\varphi)$. 
Then $\cD$ is an open cone in $\R^2$ with apex $(0,0)$ or 
equals $\R^2\setminus\{(0,0)\}$. 
Let $(k,l),(m,n)\in\cD\cap\Z^2$ with $\varphi^k\sigma^l$ expansive. 
In case $kn=lm$ and $km>0$,
using Remark 2.15 and Theorem 8.1 we see 
that $\varphi^m\sigma^n\lr\varphi^k\sigma^l$ and hence (ii) is 
proved. We must further consider 
the case (a) that $kn=lm$ and $km<0$ 
and the case (b) that $kn\neq lm$. Note that the case (a) happens
only when $\cD=\R^2\setminus\{(0,0)\}$ and that the case (b) occurs 
either when $\cD=\R^2\setminus\{(0,0)\}$  or when 
$\cD$ is an open cone in $\R^2$. 
Let $F$ be the subset of $S^1$ defined as follows: for the case (a)
$F$ is the closed half circle with endpoints $s(k,l)$ and $s(m,n)$, and 
for the case (b) $F=s(\cH[(k,l),(m,n)]\setminus\{(0,0)\})$. Then 
$F$ is a closed arc between $s(k,l)$ and $s(m,n)$ in $S^1$. 
Since $s\cD$ is a connected set in $S^1$ 
with $s(k,l), s(m,n)\in\cD$, 
$F\subset s\cD$ and hence $F\subset s\E_L(\varphi)\cup s\E_R(\varphi)$. 
Since $s\E_L(\varphi)$ and $s\E_R(\varphi)$ 
are open sets in $S^1$, each component of them is an open arc in $S^1$. 
Therefore, there exists an open cover for  
$F$ such that each set in it 
is a component of $s\E_L(\varphi)$ or of $s\E_R(\varphi)$. 
Since $F$ is compact, there exists a finite subcover
$\mathcal{K}=\{J_0,\dots,J_r\}$ of it with $r\geq 0$. 
Deleting every set that is included 
by some other set in $\mathcal{K}$
from $\mathcal{K}$ if any, 
we may assume that $\mathcal{K}$ is 
a minimal cover for $F$, i.e., 
if any one set $J_i$ is deleted from $\mathcal{K}$ then 
the collection of the remaining sets cannot cover $F$. 
Therefore, we may assume that $\mathcal{K}$ 
has the following property: $J_0\ni s(m,n)$ 
and $J_r\ni s(k,l)$; if $r\geq 1$ then 
for all $i=1,\dots,r$, $J_{i-1}\cap J_i\neq\emptyset$ but 
neither $J_i$ is a subset of $J_{i-1}$ 
nor $J_{i-1}$ is a subset of $J_i$; 
if $J_{i-1}$ is a component of $s\E_L(\varphi)$ 
then $J_i$ is a component of $s\E_R(\varphi)$, if $J_{i-1}$ 
is a component of $s\E_R(\varphi)$ then $J_i$ is 
a component of $s\E_L(\varphi)$, and hence 
$J_{i-1}\cap J_i$ is a component of of $s\E(\varphi)$. 

For any open arc $J$ in $S^1$, 
let $s^{-1}J$ denote the set-union 
of all open half lines in $\R^2$ with endpoint 
$(0,0)$ intersecting with $J$. 
Let $\mathcal{J}_i=s^{-1}J_i$ for $i=0,\dots,r$. 
If $r=0$, then $(k,l),(m,n)\in\mathcal{J}_0$. 
Since $\mathcal{J}_0$ is 
a component of $\E_L(\varphi)$ or 
of $\E_R(\varphi)$, 
it follows from Theorem 11.12 ((11.10)) that 
(ii) is valid when $r=0$. 
Therefore we assume that $r\geq 1$. 
Let $(k_i,l_i)$ be any lattice point in 
$\mathcal{J}_{i-1}\cap\mathcal{J}_i$
for $i=1,\dots,r$. 
Then, since 
$\mathcal{J}_i$ is a component of $\E_L(\varphi)$ or 
of $\E_R(\varphi)$ with 
$\mathcal{J}_{i-1}\cap\mathcal{J}_i\in \E(\varphi)$ for 
$i=1,\dots, r$, we see that 
$\varphi^{k_i}\sigma^{l_i}$ is expansive for all $i=1,\dots,r$ and 
it follows from Theorem 11.12 ((11.10)) that 
$\varphi^{m}\sigma^{n}\mm^\circ\varphi^{k_1}\sigma^{l_1}
\mm^\circ\dots\mm^\circ\varphi^{k_r}\sigma^{l_r}\mm^\circ 
\varphi^k\sigma^l$. 
Therefore, (ii) follows. 

Hence the theorem is proved by (i) and (ii). 
\end{proof} 
The following corollary is a continuation of Proposition 10.8. 
\begin{corollary} Let $\varphi$ be an automorphism of 
an infinite subshift $(X,\sigma)$. Let $K$ be the subgroup of 
$H(X)$ generated by $\{\sigma,\varphi\}$. Let $(k,l)\in\Z^2$ 
with $\varphi^k\sigma^l$ expansive. 
Let $(i_1,j_1), (i_2,j_2)\in\Z^2$ with $i_1j_2\neq i_2j_1$. 
\begin{enumerate}
\item 
If $\varphi^{i_t}\sigma^{j_t}\lr\varphi^k\sigma^l$ with 
$\varphi^{i_t}\sigma^{j_t}$ nonexpansive for $t=1,2$ 
and $\C$ is the component of 
$\E(\varphi)=\E_L(\varphi)\cap\E_R(\varphi)$ with 
$(k,l)\in\C$, then 
\begin{align*}
\C&=\cH((i_1,j_1),(i_2,j_2)), \\
C_K^\circ(\varphi^k\sigma^l)
&=\{\varphi^i\sigma^j\,|\,(i,j)\in\cH((i_1,j_1),(i_2,j_2))\cap\Z^2\}, \\
C_K(\varphi^k\sigma^l)
&=\{\varphi^i\sigma^j\,|\, (i,j)\in\cH[(i_1,j_1),(i_2,j_2)]\cap\Z^2\}.
\end{align*}
\item 
If $\varphi^{i_t}\sigma^{j_t}\mm_K^\ast\varphi^k\sigma^l$ with  
$\varphi^{i_t}\sigma^{j_t}$ not onesided-expansive 
for $t=1,2$ and $\cD$ is 
the component of $\E_L(\varphi)\cup\E_R(\varphi)$
with $(k,l)\in\cD$, then 
\begin{align*} 
\cD&=\cH((i_1,j_1),(i_2,j_2)),\\ 
(D_K^\ast)^\circ(\varphi^k\sigma^l)
&=\{\varphi^i\sigma^j\,|\,(i,j)\in\cH((i_1,j_1),(i_2,j_2))\cap\Z^2\}, \\
D_K^\ast(\varphi^k\sigma^l)
&=\{\varphi^i\sigma^j\,|\, (i,j)\in\cH[(i_1,j_1),(i_2,j_2)]\cap\Z^2\}, 
\end{align*}
\end{enumerate}
\end{corollary} 
\begin{proof} 
(1) By Propositin 10.8(2) and Theorem 11.12((11.9)). 

(2) Since
the automorphism 
$\varphi^{i_t}\sigma^{j_t}$ of $(X,\sigma)$ is not onesided-expansive 
for $t=1,2$, $\varphi^{i_t}\sigma^{j_t}$ is not 
onesided $\sigma$-expansive for $t=1,2$. 
Since $\bar{\ell}(i_1,j_1)$ and $\bar{\ell}(i_2,j_2)$ are 
two distinct lines which are not onesided $\sigma$-expansive 
for $\varphi$,  
$\cD$ is an open cone which is not an open half plane in $\R^2$. 
By Theorem 11.13(2), 
($\Lambda$ being the same as in the theorem) 
$\Lambda((D_K^\ast)^\circ(\varphi^k\sigma^l))=\cD\cap\Z^2$. 
 
Suppose that $\varphi^i\sigma^j\mm_K^\ast\varphi^k\sigma^l$ 
and $\varphi^i\sigma^j$ is not onesided $\sigma$-expansive. 
There exists $(m,n)\in\Z^2$ such that $\varphi^m\sigma^n$ is an 
expansive element in $D_K^\ast(\varphi^k\sigma^l)$ and 
$\varphi^i\sigma^j\mm\varphi^m\sigma^n$. Since 
$in\neq jm$, by Propositon 10.8(1) 
we see 
that for any $(i',j')\in\cH[(i,j),(m,n)]\cap\Z^2$, 
$\varphi^{i'}\sigma^{j'}\in D_K(\varphi^m\sigma^n)$
and that for any $(i',j')\in\cH((i,j),(m,n)]\cap\Z^2$, 
$\varphi^{i'}\sigma^{j'}\in D_K^\circ(\varphi^m\sigma^n)$.
By Proposition 10.4, 
$D_K^\ast(\varphi^m\sigma^n) 
=D_K^\ast(\varphi^k\sigma^l)$ and 
$(D_K^\ast)^\circ(\varphi^m\sigma^n) 
=(D_K^\ast)^\circ(\varphi^k\sigma^l)$. 
Hence we see that for any $(i',j')\in\cH[(i,j),(m,n)]\cap\Z^2$, 
$\varphi^{i'}\sigma^{j'}\in D_K^\ast(\varphi^k\sigma^l)$
and that for any $(i',j')\in\cH((i,j),(m,n)]\cap\Z^2$, 
$\varphi^{i'}\sigma^{j'}\in (D_K^\ast)^\circ(\varphi^k\sigma^l)$.
Since 
$\Lambda((D_K^\ast)^\circ(\varphi^k\sigma^l))=\cD\cap\Z^2$, 
$\cH((i,j),(m,n)]\cap\Z^2\subset\cD\cap\Z^2$. 
Since $(i,j)\notin \E_L(\varphi)\cup\E_R(\varphi)$, 
$(i,j)\notin\cD$. Therefore $(i,j)$ is a boundary point 
of $\cD$. We also see, by the above, that 
$\varphi^{i'}\sigma^{j'}\in D_K^\ast(\varphi^k\sigma^l)$ 
for all lattice points $(i',j')$ on 
the closed half-line with endpoint $(0,0)$ 
passing through $(i,j)$.

Therefore $(i_1,j_1)$ and $(i_2,j_2)$ are boundary points 
of the open cone $\cD$ with $i_1j_2\neq i_2j_1$. Hence 
we have $\cD=\cH((i_1,j_1),(i_2,j_2))$, and hence 
$(D_K^\ast)^\circ(\varphi^k\sigma^l)=
\{\varphi^i\sigma^j\,|\,(i,j)\in\cH((i_1,j_1),(i_2,j_2))\}$. 
We also see, by the above, that 
$\{\varphi^i\sigma^j\,|\,(i,j)\in\cH[(i_1,j_1),(i_2,j_2)]\}
\subset D_K^\ast(\varphi^k\sigma^l)$.

Assume that there were $(i',j')\in\Z^2$ 
such that $\varphi^{i'}\sigma^{j'}\in D_K^\ast(\varphi^k\sigma^l)$ 
but $(i',j')\notin \cH[(i_1,j_1),(i_2,j_2)]$. There exists 
$(m',n')\in 
\Lambda((D_K^\ast)^\circ(\varphi^k\sigma^l))$ 
such that 
$\varphi^{i'}\sigma^{j'}\mm\varphi^{m'}\sigma^{n'}$. 
We see that $\varphi^{i'}\sigma^{j'}$ is 
nonexpansive, because, otherwise, it would be 
an element of $(D_K^\ast)^\circ(\varphi^k\sigma^l)$. 
Since $\varphi^{m'}\sigma^{n'}$ is expansive and 
$\varphi^{i'}\sigma^{j'}$ is nonexpansive,  
$i'n'\neq j'm'$. Thererefore it follows from 
Propositon 10.8(1) that 
for any $(i'',j'')\in\cH((i',j'),(m',n')]\cap\Z^2$, 
$\varphi^{i''}\sigma^{j''}\in D_K^\circ(\varphi^{m'}\sigma^{n'})$.
Since either $(i_1,j_1)$ or $(i_2,j_2)$ is in $\cH((i',j'),(m',n')]$, 
either $\varphi^{i_1}\sigma^{j_1}$ or 
$\varphi^{i_2}\sigma^{j_2}$ would be 
onesided $\varphi^{m'}\sigma^{n'}$-expansive, 
contrary to the hypothesis (by definition). 
\end{proof} 

Boyle and Lind \cite[Proposition 8.3]{BoyLin} showed that 
for any expansive component $\C$ 
of 1-frames for a $\Z^d$-action on an infinite compact metric space, 
if a vector in $\C$ is ``Markov'' 
then so is every vector in $\C$. They defined a component 
containing  a marcov vector to be``Markov''. 
When $d=2$ and the space is zero-dimensional, 
this implies the following. For an automorphism $\varphi$
of an infinite subshift $(X,\sigma)$ and  
for any component of $\E(\varphi)=\E_L(\varphi)\cap\E_R(\varphi)$, 
if for some lattice point $(k,l)$ in $\C$, 
$(X,\varphi^k\sigma^l)$ is conjugate to an SFT (respectively, 
a mixing SFT) , then for every 
lattice point $(i,j)$ in $\C$, $(X,\varphi^i\sigma^j)$ is 
conjugate to an SFT (respectively, a mixing SFT).  
Moreover the statement (1) of the following corollary 
shows that a Markov component of a $\Z^2$-action 
$\alpha:(i,j)\mapsto \varphi^i\sigma^j$
is the ``same'' (in a natural sense)  as the interior of an ELR   
(essentially LR ) cone in the subgroup $K$ generated by 
$\{\sigma,\varphi\}$. This will be generalized to 
$\Z^d$-actions on zero-dimensional compact metric spaces 
for any $d$ in the next subsection. 
Recall Proposition 10.2 and Theorem 10.1. 

\begin{corollary} 
Let $\varphi$ be an automorphism of 
an infinite subshift $(X,\sigma)$. Let $K$ be 
the sugroup of $H(X)$ generated by $\{\sigma,\varphi\}$. 
For any $F\subset K$, 
let $\Lambda(F)$ denote the set 
$\{(i,j)\in\Z^2 \,|\, \varphi^i\sigma^j\in F\}$. 
\begin{enumerate}
\item[(1)]
If a component $\C$ of $\E_L(\varphi)\cap\E_R(\varphi)$ 
contains a lattice point $(m,n)$ with 
$(X,\varphi^m\sigma^n)$ conjugate to an SFT, then 
$\C\cap\Z^2
=\Lambda((C_0)_K^\circ(\varphi^m\sigma^n))$.
\item[(2)]
If a component $\cD$ of $\E_L(\varphi)\cup\E_R(\varphi)$ 
contains a lattice point $(m,n)$ with 
$(X,\varphi^m\sigma^n)$ conjugate to a mixing SFT, then 
for every $(k,l)\in\cD\cap\Z^2$ with $\varphi^k\sigma^l$ 
expansive, $(\varphi^k\sigma^l)$ is conjugate to a 
mixing SFT, and the following hold for any $(i,j)\in\Z^2$: 
\begin{enumerate}
\item [(a)]
if $\varphi^i\sigma^j\pl\varphi^k\sigma^l$ 
(respectively, $\varphi^i\sigma^j\qr\varphi^k\sigma^l$) 
and $\varphi^i\sigma^j$ is left $\varphi^k\sigma^l$-expansive 
(respectively, right $\varphi^k\sigma^l$-expansive), 
then $\varphi^i\sigma^j$ has POTP and is 
topologically mixing; 
\item[(b)]
if $\varphi^i\sigma^j\lr\varphi^k\sigma^l$ and $(i,j)\in\cD$,  
then $\varphi^i\sigma^j$ has POTP and is 
topologically mixing; 
\item[(c)]
if $\varphi^i\sigma^j\lr\varphi^k\sigma^l$ 
then $\varphi^i\sigma^j$ has POTP. 
\end{enumerate}
\end{enumerate}
\end{corollary}
\begin{proof} 
(1) By Theorem 11.12((11.9) and Proposition 10.2(1). 

(2) Suppose that 
there exists $(m,n)\in\cD\cap\Z^2$ with 
$(X,\varphi^m\sigma^n)$ conjugate to a mixing SFT.  
By Theorem 11.13(2) there exists  an extended district $D^*$ in $K$
such that $\cD\cap\Z^2=\Lambda((D^\ast)^\circ)$, where $(D^*)^\circ$ 
is the interior of $D^*$.  
Using this and Theorem 10.5, we see that 
for every $(k,l)\in\cD\cap\E(\varphi)\cap\Z^2$, 
$(X,\varphi^k\sigma^l)$ is conjugate to a mixing SFT, and 
(a) and (c) hold. 

To prove (b), suppose that 
$\varphi^i\sigma^j\lr\varphi^k\sigma^l$ and $(i,j)\in\cD$ 
with $(k,l)\in\cD\cap\E(\varphi)\cap\Z^2$. Then, since 
$(i,j)\in\E_L(\varphi)\cup\E_R(\varphi)$, by Proposition 
11.10(3) $\varphi^i\sigma^j$ is 
onesided $\varphi^k\sigma^l$-expansive. Therefore, 
since $\varphi^i\sigma^j\lr\varphi^k\sigma^l$, (b) 
follows from (a). 
\end{proof} 

We note that  
all information provided by (a),(b),(c) of (2) in Corollary 11.15
is concerned with the dynamics of $\varphi^i\sigma^j$ with 
$(i,j)$ being on the boundary of a component $\C$ of $\E(\varphi)$ 
such that $\C\subset\cD$ for $\cD$ in the corollary. 
This is seen by the following proposition. 

\begin{proposition} Let $\varphi$ be an automorphism of an infinite 
subshift $(X,\sigma)$. Let $\cD$ be a component of 
$\E_L(\varphi)\cup\E_R(\varphi)$ and let 
$(k,l)\in \cD\cap\E(\varphi)\cap\Z^2$.
If $\varphi^i\sigma^j\pl\varphi^k\sigma^l$ 
(respectively, $\varphi^i\sigma^j\qr\varphi^k\sigma^l$) 
and $\varphi^i\sigma^j$ is left $\varphi^k\sigma^l$-expansive 
(respectively, right $\varphi^k\sigma^l$-expansive), 
then there exists a component $\C$ of $\E(\varphi)$ such that 
$\C\subset\cD$ and 
$(i,j)$ is contained in the closure of $\C$ in $\R^2$. 
\end{proposition}
\begin{proof} If $il=jk$ then 
$\varphi^i\sigma^j\lr^\circ\varphi^k\sigma^l$
and hence the conclusion trivially holds. 
Suppose that $il\neq jk$. 
Since $\varphi^i\sigma^j$ is left $\varphi^k\sigma^l$-expansive, 
by Proposition 11.10(2) 
there exists a component $\hat{\C}^L$ of $\E_L(\varphi,(k,l))$ 
containing $(i,j)$. By Proposition 11.10((11.4)),  
$\hat{\C}^L$ is an open cone in $\R^2$. 
Let $\C_0=\hat{\C}^L\cap\cH((i,j),(k,l))$. Then $\C_0$ 
is an open cone in $\R^2$. 
Since $\varphi^i\sigma^j\pl\varphi^k\sigma^l$, 
it follows from \cite[Remark 7.9, Proposition 7.6]{Nasu-te} that 
$(\varphi^i\sigma^j)^r(\varphi^k\sigma^l)^s\pl^\circ\varphi^k\sigma^l$, 
i.e., $(\varphi^i\sigma^j)^r(\varphi^k\sigma^l)^s$ is an 
essentially weakly $p$-L, right $\varphi^k\sigma^l$-expansive 
automorphism of $(X,\varphi^k\sigma^l)$, for all $r\geq 0, s\geq 1$. 
Therefore,
for any lattice point $(i',j')$ in $\C_0$, $\varphi^{i'}\sigma^{j'}$
is both left $\varphi^k\sigma^l$-expansive 
and right $\varphi^k\sigma^l$-expansive, 
and hence expansive, and 
$\varphi^{i'}\sigma^{j'}\in PL_K^\circ(\varphi^k\sigma^l) 
\subset (D_K^\ast)^\circ(\varphi^k\sigma^l)$, where 
$K$ is the subgroup of $H(X)$ generated by $\{\sigma,\varphi\}$. 
Hence $\C_0\subset\E(\varphi)$, and by Theorem 11.13(2)
$\C_0\subset\cD$. Let 
$\C$ be the component of $\E(\varphi)$ with $\C\supset\C_0$.
Then, $\C\subset\cD$, and $(i,j)$ 
is contained in the closure of $\C$, 
because $(i,j)$ is on the boundary of $\C_0$. 
\end{proof} 

\begin{example}  Let $A$ be an alphabet.  
Let $B=A^{r+s}$ with $r,s\geq 1$. 
Let $f:A^{r+s+1}\to A$ be a bipermutive 
local rule (see Example 9.20(1)). 
Let $T=(p,q:G_A^{[r+s+1]}\to G_A)$
be the textile system such that 
for each arc $w=a_{-r}a_{-r+1}\dots a_s$ in $G_A^{[r+s+1]}$ 
with $a_j\in A$, $p(w)=a_0$ and $q(w)=f(a_{-r}\dots a_s)$.  
Then $T^*$ is 1-1 and LL with $(X_{T^*},\sigma_{T^*})=(X_B,\sigma_B)$,  
and $\psi=\varphi_{T^*}$ is an LL automorphism of 
$(X_B,\sigma_B)$. 
Hence $P_L(\psi)=Q_L(\psi)$=0 (by Proposition 7.10(1)). 
Since $\eta_{(T^*)^*}=\eta_T$ is not one-to-one, $\psi$ is 
not right $\sigma_B$-expansive. Therefore by Proposition 9.15 we see 
that $-p_L(\psi)=q_L(\psi)=0$. Hence the open intervals
$(-\infty,q_L(\psi))=(-\infty,0)$ and $(-p_L(\psi),\infty)=(0,\infty)$ 
are the components of 
$E_R(\psi)=(-\infty,0)\cup(0,\infty)$ and  
\begin{align*}
\E_R(\psi)
&=\Ll^+(E_R(\psi))\cup\Ll^-(E_R(\psi)))
\cup\hbar^+(0,0)\cup\hbar^-(0,0)\\ 
&=\R^2\setminus\ell_0(0,0). 
\end{align*}

Since $\xi_{(T^*)^*}=\xi_T$ is one-to-one, 
$0$ is a left $\sigma_B$-expansive direction for $\psi$. 
Therefore $p_R(\psi)\neq 0$ and $-q_R(\psi)\neq 0$, 
because by Theorem 11.8(3) 
neither $p_R(\psi)$ nor $-q_R(\psi)$ 
is left $\sigma_B$-expansive. 
For $t>0$, $t$ is a nonexpansive direction for 
the endomorphism $\varphi_T$ of $(X_T,\sigma_T)$ 
if and only if $1/t$ is a nonexpansive direction 
for $\psi=\varphi_{T^*}$ of $(X_B,\sigma_B)
=(X_{T^*},\sigma_{T^*})$. 
(For $t$ is a nonexpansive direction for the endomorphism 
$\xi_T^{-1}\eta_T$ of $(Z_T,\varsigma_T)$
if and only if $1/t$ is a nonexpansive direction for 
the automorphism $\xi_{T^*}^{-1}\eta_{T^*}$ 
of $(Z_{T^*},\varsigma_{T^*})$, 
and $(X_T,\sigma_T, \varphi_T)$ and 
$(Z_T,\varsigma_T, \xi_T^{-1}\eta_T)$ 
are topologically conjugate and so are 
$(X_B,\sigma_B,\psi)=(X_{T^*},\sigma_{T^*},\varphi_{T^*})$
and $(Z_{T^*},\varsigma_{T^*},\xi_{T^*}^{-1}\eta_{T^*})$.) 
Therefore,
since $-s$ and $r$ are all of 
the nonexpansive directions for $\varphi_T$, 
$-1/s$ and $1/r$  are all of the nonexpansive
directions for $\psi$ other than $0$. 
Therefore, since $p_L(\psi)=q_L(\psi)=0$, it follows from 
Propositions 9.9(2) and 9.12(2) that
$p_R(\psi)=-1/s$ and $-q_R(\psi)=1/r$, 
which are not left $\sigma_B$-expansive, by Theorem 11.8(3). 
Therefore,  
$E_L(\psi)=(-\infty,-1/s)\cup(-1/s,1/r)\cup(1/r,\infty)$ 
and hence
\begin{align*}
\E_L(\psi)
&=\Ll^+(E_L(\psi))
\cup\Ll^-(E_L(\psi))\cup\hbar^+(0,0)\cup\hbar^-(0,0)\\
&=\R^2\setminus(\ell_{-1/s}(0,0)\cup\ell_{1/r}(0,0)). 
\end{align*}
Therefore we have 
\[\E_L(\psi)\cup\E_R(\psi)=\R^2\setminus\{(0,0)\},\]
which consists of a unique component. Therefore, 
by virtue of Theorem 11.13(2), 
$\{\psi^i\sigma_B^j\,|\,(i,j)\in \Z^2\setminus\{(0,0)\}\}$
is the interior of a unique extended district $K$ in $K$, 
where $K$ is 
the subgroup of $H(X_B)$ generated by $\{\psi,\sigma_B\}$. 
Since $(X_B,\sigma_B)$ is a mixing topological Markov shift, 
so is $(X_B,\psi^i\sigma_B^j)$ for 
all $(i,j)\in \Z^2$ such that $\psi^i\sigma_B^j$ is expansive 
(i.e., $j/i\neq -1/s$, $j\neq0$ and $j/i\neq 1/r$), by 
Corollary 11.15(1). Moreover, by constructing 
LR textile systems directly, it was known 
in \cite[Section 10, Example 3]{Nasu-t} that 
$K$ is the union of ELR cones in $K$, hence 
it follows from Corollary 11.15(2)(b) that 
$\psi^i\sigma^j$ has POTP and is topologically mixing for all 
$(i,j)\in\Z^2$ with $(i,j)\neq (0,0)$, and hence in 
particular, $\psi^s\sigma_B^{-1}$, $\psi$ and $\psi^r\sigma_B$ 
have POTP and are topologically mixing. 
(Furthermore, in \cite[Section 10, Example 3]{Nasu-t}, it was seen that
$(X_B,\psi^i\sigma_B^j)$ is 
conjugate to a full-shift whenever $\psi^i\sigma_B^j$ is expansive, 
and more precise overall dynamics of $\psi$ was given.) 
We have 
$E(\psi)=(-\infty,-1/s)\cup(-1/s,0)\cup(0,1/r)\cup(1/r,\infty)$.
Let 
\begin{align*}
\C_1&=\Ll^+((-1/s,0)),\\
\C_2&=\Ll^+((-\infty,-1/s))\cup h^-(0,0)\cup\Ll^-((1/r,\infty))
=\bar{\C}_\psi,\\
\C_3&=\Ll^-((0,1/r)), \q\text{and}\\
\bar{\C}_i&=-\C_i, \q i=1,2,3.
\end{align*} 
Then \[\E(\psi)=\C_1\cup\C_2\cup\C_3\cup
\bar{\C}_1\cup\bar{\C}_2\cup\bar{\C}_3,\] 
and 
$\C_1,\C_2, \C_3, \bar{\C}_1,\bar{\C}_2$ and $\bar{\C}_3$  are 
the components of $\E(\psi)$ or 
the expansive components for the automorphism $\psi$. 
Let
$\D_k^\circ=\{\psi^i\sigma_B^j\,|\,(i,j)\in\C_k\}$ for $k=1,2,3$, and 
$\bar{\D}_k^\circ=\{\psi^i\sigma_B^j\,|\,(i,j)\in\bar{\C}_k\}$ for $k=1,2,3$. 
Then, by virtue of Theorem 11.13(1),  
$\D_1^\circ, \D_2^\circ, \D_3^\circ, 
\bar{\D}_1^\circ, \bar{\D}_2^\circ, \bar{\D}_3^\circ$
are the interiors of the ELR cones in $K$. 
Using Theorem 11.12, we see that
$\D_1^\circ\mm\D_2^\circ\mm\D_3^\circ\mm
\bar{\D}_1^\circ\mm\bar{\D}_2^\circ\mm\bar{\D}_3^\circ\mm\D_1^\circ$
(see the paragraph following Proof of Theorem 10.5; 
cf. \cite[Example 9.8]{Nasu-te}). 
\end{example} 
\subsection{``Resolvingness'' and `` 
expansiveness'' for $\Z^d$-actions} 

In this section, we discuss about the relation between 
``expansiveness'' and ``resolvingness'' for $\Z^d$-actions. 
The discussions will be made by using the framework provided by 
Boyle and Lind \cite{BoyLin}. The reader is referred to 
\cite{BoyLin} for the notions due to them which are used without
definitions in this section. 

Let $X$ be an infinite compact  metric space.
Let $K$ be the commutative subgroup 
generated by 
$\{\psi_1,\dots,\psi_d\}$. Let $\alpha:\Z^d\to K$ be 
the $\Z^d$-action defined by  
$\alpha^{\bfm}=\psi_1^{m_1}\cdots\psi_d^{m_d}$, 
where $\bfm=(m_1,\dots,m_d)$ is an integral vector in $\R^d$. 

Let $\|\cdot\|$ be the Euclidean norm on $\R^d$.
For any subset $F$ of $\R^d$ 
and $t\geq 0$, let $F^t$ denote 
the set of all points $\bfv\in \R^d$ with 
$\inf\{\|\bfv-\mathbf{w}\|\,|\,\mathbf{w}\in F\}\leq t$, which is 
the subset made by thickenning $F$ by $t$. 
A subset $F$ of $\R^d$ is said to be \itl{expansive for $\alpha$} 
if there exists $\delta>0$ and $t\geq 0$ such that 
for any $x,y\in X$ it holds that if 
$d_X(\alpha^{\bfm}(x),\alpha^{\bfm}(y))\leq \delta$ 
for all integral vectors $\bfm\in F^t$, then $x=y$ 
(\cite[Definition 2.2]{BoyLin}). 

Let $\E(\alpha)$ be the subset of $\R^d$ obtained from 
the set-union of all lines 
passing through $\bfzero$ (i.e., all one-dimensional subspaces) 
in $\R^d$ that are expansive for $\alpha$, by 
subtracting $\{\bfzero\}$. 
Boyle and Lind \cite{BoyLin} showed that 
if $E(K)\neq\emptyset$, then each (connected)
component $\C$ of $\E(\alpha)$ is an open cone in $R^d$ and has the 
property that $\C\cap(-\C)=\emptyset$. They
call it an \itl{expansive component of 1-frames for $\alpha$}. 

When $X$ is zero-dimensional, 
an expansive component of 1-frames for $\alpha$ which
contains an integral vector $\bfk$ with $(X,\alpha^\bfk)$ 
conjugate to an SFT is called a \itl{Markov component} for $\alpha$. 
( The definition of Boyle and Lind \cite{BoyLin} is 
more general: they defined a `` Markov vector'' 
without the hypothesis that $X$ is zero-dimensional,  
proved that if an expansive component $\C$ of 1-frames for $\alpha$
contains a Markov vector then  all vectors in $\C$ are Markov, and 
defined an expansive component containing 
a Markov vector to be ``Markov''.)
 
Theorems 11.12((11.19)), 11.13(1) and   
Corollary 11.15 are generalized to the 
following theorem. 
\begin{theorem} 
Let $X$ be an infinite zero-dimensional compact  metric space.
Let $K$ be the commutative subgroup of $H(X)$ 
generated by 
$\{\psi_1,\dots,\psi_d\}$. Let $\alpha:\Z^d\to K$ be 
the $\Z^d$-action defined by  
$\alpha^{(m_1,\dots,m_d)}=\psi_1^{m_1}\cdots\psi_d^{m_d}$.  
\begin{enumerate}
\item 
Then two integral vectors
$\bfm=(m_1,\dots, m_d)$ and $\bfn=(n_1,\dots,n_d)$ in $\R^d$ 
belong to the same 
expansive component of 1-frames for $\alpha$ 
if and only if 
$\psi^{m_1}\dots\psi^{m_d}\in C_K^\circ(\psi^{n_1}\dots\psi^{n_d})$, 
i.e., 
$\alpha^\bfm$ is an expansive 
essentially weakly LR automorphism of $(X,\alpha^\bfn)$ 
with $\bfm$ expansive, 
(or $\alpha^\bfm\lr^\circ\alpha^\bfn$). 
\item 
\begin{enumerate}
\item
A subset of integral vectors is given as 
$\C\cap\Z^d$, where $\C$ is an expansive component of 1-frames for 
$\alpha$, if and only if it is given as  
$\Lambda(C^\circ)=\{\bfn\in\Z^d\,|\,\alpha^\bfn\in C^\circ\}$, 
where $C^\circ$ is the interior of an 
essentially weakly LR cone $C$ in $K$.
\item
A subset of integral vectors is given as 
$\C\cap\Z^d$, where $\C$ is a \itl{Markov component} for 
$\alpha$ if and only if it is given as  
$\Lambda(C^\circ)=\{\bfn\in\Z^d\,|\,\alpha^\bfn\in C^\circ\}$, 
where $C^\circ$ is the interior of an ELR cone $C$ in $K$. 
\end{enumerate}
\end{enumerate}
\end{theorem} 
\begin{proof} 
(1) We may assume that
the vectors $\mathbf{m}$ and $\mathbf{n}$ are linearly independent, 
because otherwise, the theorem is clear (by Corollary 8.8(1)). 
Put $\varphi=\psi_1^{m_1}\cdots\psi_d^{m_d}=\alpha^{\mathbf{m}}$ and 
$\tau=\psi_1^{n_1}\cdots\psi_d^{n_d}=\alpha^{\mathbf{n}}$ 
and suppose that $\tau$ and $\varphi$
are expansive. 
Let $\pi:\R^2\to\R^d$ be the map 
$(a,b)\mapsto a\mathbf{m}+b\mathbf{n}$. 
Then $\pi$ is a linear imbedding of $\R^2$ into $\R^d$ 
and $\pi(\R^2)$ is a 2-plane (2-dimensional subspace) in  
$\R^d$. We note the following fact: for any line 
$\bel$ passing through $(0,0)$ (which may be horizontal) in 
$\R^2$, it holds that $\bel$ is an expansive line for 
the automorphism $\varphi$ of $(X,\tau)$ if and only if 
the line $\pi(\bel)$ in $\R^d$ is expansive 
for the $\Z^d$-action $\alpha$. This is clear 
when $\bel$ is horizontal. When $\bel$ is not horizontal, 
the ``only-if'' part is clear and 
the ``if'' part follows from the 
continuity of $\psi_i, \psi_i^{-1}, i=1,\dots,d$, because 
for any $t\geq 0$ there exist $s, e\geq 0$ 
such that for any line $\bel$ 
in $\R^2$ and for each lattice point $\mathbf{k}=(k_1,\dots,k_d)$ 
in $\R^d$ that is contained in 
the tube $(\pi(\bel))^{t}$ made by thickening $\pi(\bel)$ by $t$,
there exists a lattice point $(I,J)$ in $\R^2$ such that
$(I,J)\in\R^2[\bel-(0,s),\bel+(0,s)]$ and 
if $\mathbf{k}-\pi(I,J)=(l_1,\dots,l_d)$ 
(hence $\psi_1^{k_1}\cdots\psi_d^{k_d}(x)
=\psi_1^{l_1}\cdots\psi_d^{l_d}(\varphi^I\tau^J(x))$
for $x\in X$) then $|l_i|\leq e$ for $i=1,\dots,d$.  

Assume that $\C$ is a component of $\E(\alpha)$ and
$\mathbf{m},\mathbf{n}\in \C$. Then by 
\cite[p. 71]{BoyLin} $\C$ is an open cone in $\R^d$ 
(whose apex is the origin $\mathbf{0}$). 
Therefore $\C\cap\pi(\R^2)$
is an open cone in the 2-plane $\pi(\R^2)$ 
and contains $\mathbf{m},\mathbf{n}$. 
Let $S$ be the subset of $\R^2$ such that $\pi(S)=\C\cap\pi(\R^2)$. 
Then $S$ is an open cone in $\R^2$ containing points $(1,0),(0.1)$. 
By the fact noted above,
any line $\bel$ in $\R^2$ passing through a point in $S$ and $(0,0)$ 
is expansive for the automorphism $\varphi$ of $(X,\tau)$, because 
$\pi(\bel)$ is an expansive for $\alpha$. Hence $(1,0)$ and 
(0,1) belong to the same expansive component for the automorphism 
$\varphi$. Therefore it follows from Theorem 11.12 (see (11.9)) that 
$\varphi$ is an essentially weakly LR automorphism of $(X,\tau)$.  
Hence the ``only-if'' part of the theorem is proved. 

Assume that $\varphi\lr^\circ\tau$.
Then by Theorem 11.12 (see (11.9)), $(0,1)$ and (1,0) belong 
to the same expansive component, say $\bar{\C}$, 
for the automorphism $\varphi$
of $(X,\tau)$. Since $\bar{\C}$ is an open cone in $\R^2$,  
 $\pi(\bar{\C})$ is an open cone in $\pi(\R^2)$, which contains 
$\mathbf{m}=\pi(0,1)$ and $\mathbf{n}=\pi(1,0)$. 
Since every line in $\R^d$ passing through $\bfzero$ 
and a point in $\pi(\bar{\C})$ is 
expansive for $\alpha$ by the fact noted above, there exists 
a component $\E(\alpha)$ which contains 
$\mathbf{m}$ and $\mathbf{n}$. 

(2) By (1) and Proposition 10.2. 
\end{proof} 

Boyle and Lind \cite[Proposition 8.3]{BoyLin} showed that 
for any expansive component $\C$ 
of 1-frames for a $\Z^d$-action $\alpha$ on 
an infinite compact metric space $X$, 
if a vector in $\C$ is ``Markov'' 
then so is every vector in $\C$. 
When $X$ is zero-dimensional, 
this implies the following: 
if $(X,\alpha^\bfk)$ is conjugate to an SFT (respectively, 
a mixing SFT) for some integral vector $\bfk$in $\C$, 
then for every integral vector $\bfm$ in $\C$
$(X,\alpha^\bfm)$ is 
conjugate to an SFT (respectively, a mixing SFT). 
This is also a consequence of 
Theorem 11.18 and Propositio 10.2.

Let $X$ be an infinite compact metric space.
Let $K$ be the commutative subgroup of $H(X)$ generated by 
$\{\psi_1,\dots,\psi_d\}$. 
Let $\alpha:\Z^d\to K$ be the $\Z^d$-action defined by  
$\alpha^{(m_1,\dots,m_d)}=\psi_1^{m_1}\cdots\psi_d^{m_d}$.  
An integral vector $\mathbf{k}=(k_1,\dots,k_d)$ is called an
\itl{expansive integral vector for $\alpha$} 
if $\alpha^\mathbf{k}=\psi_1^{k_1}\cdots\psi_d^{k_d}$ is expansive. 

Let $\mathbf{k}$ be an expansive integral vector for $\alpha$ 
and $V$ a subspace of $\R^d$. Define 
\begin{gather*} 
(V,\mathbf{k})=\{\mathbf{v}+b\mathbf{k}\,|\,\mathbf{v}\in V, b\in \R\}\\
(V,\mathbf{k})_+=\{\mathbf{v}+b\mathbf{k}\,|\,
\mathbf{v}\in V, b\geq 0\},\q
(V,\mathbf{k})_-=\{\mathbf{v}+b\mathbf{k}\,|\,
\mathbf{v}\in V, b\leq 0\}.
\end{gather*}
Then $(V,\mathbf{k})$ is the 
\itl{subspace of $\R^d$ spanned by $V$ and $\mathbf{k}$}.
If $\mathbf{k}$ is not in $V$, then
$(V,\mathbf{k})_+$ is the closed half-space of 
$(V,\mathbf{k})$ containing $\mathbf{k}$ with boundary $V$, and 
$(V,\mathbf{k})_-$ is the closed half-space of 
$(V,\mathbf{k})$ containing $-\mathbf{k}$ with boundary $V$. 

We say that $V$ is \itl{right $\mathbf{k}$-expansive} for $\alpha$ 
if $(V,\mathbf{k})_+$ is an expansive subset of $\R^d$ for $\alpha$,  
and 
\itl{left $\mathbf{k}$-expansive} for $\alpha$ if $(V,\mathbf{k})_-$ 
is an expansive subset of $\R^d$ for $\alpha$. Note that $(V,\bfk)$ is 
expansive because $(V,\bfk)$ includes the expansive line 
$\bel(\bfk)$ for $\alpha$,  where for $\bfv\in\R^d$, $\bel(\bfv)$ 
denotes the line passing through $\bfzero$ and $\bfv$. 
If $V\ni\mathbf{k}$, then $V$ is expansive for $\alpha$. 
\begin{proposition} Let $\alpha$ be a $\Z^d$-action and 
$\mathbf{k}$ an expansive integral vector for $\alpha$. 
Then a subspace $V$ of $\R^d$ is expansive 
for $\alpha$ if and only if $V$ is right $\mathbf{k}$-expansive 
for $\alpha$ and left $\mathbf{k}$-expansive for $\alpha$.
\end{proposition} 
\begin{proof} 
Since $(V,\mathbf{k})_+\supset V$ and 
$(V,\mathbf{k})_-\supset V$, ``only-if'' part follows. 
The proof of the ``if''-part is given 
by \cite[Proposition 3.10]{BoyLin} 
(it is also found in the proof of the next proposition).
\end{proof} 

Let $\mathbb{G}_k$ be the Grassmann manifolds consisting 
of the $k$-dimensional subspaces of $\R^d$. 
For a $\Z^d$-action $\alpha$, an expansive integral vector 
$\mathbf{k}$ and an integer $k\geq 0$, let 
$\mathbb{E}_{R,k}(\alpha,\mathbf{k})$ (respecttively, 
$\mathbb{E}_{L,k}(\alpha,\mathbf{k})$) denote the set of 
all subspaceces 
in $\mathbb{G}_k$ that are right 
$\mathbf{k}$-expansive 
(respectively, left
$\mathbf{k}$-expansive for) $\alpha$. 

The following two propositions and a corollary generalizes   
\cite[Lemma 3.4, Theorem 3.6]{BoyLin} and a result on 
expansive components of 1-frames for a $Z^d$-action 
in \cite[p. 71]{BoyLin} 
in the case that the $\Z^d$-action $\alpha$ has an expansive 
integral vector. The proofs are given 
by mimicking the proofs 
of Boyle and Lind in \cite{BoyLin} 
for the corresponding results. 

Suppose that
$\bfk$ is an expansive integral vector for 
a $\Z^d$-action $\alpha$. 
Then the line $\bel(\bfk)$ is an expansive subspace for 
$\alpha$ and hence $\R^d$ is expansive for $\alpha$.  
Hence, we can use 
``$E$ \itl{codes} $F$'' for subsets $E,F$ of $\R^d$
in the sense of \cite[definition 3.1]{BoyLin}. Then 
for any subset $F$ of $\R^d$ it holds that 
$F$ is expansive for $\alpha$ if and only if 
there exists $t\geq 0$ such that $F^t$ codes $R^d$
(\cite{BoyLin}). 
\begin{proposition} Let $\alpha$ be a $\Z^d$-action 
on a infinite compact metric space, 
$\mathbf{k}$ an expansive integral vector for $\alpha$ 
and $k$ a nonnegative integer. Then 
$\mathbb{E}_{R,k}(\alpha,\mathbf{k})$ (respectively, 
$\mathbb{E}_{L,k}(\alpha,\mathbf{k})$) is an open subset 
of $\mathbb{G}_k$. 
\end{proposition} 
\begin{proof}  
Let $V$ be a subspace of $\R^d$. Let $\pi_V$ denote orthoganal 
projection to $V$ along its orthoganal complement $V_{\perp}$. 
Let $t,r,s\geq 0$. Then
$V^t=\{\bfv\in\R^d\,|\, \|\pi_{V_\perp}(\bfv)\|\leq t\}$. 
Define $V^t(r)=\{\bfv\in V^t\,|\, \|\pi_{V}(\bfv)\|\leq r\}$ (as in
\cite{BoyLin}) and further define 
\begin{align*}
(V,\mathbf{k})_{+,s}^{t}=\{\bfv+b\bfk\,|\,
\bfv\in V^t, 0\leq b\leq s\},\\ 
(V,\mathbf{k})_{+,s}^t(r)=\{\bfv+b\bfk\,|\,
\bfv\in V^t(r), 0\leq b\leq s\}.  
\end{align*}

Suppose that $V\in\mathbb{G}_k$ and 
$V$ is right $\bfk$-expansive for $\alpha$. 
Then there exists $t\geq 0$ such that $(V,\mathbf{k})_+^t$ 
codes $\R^d$ and hence hence $B(t+2)$, where 
$B(r)=\{\bfv\in \R^d\,|\, \|\bfv\|\leq r\}$ for $r\geq 0$.  
Therefore, 
a standard compactness argument shows that
there exists $r, s\geq t+2$  
such that $(V,\mathbf{k})^t_{+,s}(r)$ 
codes $B(t+2)$.
(Using $(V,\mathbf{k})^t_{+,s}(r)$ as a ``pattern'' 
with $B(t+2)$ as a ``head'' for ``sliding pattern coding'', 
we see that $(V,\mathbf{k})^t_{+,s}$ codes 
$(V,\bfk)_-^t$, which proves 
the ``only-if'' part of the previous proposition.) 
For $W\in\mathbb{G}_k$ sufficiently close to $V$, 
$(W,\mathbf{k})^{t+1}_{+,s+1}(r+1)$ contains 
$(V,\mathbf{k})^t_{+,s}(r)$, which codes $B(t+2)$. 
Let $t_V=t+1$, $r_V=r+1$ and $s_V=s+1$. Then we see that 
there exists a neighborhood $\mathcal{N}_V$ of $V$ in 
$\mathbb{G}_k$ such that for every $W\in\mathcal{N}_V$, 
$(W,\mathbf{k})^{t_V}_{+,s_V}(r_V)$ codes $B(t_V+1)$. 
Using $(W,\bfk)^{t_V}_{+,s_V}(r_V)$ as a ``pattern'' 
with $B(t_V+1)$ as a ``head'' for ``sliding pattern coding'', 
$(W,\mathbf{k})^{t_V}_+$ codes $(W,\mathbf{k})^{t_V}_-$
and hence $(W,\mathbf{k})^{t_V}$ and hence $\R^d$, so that 
$W$ is right $\bfk$-expansive. 
\end{proof} 
\begin{proposition} 
Let $\alpha$ be a $\Z^d$-action on a infinite compact 
metric space and 
$\mathbf{k}$ an expansive integral vector for $\alpha$. 
Suppose that $V$ has dimension $k<d-1$ 
and is not right $\bfk$-expansive 
(respectively, left  $\bfk$-expansive) for $\alpha$.
Then $V$ is a subspace of  
a $(d-1)$-dimensional subspace which is not 
right $\bfk$-expansive (respectively, left  $\bfk$-expansive) 
for $\alpha$. 
\end{proposition} 
\begin{proof} 
Let $\mathcal{K}=\{W\in\mathbb{G}_{d-1}\,|\,W\supset V\}$. 
Assume that every subspace in $\mathcal{K}$ 
is right $\bfk$-expansive for $\alpha$. Then 
by the proof of the preceding proposition 
and by the same argument as in \cite[Lemma 3.5]{BoyLin}, 
there are $t,r,s\geq 0$ such that $(W,\bfk)^t_{+,s}(r)$ 
codes $B(t+1)$. Hence by the same argument as in 
\cite[Proof of Theorem 3.6]{BoyLin}, $V$ is expansive, 
contradicting the hypothesis that  $V$ is not right 
$\bfk$-expansive. 
\end{proof} 

For a $\Z^d$-action $\alpha$ on an infinite compact metric space 
and an expansive integral vector $\mathbf{k}$ for $\alpha$, 
define 
\begin{align*}
\E_L(\alpha,\bfk)=(\cup_{\bel\in\mathbb{E}_{L,1}(\alpha,\bfk)}\bel)
\setminus\{\bfzero\},\q\q
\E_R(\alpha,\bfk)=(\cup_{\bel\in\mathbb{E}_{R,1}(\alpha,\bfk)}\bel)
\setminus\{\bfzero\}.
\end{align*}
\begin{corollary} 
Let $\alpha$ be a $\Z^d$-action on an infinite compact metric space 
and $\mathbf{k}$ an expansive integral vector for $\alpha$. 
Then $\E_R(\alpha,\bfk)$ (respectively, $\E_L(\alpha,\bfk)$)
is an open subset of $\R^d$ and each component $\C$ of it 
is an open cone in $\R^d$ with $\C\cap(-\C)=\emptyset$. 
\end{corollary}  
\begin{proof} 
Use Propositions 11.20 and 11.21 and the same argument 
as in \cite[p. 71]{BoyLin}. 
\end{proof}

\begin{remark} 
Let $\alpha$ be a $\Z^d$-action 
on an infinite compact metric space. 
\begin{enumerate} 
\item 
For any expansive integral vector $\bfk$, 
\[\E_L(\alpha,\bfk)\cap\E_R(\alpha,\bfk)=\E(\alpha).\]
\item 
For any expansive integral vectors $\bfk$ and $\bfm$ for $\alpha$, 
the component of $\E_R(\alpha,\bfk)$ 
(respectively,  $\E_L(\alpha,\bfk)$) containing $\bfm$ includes
the component of $\E(\alpha)$ containing $\bfm$. 
\item 
If $d=2$, then 
for any expansive integral vectors $\bfk$ and $\bfm$ for $\alpha$,
\[\E_L(\alpha,\bfk)\cup\E_R(\alpha,\bfk)=
\E_L(\alpha,\bfm)\cup\E_R(\alpha,\bfm).\]
\end{enumerate}
\end{remark} 
\begin{proof} 
(1) This follows from Proposition Proposition 11.19. 

(2) By (1). 

(3) Let $\bel$ be a line passing through the origin in 
the plane $\R^2$. Then if $\bfk$ and $\bfm$ belong to the same one of 
the closed half-planes with boundary $\bel$, then 
$\bel$ is right $\bfk$-expansive if and only if 
$\bel$ is right $\bfm$-expansive, and 
$\bel$ is left $\bfk$-expansive if and only if 
$\bel$ is left $\bfm$-expansive. 
If $\bfk$ and $\bfm$ belon to
different ones of the half-planes, then 
$\bel$ is right $\bfk$-expansive if and only if 
$\bel$ is left $\bfm$-expansive, and 
$\bel$ is left $\bfk$-expansive if and only if 
$\bel$ is right $\bfm$-expansive. Therefore (3) follows. 
\end{proof} 

The statement (3) of this remark proves the fact 
stated in Subsection 2.1 
that if automorphism $\varphi$ of an  
invertible dynamical system $(X,\tau)$ is 
onesided $\varphi^k\tau^l$-expansive for some $\varphi^k\tau^l$ 
which is expansive, then for all 
$\varphi^m\tau^n$ that is expansive, $\varphi$ is 
onesided $\varphi^m\tau^n$-expansive. 
Using the apploach proving (3) of the remark and using 
Proposition 11.6,  
we have easy proofs for generalizations of
Propositions 11.10 and 11.11 
to automorphisms of invertible dynamical systems 
or to $Z^2$-actions.   

The result (11.8) of Theorem 11.12 is generalized to the 
following theorem. 
 
\begin{theorem} 
Let $X$ be an infinite zero-dimensional compact  metric space.
Let $K$ be the commutative subgroup 
of $H(X)$ generated by  
$\{\psi_1,\dots,\psi_d\}$. Let $\alpha:\Z^d\to K$ be 
the $\Z^d$-action defined by  
$\alpha^{(m_1,\dots,m_d)}=\psi_1^{m_1}\cdots\psi_d^{m_d}$.  
Let $\bfk$ be an expansive integral vector for $\alpha$ 
and $\mathbf{m}$ an integral vector. 
Then $\bfm$ belongs to 
the component of $\E_R(\alpha,\bfk)$ (respectively, 
$\E_L(\alpha,\bfk)$) containing $\bfk$ if and only if 
$\psi_1^{m_1}\cdots\psi_d^{m_d}\in 
PL_K^\circ(\psi_1^{k_1}\cdots\psi_d^{k_d})$ 
(respectively, $\psi_1^{m_1}\cdots\psi_d^{m_d}
\in QR_K^\circ(\psi_1^{k_1}\cdots\psi_d^{k_d})$), 
i.e., 
$\alpha^\bfm$ is 
a right $\alpha^\bfk$-expansive, 
essentially weakly $p$-L 
(respectively, left $\alpha^\bfk$-expansive, 
essentially weakly $q$-R) automorphism of 
$(X,\alpha^\bfk)$ (or
$\alpha^\bfm\pl^\circ\alpha^\bfk$
(respectively, 
$\alpha^\bfm\qr^\circ\alpha^\bfk$)).
\end{theorem}
\begin{proof}
We may assume that
the vectors $\mathbf{m}$ and $\mathbf{k}$ are linearly independent, 
because otherwise, the theorem is clear (by Corollary 8.8(1)). 
Put $\varphi=\psi_1^{m_1}\cdots\psi_d^{m_d}=\alpha^{\mathbf{m}}$ and 
$\tau=\psi_1^{k_1}\cdots\psi_d^{k_d}=\alpha^{\mathbf{k}}$. 
Then $\tau$ is expansive.  
Let $\pi:\R^2\to\R^d$ be the map 
$(a,b)\mapsto a\mathbf{m}+b\mathbf{k}$. 
Then $\pi$ is a linear imbedding of $\R^2$ into $\R^d$ 
and $\pi(\R^2)$ is a 2-plane (2-dimensional subspace) in  
$\R^d$. Then we note the following fact: for any line 
$\bel$ passing through $(0,0)$ (which may be horizontal) in 
$\R^2$, it holds that $\bel$ is right $\tau$-expansive for 
the automorphism $\varphi$ of $(X,\tau)$ if and only if 
the one-dimensional subspace  
$\pi(\bel)$ of $\R^d$ is a right $\bfk$-expansive 
for the $\Z^d$-action $\alpha$. This is clear 
when $\bel$ is horizontal. When $\bel$ is not horizontal, 
the ``only-if'' part is clear and 
the ``if'' part follows from the 
continuity of $\psi_i, \psi_i^{-1}, i=1,\dots,d$, because 
for any $t\geq 0$ there exist $s, e\geq 0$ 
such that for any line $\bel$ 
in $\R^2$ and for each lattice point $\mathbf{n}=(n_1,\dots,n_d)$ 
in $\R^d$ that is contained in 
the subset of $\R^d$ made by thickening the closed half 2-plane 
containing $\bfk$ 
with boundary $\pi(\bel)$ by $t$,
there exists a lattice point $(I,J)$ in $\R^2$ such that
$(I,J)\in\R^2[\bel-(0,s),\cdot)$ and that 
if $\mathbf{n}-\pi(I,J)=(l_1,\dots,l_d)$ 
(hence $\psi_1^{n_1}\cdots\psi_d^{n_d}(x)
=\psi_1^{l_1}\cdots\psi_d^{l_d}(\varphi^I\tau^J(x))$
for $x\in X$) then $|l_i|\leq e$ for $i=1,\dots,d$.  

Assume that $\C$ is the component of $\E_R(\alpha,\bfk)$ 
containing $\bfk$ and that $\mathbf{m}\in \C$. 
Then by the preceding 
corollary, $\C$ is an open cone in $\R^d$. 
Therefore $\C\cap\pi(\R^2)$ is an open cone in 
the 2-plane $\pi(\R^2)$ and contains $\bfk$ and $\bfm$. 
Let $S$ be the subset of $\R^2$ such that $\pi(S)=\C\cap\pi(\R^2)$. 
Then $S$ is an open cone in $\R^2$ containing points $(0,1),(1,0)$. 
By the fact noted above,
any line $\bel$ in $\R^2$ passing through 
a point in $S$ and $(0,0)$ 
is right $\tau$-expansive for the automorphism 
$\varphi$ of $(X,\tau)$, because 
$\pi(\bel)$ is right $\bfk$-expansive for $\alpha$. Hence $(0,1)$ 
and (1,0) belong to the same right $\tau$-expansive 
component for the automorphism $\varphi$. 
Therefore it follows from Theorem 11.12 (see (11.8)) that 
$\varphi$ is an essentially weakly $p$-L, 
right $\tau$-expansive automorphism of $(X,\tau)$.  
Hence the ``only-if'' part of the theorem is proved. 

Assume that $\varphi\pl^\circ\tau$.
Then by Theorem 11.12 (see (11.8)), $(0,1)$ and (1,0) belong 
to the same expansive right $\tau$-expansive component, 
say $\bar{\C}$, for the automorphism $\varphi$
of $(X,\tau)$. Since $\bar{\C}$ is an open cone in $\R^2$,  
 $\pi(\bar{\C})$ is an open cone in $\pi(\R^2)$, which contains 
$\mathbf{m}=\pi(1,0)$ and $\mathbf{k}=\pi(0,1)$. 
Since every line in $\R^d$ passing through $\bfzero$ and 
some point in $\pi(\bar{\C})$ is 
right $\bfk$-expansive for $\alpha$ by the fact noted above,  
the component of $\E_R(\alpha,\bfk)$ containing $\bfk$ contains 
$\mathbf{m}$. 
\end{proof} 
\begin{proposition} 
Hypothesis being as in Theorem 11.24,
let $\C^R(\alpha,\bfk)$ (respectively, $\C^L(\alpha,\bfk)$)  
denote the component 
of $\E_R(\alpha,\bfk)$ (respectively, 
$\E_L(\alpha,\bfk)$) containing $\bfk$ 
for any expansive integaral vector $\bfk$  for $\alpha$. 
Then if $\bfk$ and $\bfm$ are integral vectors 
contained in the same component of $\E(\alpha)$, 
then $\C^R(\alpha,\bfk)=\C^R(\alpha,\bfm)$ 
(respectively, $\C^L(\alpha,\bfk)=\C^L(\alpha,\bfm)$). 
\end{proposition}
\begin{proof} 
Suppose that $\bfk$ and $\bfm$ belong to 
the same component of $\E(\alpha)$. Then by Theorem 11.18
$\alpha^\bfk\lr^\circ\alpha^\bfm$. 
Let $\bfn$ be any lattice point in $\C^R(\alpha,\bfk)$. Then 
by Theorem 11.24 $\alpha^\bfn\pl^\circ\alpha^\bfk$.
Hence by \cite[Proposition 9.4(2)]{Nasu-te}, we have 
$\alpha^\bfn\pl^\circ\alpha^\bfm$. Hence by Theorem 2.4 
$\alpha^\bfn\in\C_R(\alpha,\bfm)$. Therefore 
$\C^R(\alpha,\bfm)\cap\Z^d\supset\C^R(\alpha,\bfk)\cap\Z^d$. 
Since $\C^R(\alpha,\bfk)$ and 
$\C^R(\alpha,\bfm)$ are open cones with apex 
$\bfzero$ in $\R^d$, we have 
$\C^R(\alpha,\bfm)\supset\C^R(\alpha,\bfk)$. Therefore, 
sinse the converse inclusion of this is symmetrically obtained,  
we have $\C^R(\alpha,\bfm)=\C^R(\alpha,\bfk)$. 
\end{proof}

\section{Directionally essentially weakly resolving endomorphisms} 

In this section, we shall prove that if an endomorphism $\varphi$
of a subshift is ``directionally essentially weakly \st'', 
then it is ``essentially weakly \st'', 
where \rt is any resolving term. (Recall that 
a \itl{resolving term} means any one of the ten terms 
 ``$p$-L'',  
`` $p$-R'', ``$q$-R'', ``$q$-L''. 
``LR'', ``RL''
``LL'', ``RR'', 
``$p$-biresolving'' and ``$q$-biresolving''.) 
It is known that 
when $\varphi$ is an onto endomorphism of an SFT 
and  \rt is any resolving term excepting LL and RR,  
if $\varphi$ is directionary essentially \st, then 
it is essentially \st, 
 \cite[Propositions 8.5, 8.6]{Nasu-te}. 
The proof of this depends on 
the long theory of resolvable textile systems developed in 
\cite[Section 7]{Nasu-t}. We shall also present 
a much shorter proof 
proving in particular that for the important special case when  
$\varphi$ is an onto endomorphism of an SFT $(X,\sigma)$ with 
$\varphi$ one-to-one or $\sigma$ topologically transitive, if
$\varphi$ is directionary essentially \st,  
then it is essentially \st, where \rt is 
any resolving term (including  LL and RR). 

\begin{proposition} Let $\varphi$ be an  
endomorphism of a  subshift $(X,\sigma)$ and let $s\geq 1$. 
Let \rt be any resolving term. 
If $\varphi$ is an essentially 
weakly \rt endomorphism of $(X,\sigma^s)$, then 
$\varphi$ is an essentially
weakly \rt endomorphism of $(X,\sigma)$. 
\end{proposition} 
\begin{proof} 
Suppose that $\varphi$ is an essentially weakly $p$-L 
endomorphism of $(X,\sigma^s)$. Then there exists 
a onesided 1-1 half-textile-subsystem 
$\hf{U}$ of a weakly $p$-L textile system 
$T$ and a conjugacy 
$\theta:(X,\sigma^s,\varphi)\to 
(X_{\hf{U}},\sigma_{\hf{U}},\varphi_{\hf{U}})$ 
between commuting systems. Let $\tau=\theta\sigma\theta^{-1}$. 
Then $\tau$ is an expansive automorphism of 
$(X_{\hf{U}},\sigma_{\hf{U}})$ 
such that $\tau^s=\sigma_{\hf{U}}$. 
Since $\tau$ is expansive, there exists $l\geq 0$ such that 
if $(x_k)_{k\in\Z}$ is a $\tau$-orbit with
$x_k=(a_{k,j})_{j\in\Z}$, 
where $a_{k,j}\in L_1(X_{\hf{U}})$, then 
$(a_{k,-l}a_{k,-l+1}\dots a_{k,l})_{k\in\Z}$ 
uniquely determines the orbit $(x_k)_{k\in\Z}$. 
Since $T$ is weakly $p$-L, so is $T^{[2l+1]}$. Therefore, 
replacing $T^{[2l+1]}, \hf{U}^{[2l+1]}$ 
and $\rho_{X_{\hf{U}},l,l}\theta$ 
(recall that $\rho_{X_{\hf{U}},l,l}$ is the 
higher-block conjugacy of $(l,l)$-type on $X_{\hf{U}}$) by 
$T, \hf{U}$ and $\theta$, respectively, 
we see that there exists 
a onesided 1-1 half-textile-subsystem $\hf{U}$ 
of a weakly $p$-L textile system 
$T$ and a conjugacy 
$\theta:(X,\sigma^s,\varphi)\to 
(X_{\hf{U}},\sigma_{\hf{U}},\varphi_{\hf{U}})$ 
such that $\tau$ with $\tau=\theta\sigma\theta^{-1}$ is 
an automorphism of $(X_{\hf{U}},\sigma_{\hf{U}})$, 
$\tau^s=\sigma_{\hf{U}}$ and 
if for $x\in X_{\hf{U}}$ and $k,j\in\Z$, 
we define $a_{k,j}(x)\in L_1(X_{\hf{U}})$
by 
\[(a_{k,j}(x))_{j\in\Z}=\tau^k(x),\] 
then 
$(a_{k,0}(x))_{k\in\Z}$ uniquely determines $x$. 

Let $T=(p,q:\G\to G)$. Then $p$ is weakly left-resolving. 
We define a textile system 
${T_1}=({p_1},{q_1}:{\G_1}\to{G_1})$ 
as follows: 
${G_1}$ is the graph such that 
\[A_{{G_1}}=\{(a_k)_{0\leq k\leq s}\,|\, a_k\in A_G\},\q 
V_{{G_1}}=\{(a_k)_{0\leq k\leq s-1}\,|\, a_k\in A_G\}\]
(i.e., $A_{{G_1}}=A_G^{s+1}$ and $V_{{G_1}}=A_G^s$),  
and each arc $(a_k)_{0\leq k\leq s}$ goes 
from $(a_k)_{0\leq k\leq s-1}$ to 
$(a_k)_{1\leq k\leq s}$ in ${G_1}$; 
the graph ${\G_1}$ and the homomorphisms 
${p_1}$ and ${q_1}$ are given 
as follows:
\begin{align*}
A_{{\G_1}}&=\{(\alpha_k)_{0\leq k\leq s}\,|\, \alpha_k\in A_\G, 
\alpha_0\alpha_s\in L_2(\G)\}, \\
V_{{\G_1}}&=\{(\alpha_k)_{0\leq k\leq s-1}\,|\, \alpha_k\in A_\G\}
\end{align*}
and each arc $\hat{\alpha}=(\alpha_k)_{0\leq k\leq s}$ goes 
from $(\alpha_k)_{0\leq k\leq s-1}$ to 
$(\alpha_k)_{1\leq k\leq s}$ in ${\G_1}$ with 
\[({p_1})_A(\hat{\alpha})=
(p_A(\alpha_k))_{0\leq k\leq s}\q\text{and}\q
({q_1})_A(\hat{\alpha}))=
(q_A(\alpha_k))_{0\leq k\leq s}.\]

To see that ${T_1}$ is weakly $p$-L, suppose that 
$\hat{\alpha}=(\alpha_k)_{0\leq k\leq s}$ and 
$\hat{\beta}=(\beta_k)_{0\leq k\leq s}$ are in $A_{{\G_1}}$,
$t_{{\G_1}}(\hat{\alpha})=t_{{\G_1}}(\hat{\beta})$ and 
$({p_1})_A(\hat{\alpha})=({p_1})_A(\hat{\beta})$. Then 
$(\alpha_k)_{1\leq k\leq s}=(\beta_k)_{1\leq k\leq s}$
and $(p_A(\alpha_k))_{0\leq k\leq s}=
(p_A(\beta_k))_{0\leq k\leq s}$, 
and hence    
$\alpha_s=\beta_s$ and $p_A(\alpha_0)=p_A(\beta_0)$. 
Since $\alpha_0\alpha_s,\beta_0\beta_s\in L_2(\G)$ and $p$ is 
weakly left resolving, it follows that $\alpha_0=\beta_0$, 
and hence $\hat{\alpha}=\hat{\beta}$. Therefore, 
${p}_1$ is weakly left-rsolving and 
hence ${T_1}$ is weakly $p$-L.

For $x\in X_{\hf{U}}$ and $k,j\in\Z$, 
define $\alpha_{k,j}(x)\in A_\G$ and $b_{k,j}(x)\in A_G$ by  
\[(\alpha_{k,j}(x))_{j\in\Z}=\xi_{\hf{U}}^{-1}(\tau^k(x)), 
\q\q (b_{k,j}(x))_{j\in\Z}=\varphi_{\hf{U}}(\tau^k(x)).\] 
Then since $(a_{k,j}(x))_{j\in\Z}=
\xi_{\hf{U}}((\alpha_{k,j}(x))_{j\in\Z})$ 
and $(b_{k,j}(x))_{j\in\Z}=
\eta_{\hf{U}}((\alpha_{k,j}(x))_{j\in\Z})$, we have 
\[a_{k,j}(x)=p_A(\alpha_{k,j}(x)), \q  b_{k,j}(x)=
q_A(\alpha_{k,j}(x)).\]
Let $Y=\{(a_{k,0}(x))_{k\in\Z}\,|\,x\in X_{\hf{U}}\}$. 
Then we have 
a subshift $(Y,\sigma_Y)$ and a conjugacy  
$\chi:(X_{\hf{U}},\tau)\to (Y,\sigma_Y)$ by 
$\chi(x)=(a_{k,0}(x))_{k\in\Z}$, because 
$(a_{k,0}(x))_{k\in\Z}$ uniquely determines $x$. 
Let $\psi=\chi\varphi_{\hf{U}}\chi^{-1}$. 
Then $\psi((a_{k,0}(x))_{k\in\Z})=(b_{k,0}(x))_{k\in\Z}$ 
for $x\in X_{\hf{U}}$ and 
$\chi:(X_{\hf{U}},\tau,\varphi_{\hf{U}})\to (Y,\sigma_Y,\psi)$ 
is a conjugacy 
between commuting systems. 
Let $Z=\{(\alpha_{k,0}(x))_{k\in\Z}\,|\,x\in X_{\hf{U}}\}$. Then 
we have another subshift $(Z,\sigma_Z)$. 
We see that 
$(Z^{[s+1]},\sigma_Z^{[s+1]})$ is a subshift of 
$(\hf{Z}_{{T_1}},\hf{\sigma}_{{T_1}})$, 
because for $x\in X_{\hf{U}}$, 
$\alpha_{s,0}(x)=\alpha_{0,0}(\tau^s(x))=
\alpha_{0,0}(\sigma_{\hf{U}}(x))
=\alpha_{0,1}(x)$,
and hence 
\[\alpha_{0,0}(x)\alpha_{s,0}(x)\in L_2(Z_{\hf{U}})
\subset L_2(\hf{Z}_T)\subset L_2(\G).\] 
We see that
$\hf{\xi}_{{T_1}}(Z^{[s+1]})\supset\hf{\eta}_{{T_1}}(Z^{[s+1]})$, 
for the following reasons: 
$Y=\chi(X_{\hf{U}})\supset
\chi(\varphi_{\hf{U}}(X_{\hf{U}}))=\psi(Y)$;
since $Y=\{(p_A(\alpha_{k,0}(x)))_{k\in\Z}\,|\,x\in X_{\hf{U}}\}$, 
$Y^{[s+1]}=\hf{\xi}_{{T_1}}(Z^{[s+1]})$; since
$\psi(Y)=\{(q_A(\alpha_{k,0}(x)))_{k\in\Z}\,|\,x\in X_{\hf{U}}\}$,
$(\psi(Y))^{[s+1]}=\hf{\eta}_{{T_1}}(Z^{[s+1]})$. 
We also see that $\hf{\xi}_{{T_1}}\,|\,Z^{[s+1]}$ is one-to-one, 
because $(a_{k,0}(x))_{k\in\Z}$ uniquely determines $x$ and 
hence $(\alpha_{k,0}(x))_{k\in\Z}$. 

Therefore we have a onesided 1-1 
half-textile-subsystem $\hf{{U_1}}$ 
of ${T_1}$ with $Z_{\hf{{U_1}}}=Z^{[s+1]}$. 
Since ${T_1}$ is weakly $p$-L, $\varphi_{\hf{{U_1}}}$ 
is a weakly $p$-L
endomorphism of $(X_{\hf{{U_1}}},\sigma_{\hf{{U_1}}})$. 
Passing through the conjugacies 
\[(X,\sigma,\varphi)\stackrel{\theta}{\to}
(X_{\hf{U}},\tau,\varphi_{\hf{U}}) 
\stackrel{\chi}{\to}(Y,\sigma_Y,\psi)\to 
(Y^{[s+1]},\sigma_Y^{[s+1]},\psi^{[s+1]})
=(X_{\hf{{U_1}}},\sigma_{\hf{{U_1}}},\varphi_{\hf{{U_1}}}),\]
$\varphi$ is an essentially weakly $p$-L endomorphism of $(X,\sigma)$. 

To prove the proposition, 
for the case where \rt 
is ``$q$-R'' 
it suffices to show 
that if $T$ in the above is weakly $q$-R, 
then so is ${T_1}$. 
To do this, 
suppose that 
$\hat{\alpha}=(\alpha_k)_{0\leq k\leq s}$ and 
$\hat{\beta}=(\beta_k)_{0\leq k\leq s}$ are in $A_{{\G_1}}$,
$i_{{\G_1}}(\hat{\alpha})=i_{{\G_1}}(\hat{\beta})$ and 
$({q_1})_A(\hat{\alpha})=({q_1})_A(\hat{\beta})$. Then    
$\alpha_0=\beta_0$ and $q_A(\alpha_s)=q_A(\beta_s)$. 
Since $\alpha_0\alpha_s,\beta_0\beta_s\in L_2(\G)$ and $q$ is 
weakly right-resolving, it follows that $\alpha_s=\beta_s$, 
and hence $\hat{\alpha}=\hat{\beta}$. Therefore, 
${q_1}$ is weakly right-rsolving 
and hence ${T_1}$ is weakly $q$-R.

It is analogously proved that the proposition is valid 
for the other cases. 
\end{proof}

Let $\varphi$ be an endomorphism of a  
dynamical system $(X,\tau)$. 
Let \rt be any property of 
endomorphisms of dynamical systems. 
We say that $\varphi$ is \itl{directionally \rt} if 
there exist $r,s\geq 1$ such that the endomorphism
$\varphi^r$ of $(X,\tau^s)$ is \rt (\cite[Section 7]{Nasu-te}). 

\begin{theorem} Let \rt be any resolving term. 
Then an endomorphism $\varphi$
of a subshift $(X,\sigma)$ 
is directionally essentially weakly \rt 
if and only if it is essentially weakly \st. 
\end{theorem} 
\begin{proof} 
Suppose that $\varphi^r$ is an essentially weakly \rt
endomorphism of $(X,\sigma^s)$ with $r,s\geq 1$. 
Then, by Proposition 12.1 $\varphi^r$ is an 
essentially weakly \rt endomorphism of $(X,\sigma)$. 
Therefore, so is $\varphi$ 
by Theorems 
8.1(1), 8.6(1) and 8.7(1). 
\end{proof} 

In the following proposition and theorem, 
(1) and (2)(a) are part of already known results 
\cite[Propositions 8.5 and 8.6]{Nasu-te},
whose proof depends on the long theory of 
resolvable textile systems developed in 
\cite[Section 7]{Nasu-t}. Not only 
the proof of them given in the following is much simpler than that 
of \cite[Propositions 8.5 and 8.6]{Nasu-te}, 
but also (2)(b) is a new result. 

\begin{proposition} 
Let $\varphi$ be an onto endomorphism of an SFT $(X,\sigma)$. 
Let $s\geq 1$.
\begin{enumerate}
\item \cite{Nasu-t, Nasu-te} 
If $\varphi$ is an essentially  \rt endomorphism of 
$(X,\sigma^s)$, then 
$\varphi$ is an essentially  \rt endomorphism of $(X,\sigma)$, 
where \rt is any one of 
``$p$-L'', 
``$p$-R'' and  ``$p$-biresolving''.
\item 
Suppose that $\varphi$ is one-to-one or $\sigma$ is topologically 
transitive. 
\begin{enumerate} 
\item \cite{Nasu-t, Nasu-te}
If $\varphi$ is an essentially \rt endomorphism of 
$(X,\sigma^s)$, then 
$\varphi$ is an essentially \rt endomorphism of $(X,\sigma)$, 
where \rt is any one of ``$q$-R'', 
``$q$-L'', ``LR'', ``RL'',  and ``$q$-biresolving''. 
\item If $\varphi$ is an essentially LL (respectively, essentially RR) 
endomorphism of $(X,\sigma^s)$, then 
$\varphi$ is an essentially LL (respectively, essentially RR) 
endomorphism of $(X,\sigma)$.
\end{enumerate}
\end{enumerate}
\end{proposition} 
\begin{proof} 
The proposition is proved by Proposition 12.1 and Remark 2.10. 
\end{proof}
\begin{theorem} 
Let $\varphi$ be an onto endomorphism of an SFT $(X,\sigma)$. 
\begin{enumerate}
\item \cite{Nasu-t, Nasu-te}  
If $\varphi$ is directionally essentially \st, then
$\varphi$ is an essentially \st, 
where \rt is any one of 
``$p$-L'',
``$p$-R'' and  ``$p$-biresolving''.
\item 
Suppose that $\varphi$ is one-to-one or $\sigma$ is topologically 
transitive. 
\begin{enumerate} 
\item \cite{Nasu-t, Nasu-te}
If $\varphi$ is directionally essentially \st, then 
$\varphi$ is essentially \st, where \rt is any one of
``$q$-R'', ``$q$-L'', ``LR'', ``RL'',  and ``$q$-biresolving''. 
\item If $\varphi$ is directionally essentially LL (respectively, 
directionally essentially RR), then 
$\varphi$ is essentially LL (respectively, essentially RR). 
\end{enumerate}
\end{enumerate}
\end{theorem} 
\begin{proof} 
The theorem is proved by using Proposition 12.3 and 
Theorems 8.1(3), 8.6(2) and 8.7(2). 
\end{proof}

\end{document}